%% file: tesis.tex
\documentclass[a4paper,11pt]{book} 
\usepackage[utf8]{inputenc}
\usepackage{amsfonts, amsmath, amssymb}
\usepackage[inner=3cm,outer=3cm,bindingoffset=5mm]{geometry}
\usepackage{amsthm}
\usepackage{graphicx}
\usepackage[colorlinks=true,citecolor=black,linkcolor=black]{hyperref}
\usepackage{url}
\usepackage{subfigure}
\usepackage{fourier-orns}
\usepackage{indentfirst}
\usepackage{type1cm}
\usepackage{lettrine}
\usepackage[nottoc,notlot,notlof]{tocbibind}

\usepackage{fncychap}

\makeatletter
\def\@chapapp{}
\makeatother

\ChNumVar{\fontsize{50}{30}\usefont{OT1}{pzc}{m}{n}\selectfont\centering}
\ChTitleVar{\fontsize{20}{30}\usefont{OT1}{pzc}{m}{n}\selectfont\centering}

\usepackage{fancyhdr}
\newtheorem{tma}{Theorem}[chapter]
\newtheorem{maintma}{Main Theorem}
\newtheorem{prop}{Proposition}[chapter]
\newtheorem{corol}{Corollary}[chapter]
\newtheorem{defn}{Definition}[chapter]
\newtheorem{lema}{Lemma}[chapter]
\newtheorem{step}{Step}

\def\grad{\mathop{\rm grad}\nolimits}

\def\Hess{\mathop{\rm Hess}\nolimits}
\def\vol{\mathop{\rm vol}\nolimits}
\def\trace{\mathop{\rm trace}\nolimits}
\def\Ric{\mathop{\rm Ric}\nolimits}
\def\Rm{\mathop{\rm Rm}\nolimits}
\def\inj{\mathop{\rm inj}\nolimits}
\def\Vol{\mathop{\rm Vol}\nolimits}
\def\Area{\mathop{\rm Area}\nolimits}

\title{ Ricci flow on cone surfaces and a three-dimensional expanding soliton}
\author{Daniel Ramos}

\begin{document}
\renewcommand{\thepage}{\small{\Roman{page}}}
\include{portada}  \thispagestyle{empty}
\tableofcontents
\renewcommand{\thepage}{\arabic{page}}
\setcounter{page}{0} \thispagestyle{empty}
\include{intro}  \thispagestyle{empty}
\include{conesolitons}  \thispagestyle{empty}
\include{survey}  \thispagestyle{empty}
\include{cone_rf}  \thispagestyle{empty}

\include{smoothening}  \thispagestyle{empty}
\include{cuspsoliton}  \thispagestyle{empty}
\appendix
\makeatletter
\def\@chapapp{}
\makeatother
\include{compactness}  \thispagestyle{empty}

\bibliographystyle{alpha}
\bibliography{biblio}

\end{document}

%% file: portada.tex
%
%

\newgeometry{left=3cm,right=3cm,top=3cm,bottom=3cm}
\thispagestyle{empty} 
{\centering

\vfill
 
{
\Huge
\vspace{1em}
\includegraphics[width=0.8\textwidth]{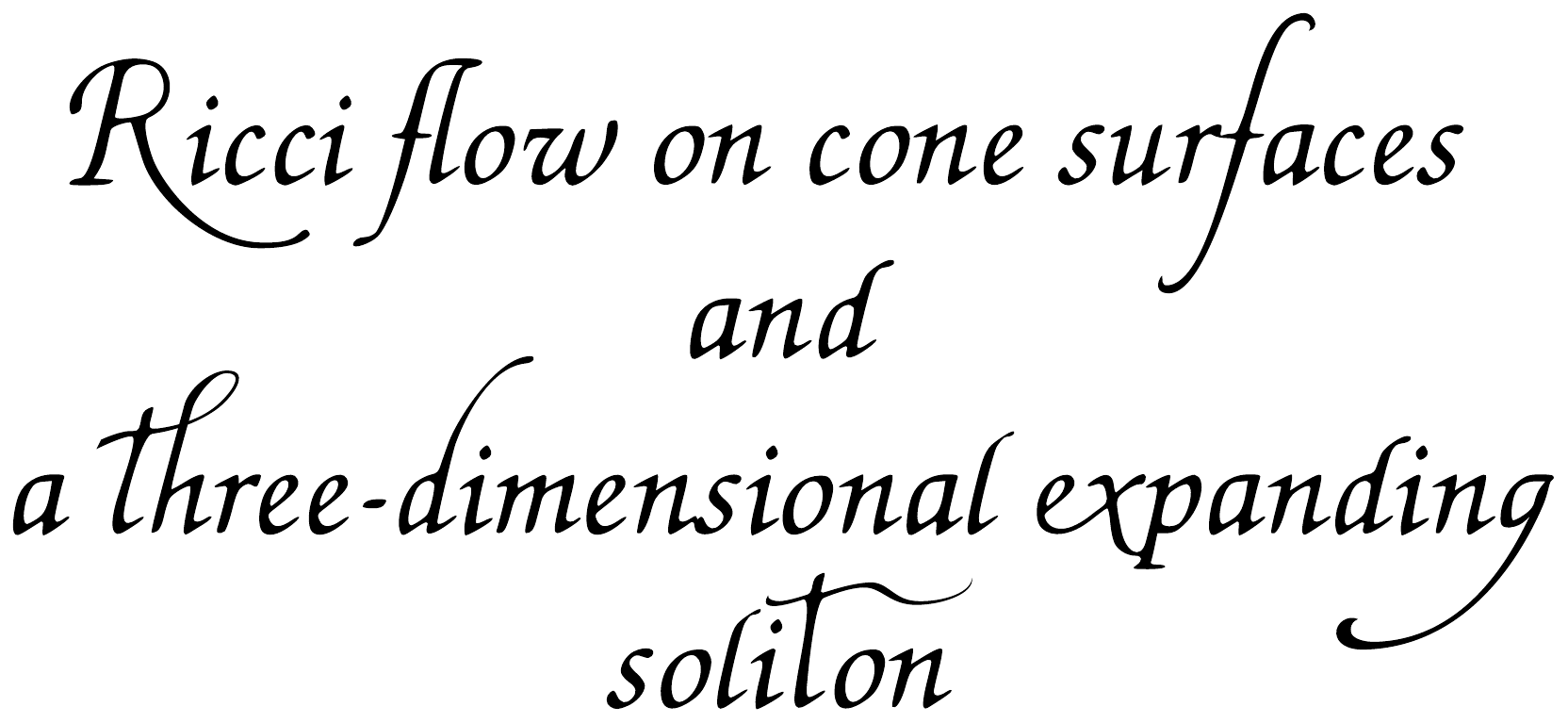}
\vspace{1em}
}
\vfill

\large
a Ph.D. dissertation to obtain the degree of \\ \textbf{Doctor in Mathematics} \\ presented by
\vfill
\LARGE
\bfseries Daniel Ramos Guallar \mdseries
\vfill
\large
advisored by Doctor
\vfill
\Large
\bfseries Joan Porti Piqué \mdseries
\vfill \vfill 
\large

Universitat Autònoma de Barcelona

December 2013

\vspace{1cm}
{\Huge \aldine}

}
\restoregeometry

\newpage \thispagestyle{empty} \null

\newpage \thispagestyle{empty} 

%

This is an archive version of the PhD thesis of the author, Daniel Ramos, defended on 28 January 2014 at Universitat Autònoma de Barcelona. This work is also available at Tesis Doctorales en Red (\url{http://hdl.handle.net/10803/133325}) since its defense. Some individual chapters are published in the referenced journals. This document is available on arXiv.
\medskip

Long-term contact: \texttt{daniel.ramosg@gmail.com}

%
%
%
%
%
%
%
%
%

\newpage \thispagestyle{empty} \null

\nocite{BBBMP}


%% file: intro.tex
\chapter{Introduction}

\lettrine{\indent E}{volution equations} are models for the change of certain quantities over time. One of
the main and earliest instances of evolution equations is the heat equation, that models the diffusion of the temperature on a medium,
$$\frac{\partial}{\partial t} u = \Delta u .$$
 Here $u(x,t)$ is the function determining the temperature at each point $x$ and
each time $t$ of a manifold ideally made of a homogeneous material. The Laplacian $\Delta= \trace\nabla^2$ is in general the Riemannian
Laplacian and depends on the metric of the manifold. This equation models two physical principles: the heat flows from points at high
temperature to points at low temperature; and
the speed of flowing is the bigger as the greater is the difference of temperatures.

Ricci flow is an evolution equation for a Riemannian metric on a manifold. That is, the quantity evolving is not scalar but
tensorial, the Riemannian metric tensor $g$. The Ricci flow equation is
\begin{equation}
 \frac{\partial}{\partial t} g = -2 \Ric(g) \label{RFequation}
\end{equation}
where $\Ric(g)$ is the Ricci curvature tensor of the time-dependent metric $g(t)$. This equation has strong similarities with the heat
equation. Both equations equal a first order time derivative with an expression involving second order space derivatives (although the heat
equation in $\mathbb R^n$ is a parabolic linear PDE and the Ricci flow is a weakly parabolic nonlinear PDE). The Laplacian is the trace of
the Hessian operator, that contains the space second derivatives of the function, and the Ricci tensor is one trace of the Riemann curvature
tensor
$$\Rm(X,Y)Z=\nabla^2_{X,Y} Z - \nabla^2_{Y,X} Z$$
that contains commutators of covariant second derivatives and actually depends on the second derivatives of the entries $g_{ij}$ of the
metric in a coordinate chart.
The evolution of the metric implies an evolution of the Riemann curvature tensor and other associated quantities, for instance the
evolution of the scalar curvature $R$ is 
$$\frac{\partial}{\partial t} R = \Delta R + 2|\Ric|^2 ,$$
which is almost a heat diffusion equation, perturbed by a reaction term. Thus, one can see the Ricci flow as an evolution equation for
the shape of a manifold, such that the several curvature quantities satisfy some reaction-diffusion equations. One would expect that,
similarly as the heat equation spreads the temperature and distributes it evenly across the manifold, the Ricci flow should smooth out the
curvature
in some sense. The Ricci flow theory is the study of the evolution of a given initial metric $g_0$ on a given manifold $\mathcal M$ under
the equation \eqref{RFequation}.

\medskip
Ricci flow was introduced in 1982 by Richard Hamilton in his seminal paper \emph{Three-manifolds with positive Ricci curvature}
\cite{Hamilton_3mfds}. There, the \emph{normalized} Ricci flow is defined as
$$ \frac{\partial}{\partial t}g = -2 \Ric + \frac{2}{n}rg ,$$
where 
$$r=\frac{\int R d\mu}{\int d\mu} $$
is the average scalar curvature, and this term serves to keep the volume constant along the flow. This flow is used there
to uniformize closed three-manifolds with $\Ric>0$. That is, any closed Riemannian three manifold with positive Ricci curvature evolves
under the Ricci flow to a metric of constant positive curvature.

In 1988 Hamilton published \emph{The Ricci flow on surfaces} \cite{Hamilton_surfaces}, where he used the two-dimensional normalized Ricci
flow 
$$\frac{\partial }{\partial t}g=(r-R)g$$
to  prove the existence and the long time behaviour of the flow on some surfaces $\mathcal M$. In dimension two, the average scalar
curvature depends only on the topology of the surface $\mathcal M$, and is constant along the flow,
$$r=\frac{4\pi\chi(\mathcal M)}{\Vol(\mathcal M)}$$
where $\chi(\mathcal M)$ is the Euler characteristic of the surface. Hamilton proved the convergence of the normalized Ricci flow to a
metric of constant curvature on the easier case $\chi(\mathcal M)\leq 0$ (i.e. $r\leq 0$), but not completely on the case $\chi(\mathcal
M)>0$ (the sphere), where he used an additional hypothesis of $R>0$. This case was hence much harder, and Hamilton developed a technique of
Harnack inequalities and entropies for surfaces with positive curvature that imply that the metric evolves to a so called soliton
metric. Soliton metrics for the normalized Ricci flow are defined to be metrics that evolve only by diffeomorphisms, and hence the
intrinsic shape remains unchanged. The unnormalized counterpart is that solitons for the Ricci flow are metrics that evolve only by
diffeomorphisms and homotheties. The initial metrics that give rise to a soliton are the same for normalized or unnormalized Ricci flow. 

Hamilton proved that the only smooth soliton on the topological sphere is the round sphere. When studying the soliton metrics for the flow
on the sphere, Hamilton showed (Theorem 10.1 in \cite{Hamilton_surfaces}) that when removing two
points (thus on a topological cylinder), there is a family  of possible soliton metrics, but when one tries to patch again the two points to
get a closed surface, only the constant curvature metrics  yield a smooth surface. Hence, the only two-dimensional soliton on the sphere is
actually the round sphere. However, if two singular  cone points were admitted, Hamilton noted that ``The other solutions we have found
exist on orbifolds''. We will bring an exhaustive and constructive enumeration of all solitons on smooth and cone surfaces. This includes Hamilton's cone solitons, and other new solitons in the noncompact cases.

In 1991, Bennet Chow completed the uniformization of smooth surfaces by Ricci flow in \emph{The Ricci flow on the 2-sphere}
\cite{Chow_sphere}, where he
modified the entropy formula by Hamilton, and was able to prove that under the normalized Ricci flow any metric on the sphere becomes
eventually a metric of positive curvature, thus Hamilton's result applies and this uniformizes the sphere. Years later, in 2006, a
brief note of X. Chen, P. Lu and G. Tian \cite{ChenLuTian}, simplified part of the proof of Hamilton by noting the remarkable fact that any
soliton solution on a smooth surface must be rotationally symmetric, thus simplifying considerably the equations in Theorem 10.1 of
\cite{Hamilton_surfaces}. This fact allowed them to
show that the Theorem of Uniformization of Surfaces can be proved by the Ricci flow (Hamilton used the uniformization in his proof).

Also in 1991, Lang-Fang Wu and Chow published three papers, \cite{Wu}, \cite{ChowWu} \cite{Chow_orbifolds}, dealing with the orbifold
solitons that Hamilton found.
These solitons occur on the teardrop and football orbifolds (the so-called bad orbifolds), the only two families of two-dimensional
orientable orbifolds that do not admit a constant curvature metric. Chow and Wu proved that any metric over a bad two-orbifold converges to
the soliton
solution under the normalized Ricci flow; first for positive curvature \cite{Wu}, and later for curvature changing sign \cite{ChowWu}, using
variations of the entropy technique \cite{Chow_orbifolds}. We will generalize this result to the broader class of cone surfaces, by means of
different techniques developed by G. Perelman.

In 2002 and 2003 Grisha Perelman posted his celebrated papers, \cite{Perelman1}, \cite{Perelman2}, \cite{Perelman3}, where he introduced
several important techniques for the
unnormalized Ricci flow, in the scope of three-dimensional smooth manifolds. One of this techniques involves the notion of
$\kappa$-noncollapsing and $\kappa$-solution, that allows one to relate and control the curvature, volume and injectivity radius of a Ricci-evolving manifold.
This is very useful for obtaining limits of sequences of manifolds and flows. In particular, every $\kappa$-solution of the Ricci flow
contains a shrinking \emph{asymptotic soliton} as $t\rightarrow -\infty$, and one has a clear classification of three-dimensional shrinking
solitons. In addition, Perelman proved that in the event of a singular time (the curvature going to infinity at some point), a blow-up
rescaling
of this point is modeled with a $\kappa$-solution. This fact allowed him to find a \emph{canonical neighbourhood} for the singular
points, that subsequently allows a surgery process, and eventually leaded to a proof of Thurston's Geometrization.

Perelman's techniques also apply in dimension two. This yields to a comparatively much simpler picture, due to several particularities of
the flow in this dimension. Firstly, the topology controls the evolution of the area. Under the unnormalized Ricci flow,
$$\frac{d}{dt}\Area(\mathcal M) = -\int_{\mathcal M} R d\mu = -4\pi\chi(\mathcal M) .$$
Hence, any initial metric on the sphere will shrink its area towards zero, and eventually will develop some type of singularity for finite
time. A blow-up rescaling technique may keep the area bounded and will bring a global model of the manifold. Secondly, all two-dimensional
asymptotic solitons are compact, and hence the only $\kappa$-solutions in dimension 2 are round spheres (this was clarified by R. Ye
\cite{Ye} in 2004), which essentially guarantees that a shrinking topological sphere rescales to a round sphere.

\bigskip
The main purpose of this thesis is the study of Ricci flow on cone surfaces. One possible and reasonable flow is the
\emph{angle-preserving flow}, that evolves the surface without changing the cone angles. This is the flow considered by Chow and Wu, and we
generalize their result by using Perelman's techniques. A second possible flow is the \emph{smoothening flow} that instantaneously removes
the cone points and makes the surface smooth. We study this flow using tecnhiques developed by Peter Topping for removing cusps with Ricci
flow \cite{Topping_revcusp}. A second purpose of this thesis is the study of some gradient Ricci solitons, in dimension two and three. We
use phase portraits and other techniques from dynamical systems to give very explicit constructions of these solitons.
\medskip

There are four main results in this thesis. Three of them are related to the theory of Ricci flow on cone surfaces, plus a fourth result
on a three-dimensional example of a soliton.

\medskip

First result, in Chapter \ref{Ch:conesolitons}, is an exhaustive enumeration of all complete gradient Ricci solitons on smooth and cone
surfaces, both compact
and noncompact, with an arbitrary lower bound on the curvature. This is achieved by exploiting the symmetry properties of solitons, that
give rise to a specific ODE in polar geodesic coordinates. This ODE can be studied by means of a phase portrait analysis. This result has
been submitted to publication in \cite{Ramos3}.

\begin{maintma}
All gradient Ricci solitons on a surface, smooth everywhere except possibly on a discrete set of cone-like singularities, complete, and
with curvature bounded below fall into one of the following families:
\begin{enumerate}
 \item Steady solitons:
       \begin{enumerate}
        \item Flat surfaces.
        \item The smooth cigar soliton.
	\item The cone-cigar solitons of angle $\alpha\in(0,+\infty)$.
       \end{enumerate}
 \item Shrinking solitons:
       \begin{enumerate}
        \item Spherical surfaces.
	\item Teardrop and football solitons, on a sphere with one or two cone points.
	\item The shrinking flat Gaussian soliton on the plane.
	\item The shrinking flat Gaussian cones.
       \end{enumerate}
 \item Expanding solitons:
       \begin{enumerate}
        \item Hyperbolic surfaces.
	\item The $\alpha\beta$-cone solitons, with a cone point of angle $\beta>0$ and an end asymptotic to a cone of angle $\alpha>0$.
	\item The smooth blunt $\alpha$-cones.
	\item The smooth cusped $\alpha$-cones in the cylinder, asymptotic to a hyperbolic cusp in one end and asymptotic to a cone of angle
$\alpha>0$ in the other end.
	\item The flat-hyperbolic solitons on the plane, that are universal coverings of the cusped cones.
	\item The expanding flat Gaussian soliton on the plane.
	\item The expanding flat Gaussian cones.
       \end{enumerate}
\end{enumerate}
\end{maintma}
Explicit descriptions of these solitons are given in the corresponding sections of Chapter \ref{Ch:conesolitons}, and some pictures
are drawn in Figures \ref{cigar_pic} to \ref{cone4_pic} in the same chapter. 

\medskip

Second result is a uniformization theorem for closed cone surfaces with cone angles less than or equal to $\pi$. This involves two
chapters and the appendix. In Chapter \ref{Ch:survey} we survey some known results, mostly from Hamilton and Perelman, and we draw a line of
argument
for proving the uniformization of smooth surfaces. This uniformization was already proven by Hamilton \cite{Hamilton_surfaces} and Chow
\cite{Chow_sphere}, but we
propose a different path for handling the case of positive Euler characteristic (the sphere) by using Perelman's $\kappa$-solutions
technique (originally developed for the three-dimensional case in \cite{Perelman1}). This technique uses several rescaling blow-ups that
need some
compactness theorems for classes of manifolds to ensure the existence and nondegeneracy of the limits of sequences of rescalings. In
Chapter \ref{Ch:cone_rf} we adapt the line of argument of Chapter \ref{Ch:survey} to the case of cone surfaces. We discuss the notion of
cone surface and the
existence theorems for the angle-preserving flow (given by the work of H. Yin \cite{Yin1}, \cite{Yin2} and R. Mazzeo, Y. Rubinstein and
N. Sesum \cite{MazRubSes}), and we adapt
some maximum principles and Harnack inequalities to work on the cone setting. Finally, in Appendix \ref{Ch:compactness} we prove the cone
version of the compactness theorems for classes of surfaces and classes of flows required to complete the proof.

\begin{maintma}
Let $(\mathcal M,(p_1,\ldots,p_n),g_0)$ be a closed cone surface, and assume that the cone angles are less than or equal to $\pi$. Then
there is an angle-preserving Ricci flow that converges, up to rescaling, to either
\begin{itemize}
 \item a constant nonpositive curvature metric, if $\hat\chi(\mathcal M)\leq 0$, or
 \item a spherical (constant positive curvature) metric, a teardrop soliton or a football soliton, if $\hat\chi(\mathcal M)> 0$,
\end{itemize}
where $\hat\chi(\mathcal M)$ is the conic Euler characteristic of the cone surface.
\end{maintma}

\medskip

Third result, in Chapter \ref{Ch:smoothening}, discusses a different flavour of Ricci flow on cone surfaces, that exposes the nonuniquenes
of the
solutions to the Ricci equation for initial cone surfaces. Whereas the flow discussed in Chapter \ref{Ch:cone_rf} keeps the angles fixed,
this is
only a boundary condition for the PDE problem, and can be altered. We discuss a smoothening Ricci flow, that has as initial condition a cone
manifold and instantaneously removes the cone points and turns a neighbourhood of the cone point into a smooth disc, although with very high
curvature. This is achieved by approximating the initial cone manifold with smooth, highly curved surfaces that resemble the vertices of the
cone points. Then the standard smooth Ricci flow is applied to the approximating surfaces, and a compactness theorem (different from the one
used in Appendix \ref{Ch:compactness}) is used to get a limit smooth flow, that has as
initial condition (limit as $t\rightarrow 0$) the
original cone surface we started with. This procedure was inspired by the work of P. Topping \cite{Topping_revcusp}, where he removes cusp
singularities on surfaces with Ricci flow. On that point of view, cusps are zero-angle cone points and we adapt the barriers and bounds on
that work to nonzero angles. This result has been accepted to publication in \cite{Ramos1}.

\begin{maintma}
Let $(\mathcal M,(p_1\ldots p_n),g_0)$ be a closed cone surface with bounded curvature. There exists a Ricci flow $g(t)$ smooth on the
whole $\mathcal M$, defined for $t\in (0,T]$ for some $T$, and such that
$$g(t) \rightarrow g_0 \quad \mathrm{ as } \quad t \rightarrow 0^+ .$$
This Ricci flow has curvature unbounded above and uniformly bounded below over time.

Furthermore, any other smooth Ricci flow $\tilde g(t)$ on $\mathcal M$, defined for $t\in (0,\delta]$ for some $\delta < T$, such
that $\tilde g(t) \rightarrow g_0$ as $t \rightarrow 0^+ $ and such that its curvature is uniformly bounded below agrees with
the stated flow $g(t)$ for $t\in (0,\delta]$.
\end{maintma}

\medskip

Fourth result, in Chapter \ref{Ch:cuspsoliton}, is a construction and study of a new three-dimensional expanding soliton. Although this is
not directly
related to cone surfaces, the techniques are very similar and were inspired by the techniques used in Chapter \ref{Ch:conesolitons} to
classify two-dimensional solitons. We first describe the expanding soliton on the topological manifold $\mathbb R \times \mathbb T^2$, the product of
the real line with a two-torus, by assuming a nice global coordinate chart and a metric in the form of a warped product. Next we prove
that this example is the unique with this topology and a lower bound on the curvature $\sec>-\frac{1}{4}$. Furthermore, our example is
critical in the following sense: it is known \cite[Lem 5.5, Rmk 5.6]{Caoetal} that any three-dimensional expanding soliton with 
$sec > c > -\frac{1}{4}$ must be topologically $\mathbb R^3$, and therefore our example shows that the bound is sharp. This result has been
submitted to publication in \cite{Ramos2}.

\begin{maintma} \label{T:cuspsoliton_intro}
There exists an expanding gradient Ricci soliton $(M,g,f)$ over the topological manifold $M = \mathbb R \times \mathbb T^2$
satisfying the
following properties:
\begin{enumerate}
 \item The metric has pinched sectional curvature $-\frac{1}{4} < sec < 0$.
 \item The soliton approaches the hyperbolic cusp expanding soliton on one end.
 \item The soliton approaches locally the flat Gaussian expanding soliton on a cone on the other end.
\end{enumerate}
Furthermore, this is the only nonflat gradient Ricci soliton over the topological manifold $M = \mathbb R \times \mathbb T^2$ with curvature
$sec > -\frac{1}{4}$.
\end{maintma}

\bigskip
We expect that the results in this thesis can bring some insights and intuition to the phenomena appearing in two-dimensional Ricci flow,
in the study of cone surfaces, and in the description of low-dimensional Ricci solitons.

%% file: conesolitons.tex
\chapter{Gradient Ricci solitons on smooth and cone surfaces} \label{Ch:conesolitons}

\lettrine{\indent G}{radient Ricci solitons} are special self-similar solutions to the Ricci flow, and were introduced in
\cite{Hamilton_surfaces}
when R. Hamilton developed the Ricci flow theory for surfaces (using there a normalized version of Ricci flow). Hamilton proved that all
closed smooth surfaces with positive
curvature converge under the Ricci flow to a gradient Ricci soliton, and proved that the only solitons on a closed smooth surface are
those of constant curvature \cite[Thm 10.1]{Hamilton_surfaces}. B. Chow subsequently was able to remove the positive curvature hypothesis
in \cite{Chow_sphere}. In the same original work of Hamilton, it is described an open gradient soliton with nonconstant curvature known as
the steady cigar soliton. Also, in the course of the proof of Theorem 10.1 of \cite{Hamilton_surfaces} Hamilton found some solitons on the
topological sphere with
cone-like singularities, in particular over orbifolds. The study of the Ricci flow
converging to these orbifold solitons was carried out by L.-F. Wu  and Chow \cite{Wu}, \cite{ChowWu}, \cite{Chow_orbifolds}. This
brings to all
2-orbifolds a natural metric, including those orbifolds that do not admit a constant curvature metric (footballs and teardrops). The
existence of the football and teardrop solitons with any cone angles (not orbifold) was also considered by H. Yin \cite{Yin1},
but he gave no explicit construction.

A very useful fact in dimension two is that all nonconstant curvature gradient solitons admit a nontrivial Killing
vector field and have a rotational symmetry (see \cite[pp 241-242]{3CY} and \cite{ChenLuTian}). This allows one to pick polar coordinates
and set a
single ODE for the soliton metric. The cigar soliton appears as an explicit solution for the steady case. See \cite{TRFTA1} for
further reference.

In this chapter we gather and re-order these results, and we apply a phase portrait analysis for the ODE associated to the soliton metric.
We
obtain a complete and unified classification of all two-dimensional gradient solitons on smooth and cone surfaces. By ``surfaces'' we will
mean open or closed
topological surfaces, endowed with a complete Riemannian metric, smooth everywhere except possibly in a discrete set
of cone-like singular points, that will arise naturally in our discussion. For geometric considerations, on the theorem statements we will
discard all
Riemannian metrics without a lower bound of the Gaussian curvature, although we will find these examples in the course of the proofs. This
solves
some questions proposed in \cite[p 51]{TRFTA1}, and describes some two-dimensional solitons (smooth and conic) that didn't exist in
the literature. For the sake of completeness, we will include some known results with their proofs, so the exposition in this chapter is
self-contained. The main theorem of the chapter is the following.

\begin{tma} \label{T:coneslt_main}
All gradient Ricci solitons on a surface, smooth everywhere except possibly on a discrete set of cone-like singularities, complete, and
with curvature bounded below fall into one of the following families:
\begin{enumerate}
 \item Steady solitons:
       \begin{enumerate}
        \item Flat surfaces.
        \item The smooth cigar soliton.
	\item The cone-cigar solitons of angle $\alpha\in(0,+\infty)$.
       \end{enumerate}
 \item Shrinking solitons:
       \begin{enumerate}
        \item Spherical surfaces.
	\item Teardrop and football solitons, on a sphere with one or two cone points.
	\item The shrinking flat Gaussian soliton on the plane.
	\item The shrinking flat Gaussian cones.
       \end{enumerate}
 \item Expanding solitons:
       \begin{enumerate}
        \item Hyperbolic surfaces.
	\item The $\alpha\beta$-cone solitons, with a cone point of angle $\beta>0$ and an end asymptotic to a cone of angle $\alpha>0$.
	\item The smooth blunt $\alpha$-cones.
	\item The smooth cusped $\alpha$-cones in the cylinder, asymptotic to a hyperbolic cusp in one end and asymptotic to a cone of angle
$\alpha>0$ in the other end.
	\item The flat-hyperbolic solitons on the plane, that are universal coverings of the cusped cones.
	\item The expanding flat Gaussian soliton on the plane.
	\item The expanding flat Gaussian cones.
       \end{enumerate}
\end{enumerate}
\end{tma}
Each family of solitons is described in the corresponding section. In Section \ref{sec_basics} we recall briefly the definition of Ricci
solitons in dimension 2, and their properties of symmetry in the case of nonconstant curvature, that lead to a first-order ODE system. In
Section \ref{sec_constcurv} we enumerate the closed solitons with
constant curvature. In Sections \ref{sec_steady}, \ref{sec_shrink} and \ref{sec_expand} we study the ODE system in the steady, shrinking and
expanding cases, respectively. These three parts combined prove Theorem \ref{T:coneslt_main}. Finally, in Section \ref{sec_gallery} we
bring a gallery
of solitons embedded into $\mathbb R^3$, drawn with Maple.

\section{Gradient solitons and rotational symmetry} \label{sec_basics}

In this section we recall the basics of Ricci solitons and their properties of symmetry. See \cite{TRFTA1} for an extended introduction.

A \emph{Ricci flow} is a PDE evolution equation for a Riemannian metric $g$ on a smooth manifold $\mathcal M$,
\begin{equation} \label{rf_eqn}
 \left\{ \begin{aligned}
          \frac{\partial}{ \partial t}g(t) &= -2 \Ric_{g(t)} \\ 
	  g(0) &=g_0 .
         \end{aligned} \right.
\end{equation}
 A \emph{Ricci soliton} is a special type of self-similar solution of the Ricci flow, in the form
 \begin{equation}
 \label{solit_defn} g(t)=c(t)\ \phi_t^* (g_0)
\end{equation}
where for each $t$, $c(t)$ is a constant and $\phi_t$ is a diffeomorphism. So, $g(0)=g_0$,
$c(0)=1$ and $\phi_0=id$. The family $\phi_t$ is the flow associated to a (maybe time-dependent) vector field $X(t)$; and
in the case when this vector field is the gradient field of a function, $X=\grad f$, the soliton is
said to be a \emph{gradient} soliton. In this case, differenciating the definition of soliton \eqref{solit_defn} gives
\begin{align*}
-2 \Ric_{g(t)} =  \frac{\partial}{\partial t} g(t) &= \dot c(t)\phi_t^*g_0 + c(t)\mathcal L_{\grad f} \phi_t^*(g_0)\\
&= \dot c(t) g_0 + 2\Hess_{g(t)} f .
\end{align*}
Since $\Ric_{g(t)} = \Ric_{c(t)\phi_t^*g_0} = \phi_t^*\Ric_{g_0}$, we get
\begin{equation}\label{solit_1}
\dot c(t)\phi_t^*g_0 + 2\Hess_{g(t)} f = -2 \phi_t^*\Ric_{g_0}
\end{equation}
and evaluating at $t=0$,
$$ \epsilon g_0 +  2\Hess_{g_0} f= -2\Ric_{g_0} $$
where $\epsilon = \dot c(0)$. The soliton is said to be \emph{shrinking}, \emph{steady} or \emph{expanding} if the constant $\epsilon$ is
negative, zero or positive respectively. This constant can be normalized to be $-1$, $0$ or $+1$ respectively, being this equivalent to
reparameterize the time $t$. Therefore, a gradient Ricci flow is a triple $(\mathcal M, g, f)$ satisfying
\begin{equation}
 \Ric + \Hess f + \frac{\epsilon}{2} g =0 . \label{solit_eqn_init}
\end{equation}
This is the soliton equation for the initial manifold $\mathcal M$. Conversely, given $(\mathcal M, g_0,f)$ and $\epsilon$
satisfying \eqref{solit_eqn_init}, we can recover the soliton in the form \eqref{solit_defn}. For we use on \eqref{solit_1} the fact
that both $\Ric$ and $\Hess f$ are invariant under rescaling and equivariant under linear transformations, so
$$
\dot c(t) \phi_t^*(g_0) = -2\Hess_{g(t)} -2 \Ric_{g(t)} = \phi_t^*(-2\Hess_{g_0}f -2\Ric_{g_0}) = \phi_t^*(\epsilon\ g_0) 
= \epsilon\phi_t^*(g_0) 
$$
and hence $\dot c(t)=\epsilon\ \forall t$, that is,
$$c(t)=\epsilon t+c(0)=\epsilon t+1 .$$
That way we can rewrite (\ref{solit_1}) as
\begin{equation}
 \frac{\epsilon}{\epsilon t+1}g(t) + 2\Hess_{g(t)} f +2\Ric_{g(t)} =0
\end{equation}
with $g(t)=(\epsilon t+1)\phi_t^*(g_0)$.
This is the dynamical soliton equation along the time. In the nonsteady cases, $\epsilon\neq 0$, also called \emph{homothetic solitons}, we
can do the parameter change $t \mapsto t-\frac{1}{\epsilon}$ and obtain
\begin{equation}
 \label{gr_shr_solit} \frac{1}{2t}g(t) + \Hess_{g(t)} f + \Ric_{g(t)} =0
\end{equation}
with $g(t)=\epsilon t\phi_{(t-1/\epsilon)}^*(g(\frac{1}{\epsilon}))$.

\bigskip
In this chapter we will seek for the initial surfaces of two-dimensional solitons. In the two-dimensional case we have $\Ric =\frac{R}{2}g$,
hence the soliton equation \eqref{solit_eqn_init} becomes
\begin{equation}
 \Hess f + \frac{1}{2} (R+\epsilon)g=0 . \label{gr_shr_solit_2d}
\end{equation}

This equation only makes sense on a smooth Riemannian surface. However, we will allow some cone-like singularities for the surface, namely
points such that admit a local coordinate chart in the form
$$ dr^2 + h(r,\theta)^2 \ d\theta^2$$
for some smooth $h:[0,\delta)\times \mathbb R /2\pi\mathbb Z \rightarrow \mathbb R$ such that $h(0,\theta)=0$ and $\frac{\partial
h}{\partial r} = \frac{\alpha}{2\pi}$ needs not to be $1$. The value $\alpha$ is the cone angle at this point. To make this point smooth, it
is required that $\alpha =2\pi$ and $\frac{\partial^{2k} h}{\partial r^{2k}} =0$ for all $k\in \mathbb N$ (see
\cite[p 450]{TRFTA1}).

We will use some properties of two-dimensional gradient solitons to turn the tensor equation \eqref{gr_shr_solit_2d} into a much
simpler first order vector ODE that will allow a subsequent qualitative analysis. The main property we will use is that there exists a
Killing vector field (given by a rotation of
$\grad f$) over the smooth part of $\mathcal M$. The associated line flow of
this field is a
one-parameter group acting globally by isometries, this group must be $\mathbb S^1$ and so $\mathcal M$ is
rotationally symmetric (this argument is from \cite{3CY} and \cite{ChenLuTian}). We define then a rotationally symmetric polar coordinate
chart on $\mathcal
M$, and the analysis of the local expression of the soliton
equation will give us the ODE system that satisfies the explicit metric over $\mathcal M$.

\medskip

Let $J: T\mathcal M \rightarrow T\mathcal M$ be an almost-complex structure on $\mathcal M$, that
is, a $90^{\circ}$ rotation on the positive orientation sense, so
$$ J^2=-Id \quad , \quad g(X,JX)=0 \quad \forall X\in \mathfrak X  \mathcal M .$$

\begin{lema} Some basic properties of $J$ are
 \begin{enumerate}
  \item $g(JY,Z)=-g(Y,JZ)$,
  \item $J$ anticommutes with $\flat$: $\flat (JX) = -J(\flat X)$,
  \item $J$ commutes with $\nabla$.
 \end{enumerate}
\end{lema}
\begin{proof}  First property is elementary,
\begin{eqnarray*}
 0=g(Y+Z,J(Y+Z))&=&g(Y,JY)+g(Z,JY)+g(Y,JZ)+g(Z,JZ)\\&=&g(Z,JY)+g(Y,JZ).
\end{eqnarray*}
A warning about the notation: if $\flat(X)=\omega$ is a 1-form then $J: T^* \mathcal M \rightarrow
T^* \mathcal M$ is defined as $\omega \mapsto J\omega$ where $J\omega(W)=\omega(J(W))$. 

So, the second statement is just
$$\flat(JX)(W)=g(JX,W)=-g(X,JW)=-\flat(X)(JW)=-J(\flat X)(W) .$$

For the last statement, let $X\in \mathfrak X\mathcal M$, then $\{X, JX\}$ form a basis of
$T\mathcal M$. Let us see that $J(\nabla X)$ and $\nabla (JX)$ have the same projections over the
basis. Since $g(X,JX)=0$, taking covariant derivatives we get
$$g(\nabla X,JX) + g(X,\nabla (JX))=0$$
which implies
$$g(J(\nabla X),-X) + g(X,\nabla (JX))=0.$$
Again, differenciating $g(X,X)=g(JX,JX)$ we obtain
$$2g(X,\nabla X) = 2g(\nabla(JX),JX) ,$$
which implies
$$g(JX,J(\nabla X)) = g(JX, \nabla(JX)) .$$
\end{proof}

This construction shows that a two-dimensional gradient soliton admits a Killing vector field ( \cite[p.
241]{3CY}, \cite[p. 11]{TRFTA1}, and
\cite{Cao_KRsolitons1})
\begin{lema}
The vector field $J(\grad f)$ is a Killing vector field.
\end{lema}
\begin{proof}
Recall that a Killing vector field $W$ is such that its line flow is by isometries, or
equivalently, the metric tensor $g$ is invariant under the line flow, that can be expressed in
terms of the Lie derivative as $\mathcal L_W g =0$. Recall also that 
$$\mathcal L_W g \ (Y,Z)=\nabla \omega (Y,Z) + \nabla \omega (Z,Y)$$
where $\omega = \flat(W) = g(W,\cdot)$. Then
\begin{eqnarray*}
 \mathcal L_{J(\grad f)} g \ (Y,Z) &=& \nabla(\flat(J(\grad f)))(Y,Z) +
\nabla(\flat(J(\grad f)))(Z,Y)\\
                                   &=& \nabla J \flat \grad f (Y,Z) + \nabla J \flat \grad f (Z,Y)\\
                                   &=& \nabla J \nabla f (Y,Z) + \nabla J \nabla f (Z,Y)\\
                                   &=& J \nabla \nabla f (Y,Z) + J \nabla \nabla f (Z,Y)\\
                                   &=& \nabla^2 f (JY,Z) + \nabla^2 f (JZ,Y)\\
                                   &=& \frac{-1}{2}(R+\epsilon)g(JY,Z) + \frac{-1}{2}(R+\epsilon)g(JZ,Y)\\
                                   &=& -\frac{1}{2}(R+\epsilon)\Big( g(JY,Z) + g(JZ,Y)\Big) =0.
\end{eqnarray*}
Note that $\nabla f = df =\flat(\grad f)$, and again $J(\nabla f)$ is the 1-form $A\mapsto \nabla
f(JA)$ and $J(\nabla^2f)$ is the 2-covariant tensor field $(A,B)\mapsto \nabla^2f(JA,B)$.
\end{proof}

The Killing vector field may be null if the $\grad f$ field itself is null, otherwise, the surface admits a symmetry.

\begin{lema}
Let $(\mathcal M, g, f)$ be a gradient Ricci soliton on a surface. Then, at least one of the following holds:
\begin{enumerate}
 \item $\mathcal M$ has constant curvature.
 \item $\mathcal M$ is rotationally symmetric (i.e. admits a $\mathbb S^1$-action by isometries).
 \item $\mathcal M$ admits a quotient that is rotationally symmetric.
\end{enumerate}
Besides, if the surface has not constant curvature, no more than two cone points may exist.

\end{lema}

\begin{proof}
We adapt the argument for the smooth closed case from \cite{ChenLuTian}. We will discuss in terms of $\grad f$. If $\grad f \equiv 0$, then
by the soliton equation \eqref{gr_shr_solit_2d} we have $R=-\epsilon$ and the
curvature is constant. Let us assume then that $f$ is not constant everywhere. Therefore $J(\grad f)$ is a nontrivial Killing vector
field and its line flow, $\phi_t$, is a one-parameter group
acting over $\mathcal M$ by isometries.

Suppose that $\grad f$ has at least one zero in a point $O\in \mathcal M$. This is the case of closed smooth surfaces. 
The point $O$ is a zero of the vector field $J(\grad f)$, so it is a fixed
point of $\phi_t$ for every $t$. Then, $\phi_t$ induces $\phi^*_t$ acting on $T_O\mathcal M$ by
isometries of the tangent plane, so we conclude that the group $\{\phi_t\}$ is $\mathbb S^1$ acting by
rotations on the tangent plane. Via the exponential map on $O$, the action is global on $\mathcal
M$ and therefore the surface is rotationally symmetric.

Suppose now that $\grad f$ has no zeroes but the surface contains a cone point $P$. Then the flowlines of $\phi_t$
cannot pass through $P$, because there is no local isometry between a cone point and a smooth one. So this point $P$ is
fixed by $\phi_t$ for every $t$ and, via the exponential map, $\phi_t$ induces $\phi^*_t$ acting
on $C_P\mathcal M$ the tangent cone (space of directions) on $P$. Again, a continuous one-parameter
subgroup of the metric cone $C_P\mathcal M$ must be the $\mathbb S^1$ group acting by rotations. Besides,
if other cone points were to exist, these should also be fixed by the already given $\mathbb S^1$ action.
This implies that no more than two cone points can exist on $\mathcal M$, for otherwise the minimal
geodesics joining $P$ with two or more conical points would be both fixed and exchanged by some
$\mathbb S^1$ group element. Note that in the case of two cone points, these need not to have equal cone
angles.

Finally suppose that $\grad f$ has no zeroes and the surface has no cone points. Then the surface is smooth and the flowlines of $\grad f$
are all of them isomorphic to $\mathbb R$ (no closed orbits can appear for the gradient of a function) and foliate the surface. The
action of $\phi_t$ exchanges the fibres of this foliation. The parameter of $\phi_t$ is $t\in \mathbb S^1$ or $t\in \mathbb R$. In the
first case, $\mathbb S^1$ is acting on $\mathcal M$ and it is rotationally symmetric. In the second case, $\mathcal M
\cong \mathbb R^2$, and the flowline $\phi_t$ of the Killing vector field induces a $\mathbb Z$-action by isometries by
$$ x \mapsto \phi_1(x)$$
that acts freely on $\mathcal M$ since no point is fixed by $\phi_t$ for any $t\neq 0$ (if $\phi_t(p)=p$, then all fibres are fixed and
every point in each fibre also is, so $\phi_t = id$). Then the quotient by this action is topologically ${\mathcal M} /_\sim
\cong \mathbb R \times \mathbb S^1 $ and is rotationally symmetric. We will find nontrivial examples of this solitons as cusped expanding
solitons and their universal coverings.
\end{proof}

The fact of being rotationally symmetric allows us to endow $\mathcal M$ with polar coordinates
$(r,\theta)$ such that the metric is given by
$$g=dr^2 + h^2(r) \ d\theta^2$$
where $r\in I \subseteq \mathbb R$ is the radial coordinate, and $\theta\in \mathbb R / 2\pi\mathbb Z$ is the periodic angular coordinate.
The function $h(r)$ does not depend on $\theta$ because of the rotational symmetry; and similarly, the potential function only depends on
the $r$ coordinate, since $\grad f$ is a radial vector field. Surfaces not rotationally symmetric but with a rotationally symmetric quotient
also admit these
coordinates, changing only $\theta \in \mathbb R$.

\begin{lema}
Given the polar coordinates $(r,\theta ) \in \mathbb R \times \mathbb R / 2\pi\mathbb Z$ and the metric in the form $g=dr^2 + h^2(r) \
d\theta^2$, 
the Gaussian curvature (which equals half the scalar curvature) is given by
$$K=\frac{R}{2}=\frac{-h''}{h} ,$$
and the Hessian of a radial function $f(r)$ is given by
$$\Hess f = f''\ dr^2 + hh'f'\ d\theta^2 .$$
\end{lema}

\begin{proof}
It is a standard computation. Covariant derivatives are given by
$$\nabla_{\partial_r} \partial_r = 0 \qquad \nabla_{\partial_r} \partial_\theta = \frac{h'}{h}
\partial_\theta \qquad \nabla_{\partial_\theta} \partial_\theta = -hh' \partial_r $$
Then we contract twice the curvature tensor $R(X,Y)Z = \nabla_X (\nabla_Y Z) - \nabla_Y (\nabla_X Z) - \nabla_{[X,Y]}Z$ for the scalar
curvature, and apply $\Hess f (X,Y)=X(Y(f))-(\nabla_X Y)(f)$ for the Hessian.
\end{proof}

On that rotationally symmetric setting, the soliton equation becomes
$$\Hess f + \frac{1}{2} (R+\epsilon)g = \left( f''-\frac{h''}{h} +\frac{\epsilon}{2} \right) dr^2 + \left( hh'f' +\left(
-\frac{h''}{h} +\frac{\epsilon}{2}
\right) h^2 \right) d\theta^2 = 0 ,$$
which is equivalent to the second order ODEs system
\begin{equation}
 \left\{ \begin{array}{rcl}
	  f''-\frac{h''}{h} +\frac{\epsilon}{2} &=& 0 \\
	  \frac{h'}{h}f' -\frac{h''}{h} +\frac{\epsilon}{2} &=& 0.
         \end{array} 
\right.\end{equation}
We combine both equations to obtain
$$\frac{f''}{f'}=\frac{h'}{h} ,$$
and integrating this equation,
$$\ln f' = \ln h + C $$
so 
\begin{equation}
f' = ah \label{2dsolit_eq_potential}
\end{equation}
for some $a>0$. Hence, substituting on the system we obtain a single ODE,
\begin{equation}
 h'' -ahh' - \frac{\epsilon}{2}h=0. \label{2dsolit_eq}
\end{equation}

We summarize the computations in the following lemma,
\begin{lema}
Let $(\mathcal M, g,f)$ be a gradient Ricci soliton on a surface with nonconstant curvature. Then $\mathcal M$ admits coordinates
$(r,\theta)$, with
$r\in I \subseteq \mathbb R$ and $\theta\in \mathbb S^1$ or $\theta \in \mathbb R$, such that the metric takes the form $g=dr^2 + h^2(r) \
d\theta^2$ for some function $h=h(r)$ satisfying \eqref{2dsolit_eq}, and the potential is $f=f(r)$ satisfying
\eqref{2dsolit_eq_potential}.
\end{lema}

Setting $h'=u$, the second order ODE \eqref{2dsolit_eq} becomes a vector first order ODE
\begin{equation}
 \left\{ \begin{array}{rcl}
   h' &=& u \\
   u' &=& (au +\frac{\epsilon}{2})h .
  \end{array} \right. \label{2dsolit_sys}
\end{equation}
The solutions to system \eqref{2dsolit_sys} (and equation \eqref{2dsolit_eq}) are functions $h(r)$ that define rotationally symmetric
metrics on the cylinder $(r,\theta)\in \mathbb R \times \mathbb S^1$. This cylinder may be pinched in one or both ends, thus changing the
topology of the surface. The pinching appears as zeros of $h$. Compactness condition of the surface is equivalent to the boundary conditions
$$h(0)=0 \quad \mbox{and} \quad h(A)=0$$
for some $A>0$ such that $h(A)=0$. In this case, one or two cone angles may appear,
$$h'(0)=\frac{\alpha_1}{2\pi} \quad \mbox{and} \quad h'(A)=-\frac{\alpha_2}{2\pi}$$
where $\alpha_1$ and $\alpha_2$ are the cone angles. Smoothness conditions would be $h'(0)=1$ and $h'(A)=-1$, plus the condition
$$h^{(2k)}=0$$
at $r=0$ and $r=A$ for all $k \geq 0$. This condition ensures $\mathcal C^\infty$ regularity (\cite{TRFTA1}, Lemma A.2). This condition
holds on our solitons, since derivating $2k$ times the equation \eqref{2dsolit_eq} we obtain
$$h^{(2+2k)}-a\left( \sum_{i+j=2k} h^{(i)} h^{(j+1)} \right) -\frac{\epsilon}{2}h^{(2k)} =0 .$$
Since $i+j$ is an even number, both $i$ and $j$ must be even or odd. In both cases, there is an even index in $h^{(i)} h^{(j+1)}$. Thus
inductively, if all even-order derivatives vanish at $r=0$ up to order $2k$, then also vanishes the $2k+2$ derivative at $r=0$ (idem at
$r=A$).

We shall study the system \eqref{2dsolit_sys} for steady, shrinking and expanding solitons to obtain a complete enumeration of gradient
Ricci solitons on surfaces of
nonconstant curvature.

\medskip
Incidentally, it is interesting to note some
geometric interpretations of the functions $f$ and $h$.

\begin{lema}
Let $(\mathcal M, g, f)$ be a rotationally symmetric gradient Ricci soliton in dimension 2, then
\begin{align*}
 \grad f &= f'(r) \partial_r = a h(r) \partial_r ,\\
 J(\grad f) &= \frac{f'(r)}{h(r)} \partial_\theta = a \partial_\theta ,\\
 K &= -\frac{h''}{h}=-\left( ah' +\frac{\epsilon}{2} \right) .
 \end{align*}
If $p\in \mathcal M$ is a (smooth or conic) center of rotation. Then, using the distance to $p$ as the $r$-coordinate,
\begin{align*}
  h(r) &= \frac{1}{2\pi}\mathrm{Perimeter}(\mathrm{Disc}(r)) ,\\
  f(r) &= \frac{a}{2\pi}\mathrm{Area}(\mathrm{Disc}(r)) + f(0) ,\\
\end{align*}
where $\mathrm{Disc}(r)$ is the disc with radius $r$ centered at $p$.
On the other hand, if there is no center of rotation, then
\begin{align*}
 h(r_0) &= \frac{1}{2\pi}\mathrm{Length}(L) ,\\
 f(r_1))-f(r_0) &= \frac{a}{2\pi}\mathrm{Area}(B) ,\\
 \end{align*}
where $L$ is the level set $\{r=r_0\}$ and $B$ is the annulus bounded by the two level sets $\{r=r_0\}$ and $\{r=r_1\}$.
\end{lema}

\section{Closed solitons of constant curvature} \label{sec_constcurv}

Before looking for the specific nontrivial steady, shrinking and expanding solitons, in this section we rule out the constant curvature
solitons that also are
rotationally symmetric. Following \cite{ChenLuTian}, if we look for rotationally symmetric closed smooth solitons ($h(0)=h(A)=0$ and
$h'(0)=-h'(A)=1$), we can show that
there is no other
function $f$ than a
constant one. For we
multiply equation \eqref{2dsolit_eq} by $h'$ to get 
$$h'h''-ah(h')^2+\frac{hh'}{2}=0$$
and integrate on $[0,A]$ to obtain
$$\frac{(h')^2}{2}\Bigg|_0^A -a \int_0^A h(h')^2 dr + \frac{h^2}{4}\Bigg|_0^A = 0 ,$$
and since $h(0)=h(A)=0$, and $h'(0)=-h'(A)$,
$$0=-a \int_0^A h(h')^2 dr$$
which is impossible unless $a=0$. This is indeed the case when $f'=0$, there is no gradient vector field,
no Killing vector field, constant curvature and the soliton is a homothetic fixed metric. Note that if there is no vector field, there is
no need to be rotationally symmetric, thus one can have constant curvature surfaces of any genus.
Therefore we have seen the following lemma,

\begin{lema}
The only solitons over a compact smooth surface are those of constant curvature.
\end{lema}
More generally, rotationally symmetric closed solitons with two equal angles satisfy $h'(0)=-h'(A)=\frac{\alpha}{2\pi}$ and the same
argument applies. In this case, equation \eqref{2dsolit_eq} turns into
$$h''-\frac{\epsilon}{2}h=0$$
that can be explicitly solved. For $\epsilon=1$ the solution is
$$h(r)=c_1 e^{r/\sqrt{2}} + c_2 e^{-r/\sqrt{2}}$$
but the closedness condition $h(0)=h(A)=0$ implies $c_1=c_2=0$. Thus there are no expanding solitons with two equal cone points besides the
constant curvature ones.

For $\epsilon=0$, the solution is $h(r)=c_1 r+c_2$, that can't have two zeroes unless $h\equiv 0$. Finally, for $\epsilon=-1$ the solution
is $h(r)=c_1 \sin(r/\sqrt{2}) + c_2 \cos(r/\sqrt{2})$, and by the closedness $c_2=0$. This metric is locally the round sphere.
We have therefore seen the following,

\begin{lema}
The only solitons over a compact surface with two equal cone points are shrinking spherical surfaces.
\end{lema}

Up to now, we have examinated all possible cases with $a=0$ (equivalently, with $f$ constant and with constant curvature). Thus, we will
assume henceforth that $a\neq 0$ and $f$ is not constant.

\section{Steady solitons} \label{sec_steady}

\begin{figure}[ht]
 \centering
\includegraphics[width=0.6\textwidth]{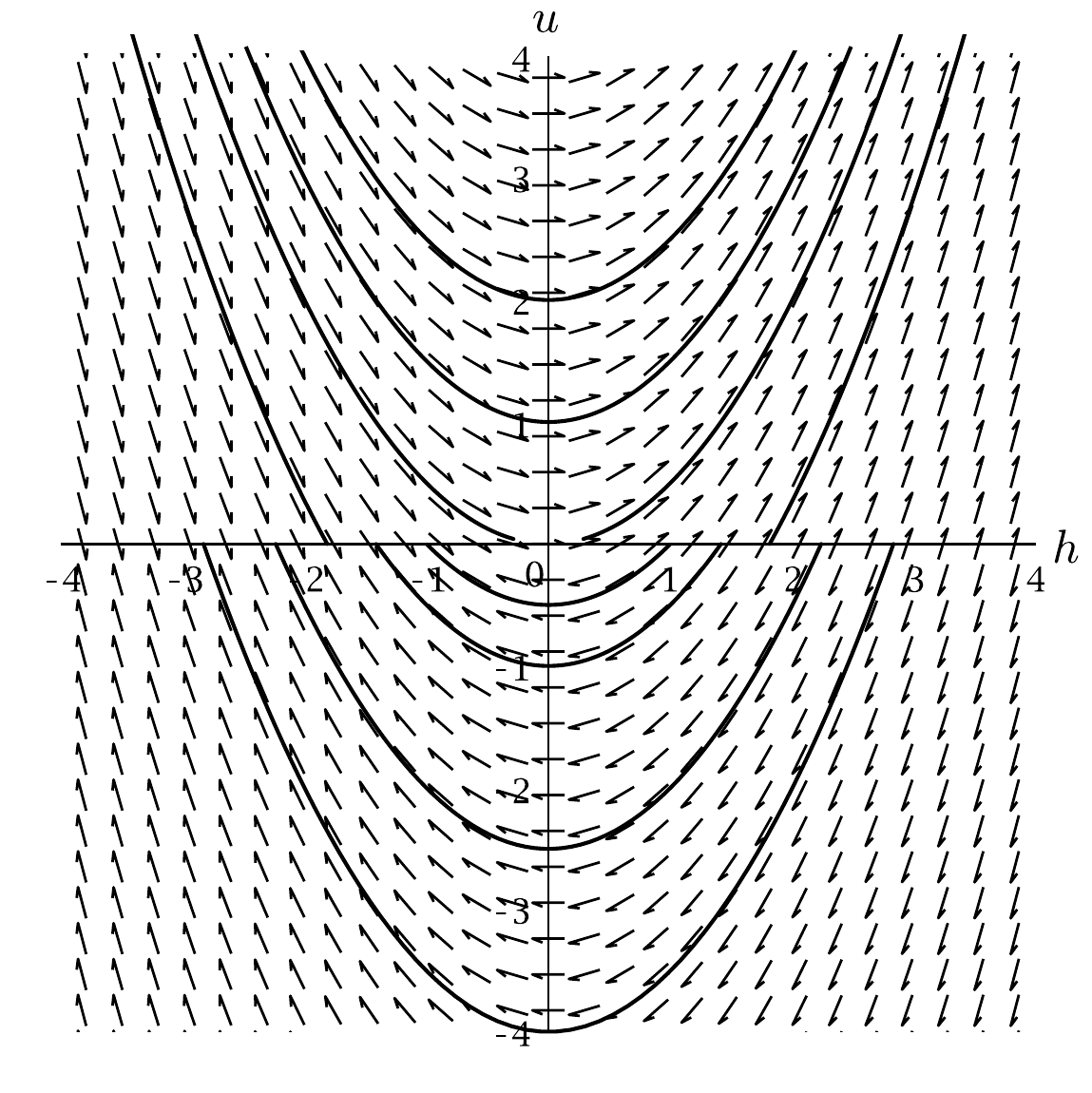}
\caption{Phase portrait of the system \eqref{2dsolit_steady_sys} with $a=1$.} \label{2dsolit_steady_pp}
\end{figure}

In this section we study the steady case ($\epsilon=0$) of rotationally symmetric solitons. The equation \eqref{2dsolit_eq} reduces to 
\begin{equation} h''-ahh'=0 \end{equation}
and the system \eqref{2dsolit_sys} to
\begin{equation}
 \left\{ \begin{array}{rcl}
   h' &=& u \\
   u' &=& auh
  \end{array} \right. \label{2dsolit_steady_sys}
\end{equation}
The phase portrait of \eqref{2dsolit_steady_sys} is shown in Figure~\ref{2dsolit_steady_pp}

This phase portrait has a line of fixed points at $\{u=0\}$, that account for the trivial steady solitons consisting on a flat cylinder of
any fixed diameter (or their universal covering, the flat plane). No other critical points are present. Only the right half-plane $\{h>0\}$
is needed, since we can take $h>0$ in the
metric definition.

Every integral curve of the system lies on a parabola. This follows from manipulating system \eqref{2dsolit_steady_sys}
$$u'=ahh'=a\left(\frac{h^2}{2}\right)'$$
and hence
$$u=a\frac{h^2}{2} +C .$$
In another terminology, the function
$$H(h,u)=a\frac{h^2}{2}-u$$
is a first integral of the system \eqref{2dsolit_steady_sys}.
Furthermore, we can finish the integration of the equation
$$h'=a\frac{h^2}{2}+C$$
by writting
$$ \frac{h'}{C+\left(\sqrt{\frac{a}{2}}h\right)^2}=1 .$$
The solution to this ODE is
\begin{equation}
 h(r)=\sqrt{\frac{2C}{a}} \tan \left( \sqrt{\frac{2}{aC}} r + D \right) \label{explodinglike}
\end{equation}
if $C>0$;
\begin{equation}
 h(r)=\sqrt{\frac{-2C}{a}} \tanh \left( \sqrt{\frac{2}{-aC}} r + D \right) \label{cigarlike}
\end{equation}
if $C<0$; and
\begin{equation}
 h(r)=\frac{1}{D-\frac{a}{2}r} \label{cigexplike}
\end{equation}
if $C=0$.

Now, let us examine each type of solution. If $C>0$, the parabola lies completely on the upper half-plane $\{u>0\}$. The equation
\eqref{explodinglike} implies that $h \rightarrow \infty$ for some finite value of $r$, and hence the metric is not complete. Furthermore,
the Gaussian curvature of the metric $K$ satisfies
$$u'=-Kh$$
and since $u$ is increasing on these solutions, the curvature is not bounded below. The case $C=1$ is sometimes called
the \emph{exploding soliton} in the literature (\cite{TRFTA1}).

If we look at $C=0$, the parabola touches the origin of coordinates, and its right hand branch defines a metric on the cylinder. The value
of $D=\frac{1}{h(0)}$ can be set so that $D=h(0)=1$ just reparameterizing $r$. With this parameterization, $r\in(-\infty, \frac{2}{a}$. For
$r\leq 0$, the function $h$ is well defined and
determines a negatively curved metric that approaches a cusp as $r\rightarrow -\infty$. However, for $r\in [0, \frac{2}{a})$ the metric is
not complete and its curvature tends to $-\infty$ as $t\rightarrow \frac{2}{a}$.

We look now at the case $C<0$, first for the solutions lying in the lower half-plane $\{u<0\}$. We can assume $h(0)=0$, $C=u(0)$, $D=0$ and
$r<0$ (this means that $-r$ is the arc-parameter of the meridians). All these arcs of parabolas join a point on the $\{h=0\}$ axis with a
point on the $\{h'=u=0\}$ axis. This means that the cylinder is pinched in one end, and approaches a constant diameter cylinder on the other
end. The metrics are complete on the cylindrical end, because from equation \eqref{cigarlike} $h\rightarrow cst$ as $r\rightarrow -\infty$.
The curvature on these metrics is bounded and positive, since $u$ and $u'<0$ are bounded. 

Some of these metrics are smooth, the particular cases of $C=u(0)=h'(0)=-1$. Note that derivating the equation $h''=ahh'$ and evaluating
at $r=0$ one sees that all even-order derivatives vanish and the surface is truly $\mathcal C^\infty$ at this point. These are the so called
\emph{cigar solitons}. There are
actually infinitely many of them, adjusting the value of $a$ and changing the diameter of the asymptotic cylinder, although all of them are
homothetic and hence it is said to exist \emph{the} cigar soliton. All the other metrics
have a cone point at $r=0$, whose angle is $-2\pi h'(0)$.

The only remaining case to inspect is the solutions with $C<0$ lying on the upper half plane $\{u>0\}$. These unbounded arcs of parabolas
rise from the axis $\{u=0\}$. We can assume (changing $D$ and reparameterizing $r$) that $r\in [0,+\infty)$. The metric is complete in
$r\rightarrow +\infty$ because of equation \eqref{cigarlike}, however, these metrics fail to be complete on $r=0$, having a metric
completion with boundary $\mathbb S^1$. The curvature on these metrics is negative and not bounded below.

We summarize the discussion in the following theorem,

\begin{tma}
The only complete steady gradient Ricci solitons on a surface with curvature bounded below are:
\begin{enumerate}
 \item Flat surfaces (possibly with cone points).
 \item The smooth cigar soliton (up to homothety).
 \item The cone-cigar solitons of angle $\alpha \in (0,+\infty)$ (up to homothety).
\end{enumerate}

\end{tma}
\emph{Remark.} There exist other steady gradient solitons, with curvature not bounded below, as described above.

Pictures of a cigar soliton and a cone-cigar soliton are shown in Figures \ref{cigar_pic} and \ref{conecigar_pic}.

\section{Shrinking solitons} \label{sec_shrink}

In this section we study the shrinking solitons ($\epsilon=-1$), besides the round sphere and the spherical footballs with two equal cone
angles found in section \ref{sec_constcurv}. When $\epsilon=-1$, the metric of $\mathcal M$ is determined by a real-valued function $h(r)$
satisfying the second order ODE
\begin{equation}
 h'' -ahh' + \frac{h}{2}=0 , \label{2dsolit_shrink_eq}
\end{equation}
or equivalently the system 
\begin{equation}
 \left\{ \begin{array}{rcl}
   h' &=& u \\
   u' &=& (au-\frac{1}{2})h .
  \end{array} \right. \label{2dsolit_shrink_sys} 
\end{equation}

\begin{figure}[ht]
 \centering
\includegraphics[width=0.5\textwidth]{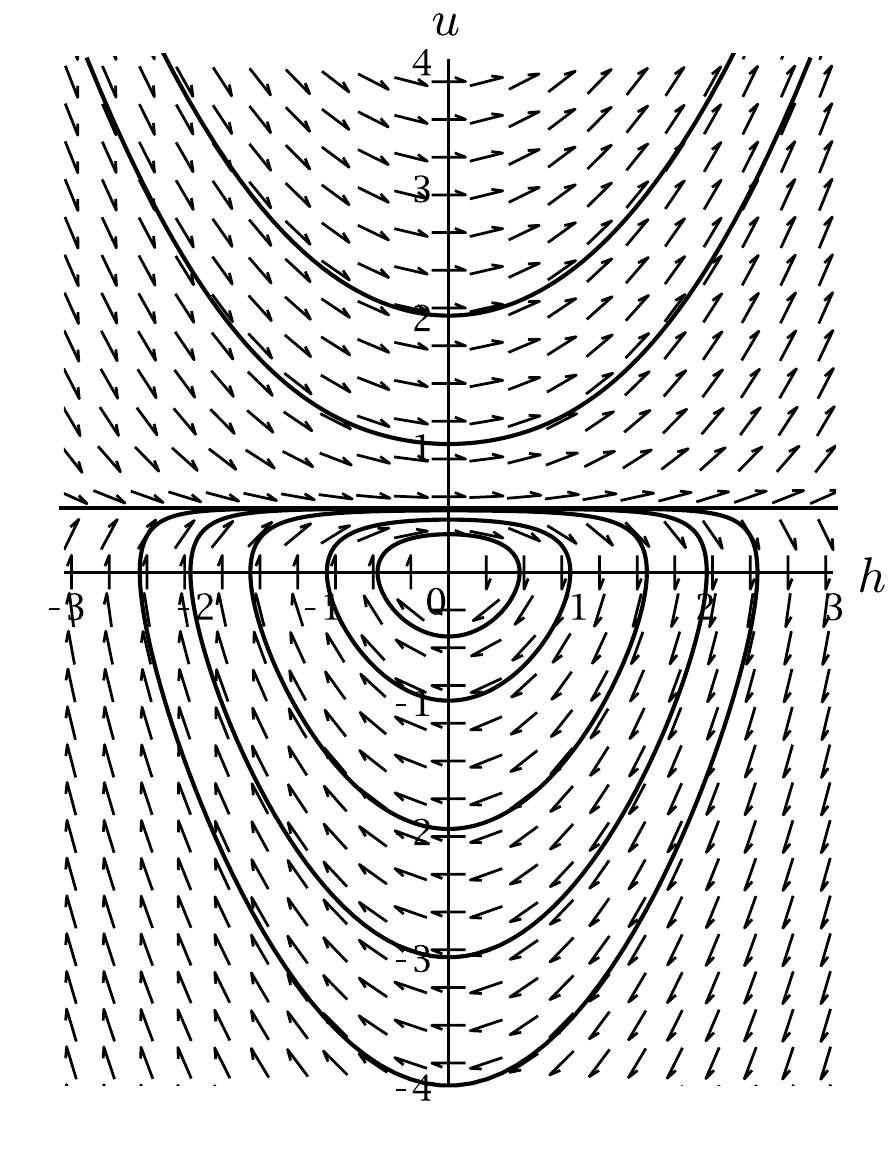}
\caption{Phase portrait of the system \eqref{2dsolit_shrink_sys} with $a=1$.} \label{2dsolit_shrink_pp}
\end{figure}

The phase portrait of this ODE system with $a=1$ is shown in Figure~\ref{2dsolit_shrink_pp}.
This phase portrait has a critical point at $(h,u)=(0,0)$ of type center, and a horizontal isocline (points such that $u'=0$) at the line
$u=\frac{1}{2a}$.
Each curve on this $hu$-plane corresponds to a solution $h$, and the intersection with the vertical axis $\{h=0\}$ are at $u(0)$ and $u(A)$,
which stand for the cone angles. Indeed, only half of each curve is enough to define the soliton, the one lying in the $\{ h>0 \}$
half-plane, since we can choose the sign of $h$ because only $h^2$ is used to define the metric.

All curves in the phase portrait represent rotationally symmetric soliton metrics over, a priori, a topological
cylinder. Closed curves (that intersect twice the axis $\{h=0\}$) are actually metrics over a doubly pinched cylinder, thus a
topological sphere with two cone points, giving the so called \emph{football solitons}. Open curves only intersect once the $\{h=0\}$ axis,
and hence are metrics over a topological plane. If the
intersection of any curve with the $\{h=0\}$ axis occurs at $u=\pm 1$, then the metric extends smoothly to this point (truly $\mathcal
C^\infty$ since derivating \eqref{2dsolit_shrink_eq} all even-order derivatives vanish at this point).  For instance, in
Figure~\ref{2dsolit_shrink_pp} there is only one curve associated to a \emph{teardrop soliton}, namely the one intersecting the vertical
axis at
some value $u(0)\in(0,\frac{1}{2})$ and at $u(A)=-1$. There is also a smooth soliton metric on $\mathbb
R^2$, namely the one associated with the
curve passing through $(h,u)=(0,1)$, and all other curves represent solitons over cone surfaces. The separatrix line,
$u=\frac{1}{2a}$,
represents the solution $h(r)=\frac{r}{2a}+c_0$, which stands for the metric $dr^2 + \frac{1}{4a^2} r^2 d\theta^2$. This is a flat
metric on the cone of angle $\frac{\pi}{a}$, a cone version of the shrinking Gaussian soliton, and we call it a \emph{shrinking Gaussian
cone soliton} (the smooth \emph{shrinking Gaussian soliton} is the case $a=\frac{1}{2}$).

Let us focus on the compact shrinking solitons.
\begin{lema} \label{lemafoot}
For every pair of values $0<\alpha_1 < \alpha_2 < \infty$, there exist a unique value $a>0$ such that the equation \eqref{2dsolit_eq} has
one solution satisfying the boundary conditions $h'(0)=\frac{\alpha_1}{2\pi}$ and $h'(A)=-\frac{\alpha_2}{2\pi}$.
\end{lema}

Equivalently, the lemma asserts that there exists a value $a$ such that the phase portrait of the system \eqref{2dsolit_sys} has one
solution curve that intersect the vertical axis $\{h=0\}$ at $u(0)=\frac{\alpha_1}{2\pi}$ and $u(A)=-\frac{\alpha_2}{2\pi}$.

\begin{proof} 
We can normalize the system by 
$$\left\{ \begin{array}{rcl} v&=&ah \\  w&=&au \end{array} \right.$$
so that on this coordinates the system becomes 
\begin{equation}
 \left\{ \begin{array}{rcl}
   v' &=& w \\
   w' &=& (w-\frac{1}{2})v
  \end{array} \right. \label{2dsolit_sys_norm}
\end{equation}
This would be the same system as \eqref{2dsolit_shrink_sys} with $a=1$, which is indeed shown on Figure~\ref{2dsolit_shrink_pp}.

The system \eqref{2dsolit_sys_norm} has the following first integral,
$$H(v,w) = v^2-2w-\ln|2w-1|$$
that is, the solution curves of the system are the level sets of $H$. Indeed, derivating $H(v(r),w(r))$ with respect to $r$,
\begin{align*}
 \frac{\partial}{\partial r} H(v,w) & = 2 v v' -2w' -\frac{2w'}{2w-1} \\
 &= 2vw -2 \left( w-\frac{1}{2} \right) v \left( 1+\frac{1}{2w-1} \right) =0 .
\end{align*}

The cone angle conditions are $\alpha_1=\frac{2\pi w(0)}{a}$, $\alpha_2 =- \frac{2\pi w(A)}{a}$, while $v(0)=v(A)=0$. Thus, the function
$w$ evaluated at $0$ and $A$ satisfies
$$H(0,w)=2w+\ln|2w-1|=C$$
for some $C\in\mathbb R$. This is equivalent, via $2w-1=-y$ and $e^{C}=k$, to the equation
\begin{equation}
 |y|=ke^{y-1}  \label{eqham}
\end{equation}
(cf. \cite{Hamilton_surfaces}). Although not expressable in terms of elementary functions, this equation has three solutions for $y$, one
for negative $y$ and
two for positive $y$ (see Figure~\ref{ploteqham}). The two positive solutions of \eqref{eqham} are the intersection of the exponential
function $e^{y-1}$ with the line
$\frac{1}{k}y$ with slope $\frac{1}{k}$. These two positive solutions are associated to a compact connected component of $H(v,w)=C$, whereas
the negative solution is associated to a noncompact component of $H$ that represent noncompact soliton surfaces.
\begin{figure}[ht]
 \centering
\includegraphics[width=0.4\textwidth]{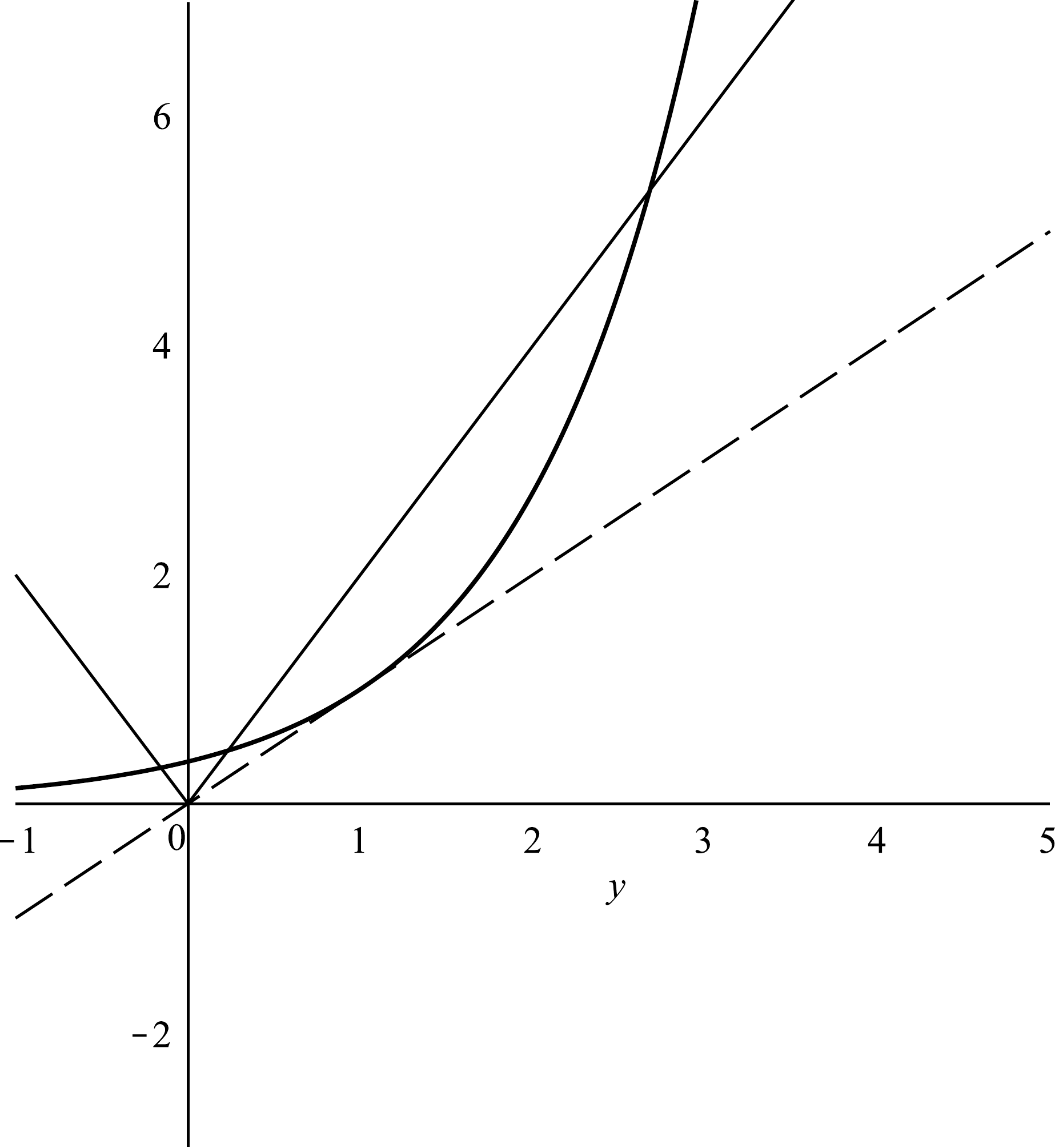}
\caption{The graphs of the exponential $e^{y-1}$ and $\frac{1}{k}|y|$. } \label{ploteqham}
\end{figure}
The two positive solutions of \eqref{eqham} exist only when $k\in(0,1)$ and actually these two solutions are equal when $k=1$ and
the line is
tangent to the exponential function at $y=1$. These two solutions $y_1$, $y_2$ of \eqref{eqham} are therefore located on $(0,1)$ and
$(1,+\infty)$ respectively, and can be expressed as 
$$y_1=1-p \quad , \quad y_2=1+q$$
with $p,q\geq 0$. 
The two cone angles, having assumed $\alpha_1<\alpha_2$, are then expressed as
$$\alpha_1 = 2\pi h'(0) = 2\pi u(0) = \frac{2\pi w(0)}{a} = \frac{2\pi}{a} \frac{1-y_1}{2} = \frac{\pi}{a} p $$
$$\alpha_2 = -2\pi h'(A) = -2\pi u(A) = -\frac{2\pi w(A)}{a} = -\frac{2\pi}{a} \frac{1-y_2}{2} = \frac{\pi}{a} q $$
and their quotient is
$$\frac{\alpha_1}{\alpha_2} = \frac{p}{q} .$$
Let $\Psi:(0,1) \rightarrow \mathbb R$ be the mapping $$k\mapsto \Psi(k) = \frac{p}{q} .$$ 
The function $\Psi$ is injective and the quotient $\Psi(k)$ ranges from $0$ to $1$ when varying $k\in(0,1)$. This is proven in
\cite[Lem 10.7]{Hamilton_surfaces}, we can visualize its
graph in Figure~\ref{football_psi}.
\begin{figure}[ht]
 \centering
\includegraphics[width=0.5\textwidth]{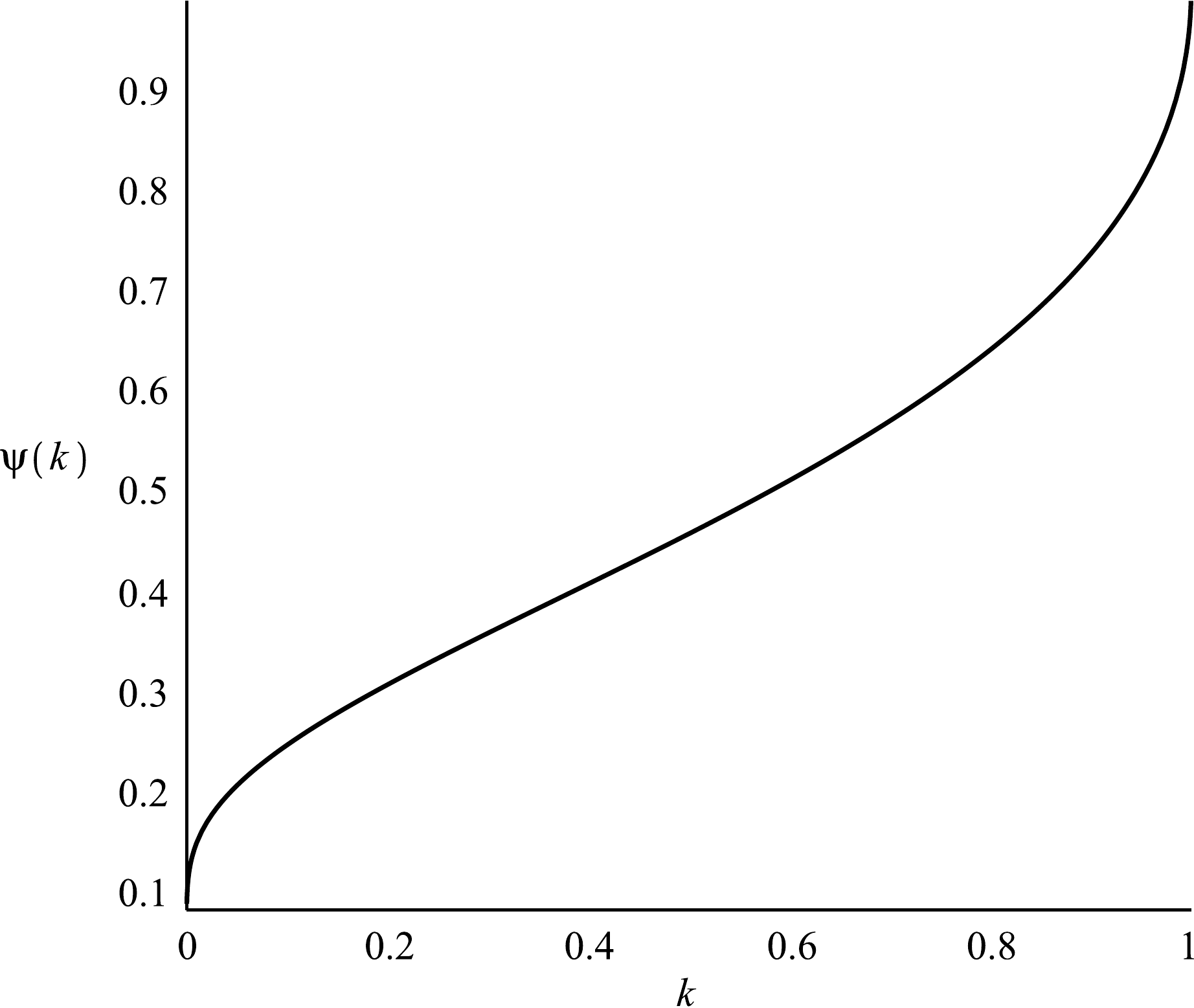}
\caption{The function $\Psi(k)$.}\label{football_psi}
\end{figure}
 Therefore, for any pair of chosen angles 
$\alpha_1 < \alpha_2$ there exist $k=\Psi^{-1}(\frac{\alpha_1}{\alpha_2})$, such that the equation \eqref{eqham} has two positive solutions
$y_1$, $y_2$. This yields two values $p=1-y_1$, $q=y_2-1$, and finally we recover 
$$a=\frac{\alpha_1}{\pi p} = \frac{\alpha_2}{\pi q} .$$
This value makes the system \eqref{2dsolit_sys} and the equation \eqref{2dsolit_eq} to have the required solutions.
\end{proof}

We summarize our discussion in the following theorem.

\begin{tma}
 The only complete shrinking gradient Ricci solitons on a surface $\mathcal M$ with curvature bounded below are:
\begin{itemize}
 \item Spherical surfaces (including the round sphere and the round football solitons of constant curvature and two equal cone angles). 
 \item The football and teardrop solitons, if $\mathcal M$ is compact and has one or two different cone points. These cone angles can be any
real positive value.
 \item The smooth shrinking flat Gaussian soliton.
 \item The shrinking flat Gaussian cones.
\end{itemize}
\end{tma}

\section{Expanding solitons} \label{sec_expand}

\begin{figure}[ht]
 \centering
\includegraphics[width=0.6\textwidth]{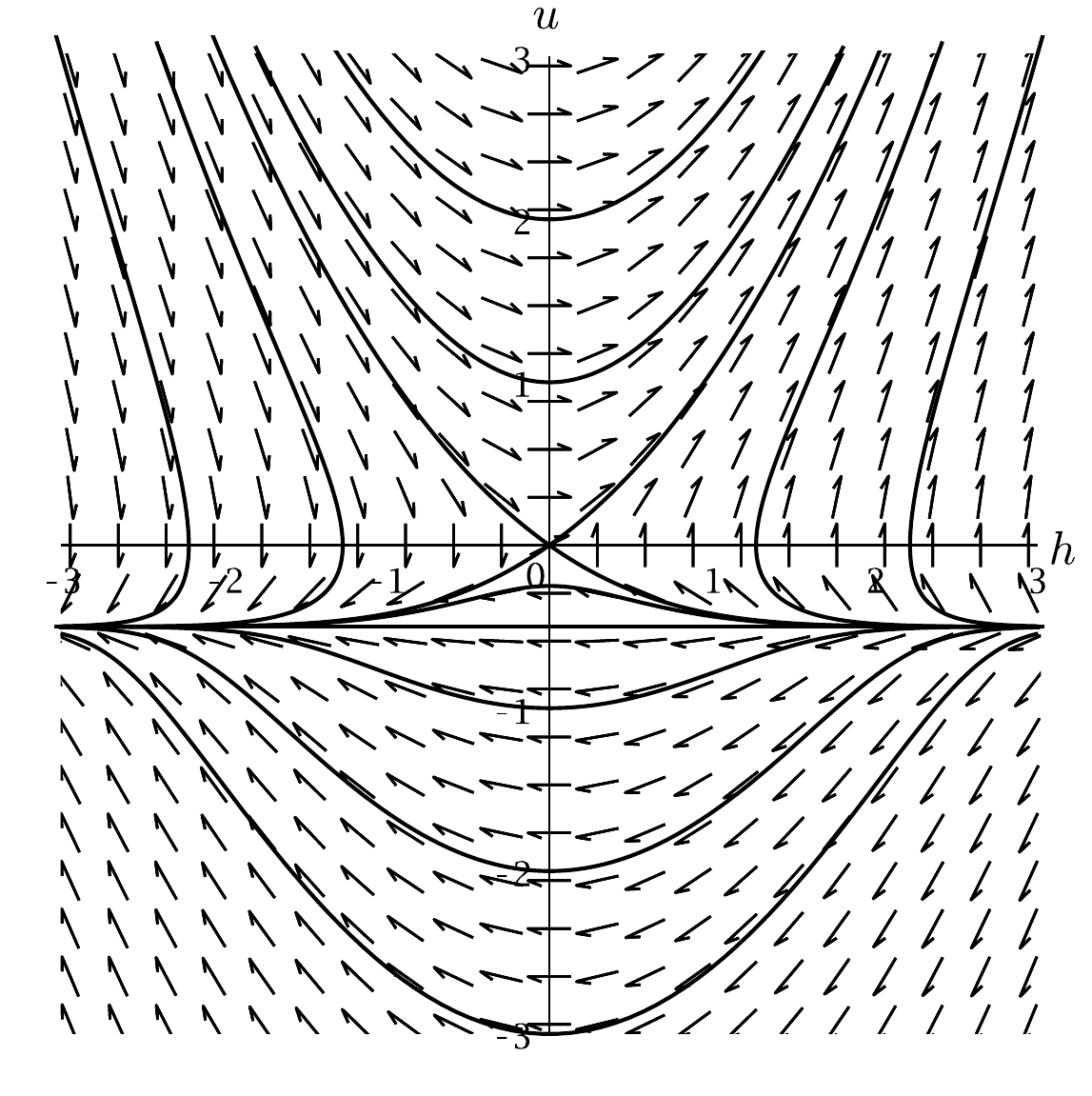}
\caption{Phase portrait of the system \eqref{2dsolit_expand_sys} with $a=1$.} \label{2dsolit_expand_pp}
\end{figure}

We end our classification with the expanding solitons ($\epsilon=1$). The equation \eqref{2dsolit_eq} and the system \eqref{2dsolit_sys}
are in this case
\begin{equation} h''-ahh'+\frac{h}{2}=0 \label{2dsolit_expand_eq}\end{equation}
and
\begin{equation}
 \left\{ \begin{array}{rcl}
   h' &=& u \\
   u' &=& \left( au+\frac{1}{2} \right) h .
  \end{array} \right. \label{2dsolit_expand_sys}
\end{equation}
The phase portrait of \eqref{2dsolit_expand_sys} is shown in Figure~\ref{2dsolit_expand_pp}.
We can rescale the system \eqref{2dsolit_expand_sys} with the change
$$\left\{ \begin{array}{rcl} v&=&ah \\  w&=&au \end{array} \right.$$
so that on this coordinates the system becomes 
\begin{equation}
 \left\{ \begin{array}{rcl}
   v' &=& w \\
   w' &=& (w+\frac{1}{2})v
  \end{array} \right. \label{2dsolit_expand_sys_norm}
\end{equation}
which is exactly the system \eqref{2dsolit_expand_sys} with $a=1$. The phase portrait of this system is indeed
shown on Figure~\ref{2dsolit_expand_pp}. We will
study the trajectories of the normalized system and next we will discuss the geometrical interpretation of each trajectory.

The system \eqref{2dsolit_expand_sys_norm} has a critical point ($v'=w'=0$) at
$(0,0)$. It has an horizontal isocline ($w'=0$) at the line $L=\{w=-\frac{1}{2}\}$, that is also an orbit solution, and hence no other
trajectory can cross it. The vertical axis $\{v=0\}$ is also an horizontal isocline. The horizontal axis $\{w=0\}$ is, on the other hand, a
vertical isocline ($v'=0$).

The linearization of the system \eqref{2dsolit_expand_sys_norm} at the critical point $(0,0)$ is
$$\left( \begin{array}{c} v' \\ w' \end{array} \right) =
\left( \begin{array}{cc} 0 & 1 \\ w+\frac{1}{2} & v \end{array} \right) 
\left( \begin{array}{c} v \\ w \end{array} \right). $$
The matrix of the linearized system at the critical point is $\left( \begin{array}{cc} 0 & 1 \\ \frac{1}{2} & 0 \end{array} \right)$, that
has determinant $-\frac{1}{2}<0$ and hence the critical point is a saddle point. The eigenvalues of this matrix are $\frac{1}{\sqrt{2}}$
and $-\frac{1}{\sqrt{2}}$, with eigenvectors respectively
$$\left( \begin{array}{c} \sqrt{2} \\ 1 \end{array} \right) \quad \mbox{and} \quad \left( \begin{array}{c} -\sqrt{2} \\ 1 \end{array}
\right).$$
These eigenvectors determine the two principal directions of the saddle point, from which four separatrix curves are emanating.

The system \eqref{2dsolit_expand_sys_norm} has the following first integral,
$$H(v,w) = v^2-2w+\ln|2w+1|$$
that is, the solution curves of the system are the level sets of $H$. Indeed, derivating $H(v(r),w(r))$ with respect to $r$,
\begin{align*}
 \frac{\partial}{\partial r} H(v,w) & = 2 v v' -2w' +\frac{2w'}{2w+1} \\
 &= 2vw -2 \left( w+\frac{1}{2} \right) v \left( 1-\frac{1}{2w+1} \right) =0 .
\end{align*}

From the system, and more apparently from the first integral, it is clear that the phase portrait is symmetric with respect to the axis
$\{v=0\}$. We will only study then the trajectories on the right-hand half-plane $\{v>0\}$. Actually this restriction agrees
with the geometric assumption of $h>0$.

Let us consider a trajectory passing through a point in the quadrant $\{v>0 , w>0\}$. Then $v'>0$ and $w'>0$ and hence the curve moves
upwards and rightwards. More carefully, as soon as $v>\delta_1>0$ and $w>\delta_2>0$, both derivatives are bounded below away from zero,
$v'>\tilde\delta_1>0$ and $w'>\tilde\delta_2>0$, and therefore $v$ and $w$ tend to $+\infty$. We can further evaluate the asymptotic
behaviour from the first integral,
$$\frac{v^2}{2w}-1 + \frac{\ln|2w+1|}{2w} = \frac{C}{2w}$$

We take the limit as $r\rightarrow +\infty$ and since $w\rightarrow +\infty$, we get that 
$$\lim_{r\rightarrow +\infty} \frac{v(r)^2}{2w(r)}= 1 ,$$ 
so the orbit approaches
a parabola (the same parabolas of the steady case and asymptotic on the shrinking case).

We now inspect the separatrix $S$ emanating (actually sinking) from the critical point at the direction $(\sqrt{2},-1)$. The associated
eigenvalue is $-\sqrt{2}$ and hence the trajectory is approaching the saddle point (hence the sinking). The curve $S$ lies in the
$w'>0$ region, and cannot cross the horizontal isocline $L$. Therefore, the
separatrix when seen backwards in $r$ must be decreasing and bounded, and hence must approach a horizontal asymptote. This asymptote must be
$L$, since if the trajectory were lying in the region $w'>\delta>0$ for infinite time, it would come from $w=-\infty$,
which is absurd since it cannot cross the isocline $L$. Therefore, over this separatrix $S$, $v\rightarrow +\infty$ and $w\rightarrow
-\frac{1}{2}$ as $r\rightarrow -\infty$.

Any trajectory lying above $S$ will eventually enter in the upper right quadrant, and hence approach asymptotically the previously
mentioned parabolas. This is clear since $w'$ is positive and bounded below away from zero.

Let us study the trajectories below $S$. All these curves intersect the axis $\{v=0\}$, and we can consider the origin of the $r$ coordinate
as such that the intersection point with the axis occurs at $r=0$. Then, the region of the curves parameterized by $r<0$ lies in the $v>0$
half-plane. Since the curves are below $S$, they lie in the
lower right quadrant and hence $v'<0$ and $v\rightarrow +\infty$ as $r \rightarrow -\infty$. If the curve lies over $L$,
then $w'>0$, and if it lies below $L$, then $w'<0$. This means that the isocline is repulsive forward in $r$ and attractive backwards in
$r$. Therefore any curve lying below $S$ will have an asymptote as $r\rightarrow -\infty$ and as before this must be the isocline $L$, that
is, $v\rightarrow +\infty$ and $w\rightarrow -\frac{1}{2}$ as $r\rightarrow -\infty$.

Let us comment about the domain of $r$. We have stated that the trajectories below $S$ are parameterized for $r\in (-\infty,0]$, although a
priori it could be $r\in (-M,0]$ for some maximal $M$ (and hence $v\rightarrow +\infty$ as $r\rightarrow -M$, and these would represent
noncomplete metrics). This is not the case, since the trajectories are approaching $v'=-\frac{1}{2}$, and hence $v'$ is bounded ($|v'|< 1$
for $r$ less than some $r_0<0$), so $v$ cannot grow to $+\infty$ for finite $r$-time.

Finally, let us study the separatrix $S$ itself. This curve is parameterized by $r\in\mathbb R$, and $(v,w)\rightarrow (0,0)$ as
$r\rightarrow +\infty$ and $(v,w)\rightarrow (+\infty,-\frac{1}{2})$ as $r\rightarrow -\infty$ (this follows from the Grobman-Hartman
theorem in the end near the saddle point, and from the asymptotic $L$ on the other end). We can give a more detailed
description of the asymptotics. As $r \rightarrow -\infty$, we know that $w\rightarrow -\frac{1}{2}$, this is
$$\lim_{r\rightarrow -\infty}\frac{v'}{-\frac{1}{2}} = 1 .$$
Then, applying the l'Hôpital rule,
$$\lim_{r\rightarrow -\infty}\frac{v}{-\frac{1}{2}r} = 1 ,$$
or $v(r)\sim -\frac{1}{2}r$ as $r\rightarrow -\infty$.
This is valid for all the trajectories asymptotic to the horizontal isocline. 
Similarly, as $r\rightarrow +\infty$, we know that
$v,w\rightarrow 0$, but furthermore we know that their quotient tends to the slope of the eigenvector determining the separatrix, i.e.
$$\lim_{r\rightarrow +\infty}\frac{v}{w} = \lim_{r\rightarrow +\infty} \frac{v}{v'} = \frac{-1}{\sqrt{2}} ,$$
which is to say
$$\lim_{r\rightarrow +\infty}(\ln{v})'= \lim_{r\rightarrow +\infty} \frac{v'}{v} = -\sqrt{2} .$$
Then, by l'Hôpital rule,
$$ \lim_{r\rightarrow +\infty}\frac{(\ln{v})'}{-\sqrt{2}} = \lim_{r\rightarrow +\infty} \frac{\ln{v}}{-\sqrt{2}r}=1 .$$
this is, $v(r)\sim e^{-\sqrt{2}r}$ as $r\rightarrow +\infty$.

At this point we have established all the important features of the phase portrait in Figure~\ref{2dsolit_expand_pp}. With the unnormalized
system, in coordinates $(h,u)$, the phase portrait is just a scaling of the one in Figure~\ref{2dsolit_expand_pp} by the factor $a$, and
thus the horizontal isocline $L$ is at $\{u=-\frac{1}{2a}\}$. We now give the geometric interpretation of each trajectory.

All the trajectories above $S$ have solutions with unbounded positive $w$. This means that the Gaussian curvature of the associated metric
$K=-(au+\frac{1}{2}) = -(w+\frac{1}{2})$
is not bounded below, and we will discard them. 

All the trajectories below the separatrix $S$ have bounded $w$ and therefore bounded curvature on the associated metric. More specifically,
the curves above the isocline $L$ will give metrics with negative curvature, and curves below $L$ will give metrics with positive
curvature. These curves will intersect the $\{h=0\}$ axis at $b<0$, and the associated metric will have a cone point of angle
$$\beta=-2\pi b$$ 
at the point of coordinate $r=0$. On the other end, the function $h(r)$ is asymptotic to $-\frac{1}{2a}r$ (recall
that the parameter is $r\in (-\infty,0 ]$) and the metric will be asymptotic to the wide part of a flat cone of angle 
$$\alpha = \frac{\pi}{a} .$$ 
We call these solitons the \emph{$\alpha \beta$-cone solitons}. These solitons have positive curvature if $\alpha<\beta$
($b<\frac{-1}{2a}$) and negative if $\alpha>\beta$ ($b>\frac{-1}{2a}$). In the case $\alpha=\beta$ ($b=\frac{-1}{2a}$) we are in the case
of the isocline $L$. This line $h' = u = -\frac{1}{2a}$, has as solution the parameterization
$$h(r)=-\frac{1}{2a}r+C$$
which represents a \emph{flat expanding Gaussian cone soliton}, with cone angle $\frac{\pi}{a}$. The special case $a=\frac{1}{2}$
yields a smooth metric at $r=0$, thus we have a flat metric on the plane known as the \emph{flat expanding Gaussian soliton}.

Other remarkable cases are those with $\beta=2\pi$ ($b=-1$), because the cone point at the apex is now blunted and the surface is smooth
(we can check from equation \eqref{2dsolit_expand_eq} that all even-order derivatives vanish at $r=0$),
we call them the \emph{blunt $\alpha$-cone solitons}. The angle $\alpha$ may be less or greater than $2\pi$ and the curvature is
positive or negative respectively. However, only the first case can be embedded symmetrically in $\mathbb R^3$). The existence of this
family was described by a different method by H.-D. Cao in \cite{Cao_KRsolitons2} in the context of Kähler-Ricci solitons.

We interpret now the separatrix $S$. This is the limiting case as the angle $\beta$ tends to zero. In this case the parameter $r$ is
not on $(-\infty,0]$ but on the whole $\mathbb R$ and thus $h(r)$ defines a smooth complete metric on the cylinder. As $r\rightarrow
+\infty$, the function $h(r)$ is asymptotic to $\frac{1}{a} e^{-\sqrt{2}r}$, that defines a hyperbolic metric of constant curvature $-2$.
This
hyperbolic metric on a cylinder is called a \emph{hyperbolic cusp}. The separatrix $L$ represents a soliton metric that approaches the
thin part of a hyperbolic cusp in one end, and the wide part of a flat cone on the other. There is still freedom to set the angle $\alpha$,
and we call these the \emph{cusped $\alpha$-cone solitons}.

Finally, there is still one more family of two-dimensional gradient solitons, namely the universal cover of the cusped $\alpha$-cones.
These solitons are metrics on $\mathbb R^2$ locally isometric to the cusped cones. These solitons are not rotationally symmetric, but
translationally symmetric, i.e. there is not a $\mathbb S^1$ group but a $\mathbb R$ group acting by isometries. The plane $\mathbb R^2$
with any of this metrics has a fixed direction (given by $\grad f$) such that a straight line following this direction (that is also a
geodesic of the soliton metric) transits gradually from a region of hyperbolic curvature on one end to a region of flat curvature on the
other. Any translation on the direction perpendicular to $\grad f$ (this is, in the direction of $J(\grad(f))$) is an isometry on these
metrics. We call these \emph{flat-hyperbolic soliton planes}.

We summarize our discussion in the following theorem.

\begin{tma}
 The only complete expanding gradient Ricci solitons on a surface $\mathcal M$ with curvature bounded below are:
\begin{itemize}
 \item Hyperbolic surfaces (possibly with cone angles).
 \item The $\alpha\beta$-cone solitons, for every pair of cone angles $\alpha,\beta>0$.
 \item The smooth complete blunt $\alpha$-cones, that are $\alpha\beta$-cones with $\beta=2\pi$.
 \item The smooth expanding flat Gaussian soliton.
 \item The expanding flat Gaussian cones.
 \item The smooth cusped $\alpha$-cone in the cylinder.
 \item The flat-hyperbolic solitons on the plane, that are universal coverings of the cusped cones.
\end{itemize}
\end{tma}

\section{Gallery of embedded solitons} \label{sec_gallery}

Just for aesthetics, we can embed some of the solitons we described into $\mathbb R^3$ and visualize them numerically as surfaces with the
inherited metric from the ambient Euclidean space. If we want to keep the rotational symmetry apparent, however, we cannot embed into
$\mathbb R^3$ a cone point of angle greater than $2\pi$, and we can't embed a rotational surface whose parallels have length $L=2 \pi\ h(R)$
if $R<\frac{L}{2\pi}$.

In order to do this, we use the metric in polar coordinates $(r,\theta) \in [0,A] \times [0,2\pi]$,
$$dr^2 + h(r)^2 \ d\theta^2 .$$
We recall that $r$ is the arc parameter of the $\{\theta=cst\}$ curves (meridians), and
that the $\{r=cst\}$ curves (parallels) are circles of radius $h(r)$ parameterized by $\theta\in[0,2\pi]$. Therefore, we can use the
$h,\theta$ as polar coordinates on the
plane, and find an appropriate third coordinate $z$ (height). When we put the stacked parallels of radius $h(r)$ at height
$z(r)$, we obtain a rotational surface whose meridians have length parameter $r$. Thus,
$$dr^2 = dh^2 + dz^2$$
or equivalently
$$\frac{dz^2}{dr^2}=1-\frac{dh^2}{dr^2}$$
which defines $z=z(r)$ as satisfying
$$(z')^2 = 1-u^2$$
with the convention that $h'=u$. Hence, to obtain an embedded surface satisfying the soliton system \eqref{2dsolit_sys} it is sufficient to
integrate the first order vector ODE
$$\left\{ \begin{array}{rcl}
   h' &=& u \\
   u' &=& (au+\frac{\epsilon}{2})h \\
   z' &=& \sqrt{1-u^2}
  \end{array} \right. $$
with initial conditions $h(0)=0$, $u(0)=b$, $z(0)=0$. Once obtained a numerical solution for $h(r)$, $u(r)$ and $z(r)$, we
can fix a value $A>0$ and then plot the set of points
$$\{(h(r)\cos\theta,h(r)\sin\theta,z(r))\in \mathbb R^3 \ \big| \ r\in[0,A] \ ,\ \theta\in [0,2\pi]\}.$$

Next we show some embedded solitons. These were obtained with the following Maple code:

\begin{verbatim}
> epsilon:=1; a:=1; b:=-1; A:=10;
> sys:= diff(h(r),r)=u(r), diff(u(r),r)=(a*u(r)+epsilon/2)*h(r), 
  diff(z(r),r)=sqrt(1-u(r)^2):
  fns:= dsolve( {sys , h(0)=0, u(0)=b, z(0)=0}, numeric, 
  output=listprocedure):
  hh:=rhs(fns[2]): uu:=rhs(fns[3]): zz:=rhs(fns[4]):
> plot3d([hh(r)*cos(theta), hh(r)*sin(theta),zz(r)],
  r=0..A, theta=0..2*Pi, scaling=constrained,grid=[40,40]);
\end{verbatim}

\newlength{\tam}
\setlength{\tam}{13em}

\begin{figure}[p]
\begin{minipage}[t]{0.45\linewidth} \centering

\includegraphics[height=\tam]{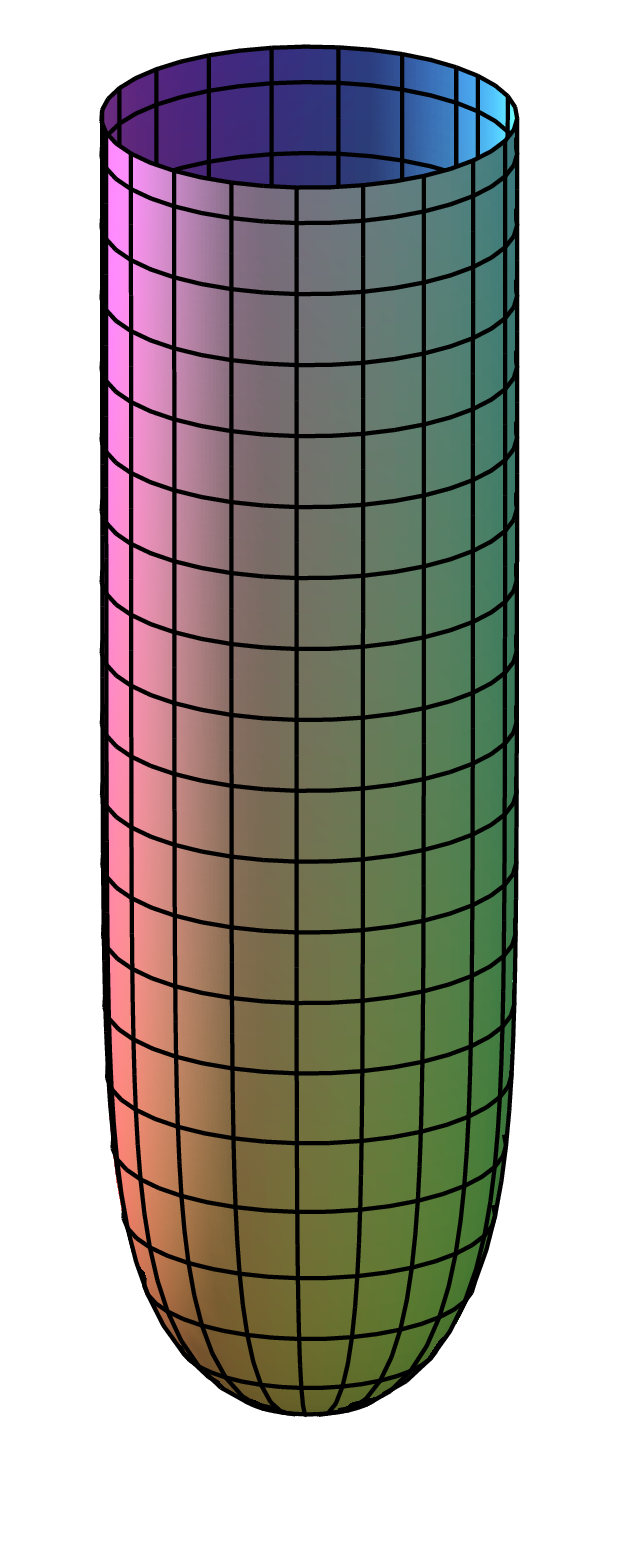}
\caption{ Hamilton's cigar soliton. \newline 
($\epsilon=0$, $a=1$, $b=-1$)} \label{cigar_pic}

\end{minipage}
\hspace{1cm}
\begin{minipage}[t]{0.45\linewidth} \centering

\includegraphics[height=\tam]{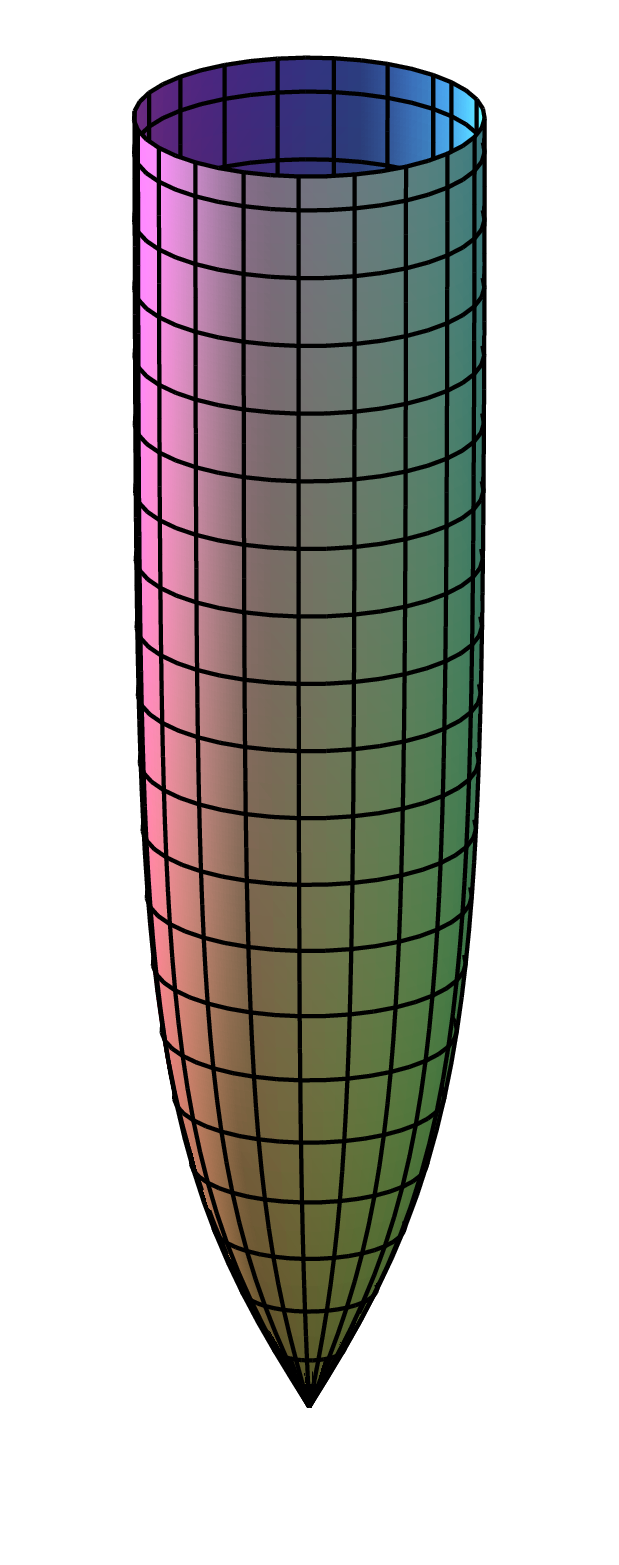}
\caption{ A cone-cigar soliton with cone angle $180^\circ$. \newline 
($\epsilon=0$, $a=1$, $b=-0.5$)} \label{conecigar_pic}
\end{minipage}
\end{figure}

\begin{figure}[p]
\begin{minipage}[t]{0.45\linewidth} \centering

\includegraphics[height=\tam]{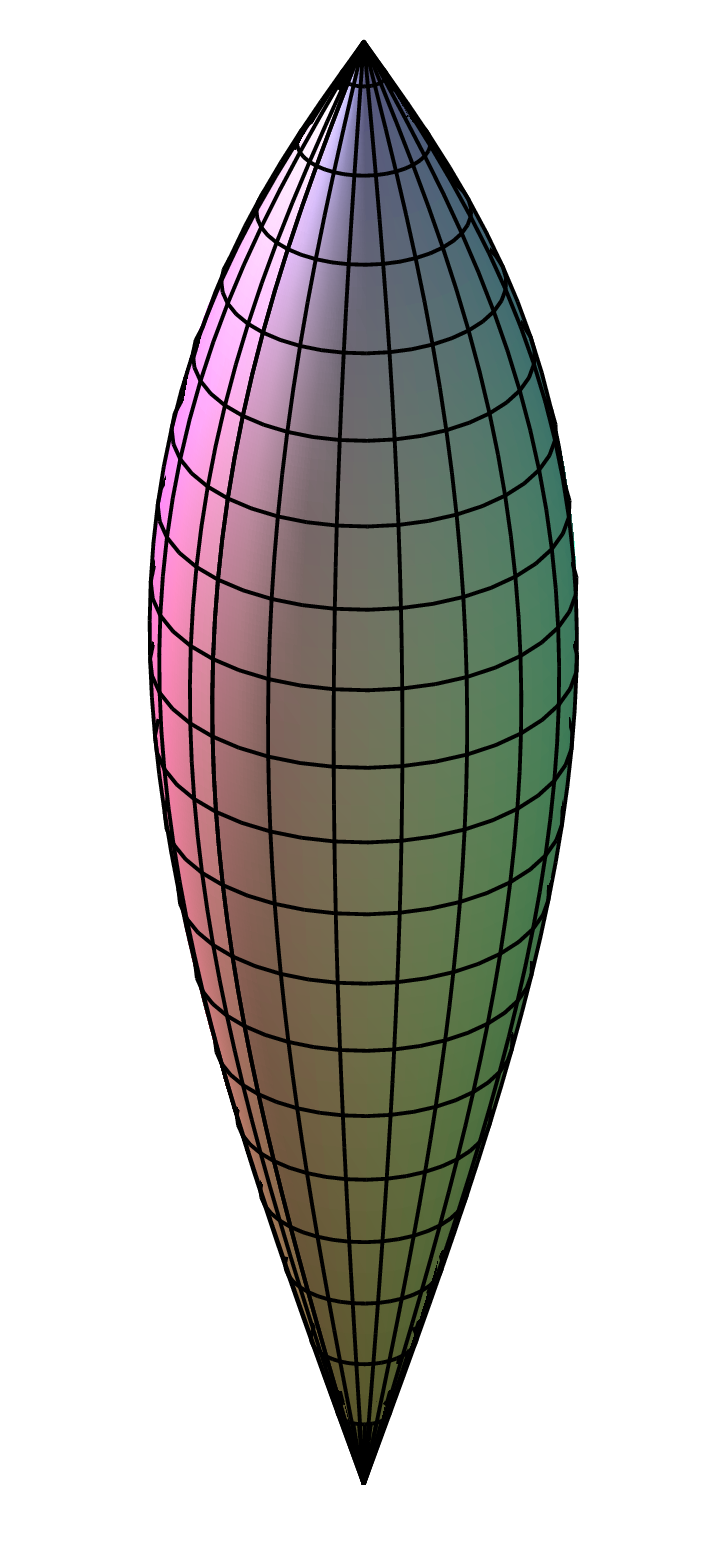}
\caption{ A football soliton with cone angles $108^\circ$ and $183.38^\circ$. \newline 
($\epsilon=-1$, $a=1$, $b=0.3$, $A=4.56$)} \label{football_pic}

\end{minipage}
\hspace{1cm}
\begin{minipage}[t]{0.45\linewidth} \centering

\includegraphics[height=\tam]{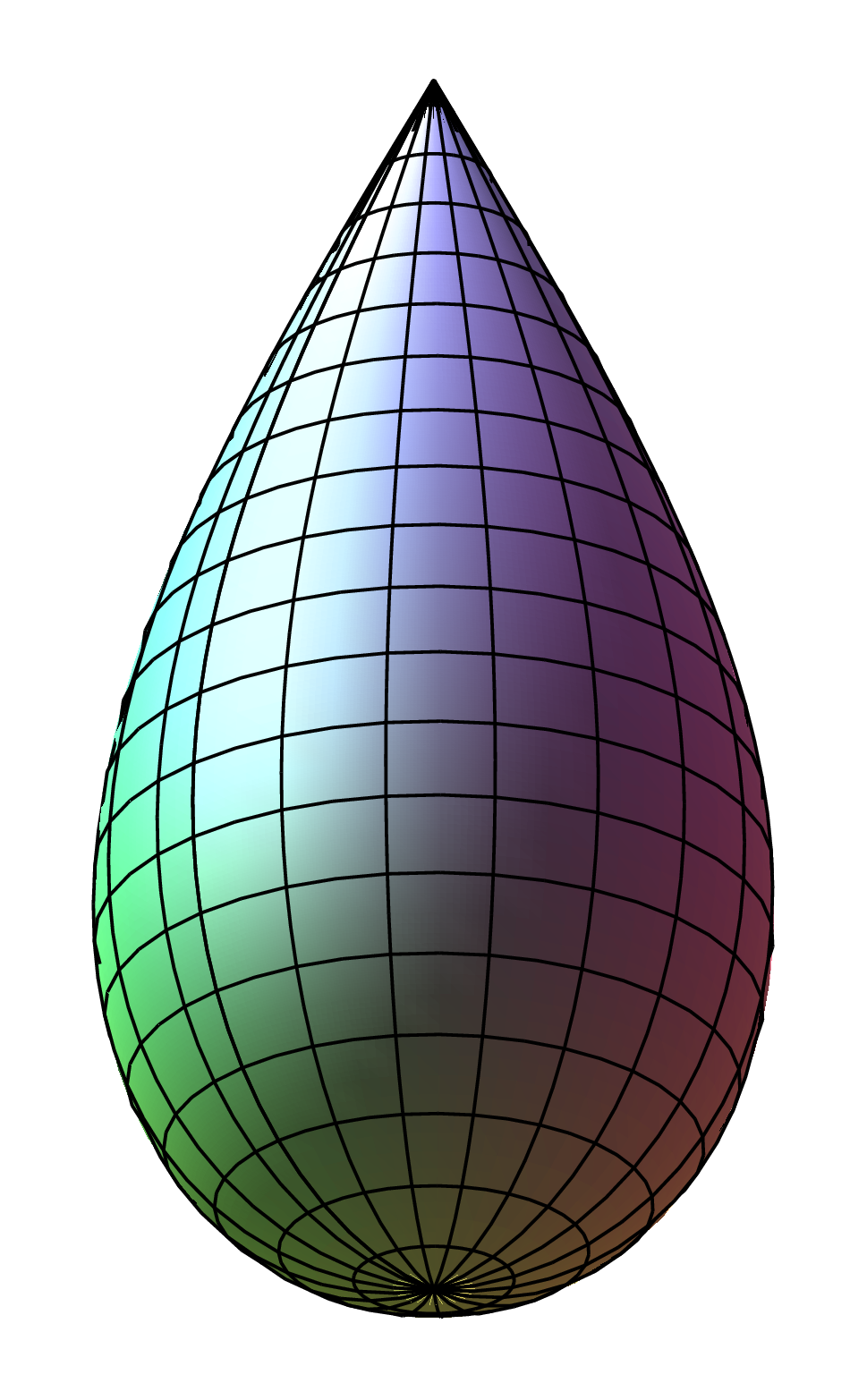}
\caption{ A teardrop soliton with cone angle $169.36^\circ$. \newline
($\epsilon=-1$, $a=0.8$, $b=-1$, $A=4.68$)} \label{teardrop_pic}

\end{minipage}
\end{figure}

\begin{figure}[p]
\begin{minipage}[t]{0.45\linewidth} \centering

\includegraphics[height=\tam]{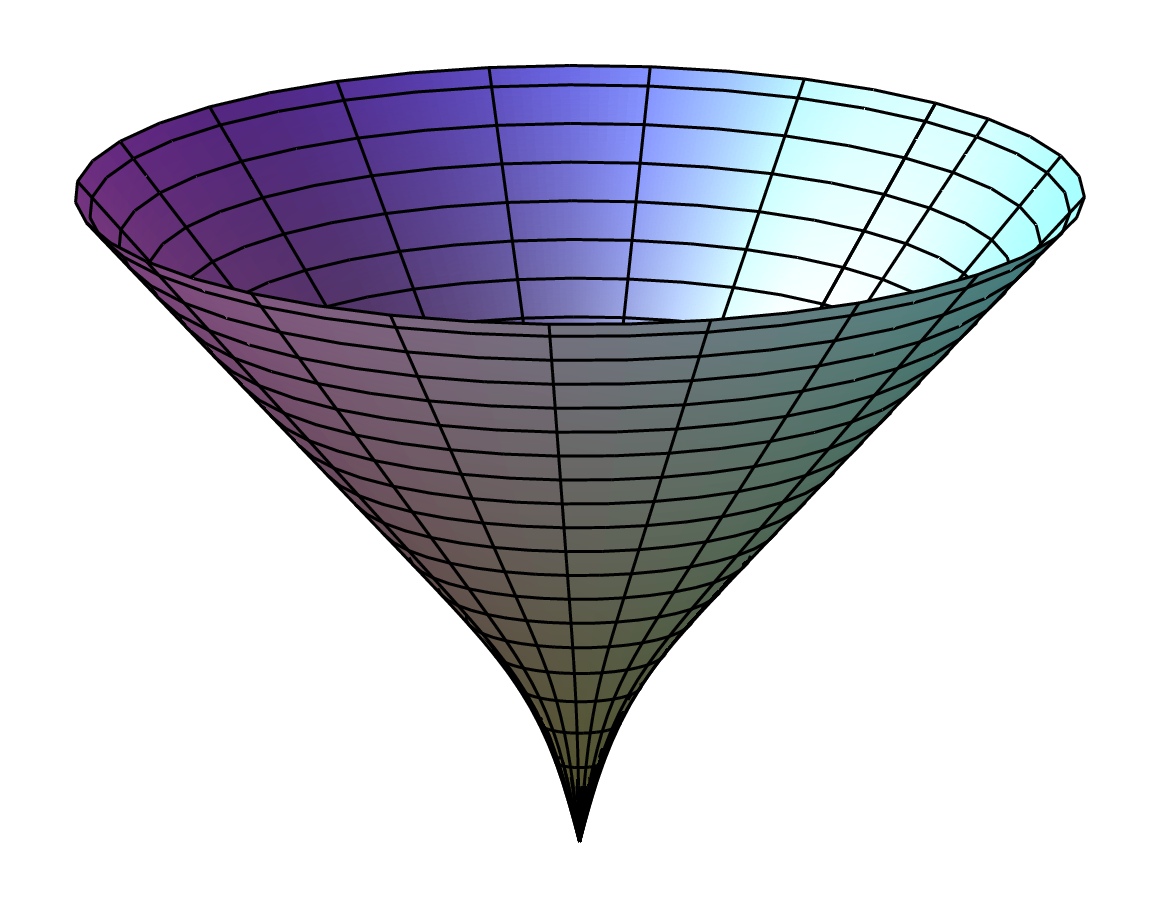}
\caption{ An $\alpha\beta$-cone soliton with asymptotic cone angle $\alpha=240^\circ$ and vertex cone angle $\beta=90^\circ$.
Note that the curvature is negative since $\alpha>\beta$. \newline
($\epsilon=1$,$a=0.75$, $b=-0.25$) }  \label{cone1_pic}

\end{minipage}
\hspace{1cm}
\begin{minipage}[t]{0.45\linewidth} \centering

\includegraphics[height=\tam]{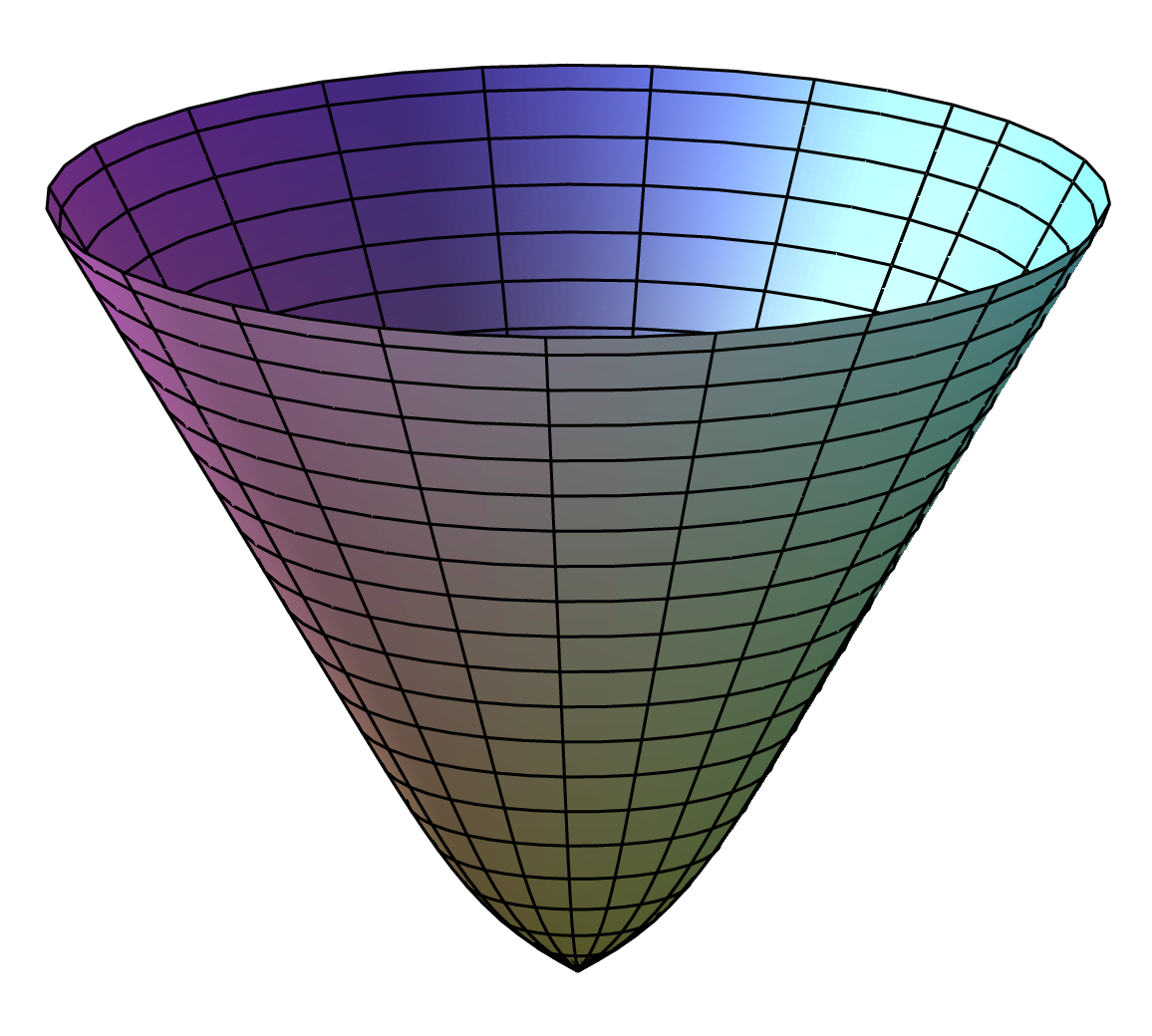}
\caption{ An $\alpha\beta$-cone soliton with asymptotic cone angle $\alpha=180^\circ$ and vertex cone angle $\beta=306^\circ$.
Note that the curvature is positive since $\alpha<\beta$. 
\newline ($\epsilon=1$, $a=1$, $b=-0.85$)} \label{cone2_pic}

\end{minipage}
\end{figure}

\begin{figure}[p]
\begin{minipage}[t]{0.45\linewidth} \centering

\includegraphics[height=\tam]{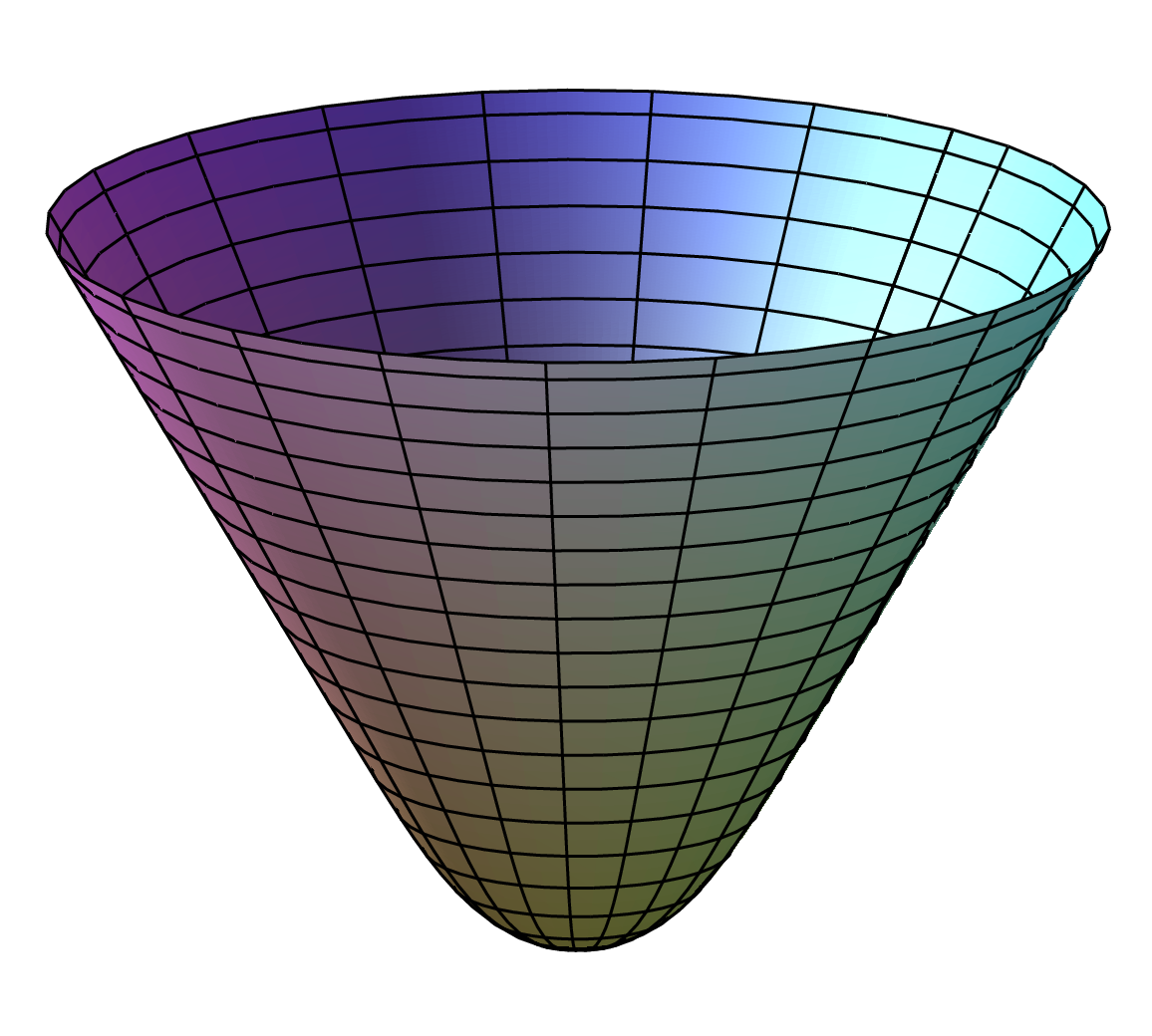}
\caption{A blunt $\alpha$-cone soliton with asymptotic cone angle $\alpha=180^\circ$. \newline
($\epsilon=1$, $a=1$, $b=-1$)}
\label{cone3_pic}

\end{minipage}
\hspace{1cm}
\begin{minipage}[t]{0.45\linewidth} \centering

\includegraphics[height=\tam]{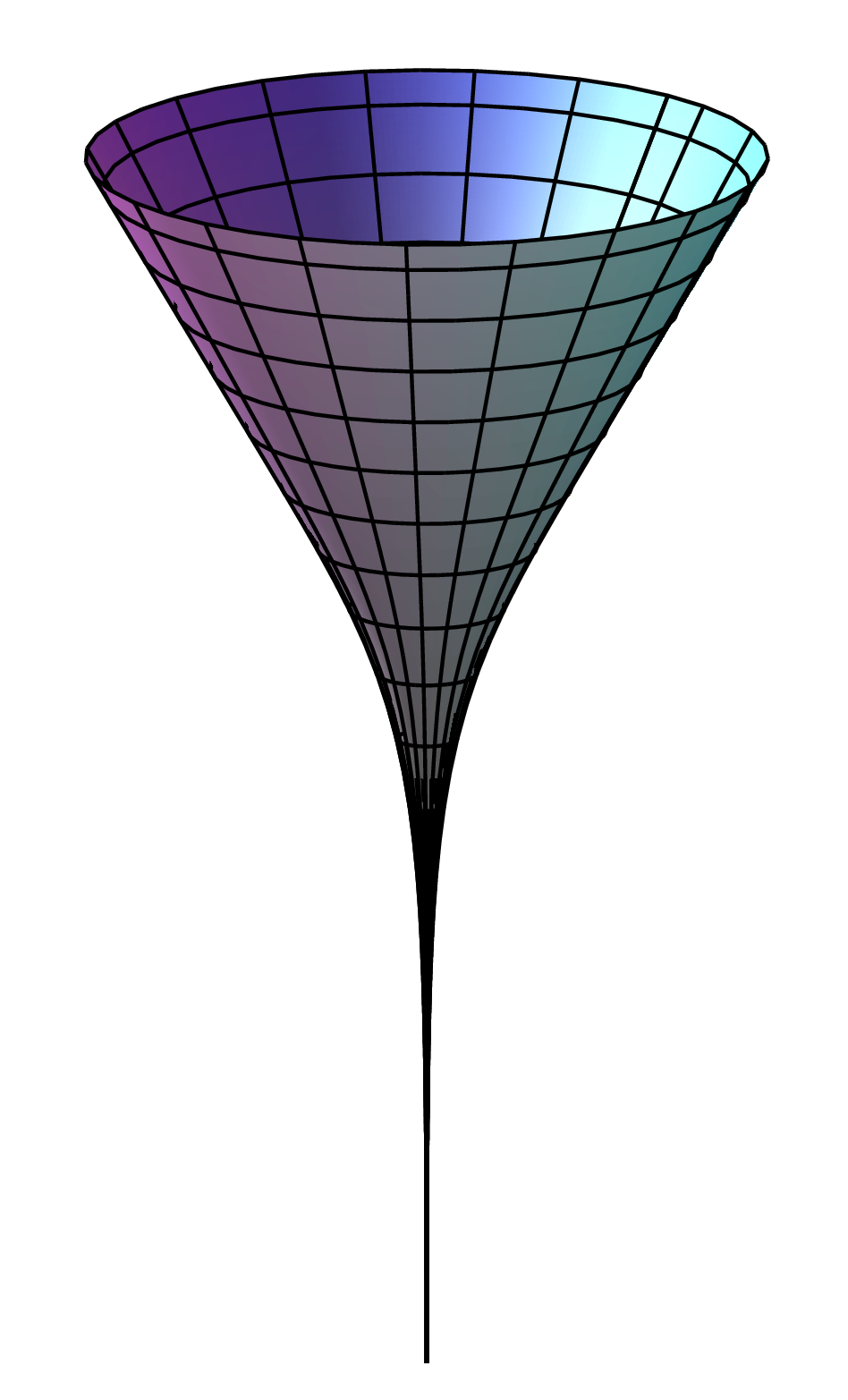}
\caption{A cusped $\alpha$-cone soliton with asymptotic cone angle $\alpha=180^\circ$. \newline
($\epsilon=1$, $a=1$, separatrix $S$ ($b\approx 0$))} \label{cone4_pic}

\end{minipage}
\end{figure}

%% file: survey.tex
\chapter{Uniformization of smooth surfaces} \label{Ch:survey}

\lettrine{\indent H}{amilton's work} on surfaces \cite{Hamilton_surfaces} is representative of several techniques and difficulties
that one finds on the Ricci flow. For the normalized Ricci flow on surfaces, Hamilton easily proved the convergence of the flow to a
metric of constant curvature in the cases where the surface satisfies $\chi \leq 0$. The main tools for this results are more or less involved computations of the evolution of the desired quantities, and applications of the \emph{maximum principle}. The case of the sphere ($\chi>0$) was more reluctant. In this case the evolution of the quantities is more subtle and several opposite effects are competing. Hamilton proved the
convergence to a metric of constant positive curvature under the additional assumption of initial positive curvature everywhere. He
developed two techniques in order to prove this: a Harnack inequality and a monotone entropy. These quantities involve higher order derivatives of the metric and hence a more precise control of the evolution. Later, Chow \cite{Chow_sphere} was able
to remove the additional assumption modifying the entropy formula. This difficulty for the convergence of the Ricci flow in this case
somehow reflects that unnormalized Ricci flow develops an infinite-curvature singularity, whereas this does not happen in non-positive Euler
characteristic.

One of the main advances in Perelman's work is the successful use of a blow-up technique for analyzing the singularities.
Perelman \cite{Perelman1} works in the scope of three-manifolds, but some techniques work in any dimension. While the blow-up technique had
been proposed by Hamilton and others for three-manifolds, Perelman's concept of $\kappa$-noncollapsing and the control of this property is
the key
to classify all singularities by means of the models called $\kappa$-solutions. These $\kappa$-solutions are limits of rescalings of the
evolving manifold at points of exploding curvature, and thus give local information of the manifold in a neighbourhood of the singularity
at the singular time. Perelman gave an enumeration of all three-dimensional $\kappa$-solutions, that allowed him to go on to a surgery
process of the manifold.

Perelman's method works in dimension two. The main difference is that the only $\kappa$-solution in dimension two is the shrinking
round sphere, that is \emph{compact}. This fact is essential, since a $\kappa$-solution arises as a limit of dilations. If the limit of
rescalings of a manifold is compact, then the limit is not only a local model, but a global one. Hence, a sphere evolving with the
unnormalized Ricci flow converges to a \emph{round point}, and therefore under the normalized Ricci flow converges to the round sphere.

In this chapter we present briefly the details of this argumental line of proving the uniformization of surfaces. Our purpose in Chapter
\ref{Ch:cone_rf} will be to modify the appropriate points of this argument to adapt to the cone surfaces setting. Therefore, we will survey
quickly most
parts, while paying more attention to those that will require modifications.

The chapter is organized as follows. In Section \ref{S:srv:HamRFsurf} we review the main arguments in Hamilton's article
\cite{Hamilton_surfaces},
especially the part dealing with the case $\chi\leq 0$. Then we stop the discussion of surfaces, and in Section \ref{S:srv:maxppl} we review
the techniques of
maximum principles and Harnack inequalities, applied to Ricci flow in two and in $n$ dimensions, due to S.-T. Yau \cite{LiYau} and Hamilton
\cite{Hamilton_surfaces}, \cite{Hamilton_Harnack}. These results are essential to Perelman's work. In Section \ref{S:srv:sings} we begin the
discussion
of the singularities of the Ricci flow, and the technique of parabolic rescalings. We also introduce the $\kappa$-noncollapsing property
of Perelman, that is the key property to successfully apply a blow-up by parabolic rescalings. In Section \ref{S:srv:noncolthm} we review
the
sophisticated Noncollapsing Theorem of Perelman, that ensures that the noncollapsing property is essentially preserved along the flow. Once
we can apply a rescaling blow-up to a singularity of the flow, the limit that we obtain is a so-called $\kappa$-solution. In Section
\ref{S:srv:ksols} we review the properties of $\kappa$-solutions. In Section \ref{S:srv:classif_ksols} we come back to the discussion of the
flow on surfaces, and we show the classification of all $\kappa$-solutions in dimension two. This is a classification of all possible models
of infinite-curva ure singularities, and the only possible case in this classification is the shrinking round sphere. On 
surfaces with $\chi>0$ the Ricci flow must develop an infinite-curvature singularity, and therefore the flow converges, up to rescaling,
to a round sphere.

\section{Hamilton's Ricci flow on smooth surfaces} \label{S:srv:HamRFsurf}
We survey in this section the original work of Hamilton for Ricci flow on surfaces, \cite{Hamilton_surfaces}. In this section we will use
the normalized Ricci flow as in Hamilton's article,
$$\frac{\partial}{\partial t}g = (r-R) g$$
with $g(0)=g_0$, where $r$ is the average scalar curvature, a constant defined by
$$r=\frac{\int R \ d\mu}{\int d\mu} = \frac{4\pi\chi(\mathcal M)}{\Area(\mathcal M)}$$
that depends on the topology by Gauss-Bonnet theorem. This normalization is a time-dependent rescaling so that the area of the compact
surface is kept invariant. As we will see, this has the advantage that the flow is defined for all time, given any initial metric, and
therefore no
singularities appear on this flow. The standard (unnormalized) Ricci flow develops singularities in the case of positive Euler
characteristic (so $r>0$), i.e. on the sphere. The normalized Ricci flow avoids this problem, but the study in the case of the sphere is
much more complicated than on the nonpositive cases. We first review the easier cases of $\chi(\mathcal M)\leq 0$, and the general arguments
of maximum principles.
Later we sketch Hamilton's argument for the case $R>0$ and Chow's completion for the case $r>0$ \cite{Chow_sphere}. We propose a
different way for the case $r>0$ using Perelman's analysis of the singularity models and therefore we will not discuss it in great detail.
However, some techniques we use in Chapter \ref{Ch:cone_rf} follow closely Hamilton's, so it is worth to review. 
\medskip

First attempt to gain some information of the flow is computing the evolution of the (scalar) curvature and applying the maximum principle
(see Section \ref{S:srv:maxppl} below).

\begin{tma}[{\cite[Sec 3]{Hamilton_surfaces}}]
The evolution of the scalar curvature for the normalized Ricci flow is
$$\frac{\partial}{\partial t} R = \Delta R + R^2 -rR$$
and therefore:
\begin{itemize}
 \item If $R\geq 0$ at $t=0$, then it remains weakly positive for all time.
 \item If $R\leq 0$ at $t=0$, then it remains weakly negative for all time.
 \item If $-C<R<-\epsilon<0$ at $t=0$, then $$re^{-\epsilon t} \leq r-R \leq C e^{rt} .$$
 \item If $R<0$, then the metric is defined for all time and converges exponentially fast to a metric of constant negative curvature.
\end{itemize}
\end{tma}

In the case of $R>0$ somewhere, the inequalities are set in the wrong direction and we cannot get useful estimates as in the case $R<0$
everywhere. This precludes the difficulty due to the positive curvature. To obtain better results, we need to introduce another ingredient:
gradient
solitons. For the normalized Ricci flow there are no distinctions between shrinking, steady or expanding solitons, since the area is fixed
and therefore no homothetic factor applies. The soliton equation for the normalized Ricci flow is
\begin{equation} \label{E:solitonNRF}
 \Hess f -\frac{1}{2}(R-r)g=0
\end{equation}
where $f$ is a potential function such that the metric evolves by diffeomorphisms induced by $\nabla f$. If the Ricci flow is not a
gradient soliton, there is no such function $f$ satisfying the soliton equation, but we can always solve the traced equation. The
\emph{potential} $f$ of a Ricci flow is the solution of
$$\Delta f =R-r$$
normalized to have mean value zero, $\int f =0$
(see Section \ref{S:srv:classif_ksols}). Then we define the soliton quantity
$$M:=\Hess f -\frac{1}{2}(R-r)g = \Hess f -\frac{1}{2}\Delta f \ g$$
which is the ``traceless Hessian'' of $f$, and vanishes iff the flow is a soliton. We also define
$$h=\Delta f + |\nabla f|^2$$
that helps controlling $R=h+|\nabla f|^2+r \leq h+r$.

\begin{tma}
The evolution of $h$ for the normalized Ricci flow is
$$\frac{\partial}{\partial t}h=\Delta h -2|M|^2 +rh$$
and therefore:
\begin{itemize}
 \item If $h\leq C$ at $t=0$, then $h\leq Ce^{rt}$ for all time.
 \item For any initial metric, there is a constant $C$ such that $$-C\leq R \leq Ce^{rt}+r .$$
 \item For any initial metric, the normalized Ricci flow has solution for all time $t\in [0,+\infty)$
 \item If $r< 0$, then the metric converges exponentially fast to a metric of constant negative curvature.
\end{itemize}
\end{tma}

This proves the uniformization theorem for the case $\chi< 0$. The case $r=\chi=0$ can be solved by similar techniques applying the maximum
principle to
the evolution of several other quantities, such as $|\nabla f|^2$ and $\nabla R$ (this argument is in the style of
Bernstein-Bando-Shi estimates, and is different from Hamilton's original, see for instance \cite[Sec 5.6]{ChowKnopf}).

The case $r>0$ turned out to be much harder. The argument, roughly speaking, is just applying once again the maximum principle to the
quantity $|M|^2$.

\begin{tma}
The evolution of $|M|^2$ for the normalized Ricci flow is
$$\frac{\partial}{\partial t} |M|^2 =\Delta |M|^2  -2|\nabla M|^2 -2R|M|^2$$
and therefore:

If $R\geq c > 0 $, then $|M|^2\leq Ce^{-ct}$ for all time.
\end{tma}

This implies that the metric converges exponentially fast to a soliton metric. Together with a classification of solitons (all
gradient solitons on a smooth compact surface have constant curvature), this yields the convergence to a round sphere. 
The main issue, however, is proving the positive lower bound for the curvature that is required.

Hamilton proved the case $R>0$, developing a Harnack inequality and an entropy quantity. This allows to prove that after some time there
exists the desired positive lower bound for the curvature. Although at the end Harnack inequalities are no more than
applications of the maximum
principle, they involve higher order derivatives of the quantities involved. Heuristically, in the case $R<0$ it suffices to watch at the
distribution of the curvature, but in the case $R>0$ one needs to look also at the distribution of the convexity of the curvature to guess
how it will evolve. This allows to compare the curvature at different points and different times. Finally, the case $r>0$ was proved by B.
Chow, \cite{Chow_sphere}, using a modified version of the Harnack inequality and a
modified entropy. 

The details for all the proofs can be found in the original articles \cite{Hamilton_surfaces}, \cite{Chow_sphere}, or in the reference
book \cite[Ch 5]{ChowKnopf} with some variations. We won't survey on these proofs, since we propose an alternative path using Perelman's
technique. However, we review the Harnack inequalities in the next section since we will need them.

\section{Maximum principles and Harnack inequalities} \label{S:srv:maxppl}
Under the name of \emph{maximum principles} (for parabolic equations) there is a whole family of theorems sharing a common ``principle'',
including the so-called Harnack inequalities. Hamilton made a strong use of these principles, and developed some new ones. It is not our
purpose to develop all the proofs, but we will give a brief survey on these theorems since there is an important modification when adapting
maximum principles from smooth to cone surfaces.

The fundamental idea is the fact that one has knowledge of space derivatives on a maximum (or minimum) point, and this restricts the time
derivative on the function on that same
point, locally. One of the easiest maximum principles is the following.
\begin{lema}
Let $\mathcal M$ be a closed Riemannian manifold, and let $f:\mathcal M  \times [0,T] \rightarrow \mathbb R$ satisfy
$$\frac{\partial}{\partial t} f \leq \Delta f + \langle X(t), \nabla f \rangle$$
and $f(x,0)=f_0$, for some function $f_0$ and some vector field $X(t)$.
Then the space maximum of $f$, $\max f(\cdot,t)$ is decreasing (nonstrictly) in $t$.
\end{lema}
\begin{proof}
Suppose that the maximum of $f$ at time $t_0$ occurs at a point $(x_0,t_0)$. Since $x_0$ is interior to $\mathcal M$ (it is
closed), the Laplacian is nonpositive, and the gradient is null. Therefore,
the derivative in time of $f$ at the maximum point $(x_0,t_0)$ is nonpositive, and hence the maximum cannot increase. 
\end{proof}

A more refined version of the principle includes a ``reaction'' term, that controls how much the function can grow.
\begin{lema}
Let $\mathcal M$ be a closed Riemannian manifold, and let $f:\mathcal M \times [0,T] \rightarrow \mathbb R$ satisfy
$$\frac{\partial}{\partial t} f \leq \Delta f + \langle X(t), \nabla f \rangle + F(f,t)$$
suppose that $\phi:[0,T]\rightarrow \mathbb R$ is the solution of
$$\left\{ \begin{array}{rcl}
          \frac{\partial \phi}{\partial t} &=& F(\phi(t),t) \\
	  \phi(0) &=&\alpha
         \end{array}
\right.$$
Suppose also that $f(x,0)\leq \alpha$, then $f(x,t)\leq \phi(t)$ for all $t\in[0,T]$.
\end{lema}

Harnack inequalities are elaborated uses of the maximum principle. Most times, these theorems are used to control the
value of the evolving function at different points and at different times. We illustrate this first with Yau's Harnack inequality
(\cite{LiYau}, see \cite[Thm 6.1]{Hamilton_surfaces}).

\begin{tma}[Harnack inequality for manifolds]
Let $\mathcal M$ be a closed Riemannian manifold with nonnegative Ricci curvature, and let $f:\mathcal M \times [0,T] \rightarrow \mathbb
R$ be a solution of the heat equation
$$\frac{\partial f}{\partial t}=\Delta f $$
such that $f>0$ everywhere. Then, 
\begin{enumerate}
 \item the function $L=\ln f$ satisfies
$$\Delta L \geq -\frac{n}{2t} .$$
\item For any two points $x_1,x_2\in \mathcal M$ and any two times $0<t_1<t_2$,
$$f(x_1,t_1) t_1^{{n}/{2}} \leq f(x_2,t_2) t_2^{{n}/{2}} e^{{A}/{4}} $$
where $A = \frac{\mathrm{dist}(x_1,x_2)^2}{t_2-t_1}$.
\end{enumerate}
\end{tma}

\begin{proof}[Sketch of the proof] 
From the evolution equation for the function $f$ we deduce the
evolution of the function $Q=\Delta L$ involving second derivatives of $f$ (roughly speaking, we compute the evolution of the
convexity of $f$):
\begin{align*}
\frac{\partial Q}{\partial t} &= \Delta Q + 2\langle \nabla Q ,\nabla L \rangle + 2\Ric(\nabla L,\nabla L) + 2 |\nabla^2 L|^2  .\\
 &\geq \Delta Q + 2\langle \nabla Q ,\nabla L \rangle +\frac{2}{n} Q^2
\end{align*}
The curvature appears when commuting covariant derivatives, and we use the curvature bound for the inequality. We apply the maximum
(minimum) principle to this new equation to obtain the inequality
$Q \geq -\frac{n}{2t}$ of the first statement.
The argument
is concluded by integrating this inequality along a geodesic $\gamma$ joining
the two points with an arc-parameter proportional to time,
$$L(x_2,t_2)-L(x_1,t_1)= \int_{t_1}^{t_2} \frac{dL(\gamma(t),t)}{dt} \ dt = \ldots \geq
\frac{-n}{2}\ln\left(\frac{t_2}{t_1}\right)-\frac{A}{4}$$
and exponentiating.
\end{proof}

Following the same pattern, Hamilton developed in \cite{Hamilton_surfaces} a Harnack inequality for Ricci flow on surfaces.
%
The main difference with respect to the previous Harnack theorem is that now the metric is evolving. This was one of the key tools for the
convergence of the flow in the case of the sphere.
A few years later, Hamilton developed a great generalization of the Harnack inequality for Ricci flow in any dimension
\cite{Hamilton_Harnack}, known as Hamilton's matrix Harnack inequality for the Ricci flow. This requires handling the evolution of tensor
quantities and applying maximum principles to them. Stated in index-notation, it is the following.

\begin{tma}[Harnack inequality for Ricci flows] \label{T:Harnack_RF}
Let $(\mathcal M,g(t))$ be a complete Ricci flow on an open or closed manifold, with 
$$\Rm \geq 0$$ 
for $t\in[0,T]$. Let
$$P_{ijk}=\nabla_i R_{jk} - \nabla_j R_{ik}$$
and
$$M_{ij}=\Delta R_{ij} -\frac{1}{2} \nabla_i\nabla_j R + 2R_{kij}^l R_l^k - R_i^k R_k^j + \frac{1}{2t}R_{ij} .$$
For any one-form $W^i$ and any 2-form $U^{ij}$, let
$$ Z = M_{ij} W^i W^j +2P_{kij} U^{ki}W^j + R_{ijkl}U^{ij}U^{kl}$$
Then,
$$ Z\geq 0 .$$
\end{tma}

\begin{proof}[Sketch of the proof]
The structure of this theorem is similar to the previous Harnack inequalities. The source for obtaining 
$Z$ is picking originally an
evolving quantity, that we choose to be
$$ S=\Ric + \Hess f + \frac{1}{2t}g $$
which is precisely the quantity that vanishes iff the Ricci flow is a gradient soliton (shrinking or expanding). This is a tensor quantity,
we compute a second covariant derivative, and after some commuting and contracting we get the quantity $Z$. Then we compute again the
evolution of $Z$ and the resulting equation is suitable to apply the maximum principle to get $Z\geq 0$. This statement is useful as it is,
but one could continue the analogous proof and integrate the quantity $Z$ along a geodesic path to get some results about $S$ and comparing
$R$ on different spacetime points \cite[Cor 1.3]{Hamilton_Harnack}.

In order to apply the maximum principle to $Z$, the argument \cite[Sec 5]{Hamilton_Harnack} consists on creating a barrier function
$\hat Z$ that is very positive if $t\rightarrow 0$ or if the point tends to infinity in space, if the manifold is not compact. The barrier
function is 
$$\hat Z = \hat M_{ij}W^iW^j +2P_{kij} U^{ki}W^j + \hat R_{ijkl}U^{ij}U^{kl}$$
where
$$\hat M_{ij}=M_{ij}+\frac{1}{t}\varphi \ g_{ij}$$
and
$$\hat R_{ijkl} = R_{ijkl} +\frac{1}{2}\psi \ (g_{ik}g_{jl}-g_{il}g_{jk})$$
for suitable chosen functions $\varphi$ and $\psi$. 
\end{proof}

These functions are given by the following lemma \cite[Lem 5.2]{Hamilton_Harnack}.
\begin{lema} \label{barriers_Hamilton}
For any $C$, $\eta>0$ and any compact set $K$ in space-time, we can find functions $\psi=\psi(t)$ and $\varphi=\varphi(x,t)$ such that
\begin{enumerate}
 \item $\delta\leq\psi\leq\eta$ for some $\delta>0$, for all $t$;
 \item $\epsilon\leq\varphi\leq\eta$ on the compact set $K$ for some $\epsilon>0$, for all $t$. Furthermore, $\varphi(x,t)\rightarrow\infty$
if
$x\rightarrow\infty$, i.e. the sets $\{x\ |\ \varphi(x,t)<M\}$ are compact for all $t$ and all $M$;
 \item $\frac{\partial \varphi}{\partial t} > \Delta \varphi + C \varphi$;
 \item $\frac{\partial \psi}{\partial t} > C \psi$;
 \item $\varphi\geq C \psi$.
\end{enumerate}
\end{lema}
When we deal with cone surfaces, this is the only point on Hamilton's argument that will require a modification: a maximum principle with a
barrier on the cone points.

The quantity $Z$ involved in Hamiton's matrix Harnack is a 3-tensor quantity applied to some 1-form and 2-form. A useful consequence is an
inequality obtained tracing two of the indices \cite[Cor 1.2]{Hamilton_Harnack}.

\begin{tma}[Traced Harnack inequality for Ricci flows] \label{T:Harnack_RF_traced}
Let $(\mathcal M,g(t))$ be a complete Ricci flow on an open or closed manifold, with 
$$\Rm \geq 0$$ 
for $t\in[0,T]$. Then, for any vector field $X$,
$$ H(X) = \frac{\partial R}{\partial t} +\frac{1}{t} R + 2\langle \nabla R, X\rangle + 2\Ric(X,X) \geq 0 .$$
\end{tma}

Let us remark that from the construction of the quantity $Z$, the inequalities in Theorems \ref{T:Harnack_RF} and \ref{T:Harnack_RF_traced}
turn into equalities when the Ricci flow is an homothetic soliton.

\section{Singularities of the Ricci flow and $\kappa$-noncollapse} \label{S:srv:sings}

Singularity formation is a phenomenon that may occur in the long time behaviour of the $n$-dimensional Ricci
flow. Roughly speaking, the only obstruction to continue the evolution of the Ricci flow further in time is the appearance of points with
infinite curvature in finite time $T$, so the Ricci curvature term in the flow equation becomes undefined and the equation has no further
solution for $t=T$. The theorem for smooth flows is the following (\cite[Thm. 8.1]{Hamilton_formsing}, see also
\cite[Sec 5.3]{Topping_lectures}).

\begin{tma}
Let $\mathcal M$ be a smooth closed $n$-manifold, and $g(t)$ a Ricci flow on a maximal time interval $[0,T)$ and $T<\infty$, then
$$\sup_{\mathcal M} |\Rm|(\cdot,t) \rightarrow \infty$$
as $t\rightarrow T$.
\end{tma}

\begin{proof}[Sketch of the proof]
By contradiction, suppose that the curvature remains bounded for all $t\in [0,T)$. We pick a sequence $t_i\rightarrow T$, and we construct
the sequence of metrics $g_i=g(t_i)$. We will see that the metric tensor and all its derivatives are bounded. Then, by the compactness
theorems (e.g Arzelà-Ascoli in the smooth case), the sequence $g_i$ converges to a smooth metric $g(T)$. Once the flow has been
extended smoothly to a metric in $g(T)$, the flow can be further extended using $g(T)$ as the initial condition by the
short time existence results, thus contradicting the assumption.

It is an immediate consequence from the flow equation, that if the Ricci curvature is bounded, then the distortion of the metric is
bounded for finite time. Namely, if
$g(t)$ is a Ricci flow on $\mathcal M \times [0,s]$ and $|\Ric|<M$, then
$$e^{-2Mt}g(0)\leq g(t) \leq e^{2Mt}g(0)$$
for all $t\in[0,s]$.
This implies that if the curvature is bounded in $\mathcal M \times [0,T)$, then the metric can be continuously extended to a metric
in $t=T$. In order to prove that this extension can be done smoothly, one needs to apply Shi's estimates (see \cite{Shi}, see also
\cite[Ch 7]{ChowKnopf}), namely, if
$|\Rm|<M$ in
$t\in[0,T]$, then 
$$\bigg| \frac{\partial}{\partial t^l} \nabla^k \Rm \bigg| \leq C=C(l,k,M,T,n)$$ 
for $t\in[T/2,T]$. It also can be applied a local Shi's estimate, where the hypothesis applies in a ball of radius $r$ and the thesis in a
ball of radius $r/2$, with a constant $C=C(l,k,M,T,n,r)$.

\end{proof}

In the case of surfaces, there is a simple argument that ensures that any sphere evolving with the Ricci flow develops a singularity in
finite time.
\begin{lema} \label{L:RF_on_sphere_has_sing}
Let $(\mathcal M^2,g(t))$, $g(0)=g_0$, be a Ricci flow on a surface with $\chi(\mathcal M)>0$, defined on a maximal time interval $[0,T)$.
Then $T<\infty$.
\end{lema}
\begin{proof}
Computing the evolution of the area of the surface, we obtain
$$\frac{d}{dt}\mathrm{Area}(\mathcal M)=\int_{\mathcal M}\frac{d}{dt} d\mu = \int_{\mathcal M} -R\ d\mu = -4\pi\chi(\mathcal M) $$
and since  $\chi(\mathcal M)>0$, the area or $\mathcal M$ is a decreasing linear function of $t$.
$$\Area(\mathcal M)=-4\pi\chi(\mathcal M)\ t + \Area(g(0)) .$$
Therefore there is a singular time when the area collapses to zero (or earlier).
\end{proof}
A priori there could be a singular time
before the area goes to zero, with an isolated point of infinite curvature. We will see that this is not the case, and at the singular time
all the surface tends to a point while approaching constant curvature (a \emph{round point}).

\medskip

A useful technique for analyzing this infinite-curvature is parabolic rescaling. A \emph{parabolic rescaling} is a transformation of the
evolving metric into the form
$$g(t) \mapsto \lambda^2 g\left(\frac{t}{\lambda^2}\right) = \tilde g(t) .$$
It is of interest because if $g(t)$ is a Ricci flow, then so is $\tilde g(t)$. Besides, distances
get multiplied by $\lambda$, time gets multiplied by $\lambda^2$ and scalar curvature gets divided
by $\lambda^2$. The idea is to pick a sequence $\{\lambda_i\}_{i\in \mathbb N}$, and construct a
sequence of rescaled pointed Ricci flows keeping the curvature bounded. The hope is to find a convergent
subsequence of these pointed flows, using techniques of classes (spaces) of manifolds, and compactness theorems for classes of
manifolds. If this process of iterated rescalings is successful and wo obtain a limit, we call this a \emph{blow-up} of the singularity. In
Appendix \ref{Ch:compactness} we will develop equivalent compactness theorems for the case of cone surfaces and flows, that we will use in
Chapter
\ref{Ch:cone_rf}.

The compactness theorems for classes of manifolds we consider usually need as a hypothesis a control of the injectivity radius of the
manifold. This allows one to consider the sequence of manifolds as a sequence of metrics on a fixed chart, with its coordinate functions
defined over a fixed open subset of $\mathbb R^n$. 

We begin by recalling the Riemannian fact of the equivalence, under bounded curvature, of controlling volume or injectivity radius. (Cf.
\cite[Thm 10.6.8]{BurBurIva}  and \cite[Thm 4.3]{CheGroTay})
\begin{prop} \label{P:srv:inj_vol}
Let $(\mathcal M^n,g)$ be a smooth Riemannian $n$-manifold and $p\in\mathcal M$. Pick an $r>0$ such that
$$|\Rm(x)|\leq r^{-2}\ \forall x\in B(p,r)$$
Then for all $c_1>0$ exists $c_2>0$ (only depending on $c_1$ and $n$) such that
$$\frac{\inj(p)}{r}\geq c_1 \qquad \Rightarrow \qquad \frac{\Vol(B(p,r))}{r^n}\geq c_2 .$$
Conversely, for all $c_2>0$ exists $c_1>0$ (only depending on $c_2$ and $n$) such that
$$\frac{\Vol(B(p,r))}{r^n}\geq c_2 \qquad \Rightarrow \qquad \frac{\inj(p)}{r}\geq c_1 .$$
\end{prop}

Note that all quantities involved on this proposition, $\frac{\Rm}{r^{-2}}$, $\frac{\inj}{r}$, $\frac{\Vol}{r^n}$,  are scale-invariant.
Using this proposition, we can focus on controlling the volume, instead of the more elusive injectivity radius.
This motivates the definition of the $\kappa$-noncollapsed property by Perelman. 

\begin{defn}[$\kappa$-noncollapsed manifold]
Let $\kappa, \rho>0$. A Riemannian manifold $(\mathcal M^n,g)$ is \emph{$\kappa$-noncollapsed at scale $\rho$} in $p\in\mathcal M$ if
$\forall
r<\rho$ it is satisfied that
$$ |\Rm|(x)\leq \frac{1}{r^2}\ \forall\ x\in B(p,r) \ \ \Rightarrow \ \ \frac{\Vol(B(p,r))}{r^n}\geq
\kappa $$
\end{defn}

This means: every ball with radius $r<\rho$ and bounded curvature has a volume of at least $\kappa r^n$. Observe that for any given
smooth compact manifold
one can find $\kappa$, $\rho$ small enough such that the manifold is $\kappa$-noncollapsed at scale $\rho$.

\begin{defn}[$\kappa$-noncollapsed Ricci flow]
Let $\kappa, \rho>0$. A Ricci flow $(\mathcal M^n,g(t))$ is \emph{$\kappa$-noncollapsed at scale $\rho$} in $p\in\mathcal M$ if $\forall
r<\rho$ it is satisfied that
$$|Rm|(x)\leq \frac{1}{r^2} \forall\ x\in B_{t}(p,r)\ \forall t\in [t_0-r^2, t_0]\ \Rightarrow \
\frac{Vol(B_{t_0}(p,r))}{r^n}\geq \kappa$$
\end{defn}
This means: every time $t$ slice of the flow is a $\kappa$-noncollapsed manifold at scale $\rho$, with the same $\kappa$ for each $t\in
[t_0-r^2, t_0]$.
Heuristically, the class of manifolds (or flows)
enjoying a $\kappa$-noncollapse at all scale, is compact, and hence we can ensure
subconvergence of our sequence of rescaled Ricci flows.

\section{Noncollapsing theorem} \label{S:srv:noncolthm}

The $\kappa$-noncollapsing property plays an essential role in the blow-up process of the singularity. It ensures that the sequence of
blow-ups has a
uniformly bounded injectivity radius around the basepoint that allows to endow the limit with a smooth structure of $n$-manifold.
The key property is that $\kappa$-noncollapse is perserved under the Ricci flow (\cite{Perelman1}, see also \cite[Thm 26.2]{KleinerLott}
)

\begin{tma}[Noncollapsing theorem (Perelman)]
Given numbers $n\in \mathbb N$, $T<\infty$, $\rho, K, c>0$, there exists $\kappa>0$ such that the following holds: Let $(M^n,g(t))$ be a
Ricci
flow defined on $[0,T)$ such that
\begin{itemize}
 \item $|Rm|$ is bounded on every compact subinterval $[0,T']\subset[0,T)$.
 \item $(M,g(0))$ is complete with $|Rm|<K$ and $\inj(M,g(0))\geq c >0$.
\end{itemize}
Then the Ricci flow is $\kappa$-noncollapsed at scale $\rho$. Furthermore, $\kappa$ is (nonstrictly) decreasing in $T$, while all other
constants fixed.
\end{tma}

This means that if the time-zero slice is $\kappa '$-noncollapsed, then the whole Ricci flow is $\kappa$-noncollapsed for $t\in[0,T)$ (with
possibly different $\kappa$ and $\kappa '$). In particular all the time-$t$ manifolds are $\kappa$-noncollapsed.

Remark that if $M$ is compact, we can always find a scale $\rho_0$ and a parameter $\kappa_0$
such that $M$ is $\kappa_0$-noncollapsed at scale $\rho_0$. The point is that $\kappa$-noncollapsing allows us to apply the convergence
theorems.

Perelman developed several proofs of the noncollapsing theorem. We review the ``comparison geometry approach'' using the $\mathcal
L$-geodesics.
\begin{proof}[Sketch of the proof]
Fix $(p,t_0)\in \mathcal M \times [0,T)$. We focus our atention on curves $\gamma (\tau)$ parameterized backwards
in time ($\tau=t_0-t$). Following Perelman, we define the $\mathcal L$-length of a curve as
 $$\mathcal L(\gamma)=\int_{\tau_0}^{\tau_1}
\sqrt{\tau}(R(\gamma(\tau))+ |\dot\gamma(\tau)|^2)d\tau .$$
Associated to this length, we define several quantities: the $L$-length to a point
$$L(q,\tau)=\min \{\ \mathcal L (\gamma) \ | \ \gamma(0)=p,\ \gamma(\tau)=q\} ,$$
the reduced length
$$l(q,\tau)=\frac{L(q,\tau)}{2\sqrt{\tau}} ,$$
and the reduced volume of the manifold 
$$ \tilde V(\tau)=\int_M \tau^{-n/2} e^{-l(q,\tau)} dq .$$

Here ``reduced'' means adimensional, or more precisely, invariant under changes of scale on the metric. 
The rationale of this construction is that $\tilde V$ can be used to measure the maximum size of balls not collapsed at a given scale.
Therefore, the reduced volume can measure the noncollapse of previous time-sheets as seen from a chosen spacetime basepoint. The key
step is proving that $\tilde V$ is increasing (nonstrictly) in $t$. Then the argument for the proof of the theorem is as follows: on the one
hand, the
hypothesis on $(M,g(0))$ imply that the initial manifold is noncollapsed and hence one can find a positive lower bound for the reduced
volume at time $t=0$. On the other hand, if the theorem were false one could find certain basepoints $(p_k,t_k)$ and radius $r_k$ such that
the parabolic regions $B_{t_k}(p_k,r_k)\times [t_k-r_k^2,t_k]$ would be collapsed (i.e., although having bounded normalized curvature
$|\Rm|r^{-2} \leq 1$, the normalized volume $\frac{\vol B_{t_k}(p_k,r_k)}{r_k^n}$ is arbitrarily small). This would imply that the reduced
volume $\tilde V$ on a time sheet close to $t=t_k$ would be arbitrarily small, and that contradicts the monotonicity of $\tilde V$ and its
uniform lower bound at $t=0$.

On the remaining of the proof, we describe the proof of the monotonicity of $\tilde V$ (see \cite[Secs 15-23]{KleinerLott}). First step is
developing a variational theory for
the $\mathcal L$-distance analogous to the classical theory of length and geodesics. More precisely, obtaining variational formulas and
estimations for the $\mathcal L$-length, $\mathcal L$-geodesics, an $\mathcal L$-exponential map, and $\mathcal L$-Jacobi
fields. All this is with respect to a fixed spacetime basepoint. 

From the first variation of $\mathcal L$ we obtain an $\mathcal L$-geodesic equation
$$ \nabla_X X - \frac{1}{2} \nabla R + \frac{1}{2\tau}X + 2\Ric(X,\cdot) =0$$
where $X$ is the tangent vector to the $\mathcal L$-geodesic. We also get formulae for the first derivatives of the $L$-distance to the
basepoint.
$$\frac{\partial L}{\partial \tau} = 2\sqrt{\tau}R -\frac{1}{2\tau}L + \frac{1}{\tau}K$$
$$|\nabla L|^2 = -4\tau R + \frac{2}{\sqrt{\tau}}L - \frac{4}{\sqrt{\tau}} K$$
where $K$ is a grouping of terms appearing on the Harnack inequality. Specifically,
$$K=\int_0^{\bar\tau} \tau^{\frac{3}{2}} H(X(\tau)) \ d\tau$$
with
$$H(X)=-\frac{\partial R}{\partial \tau} -\frac{1}{\tau} R -2\langle \nabla R, X\rangle + 2 \Ric (X,X) .$$
This is the quantity of Theorem \ref{T:Harnack_RF_traced} after the change $\tau=-t$ and choosing $-X$ for $X$. Harnack inequality states
that if $\Rm\geq 0$ then
$H\geq 0$ and hence $K\geq 0$. However, we are not assuming this positivity hypothesis.
The noncollapsing theorem is independent of the Harnack inequality.

From the second variation of $\mathcal L$, we get an expression for the $\mathcal L$-Jacobi fields (variational fields of a variation by
$\mathcal L$-geodesics),
\begin{equation}  \label{E:LJac}
\begin{split}
Jac(Y):= & -\nabla_X \nabla_X Y - \frac{1}{2\tau}\nabla_X Y +\frac{1}{2}\nabla_Y (\nabla R) \\
& - 2(\nabla_Y \Ric)(X,\cdot) -2\Ric(\nabla_Y X, \cdot) =0 .
\end{split}
\end{equation}
 
We get a formula for the Hessian of $L$,
$$\Hess_{L(\cdot,\bar \tau)}(w,w)=Q(Y,Y) = 2\sqrt{\bar\tau} \langle \nabla_X Y, Y\rangle$$
where $Y$ is an $\mathcal L$-Jacobi field satisfying 
\begin{equation} \label{E:LJac_ics}
 Y(0)=0,\quad Y(\bar\tau)=w,
\end{equation}
and where
$$Q(\tilde Y,\tilde Y):=2\int_0^{\bar\tau} \sqrt{\tau}\langle \tilde Y, Jac(\tilde Y) \rangle \ d\tau + 2\sqrt{\bar\tau} \langle \nabla_X
\tilde Y (\bar\tau), \tilde Y(\bar\tau) \rangle$$
is a quadratic form, an ``index'', that minimizes over $\mathcal L$-Jacobi fields. This means
$$Q(Y,Y)\leq Q(\tilde Y, \tilde Y)$$
if $Y$ is an $\mathcal L$-Jacobi field and $\tilde Y$ is any other variational field along the $\mathcal L$-geodesic.

We use the test vector field $\tilde Y$ such that solves
\begin{equation} \label{E:Ltestfield}
 \nabla_X\tilde Y = -\Ric(\tilde Y,\cdot) + \frac{1}{2\tau}\tilde Y
\end{equation}
with the boundary condition $\tilde Y(0)=0$ and $\tilde Y(\bar \tau)=w$. This choice is done because $\tilde Y$ is Jacobi iff the Ricci flow
is a soliton with potential function $f=l$. Indeed, if $Y$ is a vector field satisfying \eqref{E:LJac} (Jacobi field), \eqref{E:Ltestfield},
and \eqref{E:LJac_ics}, then
$$\Hess L (Y,Y)=2\sqrt{\tau} \langle -\Ric(Y,\cdot)+\frac{1}{2\tau}Y,Y\rangle$$
and hence $l=\frac{L}{2\sqrt{\tau}}$ satisfies the soliton equation
$$\Hess l = -\Ric+\frac{1}{2\tau}g .$$
This allows us to obtain an
inequality for the Hessian of $L$ comparing our Ricci flow with a soliton. Tracing, this yields an inequality for the Laplacian,
$$\Delta L \leq \frac{n}{\sqrt{\tau}}-2\sqrt{\tau}R -\frac{1}{\tau} K$$
with equality iff we have a soliton. This is the only inequality needed to prove the monotonicity of $\tilde V$. 

The computation for the reduced volume $\tilde V(\tau)$ involves pulling back the integration domain $\mathcal M$ to the tangent space
$T_p\mathcal M$ via the $\mathcal L$-exponential map.
$$\tilde V(\tau) = \int_{T_p\mathcal M} \tau^{-\frac{n}{2}} e^{-l(\mathcal L\exp_\tau(v),\tau)} \mathcal J(v,\tau) \chi_\tau(v)\ dv$$
where $\mathcal J(v,\tau)$ is the jacobian of the $\mathcal L$-exponential and $\chi_\tau(v)$ is a cutoff function for the $\mathcal
L$-cut locus. As in the Riemannian case, the jacobian $\mathcal J$ is expressed in terms of Jacobi fields and we can use the inequalities
above. We prove monotonicity by derivating the integrand and proving that is negative, i.e. the reduced volume $\tilde V$ is decreasing in
$\tau$ (increasing in $t$), and is constant iff the flow is a soliton and iff the inequalities become equalities.

\end{proof}

\section{$\kappa$-solutions} \label{S:srv:ksols}

The $\kappa$-noncollapsing property allows one to find limits of sequences of rescaled Ricci
flows. Since we will use sequences of dilations around points of high positive curvature, the limit flow will have positive curvature, and
since we will be dilating the time before the singular moment, the limit flow will be ancient. We reserve the word
\emph{$\kappa$-solution} for those Ricci flows that enjoy these good properties.

\begin{defn}[$\kappa$-solution]
A \emph{$\kappa$-solution} is a Ricci flow 
\begin{itemize}
\item ancient, $t\in (-\infty,T]$,
\item nonflat, $Rm \neq 0$,
\item with $Rm$ positive definite and bounded in each time-slice, $|Rm|<C$,
\item $\kappa$-noncollapsed at all scales.
\end{itemize}
\end{defn}

The enhaced properties of a $\kappa$-solution allow an in-depth analysis that can't be done for a
general Ricci flow. Besides the $\kappa$-noncollapse at all scales, the positiveness of the curvature ensures that the Harnack
inequality holds on the quantities considered on the proof of the noncollapsing theorem, and thus it gives us valuable extra properties.

An important feature is that one can find a soliton ``buried'' inside every
$\kappa$-solution (Perelman \cite{Perelman1}, see also \cite{KleinerLott} Prop 39.1). More specifically, the limit backwards in time to
$t=-\infty$ is, after rescaling, a gradient shrinking soliton.

\begin{tma}[Asymptotic soliton]
Let $(\mathcal M^n,g(t))$ be a $\kappa$-solution. Take a sequence $\{\bar\tau_i\}$ tending to $+\infty$ of
(backwards) times. Pick in $(\mathcal M,g(-\bar\tau_i))$ the point $q_i$ where $l(\cdot, \bar\tau_i)$
achieves its minimum. Construct a sequence of Ricci flows using a scale such that the new basepoints are $(q_i,-1)$ (rescaling time to
$-1$), i.e.
$$g_i(t)=\frac{1}{\bar\tau_i}g(\bar\tau_i t).$$
Then, after passing a subsequence, the limit flow is a nonflat gradient shrinking soliton which is also a
$\kappa$-solution. Moreover, the soliton function is the limit $l_\infty(q,-1)$.
\end{tma}

\begin{proof}[Sketch of the proof.]
Continuing with the properties of the $\mathcal L$-distance, from the estimates for $\frac{\partial L}{\partial \tau}$ and $\Delta L$, we
can write an estimate for $\bar L(q,\tau) = 2\sqrt{\tau} L(q,\tau)$,
$$\frac{\partial \bar L }{\partial t} + \Delta\bar L \leq 2n$$
and applying the maximum principle, $\min \bar L (\cdot,\tau) -2n\tau \leq 0$, or, written in terms of
$l(q,\tau)=\frac{L(q,\tau)}{2\sqrt{\tau}}$,
$$\min l(q,\tau)\leq \frac{n}{2} .$$
We exploit this property to pick the basepoints $(q_i,\bar\tau_i)$ with $q_i$ achieving the minimum of $l$ at $\bar \tau_i$, that is
bounded above by the constant $\frac{n}{2}$.

From the same estimates as above, we can also get estimates for $\frac{\partial l}{\partial \tau}$, $|\nabla l|^2$, $\Delta l$. Using the
nonnegative curvature, the Harnack quantity $K$ is nonnegative, and the formulae simplify considerably to a more compact expression,
$$|\nabla l| +R \leq \frac{ C l}{\tau} ,$$
$$\frac{\partial l}{\partial \tau} \geq -\frac{(1+C)l}{2\tau} .$$
In particular, $|\nabla l |^2\leq \frac{C}{\tau} l $ and $\tau R \leq C l$.
These expressions allow to bound $l$ on an uniform spacetime neighbourhood around the basepoints $(q_i,\bar\tau_i)$, since the value of
$l$ at the basepoint is uniformly bounded, and the space and time derivatives are bounded. Similar uniform bounds can be obtained for $\tau
R$.

Now we take the pointed limit of the sequence of rescalings. The $\kappa$-noncollapsing ensures that the sequence of rescalings have
bounded $R$ and bounded injectivity radius. The sequence of rescaled manifolds have the same functions $l$ and $\tau R$ since both are
scale-invariant quantities. Since the sequence of rescaled functions $l_i (q,\tau)= l(q_i,\bar \tau_i)$ is uniformly bounded, we can apply
the Arzelà-Ascoli compactness theorem to get a limit function $l_\infty$. This function is not defined in terms of any $\mathcal L$-length,
but satisfies appropriate PDEs that pass to the limit (as a limit, $l_\infty$ is only Lipschitz, but it can be shown to be
smooth by PDE arguments). 

This allows to define a reduced volume $\tilde V$ with $l_\infty$ on the limit manifold, that agrees with the limit of the reduced
volumes of the sequence 
$$\tilde V_\infty(\tau) = \lim_{i\rightarrow\infty} \tilde V_i(\tau) = \lim_{i\rightarrow\infty}\tilde V(\tau\bar\tau_i)$$
Since the reduced volume on the original flow is positive and decreasing in $\tau$, the limit must be a constant, and hence the limit of
the rescalings is a soliton.

The reduced volume can also be used to show that the soliton is nonflat. On a shrinking soliton, $\tilde V$ is a constant $0<c\leq
(4\pi)^\frac{n}{2}$ and $c$ achieves the maximum value $(4\pi)^\frac{n}{2}$ if and only if the soliton is flat (and hence it is the
shrinking Gaussian soliton). To see the nonflatness, one bounds $c<(4\pi)^\frac{n}{2}$.

\end{proof}

\begin{tma}
Let $(\mathcal M^2,g(t))$ be a Ricci flow on a surface, defined on $[0,T)$, which becomes singular at time $T$. Let
$\kappa,\rho>0$ and assume $(\mathcal M,g(0))$ is $\kappa$-noncollapsed at scale $\rho$. There is a
sequence of times $t_i \rightarrow T$ such that, if $Q_i=\max R(\cdot,t_i)$ and $p_i$ is the point that achieves the maximum of $R$
at time $t_i$, then the sequence of pointed Ricci flows $(\mathcal M,g_i(t),p_i)$ with
$$g_i(t)=Q_i \ g\left( \frac{t}{Q_i}+ t_i \right)$$
has a subsequence that converges to a $\kappa$-solution.
\end{tma}

\begin{proof}
The evolution of the scalar curvature on a surface is
$$\frac{\partial}{\partial t} R = \Delta R + R^2$$
By the maximum principle, $\frac{\partial}{\partial t} R_{min} \geq 0$ and hence the curvature of $g(t)$ is bounded below by the minimum of
the curvature at time zero, $R_{g(t)}>-c$. The rescaling gets the curvature divided by $Q_i$, hence on the rescaled flow
$R_{g_i(t)}>-\frac{c}{Q_i}$.

The rescaled flow $g_i(t)$ is defined for $t\in [-t_i Q_i, (T-t_i)Q_i )$. Since $(T-t_i)Q_i>0$, we can restrict the flow to the interval
$[-t_iQ_i,0]$. The rescaling shifts the time so that the basepoint is at new time $0$. Therefore, on this interval the curvature of $g_i$
is bounded above by the curvature of the basepoint, $R_i(p_i,0)=1$.

The bounded curvature, together with the $\kappa$-noncollapse at all scales and the compactness theorems, imply that there is a
limit flow. This is a noncollapsed flow at all scales (since this is a scale-invariant property). It is ancient, since $t_iQ_i\rightarrow
\infty$ and therefore the the time domain of the limit is at least $(-\infty,0]$. It is nonflat since $R_i(p_i,0)=1$. Finally, it has
positive curvature since  $R_{g_i(t)}>-\frac{c}{Q_i}\rightarrow 0$.

\end{proof}

This theorem in dimension three is the so-called Canonical Neighbourhood theorem of Perelman, \cite[Thm 12.1]{Perelman1}, and the main
difficulties are the control of the whole Riemann tensor from the control of the scalar curvature (Hamilton-Ivey pinching), and the careful
setting of the constants controlling the size of the neighbourhood.

\section{Classification of the $\kappa$-solutions in two dimensions} \label{S:srv:classif_ksols}

The purpose of this section is proving that every two-dimensional $\kappa$-solution is actually a soliton. In the case of smooth surfaces,
the only possible soliton is the shrinking round sphere. This result is Corollary 11.3 in Perelman \cite{Perelman1}, but it was pointed out
by R. Ye \cite{Ye} that the argument had a small flaw because it appealed to Hamilton's work, that \emph{assumed} the compactness of the
surface to infere that the only soliton is the shrinking sphere. Ye proved that all two-dimensional asymptotic solitons are compact.
However, we have done in Chapter \ref{Ch:conesolitons} an exhaustive enumeration of all two-dimensional solitons on smooth and cone surfaces
and the
compactness of asymptotic solitons follows inmediately. Moreover, the arguments in Hamilton and Ye use the normalized Ricci flow, that is
not really suited when talking about $\kappa$-solutions. We rewrite their arguments in a way more suited to our discussion.
\medskip

Given an arbitrary Ricci flow, if it is not a (homothetic) soliton, there is no function $f$ such that the Ricci flow moves by homotheties
and the diffeomorphisms induced by $\grad f$, i.e. satisfying the equation
$$\Ric + \Hess f +\frac{1}{2t}g=0 .$$
We can, however, define a \emph{potential} $f$ for a Ricci flow as
the appropriate candidate.

\begin{prop}
Given a Ricci flow $(\mathcal M, g(t))$ with $\chi(\mathcal M)>0$, the equation
$$\Delta f = -\left( R+\frac{1}{t} \right)$$
has always a solution $f$, called the potential function of the flow.
\end{prop}

\emph{Remark}. Note that Hamilton's potential for the normalized Ricci flow is the solution of $\Delta f = R-r$ where $r$ is the average of
scalar curvature. This equation is the trace of the soliton equation for normalized Ricci flow \eqref{E:solitonNRF}.

The existence of $f$ is a consequence of the following general lemma

\begin{lema}
Let $\mathcal M^n$ be a smooth closed manifold. Then, for any smooth function $h$, the Poisson equation
$$\Delta f = h$$
has a unique solution $f$ with mean value zero, $\int_{\mathcal M} f \ d\mu=0$, if and only if 
$$\int_{\mathcal M} h\ d\mu = 0 .$$
\end{lema}
\begin{proof}[Sketch of the proof]

This can be expressed in the language of operators saying that the Laplacian operator has an almost inverse $G$, such
that
$$\Delta \circ G = G\circ \Delta = Id - \Pi$$
where $\Pi$ is a projection onto the kernel of $\Delta$ (i.e. the harmonic functions). In the case of a closed manifold, the harmonic
functions are constants and 
$$\Pi(f) = \frac{1}{\Vol (\mathcal M)} \int_{\mathcal M} h\ d\mu .$$
For one implication, if there is a solution of $\Delta f=h$, then by Stokes theorem
$$\int h =\int \Delta f = - \int \langle \nabla f , \nabla 1\rangle = - \int\langle \nabla f,0\rangle =0$$
and hence $\Pi(h)=0$. For the converse, if $\Pi(h)=0$, then $f=G(h)$ solves the equation,
$$\Delta f = \Delta G (h) = h - \Pi(h) =h .$$
If two functions satisfy the equation, their difference is a harmonic function, $\Delta (f_1-f_2)=h-h=0$, hence a constant. For any
solution $f$, we can choose the normalization $f-\Pi(f)$, that has mean value zero.

The general theory of linear operators brings the existence of $G$. As usual in the field, the proof involves first looking for weak
solutions on a Sobolev space (the operators acting on $H^k$, the space of $L^2$ functions with $k$ weak derivatives in $L^2$) using
linear algebra on infinite dimensions. Secondly, one proves that if the data is smooth, the solution is also smooth using elliptic
regularity. See for example \cite{Aubin}. 
\end{proof}

This approach is considered in \cite{MazRubSes}, where more refined spaces of functions are considered to extend the result to cone
surfaces. We will appeal to that result in Chapter \ref{Ch:cone_rf}. 

\medskip

\begin{proof}[Proof (Proposition)]
Since $\chi(\mathcal M)>0$, by Lemma \ref{L:RF_on_sphere_has_sing} there exists a singular time $0<T<\infty$. By reparameterizing
$t\rightarrow t-T$, we can
assume that the area collapses to zero at time $t=0$, and the flow is defined for $t<0$. Hence
$$\Area(\mathcal M)=-4\pi\chi(\mathcal M)\ t $$
Therefore, the integral
$$ \int_{\mathcal M} R+\frac{1}{t} d\mu = \int_{\mathcal M} R d\mu + \frac{1}{t} \Area \mathcal M 
= 4\pi\chi(\mathcal M) + \frac{1}{t} \Area \mathcal M = 0$$
vanishes, and therefore the equation has solution.
\end{proof}

\emph{Remark.} The potential can be defined with the same formula for surfaces with $\chi(\mathcal M)<0$, but then the time is defined for
$t>0$. For the case $\chi(\mathcal M)=0$, the potential should be defined as $\Delta f =-R$ and for $t\in \mathbb R$. In other words,
shrinking, expanding and steady solitons only can occur on surfaces with Euler characteristic positive, negative and
zero, respectively.

In the case of a soliton, this potential function actually solves the soliton equation. With Hamilton \cite{Hamilton_surfaces} in mind, we
make the following
definition,
\begin{defn}
Let $(\mathcal M, g(t))$ be a Ricci flow with $\chi(\mathcal M)>0$. Let $M$ be the quantity
$$M= \Hess f + \frac{1}{2} \left( R+\frac{1}{t} \right) g$$
which vanishes iff the Ricci flow is a gradient shrinking soliton.
\end{defn}

Remark that from the definition of $f$, $\trace M =0$, or equivalently, $M$ is the traceless Hessian of $f$. The norm of this quantity,
$|M|^2$, can be seen as a measure of how much our flow differs from a soliton.

Our aim is to use this scalar quantity in order to classify $\kappa$-solutions, since we know that a rescaled backwards-in-time limit is a
soliton. The idea is that the limit of $|M|^2$ as $t\rightarrow -\infty$ is $0$, because the limit flow as $t\rightarrow -\infty$ is a
soliton. If $|M|^2$ turned out to be decreasing in $t$, then $|M|^2$ would be identically $0$ and every $\kappa$-solution would be a
soliton. The evolution in time of this quantity could be our key tool, but $|M|^2$ is not scale invariant. We can solve this considering the
quantity $t^2|M|^2$. 

\begin{prop}
The evolution of $t^2|M|^2$ is
$$\frac{\partial}{\partial t} t^2 |M|^2 = \Delta (t^2 |M|^2) - 2 t^2 |\nabla M|^2 .$$
\end{prop}
\begin{proof}
It is basically a chain of computations, following the idea of \cite{Hamilton_surfaces}, see more detailed computations in \cite[p
217]{3CY}.
The main steps are the evolution of the following quantities:
\begin{itemize}
\item Metric (Ricci flow):
  $$\frac{\partial}{\partial t}g = -R g $$
\item Area element:
  $$\frac{\partial}{\partial t} d\mu = -R d\mu $$
\item Scalar curvature:
  $$\frac{\partial}{\partial t} R = \Delta R + R^2 $$
\item Potential function:
  $$\frac{\partial}{\partial t} f = \Delta f - \frac{1}{t}f + b $$
where $\Delta b=0$. \\
\item The $M$ quantity:
  $$\frac{\partial}{\partial t} M = \Delta M -2 RM -\frac{1}{t}M $$
\item The $|M|^2$ quantity:
  $$\frac{\partial}{\partial t} |M|^2 = \Delta |M|^2 - 2 |\nabla M|^2 -\frac{2}{t} |M|^2$$
\end{itemize}
The evolution of $t^2 |M|^2$ then follows.
\end{proof}

In particular, we have the inequality
$$\frac{\partial}{\partial t} t^2 |M|^2 \leq \Delta (t^2 |M|^2)$$
which allows us to apply the maximum principle. This is straightforward in the closed smooth case, but it will require a bit more 
caution on cone surfaces in Chapter \ref{Ch:cone_rf}. However, in both cases we will deduce that $\max_{\mathcal M } t^2|M|^2$ is decreasing
in $t$.


\begin{tma}
Every $\kappa$-solution of the Ricci flow on a cone surface is a gradient shrinking soliton.
\end{tma}
\begin{proof}
We use the quantity $t^2|M|^2$, which satisfies:
\begin{itemize}
 \item $\max_{\mathcal M} t^2|M|^2$ is decreasing on $t$.
 \item It is invariant under parabolic rescalings.
\end{itemize}
On the other hand, given a $\kappa$-solution we construct its asymptotic soliton by taking a sequence of times $\bar \tau_k \rightarrow
+\infty$; then picking appropriate $q_k$ such that $l(q_k,\bar \tau_k)\leq \frac{n}{2} =1$; and then constructing the sequence of rescalings
$$g_k(t)=\frac{1}{\bar\tau_k}g(\bar\tau_k t)$$
that subconverges to a shrinking soliton.

By monotonicity, taking $t_k=-\bar\tau_k$, $\forall x_0\in\mathcal M$ and $\forall t_0\in(-\infty,T)$,
$$t_0^2|M|^2_{g(t_0)} (x_0,t_0) \leq \max_{\mathcal M} t_k^2 |M|^2_{g(t_k)} (\cdot , t_k) \qquad \forall t_k<t_0 .$$
By the rescaling invariance, 
$$t_k^2 |M|^2_{g(t_k)}(\cdot,t_k) = |M|^2_{g_k(-1)} (\cdot,-1) .$$
Using that the limit of rescalings is a soliton, for all $ \epsilon>0$ there exists $k>0$ big enough ($t_k$ negative big enough) such that
$$|M|^2_{g_k(-1)} (\cdot,-1) <\epsilon .$$
Putting all together, for all $(x_0,t_0)$ and for all $\epsilon>0$ we conclude that
$$|M|^2_{g(t_0)} (x_0,t_0) <\epsilon ,$$
so $M\equiv 0$, which proves that the $\kappa$-solution is actually a soliton.
\end{proof}

From the classification of the solitons on surfaces on Chapter \ref{Ch:conesolitons}, the only smooth, complete, gradient shrinking soliton
on a surface is
the shrinking round sphere. Hence, we have the following corollaries.

\begin{corol}
The only $\kappa$-solution of the Ricci flow on a smooth surface is the shrinking round sphere.
\end{corol}

\begin{corol}
Let $(\mathcal M, g(t))$, $g(0)=g_0$, be a Ricci flow with $\chi(\mathcal M)>0$. Then $g(t)$ develops an infinite-curvature singularity for
some finite time $T$. After a blow-up rescaling, the limit of the rescalings is a shrinking round sphere, i.e. the original flow
converges to a round point at $t=T$. The normalized Ricci flow with the same initial data converges to a round sphere at $t=+\infty$.
\end{corol}

%% file: cone_rf.tex
\chapter{Cone surfaces evolving along Ricci flow} \label{Ch:cone_rf}

\lettrine{\indent C}{one singularities} arise usually in the study of two-dimensional orbifolds. Orbifolds are spaces locally modelled as
the
quotient of a manifold by the action by isometries of a discrete group. This identifies different directions as seen from a fixed
point. In the case of two dimensions, orientable orbifolds consist locally in
the quotient of a smooth surface by perhaps the action of a cyclic group acting by rotations, leading to the rise of singular cone points at
the center of the rotations. The space of directions on a cone point is no longer a metric circle of length $2\pi$ but a metric circle of
length $\frac{2\pi}{n}$ (this is the cone angle). General two-dimensional cone points include all angles, not only submultiples of $2\pi$.

 On the other hand, as we
saw in Chapter \ref{Ch:conesolitons}, cone singularities also arise naturally on the study of solitons. Some cone solitons (in the case of
compact surfaces) were found by Hamilton \cite{Hamilton_surfaces}, and this leaded to the work of Wu \cite{Wu} and Chow \cite{ChowWu} that
proved the convergence of the normalized flow on bad orbifolds (the teardrop and the football) to the soliton metrics. That gave a
uniformization of all compact two-orbifolds. However, despite this convergence result, some existence and collapsing issues remained
unclear. 

The Ricci flow is equivariant under isometries of the manifold. Therefore, if a manifold admits a quotient by a group, this is not altered by the effect of the flow. It seems natural, therefore, to assume an equivariant definition of the flow on orbifolds. Furthermore, the Ricci flow on smooth surfaces is \emph{conformal}, that is, preserves angles. It is therefore natural to assume that the flow also preserves the cone angles at the cone points, and hence the evolution of the metric occurs only in the smooth part of the surface, keeping the cone points unaltered. These assumptions define an \emph{angle preserving flow}. However, these are not the only possible assumptions as we will discuss, and we will see a different alternative in Chapter \ref{Ch:smoothening}.

General cone surfaces admit cone points not risen from the quotient of a disc by rotations, thus any cone angle may appear (not only
$\frac{2\pi}{n}$). Although the intrinsic geometry is essentially equivalent, the analysis used for the orbifolds does not apply. The
problem requires a different point of view, based on local models for the cone points. The existence issues for the flow, even assuming an
angle-preserving flow, remained unclear until the recent works of H. Yin \cite{Yin2} and R. Mazzeo, Y. Rubinstein and N. Sesum
\cite{MazRubSes}. In Section \ref{S:conerf:cones_and_RF} we discuss in detail the possible definitions of cone surface and the existence
results for the flow on
such surfaces. 

The other unclear issue for the flow on cone surfaces is the behaviour of the singular set along the flow. Cone points might a priori
collapse all together on a single point, or an infinite-curvature singularity might get rise to a new cone point. We will rule out those
cases in Section \ref{S:conerf:infcur_sing} using the additional metric assumption that the cone angles are less than or equal to $\pi$. 

In the remainder of the chapter, our goal will be to develop a post-Perelman argument to prove the uniformization of cone surfaces,
recovering the results of Hamilton-Chow-Wu from a different point of view. This proof will apply to both smooth and cone (or orbifold)
surfaces. The main point to adapt Perelman's arguments, as described in Chapter \ref{Ch:survey}, to cone surfaces is the development of
(barrier) maximum principles and Harnack inequalities that can work despite the cone singularities. We carry on this work in Section
\ref{S:conerf:barrier_maxpple}, using barrier techniques inspired by the work of T. Jeffres \cite{Jeffres}. Finally in Section
\ref{S:conerf:unif_conesurf} we assemble all the argument
to prove the uniformization theorem for cone surfaces with angles less than or equal to $\pi$.

\section{Cone surfaces and Ricci flow} \label{S:conerf:cones_and_RF}

Cone surfaces are topological surfaces equipped with a Riemannian metric which is smooth everywhere except on some discrete set of points
(\emph{cone points}) that look like the vertex of a cone. To formalize this, we demand that the metric in a neighbourhood of the cone
point adopts a specific canonical form, in specific coordinates. We will see three coordinate systems which are suitable to describe a cone
point: polar geodesic coordinates, conformal coordinates, and conformal coordinates with respect to a fixed cone. The regularity of the
functions describing the metric on each coordinates is crucial to the geometric meaning and to the existence of the Ricci flow. We will
therefore discuss three different definitions of cone surface and the existence of a suitable Ricci flow on these surfaces.

In Chapter \ref{Ch:conesolitons}, we found naturally some singular metrics on surfaces that one must agree on qualify them as ``cone
surfaces''. This comes from the canonical representation of a metric on geodesic polar coordinates
\begin{equation} \label{E:g_gpolar_coord}
 g=d\rho^2 + h^2 d\theta^2
\end{equation}
where $h=h(\rho,\theta)$. There, $\rho$ is the arc-parameter of the geodesic curves $\{\theta=cst\}$ (the meridians), and $2\pi h$ measures
the length of the curves $\{\rho=cst\}$ (the parallels).

\begin{defn}\label{defconesurf1}
A cone surface $(\mathcal M, (p_1,\ldots,p_n),g)$ is a topological surface $\mathcal M$ and points $p_1,\ldots, p_n\in \mathcal M$ equipped
with a smooth Riemannian metric $g$ on $\mathcal M\setminus \{p_1,\ldots, p_n\}$, such that every point $p_i$ admits an open
neighbourhood $U_i$, and an homeomorphism $(U_i, p_i) \rightarrow (D,0)$ (which is a diffeomorphism on $U_i\setminus \{p_i\}$), where the
metric on the
coordinates of $D\setminus\{0\}$ is written as
$$g=d\rho^2 + h^2 d\theta^2$$
with $h=h(\rho,\theta)$ a smooth function $h: D \rightarrow \mathbb R$, satisfying
$$h(0)=0, \qquad \frac{\partial h}{\partial \rho}(0)=\frac{\alpha_i}{2\pi},\qquad \frac{\partial^{2k} h}{\partial \rho^{2k}}(0)=0 $$
for some $\alpha_i\in(0,2\pi]$ (the cone angles). Here, $\rho>0$ is the arclength parameter measuring the distance to the singular point,
and $\theta\in [0,2\pi]/\sim$ is proportional to the
angle with respect to a fixed arbitrary half-line.
\end{defn}

All the solitons found in Chapter \ref{Ch:conesolitons} are cone surfaces in this sense. The simplest example of a cone surface is the flat
cone. The
smooth metric $d\rho^2 + r^2 d\theta^2$, with $\rho\in[0,+\infty)$, $\theta\in \mathbb R / 2\pi\mathbb Z = [0,2\pi]/\sim$ is the standard
Euclidean metric in polar coordinates. If we pick two half-lines in the plane, emanating from the origin at
angle $\alpha$, and we identify these two half lines, then the region between these lines (the region we choose), constitutes a cone of
angle $\alpha$. The metric on this surface is the same  $d\rho^2 + \rho^2 d\xi^2$, but now the angle is $\xi\in \mathbb R / \alpha\mathbb Z
= [0,\alpha]/\sim$. This is the same metric as 
$$ d\rho^2 + \left(\frac{\alpha}{2\pi}\right)^2 \rho^2 d\theta^2$$
with $\theta\in
[0,2\pi]/\sim$. 

In dimension two, all smooth metrics are locally conformally equivalent, \cite{Chern}. This provides another canonical expression of any
metric
as conformal to the Euclidean metric, these are the so-called \emph{conformal} or \emph{isothermal coordinates}, 
\begin{equation} \label{E:g_conf_coord}
  g=e^{2u}(dr^2 + r^2 d\theta^2 ).
\end{equation}
where $u:U\subset\mathbb R^2\rightarrow \mathbb R$ is the \emph{conformal factor}. 

We can express some examples of cone metrics on this conformal form, an this will motivate another definition of cone surface. This second
definition will be slightly more general than the previous one.

The conformal expression for the singular metric of the flat cone can be derived from simple geometric arguments. Consider the metric
space
resulting of identification of two half lines in the plane meeting at angle $\alpha$. Assume that there exist a Riemannian
metric of the form
$$g=\phi^2(r)\ (dr^2+r^2 d\theta^2)$$
where $\phi$ could depend also on $\theta$, but we assume that $\phi=\phi(r)$ and that there is a rotational symmetry. We look for a
function $\phi$ that produces a
flat metric and a cone angle at the origin. Let us consider a circle given by the curve $r=x=const$. On the
one hand, its length is the angle times the radius (since it is a region of the Euclidean plane),
$$L=\alpha \int_0^x \phi(r) dr.$$
On the other hand, the length of the curve measured on the metric is
$$L=\int_0^{2\pi} x \phi(x) d\theta = 2\pi x \phi(x).$$
So
$$\alpha \int_0^x \phi(r) dr = 2\pi x \phi(x) ,$$
and denoting $\Phi'(r)=\phi(r)$ we obtain
$$\alpha \Phi(x) = 2\pi x \Phi'(x) .$$
Solving this ODE, $\Phi(r) = e^C r^{\frac{\alpha}{2\pi}}$. For simplify, we can fix $C=0$ (this is just a scale factor on the metric).
We get
$$\phi(r) = \Phi'(r) = \frac{\alpha}{2\pi} r^{\frac{\alpha}{2\pi}-1} .$$
Renaming 
$$\beta =\frac{\alpha}{2\pi}-1 ,$$
with $-1<\beta\leq 0$, the metric of the Euclidean cone is
\begin{align*}
 g &= \left(  \frac{\alpha}{2\pi r} r^{ \frac{\alpha}{2\pi}}   \right)^2 (dr^2+r^2d\theta^2) \\
   &= e^{\displaystyle  2\left(   \ln (\beta+1) +  \beta \ln r \right)    } (dr^2+r^2d\theta^2) .
\end{align*}

Analogously, we can find other constant curvature examples identifying a sector of the sphere or the hyperbolic space. We of course have
the metrics in the form \eqref{E:g_gpolar_coord}, for $\xi\in[0,\alpha]/\sim$,
$$d\rho^2 + \sin^2 \rho \ d\xi^2$$
for the sphere, and
$$d\rho^2 + \sinh^2 \rho \ d\xi^2$$
for the hyperbolic space. Let us write these metrics in the conformal form.

For the sphere with the metric in the form \eqref{E:g_conf_coord}, the length of a circle given by $r=x=const$ (and radius $l$) is
$$L=\alpha \sin l = \alpha \sin \int_0^x \phi(r) dr ,$$
and this equals on the metric to
$$L=\int_0^{2\pi} x \phi(x) d\theta = 2\pi x \phi(x).$$
Hence one has the ODE
$$\alpha \sin \Phi(x) = 2\pi x \Phi'(x) ,$$
for $\Phi'(x)=\phi(x)$. This ODE can be solved (e.g. using the change $\Phi=2\arctan\eta$ and some trigonometric identities) to yield the
metric of the spherical cone,
\begin{align*}
 g &= \left(  \frac{\alpha}{2\pi r} \frac{1}{\cosh\left(\frac{\alpha}{2\pi}\ln r\right)}  \right)^2 (dr^2+r^2d\theta^2) \\
   &= e^{\displaystyle  2\left(   \ln (2(\beta+1)) -\ln\left(1+r^{2(\beta+1)}\right) +  \beta \ln r \right)    } (dr^2+r^2d\theta^2) 
\end{align*}
where $\beta =\frac{\alpha}{2\pi}-1$. Remark that when $\alpha=2\pi$, the metric turns into
$$g = \left(\frac{2}{1+r^2}\right)^2 (dr^2+r^2d\theta^2) $$
which is the spherical metric of the hemisphere, parameterized by the equatorial disc in polar coordinates.

Finally, the hyperbolic case is entirely analogous, raising the ODE
$$\alpha \sinh \Phi(x) = 2\pi x \Phi'(x) ,$$
that yields the metric of the hyperbolic cone,
\begin{align*}
 g &= \left(  \frac{\alpha}{2\pi r} \frac{1}{\sinh\left(\frac{\alpha}{2\pi}\ln r\right)}  \right)^2 (dr^2+r^2d\theta^2) \\
   &= e^{\displaystyle  2\left(   \ln (2(\beta+1)) -\ln\left(1-r^{2(\beta+1)}\right) +  \beta \ln r \right)    } (dr^2+r^2d\theta^2) .
\end{align*}
where $\beta =\frac{\alpha}{2\pi}-1$. Remark that when $\alpha=2\pi$, the metric turns into
$$g = \left(\frac{2}{r^2-1}\right)^2 (dr^2+r^2d\theta^2) $$
which is the hyperbolic metric of the Poincaré disc, parameterized by polar coordinates on the disc. Also remark that if we take the limit
as $\alpha\rightarrow 0$, the metric turns into
$$g = \left(\frac{1}{r \ln(r)}\right)^2 (dr^2 + r^2 d\theta^2)$$
which is the metric of the hyperbolic cusp, or pseudosphere.

These examples suggest that the conformal factor of a cone metric has an asymptote near $r=0$ that grows as a negative logarithm, and
therefore a cone metric can be written in conformal coordinates as 
$$g=e^{2(a+\beta \ln r)} (dr^2 + r^2 d\theta^2)$$
for some $a:U\subset \mathbb R^2 \rightarrow \mathbb R$ bounded and continuous. Thus, the cone angle is codified in the $\beta \ln r$
($\beta=\frac{\alpha}{2\pi}-1 \in(-1,0]$) and the curvature (in the Riemannian sense) for neighbour points is encoded in the $a$ function.

Having seen the prototype of cone point, the following definition is justified:
\begin{defn}\label{defconesurf2}
A cone surface $(\mathcal M, (p_1,\ldots,p_n),g)$ is a topological surface $\mathcal M$ and points $p_1,\ldots, p_n\in \mathcal M$ equipped
with a smooth Riemannian metric $g$ on $\mathcal M\setminus \{p_1,\ldots, p_n\}$, such that every point $p_i$ admits an open
neighbourhood $U_i$, and an homeomorphism $(U_i, p_i) \rightarrow (D,0)$ (which is a diffeomorphism on $U_i\setminus \{p_i\}$), where the
metric on the
coordinates of $D\setminus\{0\}$ is written as
$$g=e^{2(a_i+\beta_i \ln r)} (dr^2+r^2 d\theta^2) = e^{2a_i}\left( ds^2 + \left( \frac{\alpha_i}{2\pi}\right)^2 s^2 d\eta^2\right)$$
where $a_i:D \rightarrow \mathbb R$ is a bounded and continuous function on the whole disc, and
$\beta_i=\frac{\alpha_i}{2\pi}-1$ for some given $0<\alpha_i\leq2\pi$.
\end{defn}

\emph{Remark.} The two conformal ways to write a metric are conformal coordinates with respect to the Euclidean smooth metric or with
respect to the Euclidean cone metric. The equivalence of the two systems can be made explicit with the change of variables
\begin{align*}
s &=\frac{2\pi}{\alpha_i} r^{\frac{\alpha_i}{2\pi}} ,\\
\eta &= \theta.
\end{align*}

The \emph{cone angle} at $p_i$ is $\alpha_i=2\pi(\beta_i+1)$. We say that $\mathcal M$ has bounded curvature if it has bounded
Riemannian curvature on the smooth part of $ \mathcal M$, although it has $+\infty$ curvature in the sense of Alexandrov at the cone
points (provided the angle is less than $2\pi$). Note that in the case of $\alpha=2\pi$ ($\beta=0$), the metric needs not to be smooth
across the cone point, unless we demand additional hypothesis on the function $a$.

This indeed defines a smooth metric away from $r=0$, and has a cone structure at $r=0$. However, if we require additional hypothesis such
as bounded Riemannian curvature, further restrictions will apply to the
function $a$. Besides, although $a$ is a smooth function away from the cone point, this regularity cannot be demanded on all $\mathcal M$
without losing most significant examples. Indeed, in the hyperbolic cone the function 
$$a(r)= \ln(2(\beta+1)) - \ln\left(1-r^{2(\beta+1)}\right)$$
is bounded and continuous for small $r$, and has derivative 
$$a'(r)= \frac{2(\beta+1)r^{2\beta+1} }{1-r^{2(\beta+1)}}.$$
This derivative does not exist at $r=0$ if $\beta<-\frac{1}{2}$; hence, the function $a(r)$ might not be $\mathcal C^1$ at the
singularity. In this case the function $a(r)\sim r^{2(\beta)}\in \mathcal C^{0,\gamma}$ is Hölder continuous for $\gamma \leq
2(\beta+1)$ on the whole domain, although this Hölder regularity depends on the cone angle.

\medskip

We intend to apply a Ricci flow on cone surfaces. It is reasonable to think that cone surfaces in the sense of Definition \ref{defconesurf1}
admit a
solution to the Ricci flow, whereas the Definition \ref{defconesurf2} may be too loose. We will use an existence theorem from
\cite{MazRubSes} for the Ricci flow on a certain setting of cone manifolds that we will take as a third definition.

Strictly speaking, a cone surface is a smooth surface minus a discrete set of points, together with a specific prescribed behaviour of
the metric tensor on a neighbourhood of the removed (cone) points. Thus, the Ricci flow we try to consider is defined over an open surface,
and actually it could be an ill-posed problem or possess non-unique solutions. In this chapter, we will use a particular flavour of
the Ricci flow which keeps the cone angles fixed through time (this approach is equivalent to the one followed by Chow and Wu in the case of
orbifolds). The Ricci flow is always a conformal change of the metric on the smooth part
of the surface, i.e. maintains the angles; but keeping the cone angles unchanged is an extra assumption.
 This approach is not the
only possible. In Chapter \ref{Ch:smoothening} we will develop a different flavour of the Ricci flow that instantaneously removes the cone
points and the
manifold becomes smooth everywhere. Even more, one can have an angle-changing flow \cite{MazRubSes} or a cusp-generating flow \cite{Topgie}.

For the angle-preserving flow there are two main references that consider the existence of such Ricci flow: on the one hand the work of H.
Yin
\cite{Yin1}, \cite{Yin2}; and on the other hand the work of R. Mazzeo, Y. Rubinstein and N. Sesum \cite{MazRubSes}. Both works follow the
standard approach for proving the existence of
a PDE. First, one writes the Ricci flow equation on certain conformal form such that there is only one function evolving with an appropriate
equation. Second, one writes the equation as a fixed-point problem for certain operator $\mathcal T$ acting on a certain function space.
Then one finds the linearized part of
this operator and shows that $\mathcal T$ is a contraction and hence has a fixed point. Thus, one has to study and get estimates for the
linear problem $\partial_t u = a(x,t)\Delta_g u + \langle V,\nabla u \rangle + f(x,t)$. The technical and delicate part is choosing the
appropriate function spaces to set the initial conditions and the solutions to these equations.

We will use Mazzeo, Rubinstein and Sesum approach. Their work uses the b-calculus introduced by R. B. Melrose \cite{Melrose}, that gives an
explicit description of the asymptotic behaviour of the functions near the cone point. The idea is to
``desingularize'' the point by a blow-up, substituting the cone point by the $S^1$ boundary $\{r=0,\theta\in\mathbb R/2\pi\mathbb Z\}$ one
obtains a manifold with boundary $ \tilde{\mathcal M}$. In polar coordinates, one changes the space of derivative operators to $\mathcal V_b
= \{ r\frac{\partial}{\partial r} , \frac{\partial}{\partial \theta} \}$ instead of the usual $\{ \frac{\partial}{\partial r} ,
\frac{\partial}{\partial \theta} \}$. This is equivalent to restrict to vector fields on $\tilde {\mathcal M}$ tangent to the boundary. Then
it is
constructed $\mathcal C_b^{k,\delta}(\tilde {\mathcal M})$ as the space of functions with $k$ derivatives (taken in $\mathcal V_b$), and
after taking all the derivatives, the result is on a Hölder space $\mathcal C^{0,\delta}(\tilde {\mathcal M})$. Similarly, the space
$\mathcal C_b^{k+\delta,(k+\delta)/2}([0,T]\times \tilde {\mathcal M})$ is the space of functions of space and time, with $i$
space-derivatives (taken in $\mathcal V_b$),  $j$ time-derivatives (in the usual sense), $i+2j\leq k$, and after taking all the derivatives,
the result is on a Hölder space $\mathcal C^{\delta,\delta/2}$. Finally, the Hölder-Friedrichs domain for these functions is
$$\mathcal D_b^{k+\delta,(k+\delta)/2}([0,T]\times \tilde {\mathcal M})=\{ u\in C_b^{k+\delta,(k+\delta)/2} | \Delta_g u \in
C_b^{k+\delta,(k+\delta)/2} \}.$$

The Ricci flow equation on surfaces is 
$$\frac{\partial}{\partial t} g(t) = -R g .$$
If $g(t)=e^{2\phi(t,x)}g_0$, then the flow becomes
\begin{equation}
 \frac{\partial}{\partial t}\phi(t) = e^{-2\phi}\Delta_{g_0} \phi + R_{g_0} . \label{eqRF_conftobkg}
\end{equation}

\begin{tma}[ {\cite[Prop 3.10]{MazRubSes}} ] \label{T:ExistenceMRS}
Let $(\mathcal M, (p_1,\ldots,p_n),g_0)$ be a cone surface. Let $\bar g_0$ be a background reference metric such that on a neighbourhood of
a cone point $p_i$ is flat and takes the conic form
$\bar g_0 = dr^2 + \left( \frac{\alpha}{2\pi}\right) ^2 r^2 d\theta^2$.
Let $g_0=e^{2 \phi_0}\bar g_0$ be a cone metric, where $\phi_0\in \mathcal D_b^{k,\delta}(\tilde{\mathcal M})$, and let $g(t)=e^{2
\phi(t)}g_0$. Then, there is
a unique solution $\phi\in\mathcal D^{k+\delta,(k+\delta)/2}([0,T]\times \tilde {\mathcal M})$ to \eqref{eqRF_conftobkg} with
$\phi|_{t=0}=0$, provided $T$ is sufficiently small.
\end{tma}

Consequently, unless stated otherwise, we will use henceforth the following definition of cone surface.
\begin{defn}\label{defconesurf3}
A cone surface $(\mathcal M, (p_1,\ldots,p_n),g)$ is a topological surface $\mathcal M$ and points $p_1,\ldots, p_n\in \mathcal M$ equipped
with a smooth Riemannian metric $g$ on $\mathcal M\setminus \{p_1,\ldots, p_n\}$, such that every point $p_i$ admits an open
neighbourhood $U_i$, and an homeomorphism $(U_i, p_i) \rightarrow (D,0)$ (which is a diffeomorphism on $U_i\setminus \{p_i\}$), where the
metric on the
coordinates of $D\setminus\{0\}$ is written as
$$g = e^{2 \phi_0} \left( dr^2 + \left( \frac{\alpha}{2\pi}\right) ^2 r^2 d\theta^2 \right).$$
for certain function $\phi_0\in \mathcal D_b^{0,\delta}(\tilde{\mathcal M})$.
\end{defn}

\vspace{3em}

%

Any cone surface in the sense of Definition \ref{defconesurf3} also satisfies Definition \ref{defconesurf2} (which is too general for the
flow).
Now we check that a cone surface defined with smooth polar coordinates (Definition \ref{defconesurf1}) satisfies the boundary requirements
of the b-calculus definition, for which we have existence of the flow (Definition \ref{defconesurf3}).

\begin{lema}
A cone surface in the sense of Definition \ref{defconesurf1} is also a cone surface in the sense of Definition \ref{defconesurf3}.
\end{lema}

\begin{proof}
We express a cone metric in two coordinate charts; namely conformal coordinates with respect to a cone, and polar
geodesic coordinates.
$$ g = e^{2u} (dr^2+ \left(\frac{\alpha}{2\pi}\right)^2 r^2\ d\theta^2) = d\rho^2 + h^2\ d\xi^2 .$$
For simplicity, we will assume that the functions $u,h$ are radial, i.e. $u=u(r)$, $h=h(\rho)$. The general case $u=u(r,\theta)$,
$h=h(\rho,\xi)$ follows the same structure, only involving more terms on $\frac{\partial u}{\partial \theta}$ and $\frac{\partial
h}{\partial \xi}$. 

The Gaussian curvature can be expressed as
$$K= -\Delta_g = -e^{-2u}\Delta u = -\frac{1}{h} \frac{\partial^2 h}{\partial \rho^2}$$
where $\Delta = \frac{\partial^2}{\partial r^2} + \frac{1}{r}\frac{\partial}{\partial r} + \frac{1}{r^2}\frac{\partial^2}{\partial
\theta^2} $.
The change of coordinates can be achieved by the transformation
\begin{align*}
d\rho &= e^u dr \\
\xi &= \theta \\ 
h & = \left(\frac{\alpha}{2\pi}\right) r e^u 
\end{align*}
Assume we have a metric in the sense of Definition \ref{defconesurf1}, given by $h(\rho)\in \mathcal C^\infty([0,A))$. We need to check that
$$u(r)=\ln h -\ln r -\ln \left(\frac{\alpha}{2\pi}\right)$$ 
and $\Delta_g u$ belong to $\mathcal C_b^k$ for any $k$, i.e. when applying $k$ times the operator $r\frac{\partial}{\partial r}$ the
result is in $\mathcal C^0$. Indeed,
\begin{align}
r\frac{\partial}{\partial r}  u &= r\frac{\partial}{\partial r}\left(\ln h - \ln r -\ln \frac{\alpha}{2\pi}\right) \\
&= r\frac{1}{h}\frac{\partial h}{\partial \rho}\frac{\partial \rho}{\partial r} -1 \\
&= r e^u \frac{1}{h} \frac{\partial h}{\partial \rho} -1 \\
&= \frac{2\pi}{\alpha} \frac{\partial h}{\partial \rho} -1
\end{align}
which is $\mathcal C^0$ and tends to $0$ as $r\rightarrow 0$. Since the derivative operator is
$$r\frac{\partial}{\partial r} = r \frac{\partial \rho}{\partial r} \frac{\partial}{\partial \rho} = r e^u \frac{\partial}{\partial \rho} =
\left(\frac{2\pi}{\alpha}\right) h \frac{\partial}{\partial \rho} ,$$
when applied further, results on functions with the same regularity as $h$ and the derivatives of $h$.

For the Laplacian,
\begin{align}
r\frac{\partial}{\partial r}   \Delta_g u = -  r\frac{\partial}{\partial r} K &= 
\frac{2\pi}{\alpha} h \frac{\partial}{\partial \rho } \left( \frac{1}{h} \frac{\partial^2 h }{\partial \rho^2} \right) \\
      &= \frac{2\pi}{\alpha} \left( \frac{\partial^3 h}{\partial \rho^3} + K \frac{\partial h}{\partial \rho} \right) .
\end{align}
Hence, all further derivations respect to $h\frac{\partial}{\partial \rho}$ will yield terms in $K$ and derivatives of $h$. 
Since in Definition \ref{defconesurf1} the function $h=\frac{\alpha}{2\pi}\rho + O(\rho^3)$, the curvature has a limit as $\rho\rightarrow
0$,
$$\lim_{\rho\rightarrow 0} K = \lim_{\rho\rightarrow 0} \frac{-h_{\rho\rho}}{h} = \lim_{\rho\rightarrow 0}
\frac{O(\rho)}{\frac{\alpha}{2\pi}\rho + O(\rho^3)} = C ,$$
and hence $K$ is a $\mathcal C^0$ function. Therefore, $\Delta_g u$ is in $\mathcal C_b^k$.
\end{proof}
\medskip

To end this section, we show the Ricci flow equation in (time-dependent) polar coordinates, and we show a consequence of this equation for
the angle-preserving flow.

Recall that the Ricci flow in surfaces, in conformal coordinates $g=e^{2u}(dr^2+r^2d\theta^2)$, is written as
$$\frac{\partial u}{\partial t} = e^{-2u}\Delta u = -K $$

\begin{prop} \label{P:RF_in_gpolarcoords}
The Ricci flow equation, expressed in the time-dependent geodesic polar coordinates
$$g(\tau)=d\rho^2 + h(\rho,\xi,\tau)^2 d\xi^2$$
is
$$\frac{\partial h}{\partial \tau} = \frac{\partial^2 h}{\partial \rho^2} - \frac{\partial h}{\partial \rho} \int{\frac{1}{h}
\frac{\partial^2 h}{\partial \rho^2} \ d\rho}$$
\end{prop}
\emph{Remark.} Here $\tau=t$ is the time coordinate. On the change of coordinates, $\rho$ is independent of $\tau$, but $\rho$ may depend
on $t$ in the coordinates $(r,\theta,t)$.

\begin{proof}
Let $(r,\theta,t)$ be the fixed spacetime coordinate chart on which a Ricci flow written as $g(t)=e^{2u}(dr^2+r^2d\theta^2)$ satisfies the
equation 
\begin{equation} \label{E:L_RF_in_polar:eq0}
\frac{\partial u}{\partial t} = e^{-2u}\Delta u = -K 
\end{equation}

for $u=u(r,\theta,t)$. Let
$(\rho,\xi,\tau)$ be the spacetime coordinate chart defined by
\begin{align*}
d\rho &= e^u\ dr \\
\xi &= \theta \\
\tau &= t
\end{align*}
First equation is equivalent to the integral form $\rho=\rho(r,\theta,t)=\int_0^r e^{u(s,\theta,t)} \ ds$.
From the chain rule we have
$$0=\frac{\partial \rho}{\partial \tau} = 
\frac{\partial \rho}{\partial r} \frac{\partial r}{\partial \tau} + \frac{\partial \rho}{\partial \theta} \frac{\partial \theta}{\partial
\tau} + \frac{\partial \rho}{\partial t} \frac{\partial t}{\partial \tau} $$
that is,
\begin{equation} \label{E:L_RF_in_polar:eq1}
 e^u \frac{\partial r}{\partial \tau} + \frac{\partial \rho}{\partial t} = 0 .
\end{equation}
We define
$$h=re^u .$$
The change of variables yield
\begin{align} \label{E:L_RF_in_polar:eq2}
 \frac{\partial h}{\partial \rho} 
 &= \frac{\partial r}{\partial \rho} e^u + r e^u \left( \frac{\partial u}{\partial r} \frac{\partial r}{\partial \rho} + 
 \frac{\partial u}{\partial \theta} \frac{\partial \theta}{\partial \rho} + \frac{\partial u}{\partial t} \frac{\partial t}{\partial \rho}
 \right)  \nonumber \\
 &= 1 + r \frac{\partial u}{\partial r} ,
\end{align}
and also
\begin{align} \label{E:L_RF_in_polar:eq3}
 \frac{\partial h}{\partial \tau} 
 &= \frac{\partial r}{\partial \tau} e^u + r e^u \left( \frac{\partial u}{\partial r} \frac{\partial r}{\partial \tau} + 
 \frac{\partial u}{\partial \theta} \frac{\partial \theta}{\partial \tau} + \frac{\partial u}{\partial t} \frac{\partial t}{\partial \tau}
 \right) \nonumber \\
 &= e^u \left( \frac{\partial r}{\partial \tau} \left( 1 + r \frac{\partial u}{\partial r} \right) + r\frac{\partial u}{\partial t}\right) .
\end{align}
Hence, using \eqref{E:L_RF_in_polar:eq1} and \eqref{E:L_RF_in_polar:eq1} into \eqref{E:L_RF_in_polar:eq3}, we get
\begin{align*}
 \frac{\partial h}{\partial \tau}
 &= - \frac{\partial \rho}{\partial t} \frac{\partial h}{\partial \rho} + h \frac{\partial u}{\partial t}
\end{align*}
Using now the Ricci flow equation \eqref{E:L_RF_in_polar:eq0},
$$\frac{\partial u}{\partial t} = -K = - \frac{1}{h} \frac{\partial^2 h}{\partial \rho^2}$$
and
$$\frac{\partial \rho}{\partial t} = \int{e^u \frac{\partial u}{\partial t} \ dr} = \int{\frac{\partial u}{\partial t} \ d\rho} = 
- \int K \ d\rho ,$$
so
$$ 
 \frac{\partial h}{\partial \tau} = \left( \int K \ d\rho \right) \frac{\partial h}{\partial \rho}   +  \frac{\partial^2 h}{\partial \rho^2}
$$
as stated.

\end{proof}

\emph{Remark.} Theorem \ref{T:ExistenceMRS} can be stated as follows. Let $(\mathcal M,(p_1,\ldots,p_n),g_0)$ be a cone surface such
that in a local coordinate chart the metric can be written
$$g_0=d\rho^2 + h_0^2 d\xi^2.$$
Then, the integral PDE problem on that local chart
\begin{equation}
 \left\{ \begin{array}{ll}
          h_{\tau}=h_{\rho\rho}-h_\rho \int{ \frac{h_{\rho\rho}}{h}} & \\
	  h = 0 &\mbox{ on } \rho=0 \\
	  h=h_0(\rho,\xi) &\mbox{ on } \tau=0 \\
         \end{array}
\right.
\end{equation}
has a short-time solution $h(\rho,\xi,\tau)$ for $\tau\in[0,\epsilon)$ compatible on all local charts, and the cone angle at the
cone point does not change over time,
$$\frac{\partial}{\partial \tau} \frac{\partial h}{\partial \rho} =0 .$$

In \cite{MazRubSes} it is described the existence of an angle-changing flow, where the cone angle at a cone point can be described by
any smooth function $\alpha(t)$. This would make our choice of angle-preserving flow somehow more arbitrary. However, the angle-changing
flow has some drawbacks, namely the unbounded curvature along the flow. We show next that the angle preserving flow is the unique flow that
evolves amongst the class of cone surfaces and has uniformly bounded curvature. In particular, it is the only
one that may keep the cone metrics fitting the Definition \ref{defconesurf1}, not escaping to a more general metric in the sense of
Definition \ref{defconesurf3}.

\begin{prop} \label{prop_fixedangle}
Let $(\mathcal M,p,g(t))$ be a Ricci flow with one cone point $p$ with angle $\alpha(t)$ that might not be constant (as existence asserted
in \cite{MazRubSes}). If the curvature is uniformly bounded, then the cone angle does not change along the time. 
\end{prop}

\begin{proof}

From Corollary 3.12 in \cite{MazRubSes}, the curvature along the flow has the form
$$K \sim b_0(t) + r^{\frac{2\pi}{\alpha}} B(t,\theta) + O(r^2)$$
for some functions $b_0(t)$, $B(t,\theta)$ that do not depend on $r$. Hence,
$$ \frac{\partial K}{\partial r} =  \frac{2\pi}{\alpha} r^{\frac{2\pi}{\alpha} -1} B(t,\theta) + O(r)$$
and 
$$ r \frac{\partial K}{\partial r} =  \frac{2\pi}{\alpha} r^{\frac{2\pi}{\alpha}} B(t,\theta) + O(r^2) .$$
Since $\alpha<2\pi$, we get that $r\frac{\partial K}{\partial r}\rightarrow 0 $ as $r\rightarrow 0$ for all $t>0$. 

Now, working on the geodesic polar coordinates $(\rho,\xi,\tau)$, we note that
$$r \frac{\partial K}{\partial r} = r \frac{\partial \rho}{\partial r}\frac{\partial K}{\partial \rho} =  r e^u\frac{\partial K}{\partial
\rho} = h \frac{\partial K}{\partial \rho}$$
and hence $h \frac{\partial K}{\partial \rho} \rightarrow 0 $ as $\rho\rightarrow 0$.
We also note that on these coordinates,
$$ \frac{\partial K}{\partial \rho} = -\frac{1}{h}\left(h_{\rho \rho \rho} - \frac{h_{\rho \rho}}{h} h_\rho \right) .$$
Then, we compute the time derivative of the angle on the $(\rho,\xi,\tau)$ coordinates,
\begin{align*}
 \frac{\partial}{\partial \tau} \frac{\partial h}{\partial \rho} &= 
 h_{\rho \rho \rho} + K h_\rho + h_{\rho \rho} \int K \ d\rho \\
 &= h_{\rho \rho \rho} - \frac{h_{\rho \rho}}{h} h_\rho + h_{\rho \rho} \int K \ d\rho \\
 &= h\left( -\frac{\partial K}{\partial \rho} - K \int K \ d\rho \right) . 
\end{align*}
Hence, since $K$ is bounded and $ h \frac{\partial K}{\partial \rho}\rightarrow 0$ as $\rho\rightarrow 0$, we have
$$ \frac{\partial}{\partial \tau} \frac{\partial h}{\partial \rho} \Bigg|_{\rho=0} =0$$
which means that the cone angle does not change.
\end{proof}

\section{Infinite-curvature singularities} \label{S:conerf:infcur_sing}
In the smooth case, the only obstruction to the continuation of the flow is the explosion of the curvature at some point, as we saw in
Chapter \ref{Ch:survey}. To prove this, one picks a flow
defined on $t\in[0,T)$ and assumes uniformly bounded curvature. One selects a sequence $t_i\rightarrow T$ and then the metrics $g(t_i)$ are
equivalent to $g(0)$. By a compactness
theorem, we get a limit as $t_i\rightarrow T$ and get a metric $g(T)$ that can serve as initial data for a continuation of the
flow.

In the case of cone surfaces, the same result applies. However, the compactness theorem must be examined. A priori, other phenomena
associated with the cone points might prevent a continuation of the flow, such as two cone points collapsing close together, or a limit of
cone points with certain angle that converge to a different cone angle.

Imagine an Euclidean cone of angle $\alpha$ truncated at a distance $\epsilon$ of the vertex, and replaced with a small Euclidean cone of
angle $\tilde \alpha \neq \alpha$, and imagine that the joint of the two cones is smoothed with a partition of unity, thus giving a small
nonflat annulus with some curvature around the vertex. These are cone surfaces of angle $\tilde\alpha$. Now we take a sequence of these
cones varying $\epsilon \rightarrow 0$. The Gromov-Hausdorff limit is a surface with cone angle $\alpha$, thus we may lose the cone
structure at the limit. This, however, is due to the lack of a uniform bound on the curvature on the small annulus, the curvature on the
joint tends to $\pm \infty$ as $\epsilon\rightarrow 0$ (depending on whether $\alpha>\tilde \alpha$ or the opposite).
This same phenomenon may be observed as the asymptotic cone of the $\alpha\tilde\alpha$-cone solitons depicted in Chapter
\ref{Ch:conesolitons}. These are
surfaces with cone angle $\tilde \alpha$ and have as asymptotic cone a surface with cone angle $\alpha$.

Fortunately, in Appendix \ref{Ch:compactness} we prove that, provided uniform bounds on the injectivity radius and the curvature, the cone
structure is
preserved under Gromov-Hausdorff limits. The uniform bound for the injectivity radius on $[0,t_i)$ comes automatically from the finite
distortion of the metric for finite time (and hence finite distortion of the distances). If we assume bounded curvature on $[0,t_i)$, then
we can apply Gromov's compactness theorem to get a G-H limit, this limit is smooth except on the limit of the cone points, and the cone
angle of the limit is the same as in $g(0)$. Therefore, the continuous extension is not altered by the cone points, and the
smooth extension applies (locally) to any neighbourhood of any smooth point. Let us remark that no maximum principle is required for the
continuous extension of the metric.
We have proven the following result.
\begin{prop}
Let $\mathcal M$ be a smooth closed cone surface, $\Sigma\subset \mathcal M$ a discrete set of cone points, and $g(t)$ a Ricci flow on a
maximal time interval $[0,T)$ and $T<\infty$, then
$$\sup_{\mathcal M\setminus \Sigma} |R|(\cdot,t) \rightarrow \infty$$
as $t\rightarrow T $.
\end{prop}

Now we turn our attention to the noncollapsing property.
In the case of cone surfaces, the argument with $\mathcal L$-geodesics is still valid because of the following fact: Since
the cone angles we consider are less than $2\pi$, geodesics (as a length-minimizing path)
do not pass across any cone point, because otherwise it would be possible to shorten the path. 
It is well possible to
construct a piecewise smooth path passing through a cone
point satisfying the geodesic equation ($\nabla_{\dot \gamma}\dot \gamma =0$) on each
smooth component, this is a particular case of a broken geodesic, and does not minimize the distance between its endpoints. The same
occurs with $\mathcal L$-geodesics: a spacetime path, parameterized
backwards in time and minimizing the $\mathcal L$-functional amongst its endpoints, does not touch
any cone point in any time, because otherwise the $\mathcal L$-length could be decreased.
Hence, the reduced volume of a cone surface is defined the same way as in the smooth case, up to a set of singular cone points that have
zero measure.

Since the noncollapsing theorem is independent ot Harnack inequalities and maximum principles, we have directly the same noncollapsing
theorem for cone surfaces.
Noncollapsing has an important consequence for surfaces with cone angles less or equal than $\pi$. 

\begin{tma}
Let $(\mathcal M,(p_1,\ldots,p_n),g(t))$ be an angle-preserving Ricci flow on a cone surface defined on $[0,T]$, such that $|K|<C$ for all
$t\in[0,T]$. Assume that that all the cone angles are less than or equal to $\pi$. Then, the injectivity radius of the cone points
$$\inj(\Sigma) =\min_{p\in\{p_1,\ldots,p_n\}} \inj(p) $$
is uniformly bounded below along the flow. In particular, the distance between any two cone points is uniformly bounded below.
\end{tma}

\begin{proof}
This is a consequence of $\kappa$-noncollapse and Lemma \ref{L:cptss:inj_cone} in Appendix \ref{Ch:compactness}. By the $\kappa$-noncollapse
and the bounds
on the
curvature, there is a lower bound on the volume and hence, by Proposition \ref{P:srv:inj_vol}, a lower bound on the injectivity radius, and
this bounds
are uniform in time as far as the curvature keeps uniformly bounded. Then, by Lemma \ref{L:cptss:inj_cone}, the cone points keep a uniformly
bounded
distance between them (the injectivity radius of the cone points is uniformly bounded below). 
\end{proof}

\section{Barrier maximum principles for cone manifolds} \label{S:conerf:barrier_maxpple}

We now develop some maximum principles for cone manifolds. First we develop a generic maximum principle for functions on cone surfaces.
Next, we show an ad hoc maximum principle that suits the tensor Harnack inequality for flows on a cone surface.

One of the most standard formulations of the maximum principle is the following.

\begin{tma}
Let $U$ be an open bounded set in $\mathbb R^n$. Let $U_T=U\times (0,T]$ and let $u:U\times [0,T]\rightarrow \mathbb R$ be a function $u\in
\mathcal C^{2,1}(U_T)\cap \mathcal C(\overline{U_T})$ such that 
$$\frac{\partial u}{\partial t}\leq \Delta u .$$
Let $(x_0,t_0)\in \overline{U_T}$ realize the maximum of $u$,
$$u(x_0,t_0)=\max_{\overline{U_T}}u ,$$
then
$$(x_0,t_0)\in \overline{U_T}\setminus U_T ,$$
i.e. either $t_0=0$ or $x_0\in \partial U$.
\end{tma}

The notation $\mathcal C^{2,1}(U\times(0,T])$ means functions of class $\mathcal C^2(U)$ on the space variable, and of class
$\mathcal C^1((0,T])$ on the time
variable, and with bounded $\mathcal C^2$-norm,
$$\sup_U |u| + \sup_U |\nabla u| + \sup_U |\nabla^2 u| < \infty .$$

Heuristically, an interior maximum would have gradient zero and Hessian negative defined. So it would have negative Laplacian and the
function on that point would be decreasing when fixed on that point, giving greater values backwards in time. Therefore the space maximum
cannot increase in time, and hence the maximum is at $t=0$, or at the boundary. 

This argument fails if $u$ is not at least $\mathcal C^2$ on the interior of $U$, in particular if $U$ itself contains a cone
point, because a nonsmooth maximum point no longer needs to have zero gradient or negative Laplacian. We can workaround this problem if we
are able to find a way to guarantee that the maximum cannot occur on the cone point. The way for achieving this is constructing a new
function $\bar u_\epsilon$ depending on a parameter $\epsilon>0$ such that
\begin{itemize}
 \item $\bar u_\epsilon$ tends uniformly to $u$ when $\epsilon\rightarrow 0$
 \item on a small neighbourhood, $\bar u_\epsilon$ is strictly increasing over radial lines leaving the cone point (thus not having a
maximum at the cone
point).
\end{itemize}
We call this new $\bar u_\epsilon$ a \emph{barrier function} for $u$. 

In order to construct barriers, we begin by setting some estimates for a nonflat metric on a cone
surface, and next bringing some auxiliary functions that will help us to build such barriers. Some of these computations have been done in
Chapter \ref{Ch:conesolitons} for radial metrics.

\begin{lema} \label{L:properties_g_gpolarcoords}
Let $(\mathcal M,\{p_1,\ldots,p_n\},g)$ be a cone surface such that the metric on a neighbourhood of a cone point is written as
$$g=dr^2+h(r,\theta)^2 d\theta^2$$
for some analytic $h:[0,A)\times \mathbb R/2\pi\mathbb Z \rightarrow \mathbb R^+$, such that
\begin{itemize}
 \item $h(0,\theta)=0$, $\forall \theta$ (the neighbourhood is a disc).
 \item $\frac{\partial h}{\partial r}(0, \theta)=\frac{\alpha}{2\pi}$ (where $\alpha$ is the cone angle).
 \item $\left| \frac{\partial^2 h}{\partial r^2} \right| \leq C h$ (bounded curvature).
\end{itemize}
Then,
\begin{itemize}
 \item The Gaussian curvature is
$$K=-\frac{1}{h}\frac{\partial^2 h}{\partial r^2} .$$
\item The Hessian of a radial function $f=f(r)$ is
$$\Hess f = f_{rr} dr^2 + h h_r f_r d\theta^2$$
\item The Laplacian of a radial function $f=f(r)$ is
$$\Delta f = f_{rr} + \frac{h_r}{h} f_r $$
\item For a fixed $\theta$,
$$h(r,\theta) = \frac{\alpha}{2\pi} r +O(r^3) .$$
\item The ratio $h_r/h$ behaves as
$$r\frac{h_r}{h}-1 = O(r^2) .$$
In particular, $\frac{h_r}{h}\sim \frac{1}{r}$ as $r\rightarrow 0$.
\end{itemize}
\end{lema}

\begin{proof}
First three statements are a straightforward computation. From the bound on the curvature, we have that $|h_{rr}|\leq C|h|$, hence
$h_{rr}(0)=0$ and the Taylor expansion of $h$ has only a linear term on $r$ plus terms of order $O(r^3)$. For the last statement, 
$$r\frac{h_r}{h}-1 = \frac{rh_r -h}{h} = \frac{r\frac{\alpha}{2\pi}r+O(r^3) - r\frac{\alpha}{2\pi}r+O(r^3)}{\frac{\alpha}{2\pi}r+O(r^3)} =
\frac{O(r^3)}{\frac{\alpha}{2\pi}r+O(r^3)} = O(r^2) .$$
\end{proof}

Now we define a helpful function (cf. \cite{Jeffres}) that we will use later to build the barriers.

\begin{lema} \label{L:Jeffresbarrier}
Let $U$ be a topological disk, with given polar coordinates $(r,\theta)\in (0,r_0) \times [0,2\pi)$, a
cone angle at the origin, and a smooth Riemannian metric outside
the cone point with bounded curvature. Let $0<\delta<1$. Then the function given by
$$(r,\theta)\mapsto r^\delta$$
satisfies
\begin{itemize}
 \item $\grad r^\delta$ is pointing away from the cone point, and with norm tending to $+ \infty$ as we approach the vertex.
 \item $\Delta r^\delta >0$ if $r$ small enough.
\end{itemize}
\end{lema}

\begin{proof}
The gradient vector is 
$$\grad r^\delta = \delta r^{\delta -1} \frac{\partial}{\partial r}$$
so it is clear that it points away from the origin and its norm tends to $\infty$ as $r\rightarrow 0$.
The Laplacian of $r^\delta$ is
 $$\Delta r^\delta = \delta (\delta -1)r^{\delta -2} + \frac{1}{h}\frac{\partial h}{\partial r} \delta r^{\delta-1}
 = \delta r^{\delta -2} \left( \delta -1 + \frac{1}{h}\frac{\partial h}{\partial r}r \right) .$$
Since $r\frac{1}{h}\frac{\partial h}{\partial r} \longrightarrow 1$  as $r\rightarrow 0$, then $\Delta r^\delta >0$ for $r$ small enough.
\end{proof}

This lemma allows us to construct a barrier function that proves the maximum principle on closed cone surfaces:

\begin{tma} \label{T:maxpple_cone}
Let $(\mathcal M,(p_1,\ldots,p_n),g_0)$ be a closed cone surface, and let $u\in C^{2,1}(\mathcal M\times (0,T], g_0)$
such that 
$$\frac{\partial u}{\partial t}\leq \Delta u$$
Let $(x_0,t_0)\in \mathcal M \times [0,T]$ such that realizes the maximum of $u$ over space and time,
$$u(x_0,t_0)=\max_{\mathcal M\times [0,T]}u$$
then $t_0 =0$.
\end{tma}

The notation $C^{2,1}(\mathcal M\times (0,T], g_0)$ means functions $\mathcal C^2$ in space and $\mathcal C^1$ in time with bounded
$\mathcal C^2$-norm, this norm taken with respect to the metric $g_0$.

\begin{proof}
Applying the maximum principle over the open set $\mathcal M \setminus \Sigma$, where we denote $\Sigma =\{p_1,\ldots,p_n\}$, the maximum of
$u$ is
achieved on $t=0$ or, maybe, on
$t>0$ and $p\in \Sigma$. We will rule out the latter case. Assume by contradiction that $(p,t_0)$, $p\in \Sigma$, is the maximum of $u$
over $\mathcal M \times [0,T]$. 

Let $U$ be a small neighbourhood of $p$ such that we can dispose polar coordinates $(r,\theta)$, and $\Delta r^\delta>0$ for
some $0<\delta<1$, by Lemma \ref{L:Jeffresbarrier}. Let $\epsilon>0$, and define over $U$ the function
$$\bar u = u+\epsilon r^\delta .$$
It satisfies
$$\frac{\partial \bar u}{\partial t}= \frac{\partial u}{\partial t} \leq \Delta u \leq \Delta (u+\epsilon r^\delta) = \Delta \bar u .$$

Applying the maximum principle to the open set $U\setminus\{p\}$, $\max_{\bar U\times [0,T]} \bar u$ lies on $t=0$ or on $x\in \partial U
\cup \{p\}$. We claim that the latter cannot happen. Indeed, $\bar u$ cannot have a maximum on $x=p$ (i.e. $r=0$) because $\grad \bar u$
is pointing away from $p$ with infinite norm when $r=0$, and $\grad u|$ is bounded, so $\bar u$ is strictly increasing on radial directions
leaving $p$. On the other hand, the original $u$ has no maxima on $\partial U \times (0,T]$ because they would be interior points in
$\mathcal M\setminus \Sigma \times [0,T]$. Since $\bar u \rightarrow u $ uniformly as $\epsilon \rightarrow 0$, $\bar u $ cannot either have
maxima on $\partial U \times
(0,T]$; specifically, for any $\epsilon < \epsilon_0=\frac{1}{2} \left( \max_{\bar U\times (0,T]} u - \max_{\partial U \times (0,T]} u
\right) $, the function $\bar u$ cannot have maxima on $\partial U$ because this value would be at most $\max_{\partial U} u +\epsilon$
that is less than $u(p,t_0)$. 

Therefore, $\max \bar u$ is on $t=0$ and again since $\bar u \rightarrow u$ uniformly, $\max u$ is on $t=0$.
\end{proof}

Now we look for a cone version of Harnack inequality for Ricci flows. Recall from Chapter \ref{Ch:survey} that for the smooth case, if
$(\mathcal
M,g(t))$ is a Ricci
flow with nonnegative curvature operator, then the quantity 
$$ Z = M_{ij} W^i W^j +2P_{kij} U^{ki}W^j + R_{ijkl}U^{ij}U^{kl}$$
is nonnegative. This is proven using a maximum (minimum) principle that involves creating a barrier function for the spatial infinity. We
just modify this barrier to be a barrier also at the cone points. To modify $Z$, $M_{ij}$ and $R_{ijkl}$, the barrier relies on the
functions $\psi$ and $\phi$ as in Lemma \ref{barriers_Hamilton}.

The functions in the smooth case are:
$$\varphi=\epsilon e^{At}f(x), \quad \psi = \delta e^{Bt}$$
with $\epsilon$, $\delta$ small and $A$ $B$ sufficiently large. The function $f(x)$ depends only on the position, $f(x)\rightarrow +\infty$
as $x$ goes
to $\infty$ (the distance to a fixed basepoint tends to infinity), but the derivatives of $f$ are bounded. Then, $\varphi$ is a space
barrier
for the infinity and a time barrier for $t=0$. The function $\psi$ is only a time barrier.

The only point we need is to change $\varphi$ to be a barrier also at the cone points.

\begin{lema}
Let $(\mathcal M,(p_1,\ldots,p_n),g)$ be a cone surface. There is some $C>0$ and a function
$\mu=\mu(x)$ satisfying
\begin{enumerate}
 \item $\mu\geq 1/C$
 \item $\mu \rightarrow +\infty$ as $x$ tends to a cone point.
 \item $\Delta \mu \leq C$
\end{enumerate}
\end{lema}

\begin{proof}
On a local chart around $p_i$, we can assume that $g=dr^2+h^2d\theta^2$ for $r\in[0,r_0)$, for some $r_0$ uniform on the surface.
Without loss of generality, we can assume $r_0=1$.

It suffices to use a smooth interpolation between $\mu=-\ln r$ for $r<\frac{1}{2}$ and $\mu=\ln 2$ for $r>\frac{1}{2}$ (assume that the
interpolation only affects a very small neigbourhood of $r=\frac{1}{2}$.
This function $\mu$ obviously satisfies (1) and (2). To see (3), we only need to check it for the case $\mu=-\ln r$ for $r$ small. This
gives us, by Lemma \ref{L:properties_g_gpolarcoords},
$$\Delta \mu = \mu_{rr} + \frac{h_r}{h} \mu_r = \frac{1}{r^2} \left( 1 - r \frac{h_r}{h} \right) = \frac{1}{r^2}O(r^2) = O(1)$$
and hence $\Delta \mu$ is bounded on $[0,\frac{1}{2})$.
Finally, glue all the functions defined on neigbourhoods of the cone points, and define $\mu = \ln 2$ outside these neighbourhoods.
\end{proof}

Now we construct the new barriers.
\begin{lema} \label{barriers_Hamilton_cone}
For any $C$, $\eta>0$ and any compact set $K$ in space-time not containing cone points, the functions $\psi=\psi(t)$ and
$\tilde\varphi=\tilde\varphi(x,t)$ defined as
$$\tilde\varphi=\epsilon e^{At} ( f(x) + \mu(x)) .$$
$$\psi = \delta e^{Bt}$$
satisfy
\begin{enumerate}
 \item $\delta\leq\psi\leq\eta$ for some $\delta>0$, for all $t$;
 \item $\epsilon\leq\tilde\varphi\leq\eta$ on the compact $K$ for some $\epsilon>0$, for all $t$. Furthermore,
$\tilde\varphi(x,t)\rightarrow\infty$ if
$x\rightarrow\infty$ or $x\rightarrow \Sigma = \{p_1,\ldots,p_n\}$, i.e. the sets $\{x\ |\ \tilde\varphi(x,t)<M\}$ are compact for all $t$
and all
$M$;
 \item $\frac{\partial \tilde\varphi}{\partial t} > \Delta \tilde\varphi + C \tilde\varphi$;
 \item $\frac{\partial \psi}{\partial t} > C \psi$;
 \item $\tilde\varphi\geq C \psi$.
\end{enumerate}
\end{lema}

\begin{proof}
We have defined
$$\tilde \varphi = \varphi + \epsilon e^{At}\mu$$
where $\varphi$ is the function on the smooth case on Lemma \ref{barriers_Hamilton}. Thus, items 1 and 4 have not changed. Item 2 follows
from the fact that $\mu\rightarrow+\infty$ as $x$ tends to a cone point. Item 5 is immediate, $C\psi\leq\varphi\leq \tilde\varphi$.

We check item 3:
$$\left( \frac{\partial}{\partial t}-\Delta\right) \tilde \varphi = \left( \frac{\partial}{\partial t}-\Delta\right) \varphi + \epsilon 
e^{At} \left( A\mu - \Delta \mu \right) .$$
Since $\Delta \mu \leq C' \leq (C')^2\mu$, we have $A\mu-\Delta \mu\geq (A-(C')^2)\mu \geq C'' \mu$ if $A$ is big enough. Hence,
$$\left( \frac{\partial}{\partial t}-\Delta\right) \tilde \varphi > C\varphi + C'' \epsilon e^{At}\mu > C''' \tilde \varphi ,$$
for possibly different constants $C$'s.
\end{proof}

\section{Uniformization of cone surfaces} \label{S:conerf:unif_conesurf}
We finally assemble the properties obtained on the previous sections to reconstruct a proof for the uniformization of certain
cone surfaces, as done with the smooth case.

We restrict ourselves to surfaces with cone angles less than or equal to $\pi$. This in particular covers the case of all orbifolds,
with cone angles $\frac{2\pi}{n}$ for $n\in\mathbb N$. We will see that this restriction is an important hypothesis, and that other
phenomena might happen if dropped.

From \cite{MazRubSes} we have the existence of an angle-preserving flow. There are other angle-changing flows but, as we saw, this is the
only that may keep the curvature bounded for a short time. Even this flow can develop infinite-curvature singularities at finite
time, as it happened in the smooth case, but also analogously to the smooth case, this is the only obstruction to the continuation of the
flow.

If the curvature explodes to infinity for finite time, we can perform a sequence of pointed parabolic rescalings. The $\kappa$-noncollapsing
property, that works on the cone setting with the angle restriction, allows us to
get a pointed limit flow, with the same cone structure as the original flow, that is a $\kappa$-solution.

Since the Harnack inequality holds for flows on cone surfaces, every $\kappa$-solution has an asymptotic shrinking soliton, complete and
with bounded curvature. We classified all cone solitons in Chapter \ref{Ch:conesolitons}, and all the possible solitons are compact, namely
the teardrop
and the football solitons, or the constant curvature solitons.

Since we have a maximum principle for functions on cone surfaces, we can apply it on $\kappa$-solutions to the function $u= t^2|M|^2$, over
$\mathcal M$ and nested compact time intervals $[t_1,t_2]$ with $t_1 \rightarrow -\infty$, as in the smooth case, and obtain that $|M|=0$.
This proves that every $\kappa$-solution on a cone surface is a soliton.

\begin{tma}
Let $(\mathcal M, (p_1,\ldots,p_n), g(t))$ be a $\kappa$-solution over a cone surface with cone angles less than or equal to
$\pi$. Then it is a shrinking soliton.
\end{tma}

Since all $\kappa$ solutions are therefore compact, this gives a strong restriction on which kind of surfaces can develop an
infinite-curvature singularity.

\begin{corol}
Let $(\mathcal M, (p_1,\ldots,p_n), g(t))$ be an angle-preserving Ricci flow over a cone surface with cone angles less than or
equal to $\pi$. Then, the flow is either defined for all time; or the flow converges, up to rescaling, to a round surface; or it is a sphere
with one or two cone points and converges, up to rescaling, to a teardrop or football soliton.
\end{corol}

For the general picture of the flow, it is useful to keep the analogy with the smooth case in terms of the Euler characteristic. Recall that
for cone surfaces there is a suitable modified definition of Euler characteristic,
$$\hat\chi(\mathcal M) = \chi(\mathcal M) + \sum_{i=1}^n \beta_i$$
where $\chi(\mathcal M)$ is the Euler characteristic of the underlying topological surface, and $\beta_i = \frac{\alpha_i}{2\pi}-1$ are the
angle parameters of the cone points. In the case of orbifolds, this definition makes the conic Euler characteristic multiplicative with
respect to branched coverings (i.e. if $\tilde {\mathcal M} \rightarrow \mathcal M$ is an $n$-to-one branched covering, then
$\hat\chi(\tilde{\mathcal M}) = n \hat\chi(\mathcal M)$). Furthermore, with this definition the Gauss-Bonnet formula holds,
$$\int_{\mathcal M} K \ d\mu = 2\pi\hat\chi(\mathcal M) .$$

Recall also that the evolution of the area of the surface under the Ricci flow is 
$$\frac{\partial}{\partial t} \Area(\mathcal M) = \int_{\mathcal M} \frac{\partial}{\partial t} d\mu = \int_{\mathcal M}
R d\mu = -4\pi\hat\chi(\mathcal M) .$$

Let us inspect the case $\hat\chi(\mathcal M)\leq 0$. Then the area does not tend to zero, and there are also no isolated infinite-curvature
singularities; since the rescaling blow-up of such singularities would bring as a limit a noncompact $\kappa$-solution, which is
impossible. Therefore, the flow is defined for all $t>0$. Actually, in this case the most easy way to study the flow is to use the
normalized Ricci flow, as in the smooth case explained in Chapter \ref{Ch:survey}. Substituting the maximum principle with the cone maximum
principle from Section \ref{S:conerf:barrier_maxpple}, the
result is the same, and the surface converges to a constant curvature cone metric.

Now, let us check the case $\hat\chi(\mathcal M)> 0$. Then the topological (orientable) surface must be a sphere, $\chi(\mathcal M)=2$, and
the cone angles must satisfy
$$\sum_{i=1}^n \beta_i > -2 $$ 
If we assume the angles less than or equal to $\pi$, then $-1<\beta_i\leq \frac{-1}{2}$, and hence at most three cone points can occur. The
area tends to zero for some finite time, and there must be an infinite-curvature singularity, which, as before, cannot be isolated and must
happen at the same time as the area collapses to zero. The rescaling blow up of the surface yields a $\kappa$-solution, which we have seen
it must be a soliton, that must be spherical, the teardrop or the football soliton, according to the number and magnitude of the cone
points.

This gives a uniformization of all closed cone surfaces with angles less than or equal to $\pi$, and in particular, a uniformization of all
closed two-dimensional orbifolds.

\begin{tma}
Let $(\mathcal M,(p_1,\ldots,p_n),g_0)$ be a closed cone surface (in the sense of Definition \ref{defconesurf3}), and assume that the cone
points
are less than or equal to $\pi$. Then there is an angle-preserving Ricci flow that converges, up to rescaling, to:
\begin{itemize}
 \item a constant nonpositive curvature metric, if $\hat\chi(\mathcal M)\leq 0$.
 \item a spherical (constant positive curvature) metric, a teardrop soliton or a football soliton; if $\hat\chi(\mathcal M)\leq 0$.
\end{itemize}
\end{tma}

It is interesting to note that this result agrees with the following obstruction theorem due to M. Troyanov. 
\begin{tma}[Troyanov, \cite{Troyanov}]
Let $\mathcal M,(p_1,\ldots,p_n), (\alpha_1,\ldots,\alpha_n)$ be a surface with a collection of points and angles. Then there exists a
constant curvature cone metric $g$ such that $(\mathcal M, (p_1,\ldots,p_n),g)$ is a cone surface (in the sense of Definition
\ref{defconesurf2}) if and only if one of the following holds:
\begin{itemize}
 \item $\hat\chi(\mathcal M)\leq 0$.
 \item $\hat\chi(\mathcal M)>0$ and for each $i=1..n$
$$\beta_i>\sum_{j\neq i} \beta_j $$
where $\beta_i = \frac{\alpha_i}{2\pi}-1$.
\end{itemize}
\end{tma}

In our discussion, if some cone points are greater than $\pi$, then the nonlocal collapsing is not enough to guarantee that two cone points
stay at a
uniformly bounded distance, i.e. two cone points could approach each other asymptotically while maintaining bounded curvature and area on
the surface. This phenomenon has been also observed and confirmed by Mazzeo, Rubinstein and Sesum \cite{MazRubSes}. They observed that if
the inequality for the angles in Troyanov's theorem fails, then there is only one cone point, say $p_1$, that violates the inequality.
They proved that the normalized Ricci flow on the sphere with three or more cone points, some with angle greater than $\pi$, and not
admitting a constant curvature metric,
evolves approaching a football soliton with only two cone points, when all but the cone point $p_1$ collapse all together in a single cone
point at infinite time. At the moment of writing this thesis, the details of this convergence and the limit remain to be fully described
yet.

%% file: smoothening.tex
\chapter{Smoothening cone points on surfaces} \label{Ch:smoothening}

\lettrine[nindent=-3pt]{\indent T}{he Ricci flow} considered so far in the previous chapters is the so-called \emph{angle-preserving flow}.
This is
the flow implicitly assumed by Wu and Chow (\cite{Wu}, \cite{ChowWu}) for orbifolds, which is an equivariant definition of the Ricci flow
under the action of the isotropy group of the cone points. Yin \cite{Yin2} and Mazzeo, Rubinstein, and Sesum \cite{MazRubSes} generalized
the existence theorem to cone surfaces.

The alternative consideration of the Ricci flow just acting on the smooth part of the surface leads to consider the Ricci flow on an
open, noncomplete manifold, which does not fit in the classical theory of Hamilton, so existence and uniqueness might be lost. On that general setting, P. Topping and G. Giesen \cite{Topgie} obtained an existence theorem for Ricci flow on
incomplete surfaces, which becomes instantaneously complete but might have unbounded curvature from below. This exposes the nonuniqueness of
solutions. 

In another work, Topping \cite{Topping_revcusp} considered a complete open surface with cusps of negative curvature and proved the
existence of a instantaneously smooth Ricci flow with unbounded curvature, a \emph{smoothening flow for cusps} which erases instantaneously
the cusps.
This requires a generalized notion of initial metric for a flow, namely, the flow is defined for $t\in(0,T]$ and the initial metric is
the limit of $g(t)$ (in some sense) as $t\rightarrow 0^+$. 

Topping's technique for this result consists in capping off the cusps of the original metric $g_0$ with a smooth part near the cusp point,
in an
increasing sequence of metrics, each term with a further and smaller cap. This sequence of smooth metrics gives rise to a sequence of
(classical) Ricci flows, and the work consists in proving that this sequence has a limiting Ricci flow on $\mathcal M$ which has $g_0$ as
initial condition in that generalised sense. Our work in this chapter proves that this technique applies equally
well on cone surfaces, using
truncated or ``blunt'' cones as approximations for a cone point. In our setting, cusps would be seen as a limiting case of a zero-angle
cone. This provides an instantaneously smooth Ricci flow that smooths out the cone points of a cone surface.

This technique of approximating a nonsmooth manifold by smooth ones and limiting the sequence of corresponding Ricci flows has also been
used by M. Simon \cite{Simon} in dimension three to investigate the Gromov-Hausdorff limit of sequences of three-manifolds. Another
similar development has been done by T. Richard \cite{Richard} using this technique for smoothening out a broader class of Alexandrov
surfaces.

The chapter is organized as follows.
In Section \ref{sectiondifferentRF} we recall some properties for cone surfaces and we compare the three different flows from evolving from
a cone surface. Then we proceed to the construction of the smoothening flow. In Section
\ref{sectiontruncating} we build the truncated cones that will serve us as
approximations of a cone point. In Section \ref{sectionbarriers} we build upper barriers that, applied to our truncated cones, will give
us control on the
convergence of the sequence. In Section \ref{sectionexistproof} we put together the preceding results to prove the existence theorem.
Finally in Section \ref{sectionuniqueproof}
we prove the uniqueness theorem.

\section{Cone points and different Ricci flows} \label{sectiondifferentRF}

In Section \ref{S:conerf:cones_and_RF} of Chapter \ref{Ch:cone_rf} we made an extensive discussion of different ways to define and describe
a cone surface. We recall only a few facts that we will need here.

In this chapter we will use mainly isothermal coordinates (with respect to the Euclidean metric) for the model of cone points. A smooth
Riemnnian surface admits, on a neighbourhood $U$ of each point $p$, a coordinate chart $(U,p)\rightarrow (D,0)$ (where $D$ is the unit disc)
such that the metric is written as
$$g=e^{2u}(dx^2 + dy^2)=e^{2u}(dr^2 + r^2 d\theta^2)=e^{2u}|dz|^2$$
where $(x,y)$ are the cartesian coordinates (or in polar coordinates $(r,\theta)$, or in complex notation $z=x+iy$); and
$u:D\rightarrow \mathbb R$ is smooth. If $p$ is not smooth but a cone point, then $u:D\setminus \{0\} \rightarrow \mathbb R$ is smooth and
has subtle regularity and asymptotic behaviour as $r\rightarrow 0$. As we saw in Chapter \ref{Ch:cone_rf}, the geometry of a cone requires
$u$ to have an asymptote at $r=0$ growing like a negative multiple of the logarithm of $r$.

We will use hence the broader of the three definitions of cone surface that we discussed in Chapter \ref{Ch:cone_rf}.
\begin{defn}\label{defconesurf}
A cone surface $(\mathcal M, (p_1,\ldots,p_n),g)$ is a topological surface $\mathcal M$ and points $p_1 , \ldots , p_n \in \mathcal M$
equipped with a smooth Riemannian metric $g$ on $\mathcal M\setminus \{p_1,\ldots, p_n\}$, such that every point $p_i$ admits an open
neighbourhood $U_i$, and diffeomorphism $(U_i, p_i) \rightarrow (D,0)$ (which is a diffeomorphism on $U_i\setminus \{p_i\}$),where the
metric on the
coordinates of $D\setminus\{0\}$ is written as
$$g=e^{2(a_i+\beta_i \ln r)} |dz|^2$$
where $a_i:D \rightarrow \mathbb R$ is a bounded and continuous function on the whole disc, and $-1<\beta_i\leq0$. The \emph{cone angle} at
$p_i$ is $\alpha_i:=2\pi(\beta_i+1)$.
\end{defn}
We will assume that all cone surfaces in this chapter have bounded curvature. Recall that we say that $\mathcal M$ has bounded curvature if
it has bounded Riemannian curvature on the smooth part of $ \mathcal M$. Recall also that the Ricci flow on surfaces is 
$$\frac{\partial}{\partial t}g = -2 K g$$
where $K$ is the Gauss curvature of the surface, and in isothermal coordinates $g=e^{2u}|dz|^2$ the equation becomes
$$\frac{\partial}{\partial t}u = e^{-2u}\Delta u = -K$$
where $\Delta$ is the usual Euclidean Laplacian.

We will see now three different Ricci flows that can evolve from a cone surface as in Definition \ref{defconesurf}. First is the
\emph{angle-preserving flow} from Mazzeo, Rubinstein and Sesum \cite{MazRubSes} (cf. Yin \cite{Yin1}, \cite{Yin2}), that we saw in Chapter
\ref{Ch:cone_rf}.

 \begin{tma}[cf. \cite{MazRubSes}, cf. \cite{Yin2}] \label{Yinthm}
Let $(\mathcal M, (p_1,\ldots,p_n),g_0)$ be a cone surface in the sense of Definition \ref{defconesurf3} (in particular, satisfying
Definition \ref{defconesurf}). Then there exists some $T>0$ and a
solution $g(t)$ to the Ricci flow on $t\in [0,T]$, with $g(0)=g_0$, such that on a neighbourhood of each cone point $p_i$ the (cone) metric
is written as
$g(t)=e^{2(a_i(t)+\beta_i \ln r)} |dz|^2$.
\end{tma}

Note that in \cite{MazRubSes} it is described not only the angle-preserving flow, but also an angle-changing flow for any given smooth
functions $\alpha_i(t)$ that prescribe the cone angles at the cone points $p_i$. Thus, it would be more precise to talk about the family of 
\emph{cone-preserving flows}. This shows that there are much more than three possible flows.

In the more general setting of smooth, non-complete, open surfaces, Topping and Giesen \cite{Topgie} proved that there exists a solution
which is complete for all $t>0$ and is \emph{maximally stretched}, meaning that points spread apart at the maximum possible speed.

\begin{tma}[\cite{Topgie}] \label{Topgiethm}
Let $(\mathcal M, g_0)$ be a smooth, non complete Riemannian surface without boundary. There exists a $T>0$ and a smooth Ricci flow $g(t)$
on $t\in [0,T]$ such that $g(0)=g_0$; $g(t)$ is complete for all $t>0$; and $g(t)$ is maximally stretched, i.e. if $\tilde g(t)$ is any
other Ricci flow on $M$ with $\tilde g(0)\leq g(0)$, then $\tilde g(t)\leq g(t)$ for all $t$.
\end{tma}

In the particular case of a cone surface, viewed as an open noncomplete surface, this theorem ensures that there exists a solution flow
such that after any time $t>0$ a neighbourhood of the puncture has become as the (narrow) end of a hyperbolic cusp. This provides a
\emph{maximally stretched flow}, which is a second type of flow.

In this chapter we construct a third different solution to the
same equation: a Ricci flow which is instantaneously smooth for any $t>0$, that satisfies the Ricci equation for any $t\in(0,T]$, and that
converges to the initial nonsmooth cone metric as $t\rightarrow 0$. We call this a \emph{smoothening flow}. Despite the nonuniqueness shown
by these results, there is certain
uniqueness provided we restrict to a certain class of flows.

We will use the following general definiton of initial metric:
\begin{defn}[Cf. \cite{Topping_revcusp} Definition 1.1]\label{def1}
Let $\mathcal M$ be a smooth manifold, and $p_1,\ldots,p_n\in\mathcal M$. Let $g_0$ be a Riemannian metric on $\mathcal
M\setminus\{p_1,\ldots,p_n\}$ and let $g(t)$ be a smooth Ricci flow on $\mathcal M$ for $t\in(0,T]$. We say that $g(t)$ has initial
condition $g_0$ if 
$$g(t) \longrightarrow g_0\ \mathrm{as}\ t\rightarrow 0^+$$
smoothly locally on $\mathcal M\setminus\{p_1,\ldots,p_n\}$.
\end{defn}

The two main theorems of the chapter are the following:
\begin{tma}\label{tma1}
Let $(\mathcal M,(p_1\ldots p_n),g_0)$ be a closed cone surface; with bounded curvature. There exists a Ricci flow $g(t)$ smooth on the
whole $\mathcal M$, defined for $t\in (0,T]$ for some $T$, and such that
$$g(t)\underset{t\rightarrow 0}{\longrightarrow} g_0.$$
Furthermore, this Ricci flow has curvature unbounded above and uniformly bounded below over time.
\end{tma}

\begin{tma}\label{tma2}
Let $\tilde g(t)$ be a Ricci flow on $\mathcal M$, defined for $t\in (0,\delta]$ for some $\delta < T$, such that
$$\tilde g(t)\underset{t\rightarrow 0}{\longrightarrow} g_0$$
and assume that its curvature is uniformly bounded below. Then $\tilde g(t)$ agrees with the flow $g(t)$ constructed in Theorem
\ref{tma1} for $t\in (0,\delta]$.
\end{tma}

\section{Truncating cones}\label{sectiontruncating}
We construct in this section smooth approximations of cone surfaces, by smoothly truncating the vertices. This section is analogous to
Section 3.3 of \cite{Topping_revcusp}, where we substitute the cusp points with cone points. Let $D$ denote the unit disc, and $r=|z|$. An
appropriate elimination
of the asymptote of the conformal factor at $r=0$ gives rise to a metric which is smooth, and no longer singular at the origin.
\begin{lema}\label{lema1}
Let $g_0=e^{2(a_0+\beta \ln r)}|dz|^2$ be a cone metric on the punctured disc $D\setminus\{0\}$ with curvature bounded below, $K[g_0]\geq
-\Lambda$. There exists an increasing sequence of smooth metrics $g_k= e^{2u_k}|dz|^2$ on $D$  such that 
 \begin{enumerate}
  \item $g_k=g_0$ on $D\setminus D_{1/k}$,
  \item $g_k \leq g_0$ on $D\setminus \{0\}$,
  \item $\inf_{D_{1/k}}u_k \rightarrow + \infty$ as $k\longrightarrow +\infty$, and
  \item $K[g_k]\geq \min\{e^2 K[g_0],0\}$.
 \end{enumerate}
\end{lema}

\begin{proof}
The conformal factor $u_0=a_0+\beta \ln r$ of the cone metric tends to $+\infty$ as
$r\rightarrow 0$, so for each $k\in \mathbb N$ we pick
the minimum of $u_0$ and $k$ to obtain an increasing sequence of bounded functions tending to $u_0$. This has to be done in a way such
that the functions remain smooth.

Choose a smooth function $\psi : \mathbb R \rightarrow \mathbb R$ such that 
\begin{itemize}
 \item $\psi(s)=s$ for $s\leq -1$,
 \item $\psi(s)=0$ for $s\geq 1$,
 \item $\psi'\geq 0$ and $\psi''\leq 0$.
\end{itemize}
The smoothed minimum of $u_0$ and $k$ is
$$u_k=\psi(u_0-k)+k$$
and satisfies:
\begin{itemize}
 \item If $u_0\geq k+1$ then $u_k=k$ and therefore $K[g_k]=0$.
 \item If $u_0\leq k-1$ then $u_k=u_0$ and therefore $K[g_k]=K[g_0]$.
 \item If $k-1<u_0<k+1$ then \begin{itemize}
                              \item[\textbullet] $u_k\leq u_0$,
			      \item[\textbullet] $u_k\leq k$,
			      \item[\textbullet] $u_k \geq u_0-1$.
                             \end{itemize}
\end{itemize}
So 
$$u_k\leq \min\{u_0,k\}$$
and then (2) and (3) are satisfied. We can compute
$$\Delta u_k = \psi''(u_0-k)|\nabla u_0|^2 + \psi'(u_0-k)\Delta u_0 \leq \psi'(u_0-k)\Delta u_0 .$$
Now, since $\psi'\geq 0$, we can distinguish 
$$\Delta u_k \leq \Delta(u_0)\ \mathrm{if}\ \Delta u_0 >0$$
or
$$\Delta u_k\leq 0 \ \mathrm{if}\ \Delta u_0 \leq 0 .$$
So
$$\Delta u_k \leq  \max \{\Delta u_0, 0\}$$
and then
$$K[g_k]=-e^{-2u_k}\Delta u_k \geq \min\{e^2 K[g_0],0\} ,$$
so (4) is satisfied. Finally (1) is satisfied after passing to a subsequence, since the region of points $\{z:u_0(z)>k\}$ shrinks to a point when $k\rightarrow \infty$. 
\end{proof}

\section{Upper barriers}\label{sectionbarriers}
The conformal factor of a cone surface possesses asymptotes at the coordinates of the cone points, whereas the truncated approximations have
a big but finite value on that coordinates. This section provides a ratio of how fast the maximum value of this conformal
factors decay as the Ricci flow evolves.

\begin{lema}\label{lema2}
 Let $g(t)=e^{2u(t)}|dz|^2$ be a smooth Ricci flow on $D$ and $t\in [0,\delta]$, and assume that
 $$u(0)\leq A + \beta \ln r$$
 for some $A\in \mathbb R$. Then 
 $$u(t)< B + \frac{\beta}{2(\beta+1)}\ln t$$
for some $B$ depending only on $A$ and $\beta$.
\end{lema}

\begin{proof}
We will consider the conformal factor of several different surfaces.
The function 
$$s(r):=\ln\left(\frac{2}{1+r^2}\right)$$ 
is the conformal factor of a sphere, and the functions
$$v_0(r):=\ln(2(\beta+1))+\beta\ln r$$
$$v_1(r):=\ln(2(\beta+1))+\beta\ln r - \ln\left(1-r^{2(\beta+1)}\right)$$
are the conformal factors of Euclidean and hyperbolic cones (curvature $0$ and $-1$) respectively. Note that the Euclidean and hyperbolic cones become indistinguishable as $r\rightarrow 0$.

Considering the Ricci flow ($\frac{\partial u}{\partial t} = e^{-2u}\Delta u=-K$) on the hyperbolic cone, it evolves as
$$V_1(t)=v_1 + c(t)$$
with $c(t)$ an increasing function, so comparing with say $t=1$, we have
$$V_1(t)<v_1 + C$$
for some constant $C$ and for all $0<t<1$.

The function $s\left(\frac{r}{c_1}\right) + c_2$
is the conformal factor of a rescaled sphere (in parameter and in metric). We define
$$U(r,t):=\left\{ 
\begin{array}{lc} 
S(r,\lambda(t)):=s\left(\frac{r}{\lambda}\right) + v_1(\lambda) + C & \mathrm{if}\ 0<r\leq \lambda \\
v_1(r) + C & \mathrm{if}\ \lambda<r<1 
\end{array}\right.
$$
where $\lambda= \lambda(t)$ is a function of $t$ to be determined. Geometrically, $U$ is the conformal factor of a piecewise smooth
metric, a hemisphere near the origin and a cone with constant negative curvature away from it. It is a kind of ``blunt cone'', the transition being at
coordinate
$r=\lambda(t)$. We
still have to determine $\lambda(t)$, but we will require it to tend to $0$ as $t\rightarrow 0$.
In order to prove the lemma we will see (a) $u < U$ and (b) $\sup U(\cdot,t) \leq B + \frac{\beta}{2(\beta+1)}\ln t$.

We prove (a). We can assume that $u(0)<v_1+C$, and we know that at $r=0$ the value of $u$ is finite. Since the capping of $S(r,\lambda)$ occurs at arbitrarily big values, it is also true that $u(0)<U(0)$. Indeed, for $0<r\leq \lambda$, we have $S(r,\lambda)\geq v_1(\lambda)+C \rightarrow +\infty$ as $t\rightarrow 0$ since
$\lambda\rightarrow 0$. So $u< U$ for small positive $t$.

Suppose that for some $t_0$ there is a $0<r_0<1$ such that $u(r_0,t_0)=U(r_0,t_0)$. We can assume $t<1$. Note that the asymptote of $U$ at $r=1$ avoids the case of $r_0=1$.
If the point occurs at $\lambda\leq r_0 <1$, then $u$ would be touching the upper barrier of $V_1(t)$, which is impossible since by the maximum principle $u$ cannot pass over $V_1$.

Assume then that $0<r_0<\lambda$. We have $U-u\geq 0$ for $0\geq t \geq t_0$ and
$$u(r_0,t_0)=U(r_0,t_0),\qquad \frac{\partial}{\partial t}(U-u)\bigg|_{r_0,t_0}\leq 0,\qquad \Delta(U-u)\bigg|_{r_0,t_0}\geq 0$$
so at $(r_0,t_0)$
$$0\geq \frac{\partial U}{\partial t} - \frac{\partial u}{\partial t} = \frac{\partial U}{\partial t} - e^{-2u}\Delta u = \frac{\partial
U}{\partial t} + e^{-2U}(\Delta(U-u)-\Delta U) \geq \frac{\partial U}{\partial t} - e^{-2U}\Delta U$$
so
$$\frac{\partial U}{\partial t} \leq e^{-2U}\Delta U.$$
We now choose $\lambda(t)$ properly to contradict this assertion. 
On the one hand, at $(r_0,t_0)$
$$\frac{\partial U}{\partial t} = \frac{\partial S}{\partial t} = \left( - s'\left(\frac{r}{\lambda}\right) \frac{r}{\lambda^2} +
\frac{\beta}{\lambda} + \frac{2(\beta+1)\lambda^{2(\beta+1)}}{(1-\lambda^{2(\beta+1)})\lambda}\right) \frac{\partial \lambda}{\partial t} \geq \frac{\beta}{\lambda}\frac{\partial \lambda}{\partial t}$$
since $s'(r)<0$. 
On the other hand, one can compute
$$e^{-2U}\Delta U = e^{-2S}\Delta S = -\frac{\lambda^{-2(\beta+1)}}{(\beta+1)^2} \frac{e^{-2C}}{4} (1-\lambda^{2(\beta+1)})^2.$$
Ignoring the negligible term $\lambda^{2(\beta+1)}$ tending to zero (geometrically, assuming a flat cone), one can guess a critical value
of $\lambda$ by solving
$$\frac{\beta}{\lambda}\frac{\partial \lambda}{\partial t} = -\frac{\lambda^{-2(\beta+1)}}{(\beta+1)^2} \frac{e^{-2C}}{4}$$
e.g. with the solution
$$\lambda(t)=\left(\frac{-te^{-2C}}{2\beta(\beta+1)}\right)^{\frac{1}{2(\beta+1)}}.$$
A slight modification, say
$$\bar\lambda(t)=\left(\frac{-te^{-2C}}{4\beta(\beta+1)}\right)^{\frac{1}{2(\beta+1)}},$$
gives
$$\frac{\partial S}{\partial t}(r,\bar\lambda) \geq
\frac{\beta}{\bar\lambda}\frac{\partial \bar\lambda}{\partial t} =
\frac{1}{2}\frac{\beta}{(\beta+1)t} >
\frac{\beta}{(\beta+1)t} \left(1+\frac{te^{-2C}}{4\beta(\beta+1)}\right)^2=
e^{-2S(r,\bar\lambda)}\Delta S(r,\bar\lambda),$$
giving a contradiction as long as $C$ is big enough. Therefore there is no such time $t_0$ and so $u\leq U$.

Now we prove (b). We use the $\lambda=\bar \lambda$ just found. It is easy to check that $S(r,t)$ is nonincreasing and has a maximum at
$r=0$. Its value is
\begin{align*}
 S(0,\bar\lambda(t)) &=\ln(4(\beta+1)) + \beta \ln\left(\left(\frac{-t}{4(\beta+1)\beta}\right)^{\frac{1}{2(\beta+1)}}\right)
-\ln\left( 1+ \frac{te^{-2C}}{4\beta(\beta+1)}\right) \\
& \leq B + \frac{\beta}{2(\beta+1)}\ln t.
\end{align*}
\end{proof}

\section{Existence of the smoothening flow}\label{sectionexistproof}
We now prove the Theorem \ref{tma1}.

\begin{proof}
For simplicity assume there is just one cone point $p$. We take isothermal coordinates $z$ on a neighbourhood of $p$ such that $p$
corresponds to $z=0$, and $z\in D$ the unit disc (rescaling parameter and metric if necessary), so the metric on this chart has the form
$$g_0=e^{2(a+\beta\ln r)}|dz|^2$$
with $a:D\rightarrow \mathbb R$ a bounded and continuous function, smooth away from the origin.

We truncate the metric $g_0$ as in Lemma \ref{lema1} and we obtain an increasing sequence of smooth metrics $g_k$ on $\mathcal M$ such that:
\begin{enumerate}
  \item $g_k=g_0$ on $D\setminus D_{1/k}$,
  \item $g_k \leq g_0$ on $D\setminus \{0\}$,
  \item $\inf_{D_{1/k}}u_k \rightarrow + \infty$ as $k\longrightarrow +\infty$,
  \item $K[g_k]\geq \min\{e^2 K[g_0],0\}$ and
 \item $g_k\leq g_{k+1}$.
 \end{enumerate}
We apply Ricci flow to each initial metric $g_k$ and obtain a sequence of flows $g_k(t)$. There exist a uniform $T>0$ such that all flows
$g_k(t)$ are defined for $t\in[0,T]$. Indeed, by \cite{ChowKnopf} in dimension 2, if $\chi(\mathcal M)<2$ the flow is defined for
$t\in[0,\infty)$, and if $\chi(\mathcal M)=2$ the flow is defined for $t\in[0,\frac{\mathrm{Area}(\mathcal M)}{8\pi})$, and as $g_k\leq g_{k+1}$, then $\mathrm{Area}_{k}\leq\mathrm{Area}_{k+1}$.
So in any case the area does not tend to zero.

By the maximum principle, the initial $g_k(0)\leq g_{k+1}(0)$ implies $g_k(t)\leq g_{k+1}(t)$
and again by the maximum principle, $K[g_k(0)]\geq -\Lambda$ implies $K[g_k(t)]\geq -\Lambda $.

There exists also $\hat g(t)$, the instantaneously complete Ricci flow given by Theorem \ref{Topgiethm}. This flow is maximally stretched,
so any other flow $g_k(t)$ with the same initial condition $g_0$ satisfies
$$g_k(t) \leq \hat g(t) .$$


Now we have the sequence $g_k(t)$ satisfying
$$g_k(t) \leq g_{k+1}(t) \leq \hat g(t) ,$$
so we can define the limit flow
$$G(t)=\lim_{k\rightarrow \infty} g_k(t).$$

On any chart not containing $p$, the flow $G(t)$ is smooth by the uniform bounds of $g_k$ and the parabolic regularity theory. We need to
ensure that $G(t)$ extends smoothly to $p$, and in particular our constructed solution is different from the maximally stretched
flow (given by
Theorem \ref{Topgiethm}) and from the angle-preserving flow (given by Theorem \ref{Yinthm}), since none of them can be
extended to $p$. It is enough
to show that the conformal factor of $G(t)$ in a
neighbourhood of $p$ does not tend to $\infty$ for $t>0$. We use the Lemma \ref{lema2}. Say $G(t)=e^{2v(t)}|dz|^2$, then
$$v(t)=\lim_{k\rightarrow \infty} u_k(t) \leq C + \frac{\beta}{2(\beta+1)}\ln t$$
so $v(t)<+\infty$ for all $t>0$.
Furthermore, the uniform lower bound of the curvature on the approximative terms $g_k(t)$ also passes to the limit, so $K[G(t)] > -\Lambda$.
\end{proof}

\section{Uniqueness}\label{sectionuniqueproof}
The uniqueness issue is parallel to Topping's cusps, so we will sketch the proof and refer to \cite{Topping_revcusp} for a detailed
completion.
Although there are at least three Ricci flows with a cone surface as initial metric, say Topping's instantaneously complete flow,
Mazzeo et al./Yin's
cone flow, and our constructed smoothening flow, Topping's flow is unique amongst maximally stretched, unbounded curvature flows; and
our flow is unique amongst the lower-bounded
curvature, smoothly-extended flows.

\begin{proof}(Theorem \ref{tma2}) Recall that $\tilde g(t)$ is a Ricci flow defined on $\mathcal M$ for $t\in(0,\delta]$, with curvature uniformly bounded below, and such that $\tilde g(t)\rightarrow g_0$ as $t\rightarrow 0$. We want to show that it is unique. The proof consists in 4 steps:
\begin{step}
There exists a neighbourhood $\Omega$ of $p_i$, where the metric is written $\tilde g(t)=e^{2u}|dz|^2$, and there exists
$m\in \mathbb R$ such that $$u\geq m$$
in $\Omega$ for $t\in (0,\frac{\delta}{2}]$.
\end{step}
This step makes use of the lower curvature bound. Since $\frac{\partial u}{\partial t} = e^{-2u}\Delta u = -K[\tilde g] < \Lambda$, we have
$$u(z,t)\geq u(z,\frac{\delta}{2})-\Lambda\left(\frac{\delta}{2}-t\right)\geq \inf_{\Omega}u(\cdot,\frac{\delta}{2})-\Lambda\frac{\delta}{2}=:m.$$

\begin{step}
Actually, for every $M<\infty$, there is a small enough neighbourhood $\Omega_1$ and a small enough time $\delta_1$ such that
$$u\geq M$$
in $\Omega_1$ for $t\in (0,\delta_1)$.
\end{step}
This bound is obviously true for the conformal factor $u_0$ of the metric $g_0$, since $u_0=a+\beta\ln r$ has an asymptote on $r=0$. However, it is not clear that the factors $u(t)$ of the metrics $\tilde g(t)$ remain bounded by an arbitrary constant on a small neighbourhood for small $t$. It might happen that the functions $u(t)\rightarrow u_0$ as $t\rightarrow 0$ with $u(t)$ fixed at $r=0$ (that is, non-uniform convergence); but this case would contradict the uniform bounded below curvature. The sketch of the proof is as follows.

Define the family of functions $h(t)=\max\{M-u(t),0\}$, and the goal is proving that $h(t)\equiv 0$ for all $t<\delta_1$. We do that by showing that its $L^1$ norm on some small disc, $||h(t)||=\int_{D_\epsilon} |h(t)| d\mu$, vanishes. For, on the one hand $||h(t)||\rightarrow 0$ as $t\rightarrow 0$, since
$$||h(t)||=\max\{M-u(t),0\}\rightarrow \max\{M-u(0),0\} =0$$
because $u_0>M$. On the other hand, we claim that $\frac{d}{dt}||h(t)||\leq 0$, what proves the result. In order to prove that claim, we
change the functions $h(t)$ by a smoothed version of the maximum, in a similar fashion we did in the proof of Lemma \ref{lema1}, that is
$\hat h_\rho (t)=\Psi_\rho(M-u)\rightarrow h(t)$ as $\rho\rightarrow 0$. This allows us to compute $\frac{d}{dt}||\hat h_\rho(t)||$ in terms
of the derivatives of the controlled function $\Psi_\rho$, the lower bound on $u(t)$ given by the previous step, and the lower curvature
bound. See \cite{Topping_revcusp} for the details.

\begin{step}
With the lower bound of $u$, we can compare the flow $\tilde g(t)$ (which is conical at $t\rightarrow 0$) with
any Ricci flow smooth at $t=0$. Let $\sigma(t)$ be a smooth Ricci flow on $\mathcal M$ and $t\in[0,\delta]$. If 
$\sigma(0) <  g_0\ \mathrm{on}\ \mathcal M\setminus\{p_1,\ldots,p_n\}$, then $\sigma(t)\leq \tilde g(t)\ \mathrm{on}\ \mathcal M \ \forall t\in(0,\delta]$.
\end{step}
This step is essentially an application of the maximum principle. Let $s$ be the conformal factor of $\sigma(0)$. Since it is bounded, there
exists an $M$ and (by the previous step) a neighbourhood $\Omega$ of the cone points such that $s\leq M\leq u$ for a small time $t<t_1$ on
$\Omega$. But since $\tilde g(t)\rightarrow g_0$ and $\sigma(0)<g_0$, for an even smaller time $t<t_2$ we have $\sigma \leq \tilde g$ on the
whole $\mathcal M$. Having established the inequality for a positive time the maximum principle gives it for any time $t\in(0,\delta)$.

\begin{step}
Comparing two smoothening Ricci flows $\tilde g_1(t)$, $\tilde g_2(t)$ on $t\in(0,\delta]$ with initial metric $g_0$ and
curvature uniformly bounded below, a parabolic rescaling of one of them makes it a smooth Ricci flow even at $t=0$, so it is smaller or
equal than the other. By symmetry, also the other is smaller or equal than the one, so they are identical. 
\end{step}
The point is picking a small $t_0>0$ and the bound $K[\tilde g_1(t)]\geq -\Lambda$. We define a rescaling of $\tilde g_1(t)$ as
$$\sigma(t):=e^{-4\Lambda t_0}\ \tilde g_1(e^{4\Lambda t_0}\ t +t_0)$$
for $t\in [0,(\delta-t_0)e^{-2\Lambda t_0})$. This is a smooth Ricci flow even at $t=0$, and by the lower curvature bound it satisfies
$\sigma(0) < g_0$, so by the previous step
$\sigma(t)\leq \tilde g_2(t)$. But moving $t_0\rightarrow 0$ one gets
$\tilde g_1(t) \leq \tilde g_2(t)$
and by symmetry, also 
$\tilde g_2(t) \leq \tilde g_1(t)$.
\end{proof}

%% file: cuspsoliton.tex
\chapter{A three-dimensional expanding soliton}
\label{Ch:cuspsoliton}

\lettrine{\indent R}{ecall that} a gradient Ricci soliton is a smooth Riemannian metric $g$ on a manifold $M$ together with a potential
function
$f:M\rightarrow \mathbb R$
such that
\begin{equation}
\Ric + \Hess f +\frac{\epsilon}{2}g=0 \label{solitoneq} 
\end{equation}
for some $\epsilon\in \mathbb R$. As we saw in Chapter \ref{Ch:conesolitons}, solitons provide special examples of self-similar solutions of
Ricci flow, $\frac{\partial}{\partial t} g
= -2 \Ric$, evolving by homotheties and diffeomorphisms generated by the flow $\phi(\cdot,t)$ of the vector field $\grad f$, this is
$g(t)=(\epsilon t +1)\phi^*_t g_0$. The constant $\epsilon$ can be normalized to be $-1,0,1$ according the soliton being shrinking, steady
or expanding respectively (see \cite{TRFTA1} for a general reference). Solitons play an important role in the
classification of singular models for Ricci flow despite of (or actually due to) existing only a limited number of examples. 

In dimension
3, the only closed gradient solitons are those of constant curvature. Furthermore, by the results of Hamilton-Ivey \cite{Hamilton_surfaces},
\cite{Ivey}
and Perelman \cite{Perelman1}, the only
three-dimensional, possibly open, gradient shrinking
solitons with bounded curvature are $\mathbb S^3$, $\mathbb R^3$ and $\mathbb S^2 \times \mathbb R$ with their standard metrics, and their
quotients. Notice that in all these examples the gradient vector field is null. Examples with nontrivial
potential function in dimension three include the Gaussian flat soliton, the steady Bryant soliton \cite{Bryant}, the product of the
two-dimensional steady cigar soliton with $\mathbb R$ due to Hamilton, and a continuous family of rotationally symmetric expanding gradient
solitons due to Bryant \cite{Bryant}. In summary, shrinking and steady solitons are very few, and these are useful in the
analysis of high curvature regions of the Ricci flow. Expanding solitons are less understood and there is much more variety of them.

Bryant solitons are constructed using a dynamical system approach to study an associated ODE, by means of a phase portrait. Using similar
analysis, P. Baird \cite{Baird} found some nongradient solitons. We also used this technique in Chapter \ref{Ch:conesolitons} to get a
classification of all two-dimensional smooth and cone solitons. In this chapter we exploit this technique even further to find a new example
of soliton.

\vspace{1em}
A couple of motivating examples are the following. The hyperbolic metric $g=dr^2 + e^{-2r}(dx^2+dy^2)$ on $\mathbb R^3$ together with a
trivial potential $f=cst$ fits into the soliton
equation (\ref{solitoneq}) with $\epsilon=1$, so any quotient (hyperbolic manifold) yields an expanding soliton. An open quotient of
the hyperbolic space, for instance the cusp $\mathbb R \times \mathbb T^2$ obtained as a quotient by parabolic isometries (represented by
Euclidean translations on the $xy$-plane) yields also a soliton. 

Another interesting example occurs in $\mathbb R^n$ endowed with the Euclidean metric: it fits into the soliton equation together with a
potential function $f(p)=-\epsilon \frac{|p|^2}{4}$ for any $\epsilon\in\mathbb R$. The nonzero cases are the so called Gaussian solitons,
and even when the metric is constant in all cases (hence there is a unique solution with a given initial condition), there is more than one
soliton structure on it.

We will restrict ourselves to metrics with bounded curvature. This restriction is very natural for several reasons. First, if we consider
a smooth complete metric as initial condition for the Ricci flow, the bounded curvature ensures short time existence and uniqueness of
solutions of the
flow, both in the compact case (\cite{Hamilton_3mfds}, \cite{DeTur}) as well as on the complete 
noncompact one (\cite{Shi}, \cite{ChenZhu}). If the boundedness condition of the curvature is dropped, we can loose the
uniqueness; for instance approximating a cusp end by a high-curvature cap (cf. with Chapter \ref{Ch:smoothening}). A second reason is that
from the metric point of view, curvature  unbounded  below may  be very hard to control, especially when
related with
sequences of homotheties. A lower bound on the curvature is an essential requirement for the theory of Gromov-Hausdorff convergence and 
Alexandrov spaces, which is also underlying into the theorems for the long time existence of the Ricci flow. Even with the assumption of
bounded
curvature at the initial manifold, convergence under rescalings may be an issue if the bound is not uniform. For example, it is not clear 
what an open expanding hyperbolic manifold is at the birth time, namely $t\rightarrow -1$ when the evolution is
$g(t)=(t+1)g_{\mathrm{hyp}}$. A sequence of shrinked negatively curved manifolds does not need
to have a limit in the Gromov-Hausdorff sense, since the curvature is not uniformly bounded below. A hyperbolic cusp on $\mathbb R \times
\mathbb T^2$ as an expanding soliton (with a fixed basepoint to get pointed convergence) tends to a line while its curvature tends to
$-\infty$, as $t\rightarrow -1$ the birth time, which does not induce a natural manifold structure.

\vspace{1em}
The aim of this chapter is constructing a particular example of expanding gradient Ricci soliton on $\mathbb R \times
\mathbb T^2$, different from the constant curvature examples. Furthermore, we prove that it is the only possible nonhomogeneous soliton
on this manifold provided there is a lower sectional curvature bound equal to $-\frac{1}{4}$.

\vspace{1em}
In Section \ref{secexist} we consider the generic warped product metric $g=dr^2 + e^{2h} (dx^2 + dy^2)$ over $M=\mathbb R \times \mathbb
T^2$, where
$h=h(r)$ is
a function determining the size of the foliating flat tori, and a potential function $f=f(r)$ constant over these tori. We find a suitable
choice of $h$ and $f$ that makes the triple $(M,g,f)$ a soliton solution for the Ricci flow with bounded curvature, by means of the phase
portrait analysis of the soliton ODEs.

\begin{tma} \label{tma_exist}
There exists an expanding gradient Ricci soliton $(M,g,f)$ over the topological manifold $M = \mathbb R \times \mathbb T^2$ satisfying the
following properties:
\begin{enumerate}
 \item The metric has pinched sectional curvature $-\frac{1}{4} < sec < 0$.
 \item The soliton approaches the hyperbolic cusp expanding soliton at one end.
 \item The soliton approaches locally the flat Gaussian expanding soliton on a cone at the other end.
\end{enumerate}
More precisely, $M$ admits a metric
$$g=dr^2 + e^{2h(r)} (dx^2 + dy^2)$$
where $(r,x,y)\in \mathbb R \times \mathbb S^1 \times \mathbb S^1$; and a potential function $f=f(r)$, satisfying the soliton equation and
with the stated bounds on the curvature, such that
$$h\sim \frac{r}{2} \quad \mathrm{and} \quad f \rightarrow cst \qquad  \mathrm{as} \quad r\rightarrow -\infty$$
and
$$h\sim \ln r \quad \mathrm{and} \quad f \sim -\frac{r^2}{4} \qquad  \mathrm{as} \quad r\rightarrow +\infty .$$
\end{tma}
For the asymptotical notation ``$\sim$'', we write
$$\phi(r) \sim \psi(r) \quad \mathrm{as} \quad r\rightarrow \infty$$
if$$\lim_{r\rightarrow \infty} \frac{\phi(r)}{\psi(r)} = 1.$$

Remark that when $r\rightarrow +\infty$ the theorem states that the metric approaches $g=dr^2 + r^2 (dx^2 + dy^2)$. This is a
nonflat cone over the torus, namely its curvatures are $sec_{rx}=sec_{ry}=0$ and $sec_{xy}=-\frac{1}{r^2}$, but it indeed
approaches a flat metric when $r\rightarrow +\infty$.

\vspace{1em}
In Section \ref{secuniq} we consider the general case of a metric over $\mathbb R \times \mathbb T^2$ with $sec>-\frac{1}{4}$, and we
prove that the only nonflat solution is the example previously constructed. The lower bound on the curvature implies a concavity property
for the potential function. This leads together with the prescribed topology to a general form of the coordinate expression of the metric,
that can be subsequently computed as the example. 

\begin{tma} \label{tma_uniq}
Let $(M,g,f)$ be a nonflat gradient Ricci soliton over the topological manifold $M = \mathbb R \times \mathbb T^2$ with bounded curvature
$sec > -\frac{1}{4}$. Then it is the expanding gradient soliton depicted on Theorem \ref{tma_exist}.
\end{tma}

Let us remark that our example is also a critical case for the bound on
the curvature, in regard of the following known result, see \cite[Lem 5.5, Rmk 5.6]{Caoetal}.
\begin{prop} \label{propcao}
Let $(M^n, g,f)$ with $n\geq 3$ be a complete noncompact gradient expanding soliton with $\Ric\geq -\frac{1}{2} + \delta$ for some
$\delta>0$. Then $f$ is a strictly concave exhaustion function, that achieves one maximum, and the underlying manifold $M^n$ is
diffeomorphic to $\mathbb R^n$. 
\end{prop}

Our example proves that the proposition fails if one only assumes $\Ric>-\frac{1}{2}$ and thus the case $\Ric =-\frac{1}{2}$ is critical. In
our
critical case the soliton also has a strictly concave potential, but $f$ has no maximum, and
actually our solution admits a different topology for the manifold, namely $\mathbb R \times \mathbb T^2$.

\vspace{1em}
Most tedious computations thorough the chapter can be performed and checked using Maple or other similar software, therefore no step-by-step
computations will be shown. Pictures were drawn using Maple, and using the P4 program \cite{progP4} for the compactified portrait technique,
that we became aware thanks to J. Torregrosa.


\section{The asymptotically cusped soliton} \label{secexist}
Let us consider the metric
\begin{equation}
g=dr^2 + e^{2h(r)} (dx^2 + dy^2) \label{metrica} 
\end{equation}
where $h=h(r)$ is a one-variable real function, and a potential function $f=f(r)$ depending also only on the $r$-coordinate. The underlying
topological manifold can be taken $(r,x,y)\in 
\mathbb R \times \mathbb S^1\times \mathbb S^1 = \mathbb R \times \mathbb T^2$ since $\mathbb R^3$ with this metric admits the appropriate
quotient on the $x,y$
variables.
Standard Riemannian computations yield the following equalities.

\begin{lema}\label{lema_qtties}
The metric in the form (\ref{metrica}) associates the following geometric quantities:
\begin{align*} 
 \Ric &= -2((h')^2 + h'')dr^2 -e^{2h}(h''+2(h')^2)(dx^2 + dy^2) ,\\
 R &= -4 h'' -6 (h')^2 ,\\
 sec_{xy} &= -(h')^2  ,\\
 sec_{rx}=sec_{ry} &= -((h')^2 + h'') ,\\
 \Hess f &= f'' dr^2 + e^{2h}f'h'(dx^2 + dy^2) ,\\
 \Delta f &= 2h'f' + f'' ,\\
 \grad f &= f' \textstyle \frac{\partial}{\partial r} ,\\
 \nabla f &= f' dr ,\\
 |\grad f|^2 &= (f')^2 .\\
\end{align*}
\end{lema}

Hence the soliton equation (\ref{solitoneq}) for this metric turns into
$$\left( \frac{\epsilon}{2} + f'' -2(h')^2 -2 h''\right) dr^2 + e^{2h}\left( \frac{\epsilon}{2}+h'f'-2(h')^2 - h''\right) (dx^2 + dy^2)
=0  .$$
This tensor equation is equivalent to the ODEs system
\begin{equation}
 \left\{ \begin{array}{rcl}
	  \frac{\epsilon}{2} + f'' -2(h')^2 -2 h'' &=& 0 \\
	  \frac{\epsilon}{2}+h'f'-2(h')^2 - h''&=& 0 .
         \end{array} \label{eqssolcusp}
\right.\end{equation}
Let us remark that this system would be of second-order in most coordinate systems, but in ours we can just change variables $H=h'$ and
$F=f'$, and rearrange to get a first-order system
\begin{equation*}
\left\{ \begin{array}{rcl}
          H' & = & HF -2 H^2 + \frac{\epsilon}{2} \\
          F' & = & 2HF -2 H^2 + \frac{\epsilon}{2} .	 
         \end{array}
\right.
\end{equation*}

\begin{figure}[ht]
\centering
\includegraphics[width=0.6\textwidth]{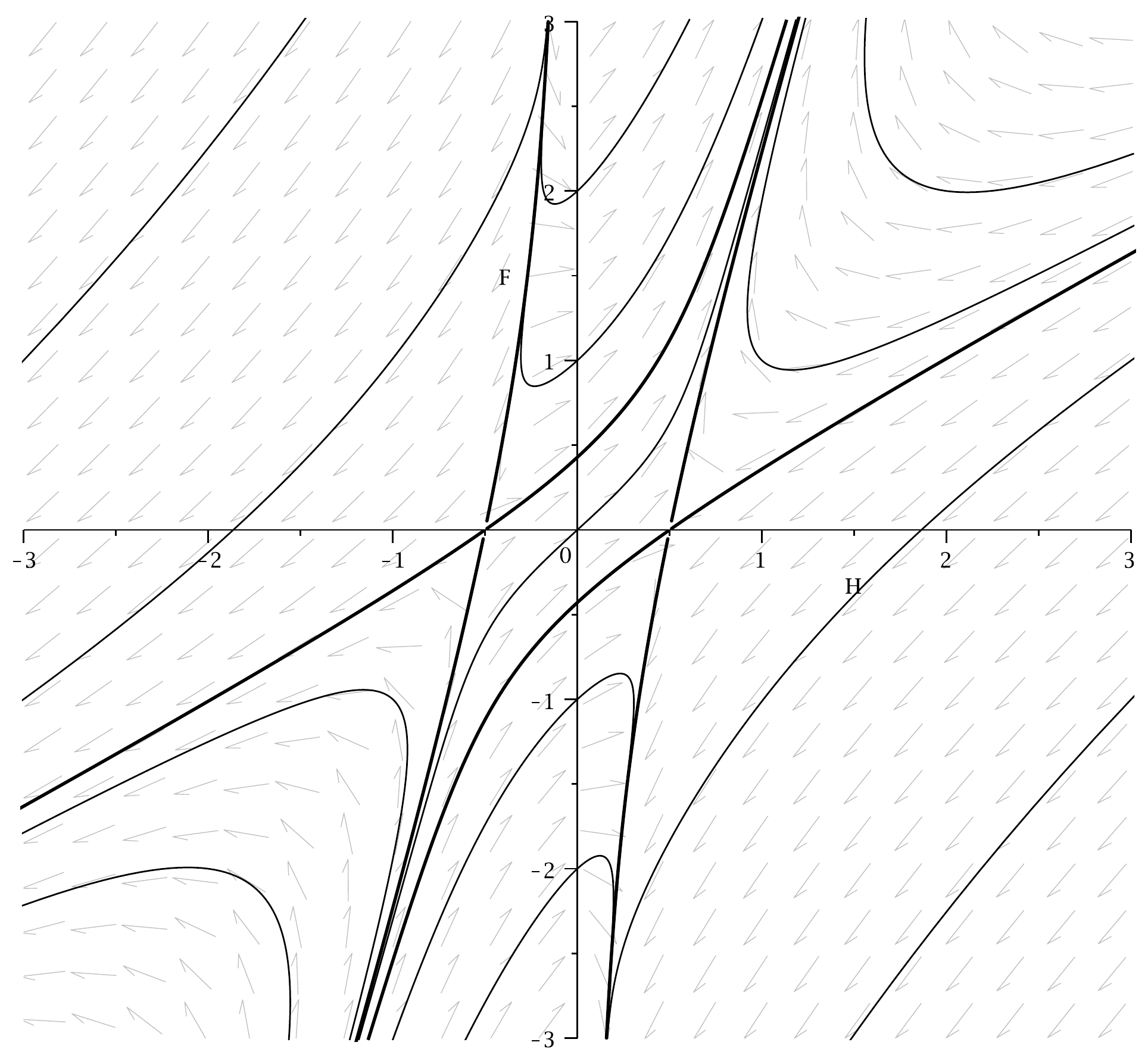}
\caption{Phase portrait of the system (\ref{eqsphase}).}
\label{retrato}
\end{figure}

We can solve qualitatively this system using a phase portrait analysis (see Figure \ref{retrato}). Every trajectory on the phase portrait
represents a soliton, but will not have in general bounded curvature. Actually, bounded curvature is achieved if and only if both $H$
and
$H'$ are bounded on the trajectory.

The critical points (stationary solutions) of the system are found by solving $\{ H'=0 , F'=0 \}$. If the soliton is shrinking
($\epsilon=-1$), there are no critical points and no trajectories with bounded curvature, agreeing with Perelman's classification. If the
soliton is steady ($\epsilon=0$), there is a whole straight line $\{H=0\}$ of fixed
points representing all of them the flat steady soliton. In this case there are neither trajectories with bounded $H$, hence all
solutions have unbounded negative curvature at least at one end. Let us assume henceforth that the soliton is
expanding ($\epsilon=1$), so our system is
\begin{equation}
\left\{ \begin{array}{rcl}
          H' & = & HF -2 H^2 + \frac{1}{2} \\
          F' & = & 2HF -2 H^2 + \frac{1}{2} .	 
         \end{array} \label{eqsphase}
\right.
\end{equation}
There are two critical points, $$(H,F)=(\pm\frac{1}{2},0).$$
The critical point $(\frac{1}{2},0)$ corresponds to a soliton with $h(r)=\frac{r}{2}+c_1$ and $f(r)=c_2$, the gradient
vector field is null, and the metric is $g_0=dr^2 + e^{r+c_1}(dx^2+dy^2)$, which is a complete hyperbolic metric, with constant
sectional curvature equal to $- \frac{1}{4}$, and possesses a cusp at $r\rightarrow -\infty$. As a Ricci flow it is $g(t)=(t+1)g_0$, it
evolves only by homotheties, and it is born at $t=-1$. The symmetric critical point $(-\frac{1}{2},0)$ represents the
same soliton, just reparameterizing $r\rightarrow -r$.

The phase portrait of the system (\ref{eqsphase}) has a central symmetry, that is, the whole phase portrait is invariant
under the change $(H,F,r)\rightarrow (-H,-F,-r)$, so it is enough to analyze one critical point and half the trajectories.

We shall see that the critical points are saddle points, and there is a separatrix trajectory emanating from each one of them that
represents the soliton metric we are looking for. Both trajectories represent actually the same soliton up to reparameterization.

\begin{lema}
 Besides the stationary solutions, and up to the central symmetry, there is only one trajectory $S$ with bounded $H$. This trajectory is a
separatrix joining a critical point and a point in the infinity on a vertical asymptote.
\end{lema}

\begin{proof}
The linearization of the system is 
$$\left(\begin{array}{c}
        H' \\ F'
        \end{array} \right)
= \left(\begin{array}{cc}
	F-4H & H \\
	2F-4H & 2H
        \end{array} \right)
 \left(\begin{array}{c}
        H \\ F
        \end{array} \right)
.$$
The matrix of the linearized system has determinant $-4H^2 \leq 0$, so the critical points are saddle points. For each one, there are
two eigenvectors determining four separatrix trajectories; being two of them attractive, two of them repulsive, according to the sign of the
eigenvalue. 

We are interested in one of the two repulsive separatrix emanating from the critical point $(H,F)=(\frac{1}{2},0)$, pointing towards the
region $\{H<\frac{1}{2}, F<0\}$. We shall see that has this is the only solution curve (together with its symmetrical) with bounded $H$
along its trajectory, so it represents a metric with bounded curvature.

In order to analyze the asymptotic behaviour of the trajectories, we perform a projective compactification of the plane, as
explained for instance in \cite{DumLliArt}, Ch. 5. The compactified plane maps into a disc where pairs of antipodal points on the boundary
represent
the asymptotic directions, Figure \ref{retrato_compact} shows the compactified phase portrait of (\ref{eqsphase}). 
\begin{figure}[ht]
\centering
\includegraphics[width=0.4\textwidth]{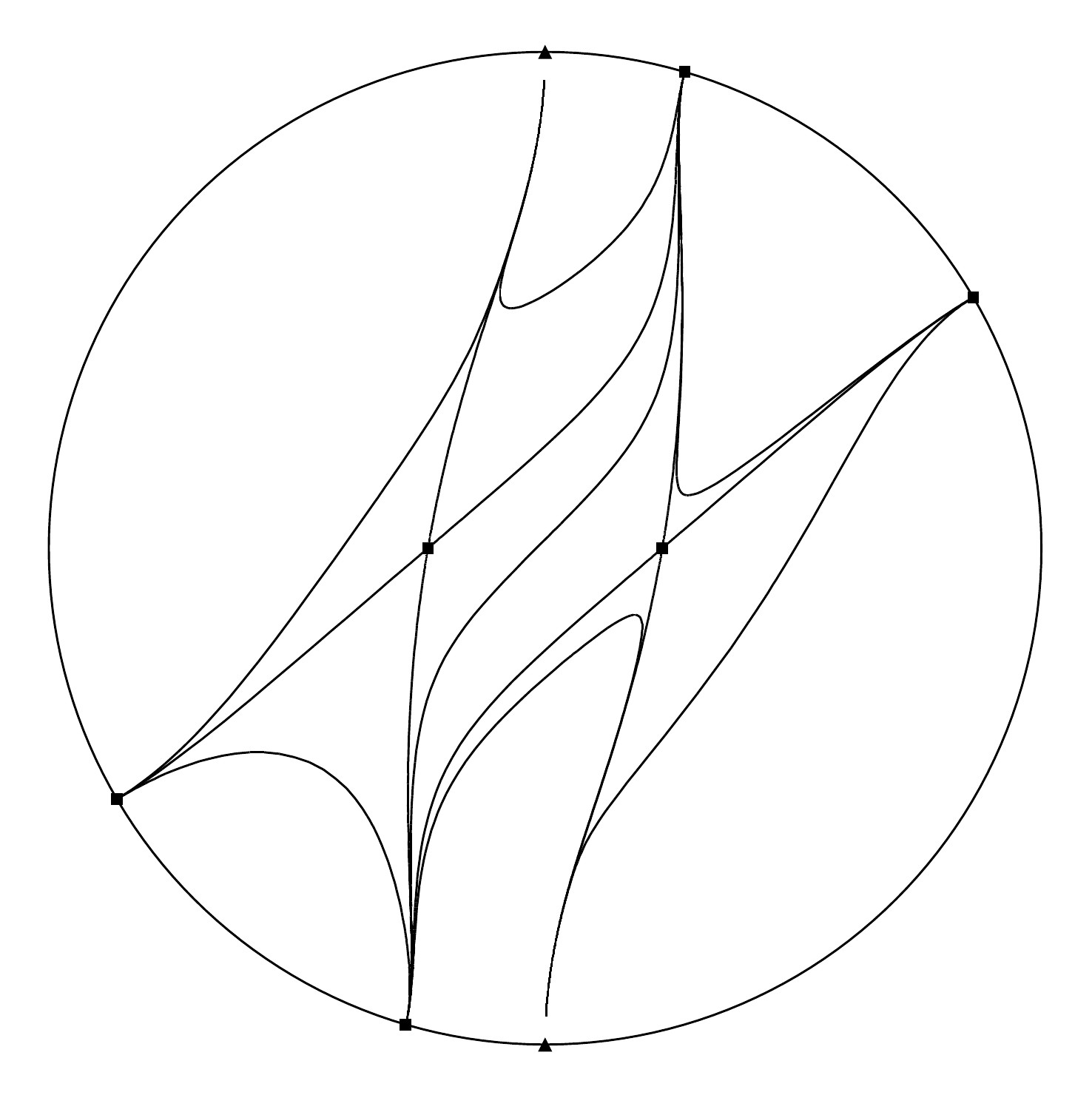}
\caption{Compactified phase portrait of the system (\ref{eqsphase}).}
\label{retrato_compact}
\end{figure}
A standard technique for polynomial systems is to perform a change of charts on the projective plane so that critical points at infinity can
be studied. A sketch is as follows: a polynomial system 
$$\left\{\begin{array}{l}\dot x = P(x,y) \\ \dot y = Q(x,y)  \end{array} \right. $$ 
can be thought as lying on the $\{z=1\}$ plane in the $xyz$-space. By a central projection this maps to a vector field and a phase portrait
on the unit sphere, or in the projective plane after antipodal identification. In order to do this, it may be necessary to resize
the vector field as 
$$\left\{\begin{array}{l}\dot x = \rho(x,y) P(x,y) \\ \dot y = \rho(x,y) Q(x,y)  \end{array} \right. $$
so that the vector field keeps bounded norm on the equator. However, this change only reparameterizes the trajectories. A global picture can
be obtained by orthographic projection of the sphere on the equatorial disc, as in Figure \ref{retrato_compact}, or it can be projected
further to a plane $\{x=1\}$ or $\{y=1\}$ in order to study the critical points at the infinity.
Let us remark that this technique works only for polynomial systems since the polynomial growth ratio suits the algebraic change of
variables.

In our system, this analysis yields that for every trajectory the ratio $H/F$ tends to either $0$, $1+\frac{\sqrt{2}}{2} $ or
$1-\frac{\sqrt{2}}{2} $ as $r\rightarrow \pm \infty$;
represented by the pairs of antipodal critical points (of type node) at infinity. The knowledge of the finite and infinite critical points,
together with their type, determines qualitatively the phase portrait of Figure \ref{retrato_compact} by the Poincaré-Bendixon theorem.
Thus, a trajectory with bounded $H$ on the $\mathbb R^2$ portrait, when seen on the $\mathbb RP^2$ portrait must have their ends either on
the finite saddle points or on the infinity node with $H/F$ ratio equal to $0$ (meaning a vertical asymptote). The only trajectory
satisfying this condition is the claimed separatrix and its symmetrical. 
\end{proof}

We shall see that this trajectory $S$ is parameterized by $r\in(-\infty,+\infty)$,
and when $r\rightarrow -\infty$ the function $h(r)$ behaves as $\frac{r}{2}$ and then the solution is asymptotically a cusp. Similarly, we
will see that $h', h'' \rightarrow 0$ when $r\rightarrow +\infty$ and then the solution is asymptotically flat.

To better understand the phase portrait it is useful to consider some isoclinic lines. This will give us the limit values for $H$,
$H'$ and the range of the parameter.
\begin{lema} \label{lemaisocl}
The vertical asymptote for the trajectory $S$ occurs at $H=0$. Furthermore, it is parameterized by $r\in(-\infty,+\infty)$ and
$H,H'\rightarrow 0$ as $r\rightarrow +\infty$.
\end{lema}

\begin{proof}

The vertical isocline $\{H'=0\}$ is the
hyperbola
$$F=2H-\frac{1}{2H}$$
and the trajectories cross it with vertical tangent vector. The horizontal isocline $\{F'=0\}$ is the hyperbola
$$F=\frac{1}{2}\left( 2H-\frac{1}{2H}\right)$$
and the trajectories cross it with horizontal tangent vector (see Figure \ref{retrato_closeup}). An oblique isocline is the hyperbola
$$F=2\left( 2H-\frac{1}{2H}\right) ,$$
since over this curve the vector field has constant direction:
$$(H',F')\big|_{(H,4H-\frac{1}{H})} = \left(2H^2-\frac{1}{2},6H^2-\frac{3}{2}\right) =
\left(2H^2-\frac{1}{2}\right) (1,3) .$$

\begin{figure}[ht]
\centering
\includegraphics[width=0.5\textwidth]{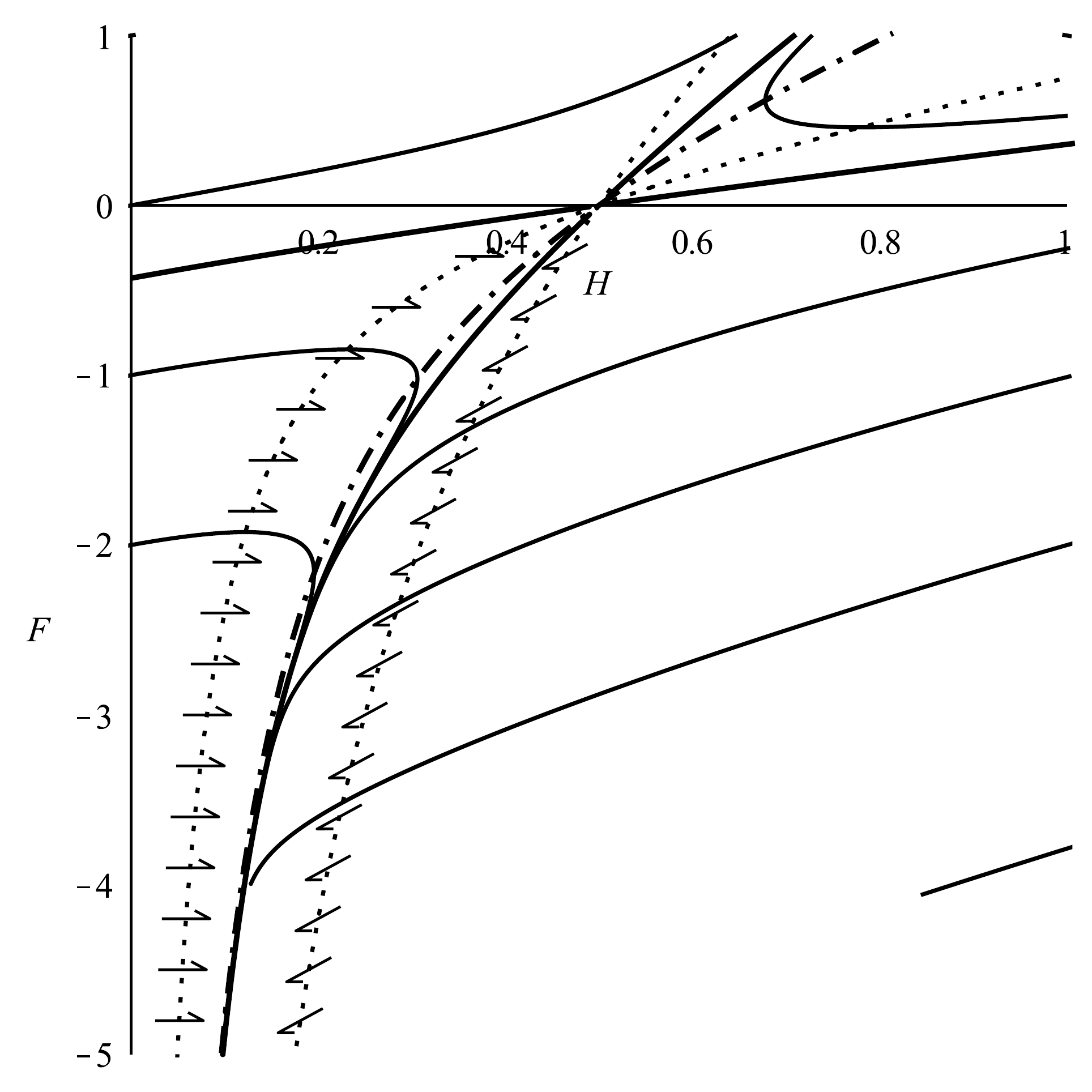}
\caption{Close-up of the separatrix trajectories (bold lines); the vertical isocline (dashes and dots); and the horizontal and oblique
isoclines (dots and arrows).}
\label{retrato_closeup}
\end{figure}

All three isoclines intersect at the critical points. Furthermore, the tangent directions at the critical point $(\frac{1}{2},0)$
have
slope $\frac{dF}{dH}\big |_{H=\frac{1}{2}} = 4$ for the vertical isocline, $\frac{dF}{dH}\big |_{H=\frac{1}{2}} = 2$ for the horizontal
one, and $\frac{dF}{dH}\big |_{H=\frac{1}{2}} = 8$ for the oblique one. The separatrix lines emanating from the critical point follow the
directions given by the eigenvectors of the matrix
of the linearized system (evaluated at the point), that is, the matrix

$$ \left(\begin{array}{cc}
	-2 & \frac{1}{2} \\
	-2 & 1
        \end{array} \right)
$$
whose eigenvalues are $\frac{-1+\sqrt{5}}{2}$ and $\frac{-1-\sqrt{5}}{2}$ with eigenvectors
$$ \left(\begin{array}{c}
        1 \\ 3+\sqrt{5}
        \end{array} \right)
\quad , \quad  
   \left(\begin{array}{c}
        1 \\ 3-\sqrt{5}
        \end{array} \right)
$$
respectively, so the separatrix lines have tangent directions with slope $3+\sqrt{5}\simeq  5.24$ and $3-\sqrt{5} \simeq 0.76$
respectively.
The repulsive separatrices are the ones associated with the positive eigenvalue, that is $\frac{-1+\sqrt{5}}{2}$, and the slope is $5.24$.

The repulsive separatrix emanating towards $\{H<\frac{1}{2} , F<0 \}$ initially lies below both the vertical and horizontal
isoclines, so it moves downwards and leftwards; and above the oblique one. The horizontal and oblique isoclines form two barriers for the
separatrix, this is, the separatrix cannot cross any of them. This is obvious for the horizontal one, since the flow is rightwards and the
trajectory is on the right. For the oblique isocline, we just check that any generic point on the isocline $(H,4H-\frac{1}{H})$ has
tangent vector $(1,4+\frac{1}{H^2})$ and a normal vector $\nu = (-4-\frac{1}{H^2},1)$ pointing leftwards and upwards for $0<H<\frac{1}{2}$.
The scalar product of the normal vector $\nu$ and the vector field $(H',F')$ over this isocline is
$$\langle \nu, (H',F')\rangle = \left(2H^2-\frac{1}{2}\right) \left(-4-\frac{1}{H^2}+3\right)=-2H^2+\frac{1}{2H^2}-\frac{3}{2}>0$$
whenever $0<H<\frac{1}{2}$. This means that the flow is always pointing to the left-hand side of the isocline branch and therefore is a
barrier. This proves that the separatrix moves downwards between the two barriers and therefore $H\rightarrow 0$. 

Actually, the vertical isocline is also a barrier for the separatrix. For, if at some point it touched the vertical
isocline, it would then move vertically downwards, keeping the trajectory on the right-hand side of the isocline. There would be then a
tangency, but it is impossible since the tangent vector should be vertical. Since the vertical isocline is becoming itself vertical, this
means that the vertical isocline acts as an attractor for the trajectories. Indeed, the trajectory lies initially in the region
$\{H>0,H'<0\}$ but $H$ must remain positive since it cannot cross the vertical isocline. Therefore $H$ is positive and decreasing, so $H'$
must tend to $0$. This implies that the trajectory tends to the vertical isocline. It is important to
note that both the vertical and horizontal isoclines come close together when $H\rightarrow 0$, but the trajectories stick to the vertical
one much faster than to the horizontal one.

We now see that $r\in(-\infty,+\infty)$. This follows immediately from the Hartman-Grobman theorem for the case $r\rightarrow -\infty$, but
the trajectory might, a priori, escape to infinity in finite time. This would require that the velocity tangent vector tends to infinity in
finite time, but this is impossible, since $H$ and $H'$ are bounded, thus $F'$ is bounded and hence the tangent vector $(H',F')$ is bounded.
\end{proof}

We have seen that not only $H\rightarrow 0$ as $r\rightarrow +\infty$, but also that $H' \rightarrow 0$. That is,
$h',h''\rightarrow 0$ as $r\rightarrow +\infty$, so all the sectional, scalar and Ricci curvatures tend to zero, the metric becoming
asymptotically flat.

At this point we have seen the existence of the soliton asserted in Theorem \ref{tma_exist}. We now give some more detailed information
about the
asymptotic behaviour of $f$ and $h$ at the ends of the manifold.

\begin{lema} \label{lemaasymp}
The asymptotic behaviour of $f$ and $h$ is
$$ h \sim \frac{r}{2} \quad \mathrm{and} \quad f \rightarrow cst \qquad  \mathrm{as} \quad r\rightarrow -\infty$$
and
$$h\sim \ln r \quad \mathrm{and} \quad f \sim -\frac{r^2}{4} \qquad  \mathrm{as} \quad r\rightarrow +\infty .$$
\end{lema}

\begin{proof}
Recall a version of the l'Hôpital rule: if $$\lim_{x\rightarrow \infty} \phi(x) = \lim_{x\rightarrow \infty} \psi(x) = 0, \pm \infty \quad
\mathrm{and} \quad \lim_{x\rightarrow \infty} \frac{\phi'(x)}{\psi'(x)} = c$$ then $$\lim_{x\rightarrow \infty} \frac{\phi(x)}{\psi(x)} =
c.$$

The case when $r\rightarrow -\infty$ follows from Hartman-Grobman theorem: the phase portrait in a small neighbourhood of a saddle critical
point has a flow that is Hölder conjugate to the flow of a standard linear saddle point 
$$\left\{ \begin{array}{l}
   \dot x = x \\
   \dot y = -y ,
  \end{array}\right.
$$
with solution $x(t) = k_1 e^t$, $y(t)=k_2 e^{-t}$.
This means that $r$ is defined from $-\infty$ onwards, that $H\rightarrow
\frac{1}{2}$ and $F,H',F'\rightarrow 0$ and $\frac{H-1/2}{F} \rightarrow 3+\sqrt{5}$ as $r\rightarrow -\infty$. Since
$H=h'\rightarrow \frac{1}{2}$ then $h\rightarrow -\infty$ and $\frac{H}{1/2}\rightarrow 1$. By l'Hôpital,
$\frac{h}{r/2}\rightarrow 1$.
Using more accurately the Hartman-Grobman theorem, there exists a Hölder function $\eta:U\subset \mathbb R \rightarrow \mathbb R$ defined on
a neighbourhood of zero, and constants $\alpha, C >0$ such that $F(r)=\eta(k_1 e^{r})$ and 
$$ | F(r) - F(r_0) | = |\eta(k_1 e^{r}) - \eta(k_1 e^{r_0}) | \leq C |k_1 e^{r} - k_1 e^{r_0} |^\alpha .$$
When $r_0 \rightarrow -\infty$, we obtain
$$ |F(r)| \leq \tilde C e^{\alpha r} ,$$
thus, $F$ is integrable on an interval $(-\infty,c]$ and thus $f\rightarrow f_0=cst$ as $r\rightarrow -\infty$. The constant $f_0$ is
actually $\sup f$ and can be chosen since it bears no geometric meaning. A more accurate description of $f$ using l'Hôpital tells
$$\lim_{r\rightarrow -\infty} \frac{h-r/2}{f-f_0}=3+\sqrt{5} .$$

For the case when $r \rightarrow +\infty$, we know from the trajectories that $H,H' \rightarrow 0$. Letting $r\rightarrow + \infty$ in the
first equation of (\ref{eqsphase}) we deduce that $HF\rightarrow -\frac{1}{2}$. Using this and letting $r\rightarrow + \infty$ in the second
equation of (\ref{eqsphase}) we conclude that $F'\rightarrow -\frac{1}{2}$, that is, $F' \sim -\frac{1}{2}$ and by l'Hôpital, $F\sim
-\frac{r}{2}$ and $f\sim -\frac{r^2}{4}$ as $r\rightarrow +\infty$.

Now, since $\lim_{r\rightarrow +\infty} \frac{F}{r} = - \frac{1}{2}$, we have
$$-\frac{1}{2} = \lim_{r\rightarrow +\infty} HF = \lim_{r\rightarrow +\infty} \frac{H}{r^{-1}} \frac{F}{r} = -\frac{1}{2} \lim_{r\rightarrow
+\infty} \frac{H}{r^{-1}}$$
thus $H\sim \frac{1}{r}$ and therefore $h\sim \ln r$ as $r\rightarrow +\infty$.
\end{proof}

It remains only to check the bounds on the sectional curvatures.

\begin{lema}
The metric (\ref{metrica}) with the function $f$ obtained as solution of the system (\ref{eqsphase}) has bounded sectional curvature
$$-\frac{1}{4}<sec<0 .$$
\end{lema}
\begin{proof}
The expression for the sectional curvatures is given in Lemma \ref{lema_qtties}. The case 
$$sec_{xy}=-(h')^2 = -H^2$$
is trivial since $0<H<\frac{1}{2}$ and therefore $-\frac{1}{4}<sec_{xy}<0$, tending to $-\frac{1}{4}$ at the cusp end, and to $0$ at the
wide end. The other sectional curvatures are
\begin{align*}
sec_{rx}=sec_{ry} & = -((h')^2+h'') \\
		  & = -H^2 - H' \\
		  & = H^2 - HF -\frac{1}{2} \\
		  & = -\frac{1}{2}\left( F'+\frac{1}{2} \right)
\end{align*}
We saw in Lemma \ref{lemaisocl} that $\{F'=0\}=\{2HF-2H^2+\frac{1}{2}=0\}$ is a barrier for the separatrix $S$. Hence, $F'<0$ along $S$
and therefore $sec_{rx}, sec_{ry}>-\frac{1}{4}$. Similarly, the set $\{H^2-HF-\frac{1}{2}=0\}$ can also be checked to be a barrier for
$S$ (actually a barrier on the opposite side), and hence $sec_{rx}, sec_{ry}<0$. We also saw in Lemma \ref{lemaasymp} the asymptotics
of $F'$, therefore we have $-\frac{1}{2}<F'<0$ and $sec_{rx}=sec_{ry}$ tend to $-\frac{1}{4}$ at the cusp end, and to $0$ at the wide end.
\end{proof}
This finishes the description of the soliton stated on Theorem \ref{tma_exist}.

\section{Uniqueness}\label{secuniq}
Let $(M,g,f)$ be a gradient expanding Ricci soliton over $\mathbb R \times \mathbb T^2$ such that $ sec > -1/4$. Then
$$\Ric + \Hess f + \frac{1}{2} g =0 ,$$
$$sec > 1/4,$$
$$\Ric > -1/2,$$
$$R > -3/2 .$$
Recall a basic lemma about solitons, that can be proven just derivating, contracting and commuting covariant derivatives on the soliton
equation, see \cite{TRFTA1}.
\begin{lema}
It is satisfied
$$ R + \Delta f + 3/2 =0 ,$$
$$ g(\grad R , \cdot ) = 2 \Ric(\grad f , \cdot) ,$$
$$ R + |\grad f |^2 + f = C .$$
\end{lema}

Since the soliton is defined in terms of the gradient of $f$, we can arbitrarily add a constant to $f$ without effect. We use this to set
$C=-3/2$ above so that we have
$$ \Delta f = |\grad f|^2 + f .$$
The bound on the curvature implies
$$\Hess f < 0 ,$$
$$\Delta f < 0 ,$$
$$ \langle \grad R , \grad f \rangle > -|\grad f|^2 .$$
First equation means that $f$ is a strictly concave function ($-f$ is a strictly convex function), i.e. $-f\circ \gamma$ is a strictly
convex real function for every (unit speed) geodesic
$\gamma$. This is a strong condition, since then the superlevel sets $A_c=\{f \geq c\}$ are totally convex sets, i.e. every geodesic
segment joining two points on $A_c$ lies entirely on $A_c$. Second equation is just a weaker convexity condition. This concavity on
this topology implies that $f$ has no maximum.

\begin{lema}
The function $f$ is negative and has no maximum.
\end{lema}
\begin{proof}
Note that $f$ is bounded above since $f=\Delta f -|\grad f|^2 < 0 $. 
Now suppose by contradiction that the maximum of $f$ is attained at some point of $\mathbb R \times \mathbb T^2$, then we can
lift this point, the metric and the potential function to the universal cover $ \mathbb R \times \mathbb R^2$. There is then a lattice of
points in the cover where the lifted function $\tilde f$ attains its maximum. But this is impossible since a strictly concave function
cannot have more than one maximum (the function restricted to a geodesic segment joining two maxima would not be strictly concave).

\end{proof}

\emph{Remark.} As stated in Proposition \ref{propcao}, if $\Ric > -\frac{1}{2} + \delta$ for any $\delta>0$, then $f$ has a maximum, and the
set
$A=\{f>f_{max} -\mu\}$ is compact and homeomorphic to a ball for small $\mu$. The function $f$ is then an exhaustion function, this is, the
whole manifold retracts onto $A$ via the flowline of $f$ and therefore $M\cong \mathbb R^3$. Thus this stronger bound on the curvature is
not
compatible with $M\cong \mathbb R \times \mathbb T^2$.

Now we prove that level sets of $f$ are compact.
\begin{lema}
The function $f$ is not bounded below and the level sets $\{f=c\}$ are compact.
\end{lema}
\begin{proof}
Consider $M=\mathbb T^2\times \mathbb R$ as splitted into $\mathbb T^2 \times (-\infty, 0] \cup \mathbb T^2\times [0,+\infty)$, each
component containing one of the two ends. Since $f$ has no maximum, there is a sequence of points $\{x_i\}$ tending to one end such that
$f(x_i)\rightarrow \sup f$. Let us assume that this end is $\mathbb T^2 \times (-\infty,0]$. Then when approaching the opposite end $f$ is
unbounded. Indeed, suppose by contradiction that there is a sequence of points $\{y_i\}$ tending to the $+\infty$ end such that
$f(y_i)\rightarrow -K
>-\infty$. There is a minimizing geodesic segment $\gamma_i$ joining $x_i$ with $y_i$. This gives us a sequence of geodesic paths (whose
length tends to infinity), each one crossing the central torus $\mathbb T^2\times \{0\}$. Since both the torus and the
space of directions of a point are compact, there is a converging subsequence of crossing points together with direction vectors that
determine a sequence of geodesic segments with limit a geodesic line $\gamma$. Now we look
at $f$ restricted to $\gamma$, this is $f\circ \gamma : \mathbb R \rightarrow \mathbb R$ such that $f\circ \gamma (t) \rightarrow \sup f$ as
$t\rightarrow -\infty$ and  $f\circ \gamma (t) \rightarrow -K > -\infty$ as $t\rightarrow +\infty$. But this is impossible since $f\circ
\gamma$
must be strictly concave. This proves that $f$ is not bounded below and that $f$ is proper when restricted to $\mathbb T^2 \times
[0,+\infty)$.

Now we consider $C_1=\min_{\mathbb T^2 \times \{0\}} f$ and $C_2<C_1<0$ such that the level set $\{f=C_2\}$ has at least one connected
component $S$ in $\mathbb T^2 \times (0,+\infty)$. Then $S$ is closed and bounded since no sequence of points with bounded $f$ can escape
to infinity. Therefore $S$ is compact. More explicitly, all level sets $\{f=C_3\}$ with $C_3<C_2$ contained in $\mathbb T^2 \times
(0,+\infty)$ are compact.

Now we push the level set $S$ to all other level sets by following the flowline $\varphi (x,t)$ of the vector field $\frac{\grad f}{|\grad
f|^2}$. Firstly, the diffeomorphism $\varphi(\cdot,t)$ brings the level set $S \subseteq \{f=C_2\}$ to the level set $\{f=C_2+t\}$,
\begin{align*}
f( \varphi (x,t)) &= f(\varphi(x,0)) + \int_0^t \frac{d}{ds} f(\varphi(x,s)) \ ds = f(x) + \int_0^t \langle \grad f
,\frac{d}{ds}\varphi(x,s) \rangle \ ds  \\
&= f(x) + \int_0^t \langle \grad f , \frac{\grad f }{|\grad f|^2} \rangle \ ds = f(x) + t .
\end{align*}
Secondly, the diameter distortion between these two level sets is bounded. If $\gamma :[0,1]\rightarrow \{0\}\times \mathbb T^2$ is a
curve on a torus,
\begin{align*}
g_{\varphi_t(x)}(\dot \gamma , \dot \gamma) &= |\dot \gamma|_x^2 + \int_0^t \frac{d}{ds} g_{\varphi_s(x)} (\dot\gamma , \dot\gamma) \ ds \\
&= |\dot \gamma|_x^2 + \int_0^t \mathcal L_{\frac{\grad f}{|\grad f|^2}} (\dot\gamma , \dot\gamma) \ ds \\
&= |\dot \gamma|_x^2 + \int_0^t \frac{2 \Hess f (\dot\gamma , \dot\gamma)}{|\grad f|^2} \ ds .
\end{align*}
Since $\Hess f <0$, this implies that $|\dot \gamma|_{\varphi_t(x)}^2 < |\dot \gamma|_x^2$ so all level sets $\{f=C_4\}$ with $C_4>C_2$
have bounded diameter, and hence are compact and diffeomorphic to $S$.
\end{proof}

Now, the level sets of $f$ are all of them compact and diffeomorphic, thus $M \cong \mathbb R \times \{f=c\} \cong \mathbb
R\times \mathbb T^2$ and therefore the level sets of $f$ are tori. This allows us to set up a coordinate system $(r,x,y)\in \mathbb R
\times \mathbb S^1 \times \mathbb S^1$ such that the potential function $f$ depends only on the $r$-coordinate. Furthermore, the gradient
of $f$ is orthogonal to its level sets, so the metric can be chosen not to contain terms on $dr\otimes dx$ nor $dr\otimes dy$. Thus the
metric can be written $g=u^2 dr^2 + \tilde g$ where $u=u(r,x,y)$ and $\tilde g$ is a family of metrics on the torus with
coordinates $(x,y)$ parameterized by $r$. Using isothermal coordinates, every metric on $\mathbb T^2$ is (globally) conformally
equivalent to the Euclidean one, thus $\tilde g = e^{2h}(dx^2 + dy^2)$ where $h=h(r,x,y)$. These conditions allow us to perform computations
that reduce to the particular case we studied in Section \ref{secexist}.

\begin{lema}
Consider the metric over $\mathbb R \times \mathbb T^2$
$$g=u^2 dr^2 + e^{2h} (dx^2 + dy^2)$$
where $u=u(r,x,y)$, $h=h(r,x,y)$, and a function $f=f(r)$. Assume that $g$, $f$ satisfy the soliton equation (\ref{solitoneq}) and that
$g$ has bounded nonconstant curvature. Then $g$ and $f$ are the ones described on the cusped soliton example of Section \ref{secexist}.
\end{lema}

\begin{proof}
The same Riemannian computations as before lead us to the soliton equation

\begin{align*}
0 =& \Ric + \Hess f +\frac{\epsilon}{2}g \\
=& \frac{1}{u} E_{11} \ dr^2 + \frac{1}{u^3}e^{2h} E_{22} \ dx^2 + \frac{1}{u^3}e^{2h} E_{33} \ dy^2  \\ 
 & + \frac{1}{u} E_{12} \ dr \ dx +\frac{1}{u} E_{13} \ dr \ dy +\frac{1}{u} E_{23} \ dx \ dy
\end{align*}

where

\begin{align*} 
E_{11} =& \textstyle
	 -u^2 \left( \frac{\partial^2 u}{\partial x^2} + \frac{\partial^2 u}{\partial y^2}\right) e^{-2h} +
	 \frac{\epsilon}{2}u^3 + f'' u -2 \left(\frac{\partial h}{\partial r}\right)^2 u -2 \frac{\partial^2 h}{\partial r^2} u +
	 \frac{\partial u}{\partial r}\left( 2\frac{\partial h}{\partial r} -f'\right) ,
\\
E_{22} =& \textstyle
	 -u^3\left(\frac{\partial^2 h}{\partial x^2} + \frac{\partial^2 h}{\partial y^2} \right) e^{-2h} 
	 -u^2\left(\frac{\partial^2 u}{\partial x^2} + \frac{\partial u }{\partial y } \frac{\partial h}{\partial y} - \frac{\partial
h}{\partial x} \frac{\partial u}{\partial x} \right) e^{-2h} \\
	  & \textstyle
	 + \frac{\epsilon}{2}u^3 + \frac{\partial h}{\partial r} f' u -2 \left(\frac{\partial h}{\partial r}\right)^2 u - \frac{\partial^2
h}{\partial r^2} u 
	 + \frac{\partial u}{\partial r} \frac{\partial h}{\partial r} ,
\\
E_{33} =& \textstyle
	 -u^3\left(\frac{\partial^2 h}{\partial x^2} + \frac{\partial^2 h}{\partial y^2} \right) e^{-2h} 
	 -u^2\left(\frac{\partial^2 u}{\partial y^2} + \frac{\partial u }{\partial x } \frac{\partial h}{\partial x} - \frac{\partial
h}{\partial y} \frac{\partial u}{\partial y} \right) e^{-2h} \\
	  & \textstyle
	 + \frac{\epsilon}{2}u^3 + \frac{\partial h}{\partial r} f' u -2 \left(\frac{\partial h}{\partial r}\right)^2 u - \frac{\partial^2
h}{\partial r^2} u 
	 + \frac{\partial u}{\partial r} \frac{\partial h}{\partial r} ,
\\
E_{12} =& \textstyle
	 \frac{\partial u}{\partial x} \left( \frac{\partial h}{\partial r} -f' \right) - u \frac{\partial^2 h}{\partial x \partial r} ,
\\
E_{13} =& \textstyle
	 \frac{\partial u}{\partial y} \left( \frac{\partial h}{\partial r} -f' \right) - u \frac{\partial^2 h}{\partial y \partial r} ,
\\
E_{23} =& \textstyle
	 \frac{\partial u}{\partial y}\frac{\partial h}{\partial x} + \frac{\partial u}{\partial x}\frac{\partial h}{\partial y} 
	 - \frac{\partial^2 u}{\partial x \partial y} .
\end{align*} 
Since the function $u$ never vanishes, nor the exponential does, the soliton equation is the PDE system
$\{E_{11}=E_{22}=E_{33}=E_{12}=E_{13}=E_{23}=0\}$. It is convenient to substitute the equations $E_{22}=0$ and $E_{33}=0$ with the linearly
equivalent $-\frac{1}{2}(E_{22}+E_{33})=0$ (equation (\ref{eq23a}) below) and $E_{22}-E_{33}=0$ (equation (\ref{eq23b})). Then, we get the
system
\begin{align} 
&	 -u^2 e^{-2h} \bigtriangleup u +
	 \frac{\epsilon}{2}u^3 + f'' u -2 \left(\frac{\partial h}{\partial r}\right)^2 u -2 \frac{\partial^2 h}{\partial r^2} u +
	 \frac{\partial u}{\partial r}\left( 2\frac{\partial h}{\partial r} -f'\right) =0 \label{eq1} ,
\\
&	 -u^3 e^{-2h} \bigtriangleup h -\frac{1}{2} u^2 e^{-2h} \bigtriangleup u 
	 + \frac{\epsilon}{2}u^3 + \frac{\partial h}{\partial r} f' u -2 \left(\frac{\partial h}{\partial r}\right)^2 u - \frac{\partial^2
h}{\partial r^2} u 
	 + \frac{\partial u}{\partial r} \frac{\partial h}{\partial r} =0 \label{eq23a} ,
\\
&	\frac{\partial u}{\partial x} \left( \frac{\partial h}{\partial r} -f' \right) - u \frac{\partial^2 h}{\partial x 		
\partial r} =0 \label{eq12} ,
\\
&	\frac{\partial u}{\partial y} \left( \frac{\partial h}{\partial r} -f' \right) - u \frac{\partial^2 h}{\partial y \partial 	r}
=0 \label{eq13} ,
\\
&	\frac{\partial u}{\partial y}\frac{\partial h}{\partial x} + \frac{\partial u}{\partial x}\frac{\partial h}{\partial y} 
	 - \frac{\partial^2 u}{\partial x \partial y} =0 \label{eq23} ,
\\
&	 2 \frac{\partial u}{\partial x}\frac{\partial h}{\partial x} -2 \frac{\partial u}{\partial y}\frac{\partial h}{\partial 	
y}  + \frac{\partial^2 u}{\partial x^2} - \frac{\partial^2 u}{\partial y^2} =0 \label{eq23b} ,
\end{align}
where $\bigtriangleup = \frac{\partial^2}{\partial x^2} + \frac{\partial^2}{\partial y^2} $ is the Euclidean Laplacian on the $xy$-surface.
We will recover our cusped soliton proving that $u\equiv 1$ and that $h(r,x,y)$ actually does not depend on $(x,y)$.

We consider first the equations (\ref{eq23}) and (\ref{eq23b}).
Since no derivatives on $r$ are present, we can consider the problem for $r$ fixed, so $u =u(r,\cdot,\cdot)$ is a function on the $xy$-torus with metric $e^{2h(r,\cdot,\cdot)}(dx^2+dy^2)$. The function $u$ must have extrema over the torus, since it is compact, so there are some critical points $(x_i,y_i)$ such that $\frac{\partial u}{\partial x}\big|_{(x_i,y_i)} = \frac{\partial u}{\partial x}\big|_{(x_i,y_i)}=0$. From the equations evaluated on a critical point, $\frac{\partial^2 u}{\partial x \partial y}\big|_{(x_i,y_i)}=0$ and $\frac{\partial^2 u}{\partial x^2}\big|_{(x_i,y_i)} =\frac{\partial^2 u}{\partial y^2}\big|_{(x_i,y_i)} = \lambda_i$ so the Hessian matrix (on the $xy$-plane) is 
$$\left(\begin{array}{cc} \lambda_i & 0 \\ 0 & \lambda_i  \end{array}\right) .$$
Suppose that every critical point is nondegenerate, that is, the Hessian matrix is nonsingular with $\lambda_i\neq 0$. Then the set of
critical points is discrete and $u$ is a Morse function for the torus. But then the Morse index on every critical point (the number of
negative eigenvalues of the Hessian) is either $0$ or $2$, meaning that every critical point is either a minimum or a maximum, never a
saddle point. Then Morse theory implies that the topology of the $xy$-surface cannot be a torus (being actually a sphere, see \cite{Mil}).
This
contradicts that every critical point is nondegenerate, so there is some point $(x_0,y_0)$ such that first and second derivatives vanish. 

We now proceed to derivate the two equations. Equations (\ref{eq23}) and (\ref{eq23b}) can be written
\begin{align}
u_{xy} &= u_y h_x + u_x h_y    \label{eq23_p} ,\\
u_{xx} - u_{yy} &= -2 u_x h_x +2 u_y h_y   \label{eq23b_p} ,
\end{align}
using subscripts for denoting partial derivation. Their derivatives are
\begin{align*}
u_{xxy} &= u_{xy}h_x + u_y h_{xx} + u_{xx}h_y + u_x h_{xy} ,\\
u_{xyy} &= u_{yy}h_x + u_y h_{xy} + u_{xy}h_y + u_x h_{yy} ,\\
u_{xxx} - u_{xyy} &= -2 u_{xx}h_x -2 u_x h_{xx} +2 u_{xy} h_y +2 u_y h_{xy} ,\\
u_{xxy} - u_{yyy} &= -2 u_{xy}h_x -2 u_x h_{xy} +2 u_{yy} h_y +2 u_y h_{yy} ,
\end{align*}
using the same notation. Evaluated at the point $(x_0,y_0)$, where all first and second order derivatives of $u$ vanish, the right-hand side
of these equations vanish and therefore all third derivatives vanish. Inductively, if all $n$-th order derivatives vanish at $(x_0,y_0)$,
then the $(n-1)$-th derivative of the equation (\ref{eq23_p}) implies that all mixed $(n+1)$-th order derivatives (derivating at least once
in each variable) vanish, then the $(n-1)$-th derivative of the equation (\ref{eq23b_p}) implies that all pure $(n+1)$-th order derivatives
(derivating only in one variable) also do; so all derivatives of all orders of $u$ vanish at $(x_0,y_0)$. Because $u(r,x,y)$ is a component
of a solution of the Ricci flow, it is an analytical function (see \cite{TRFTA2}, Ch. 13), so it must be identically constant in $(x,y)$.

At this point, we can reduce our metric to be $g=u(r)^2 dr^2 + e^{2h}(dx^2+dy^2)$ with $h=h(r,x,y)$. It is just a matter of reparameterizing the variable $r$ to get a new variable, $\bar r = \int u(r) \ dr$, such that $u(r)^2dr^2 = d\bar r^2$, so we rename $\bar r$ as $r$ and we can assume that the metric is $g=dr^2 + e^{2h}(dx^2+dy^2)$ with $h=h(r,x,y)$. 

We now look at the equations (\ref{eq12}) and (\ref{eq13}) when $u\equiv 1$, they imply
$$\frac{\partial^2 h}{\partial x \partial r} = \frac{\partial^2 h}{\partial y \partial r} = 0 ,$$
meaning that $\frac{\partial h}{\partial r}$ does not depend on $x$, $y$. 
Finally, looking at equations (\ref{eq1}), (\ref{eq23a}) when $u\equiv 1$, we get
\begin{align*}
\frac{\epsilon}{2} + f'' -2 \left(\frac{\partial h}{\partial r}\right)^2 -2 \frac{\partial^2 h}{\partial r^2} &=0 
,\\ 
\frac{\epsilon}{2} + \frac{\partial h}{\partial r} f' -2 \left(\frac{\partial h}{\partial r}\right)^2 - \frac{\partial^2 h}{\partial r^2} &=
e^{-2h} \bigtriangleup h .
\end{align*}
Since the left-hand side does not depend on $(x,y)$, nor does the term $e^{-2h} \bigtriangleup h$. Recall that a two-dimensional metric
written as $e^{2h(x,y)}(dx^2+dy^2)$ has Gaussian curvature $K=-e^{-2h} \bigtriangleup h$. So the $xy$-tori have each one constant curvature,
and the only admitted one for a torus is $K=0$. Hence $h$ only depends on $r$ and the equations
turn into the system (\ref{eqssolcusp}), that we already studied for the example of the cusp soliton. The rest of the uniqueness follows
from the discussion on Section \ref{secexist}.
\end{proof}

%
%

%% file: compactness.tex
\chapter[Compactness Theorems]{ Compactness theorems for classes of cone surfaces} \label{Ch:compactness}

\lettrine[slope=-4pt,nindent=-4pt]{\indent W}{hen studying} a sequence of manifolds, such as the sequences of rescaled flows that we use to
analyze
the
Ricci flow, we make use of the notions of convergence of sequences of manifolds, typically in the (weaker) Gromov-Hausdorff topology, or in
the (stronger) $\mathcal C^\infty$ topology. 

The theory of compactness of classes of manifolds is a key step in the Ricci flow theory (cf. \cite[Ch 3, Ch 4]{TRFTA1}). This
theory traces back to Cheeger \cite{Cheeger_finiteness}, Gromov \cite{Gromov}, Greene and Wu \cite{GreeneWu}, and Peters \cite{Peters}.
Hamilton adapted this work with
stronger hypothesis on the regularity of the curvature tensor, and proved a specific version for solutions to the Ricci flow that is the
appropriate theorem needed to perform sequences of rescalings on the flow \cite{Hamilton_compactness}.

If our manifolds have cone-like singularities, the cone structure after passing to a limit a priori might be very different from the
structure of the terms of the sequence. In this appendix we elaborate some theorems about the behaviour of
cone surfaces and Ricci flows over cone surfaces with respect to the limit of sequences. These theorems are used in Chapter \ref{Ch:cone_rf}
to
obtain limits of sequences of rescalings of the Ricci flow on cone surfaces.

First, in Section \ref{S:cptss:cptss_smooth} we review some classic compactness theorems for classes of smooth manifolds. In Section
\ref{S:cptss:conesurfs} we state the
first main theorem of the appendix, the compactness of a class of cone surfaces. The proof requires three lemmata that allow
a control of a cone surface obtained as a limit of a sequence: a control on the
number of cone points and the magnitude of their cone angles (given in Section \ref{S:cptss:conestruct}), a control of the injectivity
radius for the cone
points (given in Section \ref{S:cptss:inj_cone}), and a control of the injectivity radius for smooth points away from the cone points(given
in Section
\ref{S:cptss:inj_smooth}). These results are assembled in Section \ref{S:cptss:cptss_cone_proof} into the proof of the theorem. Finally, in
Section
\ref{S:cptss:cptss_flows} we use these compactness theorems for cone surfaces to obtain a compactness theorem for Ricci flows on cone
surfaces.

\section{Compactness theorems for smooth manifolds} \label{S:cptss:cptss_smooth}
A class of manifolds or metric spaces is simply a set of these spaces sharing certain properties. Such a class is compact if every sequence
has a subsequence converging to a limit
space, in the topology of a chosen distance between spaces. It is precompact if the limit does not belong to the class, this is, belongs
to the topological adherence. 

From a weak point of view, for metric spaces with the Gromov-Hausdorff distance, we have theorems such as
Gromov's compactness theorems, \cite{Gromov} , cf. \cite[Thm 7.4.15, Thm 10.7.2]{BurBurIva}.

\begin{tma}
Let $\mathfrak X$ be the class of metric spaces that are
\begin{itemize}
\item compact,
\item with diameter  $\mathrm{diam}(X) \leq D\quad \forall X\in\mathfrak X$,
\item $\forall \ \epsilon>0 \ \ \exists \ N(\epsilon)\in\mathbb N$ such that $\forall X\in\mathfrak
X$, $X$ has an $\epsilon$-net with $N(\epsilon)$ vertices.
\end{itemize}
Then $\mathfrak X$ is precompact in the Gromov-Hausdorff topology.
\end{tma}

\begin{tma}[Gromov's compactness]
Let $\mathfrak X$ be the class of metric spaces
\begin{itemize}
\item with dimension  $\mathrm{dim}(X)=n \quad \forall X\in\mathfrak X$,
\item with diameter  $\mathrm{diam}(X) \leq D\quad \forall X\in\mathfrak X$,
\item with (sectional) curvature  $\mathrm{sec}\geq \Lambda$.
\end{itemize}
Then $\mathfrak X$ is precompact in the Gromov-Hausdorff topology.
\end{tma}

On the other hand, in analysis there is an important compactness theorem for function spaces: the Arzelà-Ascoli theorem.

\begin{tma}[Arzelà-Ascoli]
Let $\mathfrak X$ denote the class of functions $f:\mathbb U\subset R^n \rightarrow \mathbb R$ being
\begin{itemize}
\item uniformly bounded: 

$\forall\ K$ compact, $\sup_{x\in K} |f(x)| < A$, $\forall f\in\mathfrak X ,$

\item equicontinuous:

$\forall\ K$ compact $\forall \ x_0\in K,\ \forall \epsilon > 0\ \exists\ \delta>0$ such that if
$|x-x_0|<\delta$ then $|f(x)-f(x_0)|< \epsilon \ \ \forall\ f\in \mathfrak X .$
\end{itemize}
Then $\mathfrak X$ is precompact on the topology of uniform convergence over compact sets.
\end{tma}
Furthermore, we can substitute the equicontinuity by the stronger hypothesis of $\mathcal C^1$ with uniformly bounded derivative, and apply
on degrees of regularity.

\begin{tma}[Arzelà-Ascoli, strong version]
The class of functions $\{ f\in \mathcal C^k(U)\ :\ \nabla^i f \mbox{ unif. bounded },\ i=0\ldots k \}$ is
relatively compact in the class $\{ f\in \mathcal C^{k-1}(U)\ :\ \nabla^i f \mbox{ unif. bounded },\ i=0\ldots k-1 \}$.
\end{tma}
That is, a sequence of functions $\mathcal C^k$ with uniformly bounded derivatives will have a subsequence converging to a limit
in $\mathcal C^{k-1}$.

Thus, from a stronger point of view, if our spaces are Riemannian manifolds, we can descend the metric tensor to a coordinate chart (open
subset of $\mathbb R^n$) and apply there the Arzelà-Ascoli theorem to its defining functions. For instance, bounding the curvature (now in
a Riemannian sense) we are bounding the second derivatives of the metric tensor (we have regularity up to second order), and hence we can
expect convergence at least in $\mathcal C^1$. We have then Cheeger's compactness theorem \cite{Cheeger_finiteness}, \cite[Thm
1.6]{Peters}

\begin{tma}[Cheeger's compactness]
Let $\mathfrak M$ denote the class of Riemannian $n$-manifolds with $\mathcal C^2$ regularity, satisfying
\begin{itemize}
\item $\mathrm{diam}\leq D$,
\item $\mathrm{vol}\geq V$ (equivalent to $\mathrm{inj}>i_0 $),
\item $|\mathrm{sec}|< \Lambda$.
\end{itemize}
Then $\mathfrak M$ is precompact in the biLipschitz topology, furthermore it is relatively compact on the class of $n$-manifolds with 
$\mathcal C^{1,1}$ regularity and $\mathcal C^0$ metric.
\end{tma}

Let us observe that we can eliminate the condition  $\mathrm{diam} \leq D$ if we use a pointed convergence; and that the second condition
can be chosen either $\mathrm{vol}\geq V$ or $\mathrm{inj}>i_0 $ thanks to Proposition \ref{P:srv:inj_vol} (this is equivalent to the
``Propeller lemma'' \cite[Thm 5.8]{CheegerEbin} on the global case ($\mathrm{diam}\leq D$) or the lemma of injectivity radius decay
with distance \cite[Thm 4.3]{CheGroTay} in the local case). 

From an even stronger point of view, (i.e. demanding more restrictive conditions), we can ask a control of more (or all) derivatives of
the metric tensor $\mathcal C^\infty$, and hence we have the following convergence theorem for sequences of metrics.

\begin{prop}[{\cite[Cor 3.15]{TRFTA1}}] \label{prop_ArzAsc_metrictensor}
Let $(\mathcal M,g)$ be a Riemannian $n$-manifold, and let $K\subset \mathcal M$ be a compact subset. Let $g_k$ be a sequence of
metrics on $K$ such that
\begin{itemize}
 \item $\sup_{0\leq \alpha\leq p+1} \sup_{x\in K} |\nabla^\alpha g_k | \leq C\leq \infty$.
 \item $g_k\geq \delta g$ for some $\delta>0$.
\end{itemize}
Then, there exists a convergent subsequence to a limit metric $g_\infty$ in the $\mathcal C^p$ sense.
\end{prop}

Alternatively, controlling all covariant derivatives of the curvature tensor $\Rm$, we control all derivatives of the metric tensor and the
previous proposition applies. This is Hamilton's version of the compactness theorem \cite{Hamilton_compactness}, see also \cite[Thm
3.9]{TRFTA1}.


\begin{tma}[Compactness of manifolds] \label{thm_compactness_manifolds}
Let $(M_k, O_k, g_k)$ be a sequence of complete pointed Riemannian manifolds such that:
\begin{enumerate}
\item $|\nabla^pRm|\leq C_p$, $\forall k, \ \forall p$
\item $inj(O_k)\geq i_0$, $\forall k$
\end{enumerate}
then there exists a converging subsequence in the $\mathcal C^\infty$ sense.
\end{tma}

\section{The compactness theorem for cone surfaces} \label{S:cptss:conesurfs}
Let us situate on the context of cone surfaces. We refer to Chapter \ref{Ch:cone_rf} for a discussion of the Riemannian definition of cone
surfaces;
however, in this appendix we will focus on the metric aspects of cone surfaces. In principle we work on the class of cone surfaces with
bounded $|\nabla^p \Rm|$ on all point on the smooth part, and with bounded $\inj$ on a nonsingular base point. We will need to impose some
conditions about the magnitude of the cone angles. We will restrict ourselves at first to cone angles less than $2\pi$ because these have
curvature (in the sense of Alexandrov) $+\infty$ (later we will restrict to angles less than $\pi$). Then, our cone surfaces are Alexandrov
spaces (curvature
$\geq C$). By the Gromov's compactness theorem, this class is precompact (in the pointed Gromov-Hausdorff topology), hence a sequence of
these surfaces has a converging subsequence in the weak sense (pointed G-H) towards a limit which a priori is only a metric space. To
improve this, our intention is to apply the compactness theorem for $\mathcal C^\infty$ manifolds to the smooth part, and on the other side
check that the result thereby obtained is indeed a cone surface (i.e. checking the local model around a cone point).

The first main theorem of this appendix is the following:
\begin{tma} \label{T:cptss:conesurfs}
Let $\mathfrak M$ denote the class of pointed cone surfaces satisfying
\begin{enumerate}
 \item cone points with angles $\leq \pi$,
 \item $|\nabla^p Rm_x| \leq C_p \quad \forall \ x\notin \Sigma$ where $\Sigma$ is the singular set of $M\in \mathfrak M$, for all $p\geq
0$,
 \item $\inj(O) \geq i_0 \quad \forall \ (M,O)\in \mathfrak M$, if $O\notin \Sigma$,
 \item Alternativelly, $\inj(O) \geq i_1$ and $\alpha>\alpha_0>0$, $\forall \ (M,O)\in \mathfrak M$; if $O\in\Sigma$.
\end{enumerate}
Then $\mathfrak M$ is compact in the topology of the (pointed) $\mathcal C^\infty$ convergence on the smooth part and Lipschitz on
the singular points.
\end{tma}

Let us remark several points.

\emph{Remark 1.} Recall that the Lipschitz convergence for metric spaces is those such that $d_L(M_k , M_\infty) \rightarrow 0$, where
the (bi)Lipschitz distance between two metric spaces is 
$$d_L(M_1,M_2)=\inf_f \left\{ \sup_{x,y} \ln \left| \frac{d(f(x),f(y))}{d(x,y)} \right| \right\} $$ 
where $f$ ranges over all Lipschitz homeomorphisms between $M_1$ and $M_2$ (and vice versa). This is not the same as Lipschitz convergence
of the functions $g_{ij}$ of the metric tensor. We will actually prove a stronger convergence ($\mathcal C^1$) of the metric tensor in a
certain coordinate chart. This will be enough for the Lipschitz convergence of the (metric) surfaces, but this sequence itself cannot enjoy
$\mathcal C^1$ convergence since at the singular point there are no derivatives in the smooth Riemannian sense.

\emph{Remark 2.} Lipschitz convergence implies Gromov-Hausdorff convergence. We will see below, Lemma \ref{L:cptss:conestruct}, that this
convergence is enough to guarantee that the number of cone points at the limit is the same as the number of cone points on the terms of
the approximating sequence, and the magnitude of the cone angles form convergent sequences for each cone point. Therefore, the cone
structure is preserved at the limit.

\emph{Remark 3.} Recall that the injectivity radius of a point $p$ in a smooth Riemannian
manifold is the
maximum value $\rho$ such that the exponential map $exp_p:T_pM \rightarrow M$ is a diffeomorphism of the ball $B(0,\rho)\in
T_pM$ onto its image. The injectivity radius of the whole manifold is defined as $\inj M = \inf_p \inj(p)$. With this definition applied
straightforward to the case of cone surfaces (with cone angles less than $2\pi$), the injectivity radius of any surface is zero, since any
point at a distance $\delta$ of a cone point would have injectivity radius at most $\delta$. However, the injectivity radius of the cone
point itself (defined by the exponential map from the tangent cone) is not zero. Thus, the appropriate framework is the use
of two injectivity radius: one $i_1$ for the cone points and one $i_0$ for the smooth points away from a certain (uniform) distance of the
cone points.

An outline of the proof of the theorem is the following. 

Let $\{ (M_k, g_k, O_k)\}_k$ be a sequence on the class $\mathfrak M$. We will actually apply all the argument to balls $B(O_k,R)$ in order
to get pointed convergence. First step is to ensure that the set of cone points is stable on a
subsequence, so we can locate and isolate them; and further, on each of these cone points the sequence of cone angles is convergent to a
limit cone angle. In other words, given the sequence $\{ (M_k,O_k)\}_k$ in $\mathfrak M$, there exists a subsequence such that all the
surfaces have the same number of cone points, and if $p_k\in M_k$ is a cone point of angle $\alpha_k$, then the sequence $\{\alpha_k\}$
converges $\alpha_k\to \alpha$.

Once the number of cone points is controlled, second step is to ensure that they are isolated, that is, the distance
between two of them does not tend to zero. This is true for the case when the cone angles are less than $\pi$, and for two cone angles of
angle exactly $\pi$ the points only can come close together if they also tend to infinity (get far away from the base point).

At this point, we can eliminate on each surface $\mathcal M_k$ a neighbourhood of the cone points of radius \emph{uniform} in $k$. We obtain
hence a sequence of surfaces $M_k \setminus \mathcal U_\epsilon (\Sigma_k)$ not complete but smooth, to which we can apply Hamilton's
compactness theorem. If we check carefully its proof, this starts picking on each manifold $M_k$ a net of points $\{k_k^i\}_{i=0}^N$ and
radii $\{r_k^i\}_{i=0}^N$ such that the balls $\{B(x_k^i, r_k^i)\}_{i=0}^N$ form a regular covering of $M_k$ with good properties, that
constitute an atlas (via normal coordinates by the exponential map) of the manifold $M_k$. The election of the points $x_k^i$ and the radii
$r_k^i$ relies on a control of the injectivity radius on each point $z\in M_k$ given its distance to a basepoint $O_k=x_k^0$, where by
hypothesis $\mathrm{inj} \ x_k^0 >i_0$. This result is known as the injectivity radius decay in the smooth case. We will show that its proof
adapts and is valid in spite of the cone points.

We shall use these results to prove the convergence $\mathcal C^\infty$ of the smooth part. Last step will be to check that the convergence
on a neighbourhood of a cone point is uniform for the metric on the model coordinate chart of a cone point.

\section{Cone structure is the same for GH-close cone surfaces} \label{S:cptss:conestruct}
In this section we show that two compact cone surfaces that are close enough in the Gromov-Hausdorff sense must have the same number of cone
points, and their respective cone angles must also be close. The Lemma will apply on compact balls centered at the base points.

\begin{lema} \label{L:cptss:conestruct}
For all $\epsilon, i_1, i_0, \Lambda>0$ and $\omega_0<2\pi$ there exists $\delta>0$ such that the following holds. Let $\mathcal M, \bar
{\mathcal M} $ be two compact cone surfaces with $\Sigma, \bar\Sigma$ their singular sets, satisfying
\begin{itemize}
 \item cone angles $\leq \omega_0 <2\pi$.
 \item $\inj x \geq i_1 \quad \forall x\in \Sigma$ (resp. $\bar \Sigma$).
 \item $\inj x \geq i_0 \quad \forall x\in \mathcal M \setminus \mathcal U_{\frac{i_1}{2}}(\Sigma)$ (resp. for $\bar{\mathcal M} $).
 \item $|\sec_x|\leq \Lambda  $ for all $x\notin\Sigma$ (resp. $\bar\Sigma$).
\end{itemize}
If the Gromov-Hausdorff distance between them is 
$$d_{GH}(\mathcal M, \bar{ \mathcal M})<\delta ,$$
then there exists a $2\delta$-isometry $f:\mathcal M \rightarrow \bar{ \mathcal M}$ that sends cone points to cone points, and if
$p\in\Sigma$ and $f(p) \in \bar \Sigma$ have cone angles $\omega$, $\bar\omega$, respectively, then 
$$|\omega-\bar\omega|<\epsilon .$$
\end{lema}

\begin{proof}
From the properties of Gromov-Hausdorff distance, we have \cite[Cor 7.3.28]{BurBurIva} that $d_{GH}(\mathcal M, \bar{\mathcal M})\leq
\delta$ implies that there exists $f:\mathcal M \rightarrow \bar{\mathcal M}$ a $2\delta$-isometry, that is, a possibly noncontinuous
function such that
\begin{enumerate}
 \item $dis\ f =\sup_{x,x'} |d_{\bar{\mathcal M}}(f(x),f(x'))- d_{\mathcal M}(x,x')|\leq 2\delta$.
 \item $f(\mathcal M)$ is a $2\delta$-net in $\bar{\mathcal M}$.
\end{enumerate}
Our goal is to modify the map $f$ so it sends cone points to cone points. Hence, we need to show that if $p$ is a cone point in $\mathcal
M$, then $f(p)$ is arbitrarily close to a cone point in $\bar{\mathcal M}$, choosing $\delta$ small enough.

Let $p\in\mathcal M$ be a cone point with angle $\omega<\omega_0<2\pi$. We launch four small geodesic rays spreading from $p$, of length
$d<i_1$,
on directions separated by an angle of $\alpha:=\frac{\omega}{4}$. We join the four endpoints of the rays with geodesic paths to form a
quadrilateral.

The idea is the following. The image by $f$ of the vertices of this quadrilateral defines a new quadrilateral in $\bar{\mathcal M}$. We can
compare the triangles formed by the rays and the sides of the quadrilaterals. Since the sides of the triangles are almost the same ($\pm
2\delta$), and the curvature is bounded, the angles also must be almost the same, and hence the angles around $f(p)$ will also add up less
than $2\pi$.

We proceed by contradiction, and we suppose that the image of the quadrilateral is contained in an open set of $\bar{\mathcal M}$ with
no cone points. We consider the triangle in $\mathcal M$ defined by $p$ and two rays of length $d$ forming an angle $\alpha$; with a third
side of length $l$. The image of the vertices of this triangle by $f$ defines a new triangle in $\bar{\mathcal M}$ of corresponding sides
$\bar d_1$, $\bar d_2$ and $\bar l$, with an angle in $f(p)$ of $\alpha'$. Since $f$ is a $2\delta$-isometry, we have
$$|l-\bar l|<2\delta ,$$
$$|d-d_j|<2\delta\quad j=1,2 .$$
We look for an upper bound of $\alpha'$ in terms of $\alpha$. We compare these triangles with their constant curvature models. We can
assume without loss of generality (by a dilation) that the bound on the curvature is $\Lambda =1$.

\begin{figure}[ht]
\begin{center}
\def\svgwidth{0.6\textwidth}
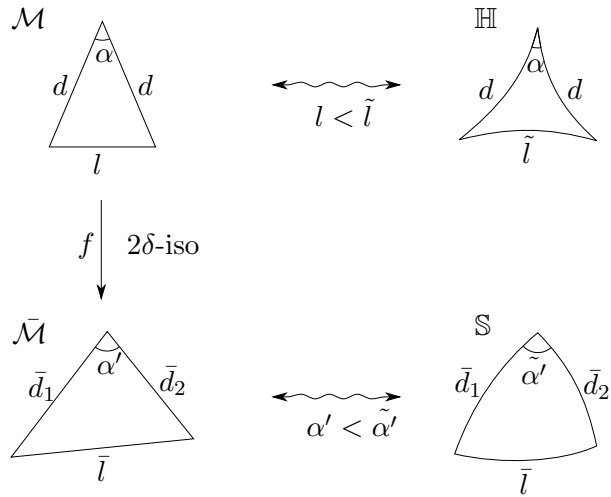
\caption{Comparison between the different triangles.}
\end{center}
\end{figure}

We compare the triangle in $\mathcal M$ with a hyperbolic triangle keeping the data side-angle-side fixed (a hinge). We compare the triangle
in $\bar{\mathcal
M}$ with a spherical triangle keeping the data side-side-side fixed (side lengths). By Toponogov's comparison theorems, we have
\begin{align}
l &< \tilde l  \label{compare1} \\
\alpha' &< \tilde{\alpha'}. \label{compare2}
\end{align}

On the other hand, by spherical/hyperbolic trigonometry (e.g. \cite[p 83]{Ratcliffe}),
\begin{align}
\cos \alpha &= \frac{\cosh d \cosh d - \cosh \tilde l}{\sinh d \sinh d} \\
\cos \tilde{\alpha'} &= - \frac{\cos \bar d_1 \cos \bar d_2 - \cos
\bar l}{\sin \bar d_1 \sin \bar d_2} 
\end{align}
and (since $\alpha'<\frac{\pi}{2}$) we have by \eqref{compare1} and \eqref{compare2},
\begin{align}
\cosh l &< \cosh \tilde l  \\
\cos \alpha' &> \cos \tilde{\alpha'}.
\end{align}
Hence,
\begin{align}
\cos \alpha &< \frac{\cosh d \cosh d - \cosh l}{\sinh d \sinh d} =: A\\
\cos \alpha' &> \cos \tilde{\alpha'} =: B.
\end{align}
We compare $A$ and $B$. Although these cannot be ordered in general, it is enough for us to show that
$$\lim_{d\rightarrow 0}\frac{A}{B}=1.$$
Indeed, if $\delta<d^3$ and since $l<2d$,
\begin{align*}
\lim_{d\rightarrow 0}\frac{A}{B} &= 
\lim_{d\rightarrow 0} -\frac{\cosh d \cosh d - \cosh l}{\cos \bar d_1 \cos \bar d_2 - \cos \bar l} \ 
\frac{\sin d_1}{\sinh d} \ \frac{\sin \bar d_2}{\sinh d} \\
&= \lim_{d\rightarrow 0} \frac{(1+\frac{d^2}{2}+O(d^4))^2 - (1+\frac{l^2}{2}+O(d^4))}{-(1-\frac{\bar d_1^2}{2}+O(d^4)) (1-\frac{\bar
d_2^2}{2}+O(d^4)) + (1-\frac{l^2}{2}+O(d^4))} \\
&= \lim_{d\rightarrow 0} \frac{1+d^2 +O(d^4) -1 -\frac{l^2}{2} +O(d^4)}{-(1-d^2 + O(d^4))+1-\frac{l^2}{2}+O(d^4)} \\
&= \lim_{d\rightarrow 0} \frac{d^2-\frac{l^2}{2} + O(d^4)}{d^2-\frac{l^2}{2} + O(d^4)} = 1
\end{align*}
Therefore, for all $\eta>0$,
$$1-\eta < \frac{A}{B} < 1+\eta$$
if we choose $d$ small enough and $\delta < d^3$. Then,
$$|A-B|<\eta |B|<\eta .$$
Finally,
$$\cos \alpha' >B > A-\eta >\cos \alpha -\eta$$ 
Given $0<\alpha<\frac{\pi}{2}$ we have $\cos\alpha>0$. Let $\eta$ be such that $\cos \alpha -\eta >0$. Then $\cos \alpha' > \cos
\alpha-\eta>0$ and therefore $\alpha'<\frac{\pi}{2}$.

Summarizing, if the image of the quadrilateral is contained in an open set of $\bar{\mathcal M}$ with no cone points, then we can apply the
comparison argument with the four triangles forming the quadrilateral in $\mathcal M$, but then the central point $f(p)$ would be conical.

Thus, we have seen that if $\delta$ is small enough and $p$ is a cone point on $\mathcal M$, then $f(p)$ is arbitrarily close to a cone
point on $\bar{\mathcal M}$. Only one cone point, since on $\bar{\mathcal M}$ there are no cone points arbitrarily close (the injectivity
radius is bounded). It is then easy now to define $\tilde f$ as $\tilde f(x)=f(x)$ if $x$ is a smooth point, and $\tilde f(p)=\bar p$ if
$p$ is a cone point, where $\bar p$ is the closest cone point on $\bar{\mathcal M}$ to $f(p)$. From our argument above, if $f$ is a
$2\delta$-isometry, then $\tilde f$ is a $(2\delta+\mu)$-isometry for any $\mu$ arbitrarily small. Therefore we can construct a
$2\delta$-isometry that sends cone points to cone points.

We finish the proof of the lemma by comparing the cone angles at $p$, $\bar p$. If $\alpha=\frac{\omega}{4}$ and
$\bar\alpha=\frac{\bar\omega}{4}$, we have seen that
$$\cos \bar\alpha >\cos \alpha -\eta_1$$
By symmetry, there is also a $2\delta$-isometry $g:\bar{\mathcal M}\rightarrow \mathcal M$ and hence
$$\cos \alpha >\cos \bar\alpha -\eta_2$$
Therefore
$$|\cos \alpha - \cos \bar\alpha |<\max \{\eta_1, \eta_2\}$$
and then $|\alpha-\bar\alpha|<\epsilon$ if $\eta_1,\eta_2$ small enough.
\end{proof}

\section{Injectivity radius bound for cone points} \label{S:cptss:inj_cone}

In this section we prove that two cone points with cone angles less than or equal to $\pi$ cannot be close together on a surface with
bounded curvature. Heuristically, forcing two cone points to be close each other would force the curvature to descend towards $-\infty$
on the region near the geodesic joining the two cone points.

\begin{lema} \label{L:cptss:inj_cone}
For all $i_0, \Lambda, D>0$ there exists $C=C(i_0, \Lambda, D)>0$ such that the following holds. Let $(\mathcal M,x_0)$ be a cone surface
with smooth base
point $x_0$, satisfying
\begin{itemize}
 \item cone angles $\leq \pi$.
 \item $\inj(x_0)>i_0$.
 \item $|\sec|<\Lambda$.
\end{itemize}
If $p,q\in B_D(x_0)$ are two cone points, then
$$d(p,q)>C$$
and
$$ \inj(p)>C .$$
\end{lema}

\begin{proof}
Let $p,q\in \mathcal M$ be two cone points, and let $\sigma$ be the shortest geodesic arc joining them, and let $|\sigma|$ be its length.
Let
$$B_R(\sigma)=\{x\in \mathcal M \ : \ d(x,\sigma)\leq R\}$$
be a neighbourhood of $\sigma$ of radius $R$. Clearly, as $R\rightarrow \infty$, $B_R(\sigma)$ is exhausting all $\mathcal M$. Let
$$N_R(\sigma)=\{x\in \mathcal M \ : \ x=\exp_y(v),\ y\in\sigma,\ v\perp\pm\dot\sigma(y),\ ||v||\leq R\}$$
be a normal neighbourhood of $\sigma$.

First step in the proof is that since the cone angles at $p,q$ are $\leq \pi$, 
$$B_R(\sigma)=N_R(\sigma).$$
Indeed, if $x\in\mathcal M$ and $y\in\sigma$ is the point that realizes the distance, $d(x,\sigma)=d(x,y)$, then $y$ must be attained
perpendicularly, and hence $y$ is the foot of the perpendicular from which the exponential emanates. To see this, if $y$ is not an endpoint
of $\sigma$, then if the angle $\theta$ between $\sigma$ and the geodesic joining with $x$ were less than $\frac{\pi}{2}$ on either side,
the distance could be shortened towards that side. If $y$ is one of the endpoints of $\sigma$, the angle $\theta$ must necessarily be
$\leq\frac{\pi}{2}$, since the space of directions at $p$ and $q$ measures $\leq\pi$, and hence the angle forming two geodesics at a cone
point must form an angle $\leq\frac{\pi}{2}$.

Let now $V$ be the construction of $N_R(\sigma)$ ported to the hyperbolic space $\mathbb H_{-\Lambda}^2$ of constant negative curvature
$-\Lambda$.  That is: first draw a geodesic segment $\tilde\sigma(t)[0,|\sigma|] \rightarrow \mathbb H_{-\Lambda}^2$ of length $|\sigma|$;
then for each $x = \exp_{\sigma(t)}(v)\in N_R(\sigma)$, add the point $\tilde x=\exp_{\tilde \sigma (t)}\tilde v \in \mathbb
H_{-\Lambda}^2$, where $\tilde v \perp  \dot{\tilde\sigma}(t)$, with the same orientation, and $|\tilde v|=|v|$.

\begin{figure}[ht]
\begin{center}
\def\svgwidth{0.9\textwidth}
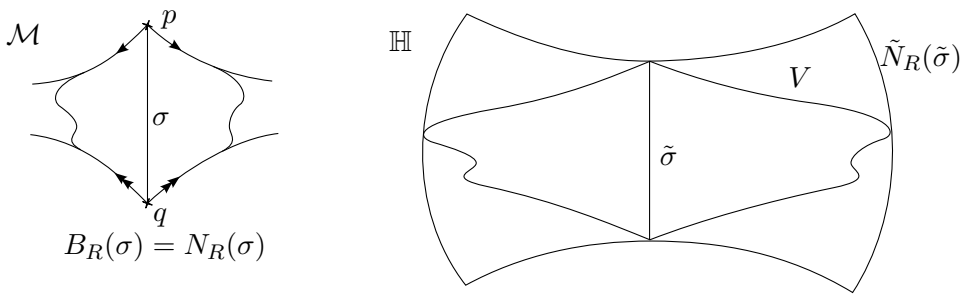
\caption{Comparison of normal neighbourhoods.}
\end{center}
\end{figure}

Let $\tilde N_R(\tilde\sigma)\subset \mathbb H_{-\Lambda}^2$ be a normal neighbourhood of $\tilde \sigma$. Then it must contain $V$, that
is, $V\subset \tilde
N_R(\tilde
\sigma)$. Choosing appropriate coordinates on $\mathbb H_{-\Lambda}^2$, the hyperbolic metric (of curvature $-\Lambda$) can be written
$$dx^2 + \cosh^2(\sqrt{\Lambda}x)\ dy^2$$
and $\tilde\sigma$ is the curve $\{y=0\}$ with arc-parameter $x$. We can compute the area of $\tilde N_R(\tilde\sigma)$,

$$Area(\tilde N_R(\tilde\sigma)) = \int_{-R}^R \int_0^{|\sigma|} \cosh(\sqrt{\Lambda}x)\ dx\ dy = 2 |\sigma|
\frac{\sinh(\sqrt{\Lambda}R)}{\sqrt{\Lambda}} .$$

Now we apply a comparison theorem. For $R=D+i_0$, $N_R(\sigma)$ contains $B(x_0,i_0)$ a smooth regular ball, whose area is bounded below by
$C(i_0,\Lambda)$, the area of a ball with same radius in the spherical space of curvature $+\Lambda$. Then,

\begin{align*}
C(i_0,\Lambda)&\leq \Area(B(x_0,i_0)) 
\leq \Area (N_R(\sigma)) \leq \Area(V)  \\
&\leq \Area (\tilde N_R(\tilde\sigma)) = 2 |\sigma| \frac{\sinh(\sqrt{\Lambda}R)}{\sqrt{\Lambda}} 
\leq 2 |\sigma| \frac{\sinh(\sqrt{\Lambda}(D+i_0))}{\sqrt{\Lambda}}  
\end{align*}
Thus, 
$$|\sigma|\geq \frac{C(i_0,\Lambda) \sqrt{\Lambda}}{2\sinh(\sqrt{\Lambda}(D+i_0))}$$
and this bounds below the length of $\sigma$.

In order to bound $\inj(p)$, it suffices to consider $p=q$ and $\sigma$ a geodesic loop based on $p$, and the argument above applies. 
\end{proof}

\section{Injectivity radius bound for smooth points} \label{S:cptss:inj_smooth}

In this section we adapt the Injectivity radius decay lemma to the case of cone surfaces. Recall the lemma for the smooth case
(\cite{Cheeger_finiteness}, \cite{CheGroTay}),

\begin{lema}
For all $i_0, R,\Lambda>0$ there exists $C>0$ such that the following holds. Let $\mathcal M$ be a smooth Riemannian manifold, and let
$x_0,z\in\mathcal M$ such that
\begin{itemize}
 \item $|\sec|<\Lambda$.
 \item $\inj (x_0) >i_0$.
 \item $d(x_0,z)<R$.
\end{itemize}
Then
$$\inj(z)>C.$$
Furthermore, while $i_0$ and $\Lambda$ are fixed, $C=C(R)$ is a continuous and decreasing function.
\end{lema}

In our case, we must check that the result is valid on a cone surface, for smooth points away from the cone points.

\begin{lema} \label{L:cptss:inj_smooth}
For all $i_0, \delta_0,R,\Lambda>0$ there exists $C>0$ such that the following holds. Let $(\mathcal M,x_0)$ be a cone surface with smooth
base point $x_0$, and let $z\in\mathcal M$ smooth such that
\begin{itemize}
 \item cone angles $\leq \pi$.
 \item $\inj(x_0)>i_0$.
 \item $|\sec|<\Lambda$ on the smooth part.
 \item $z\in B_R(x_0)$ and $d(z,\Sigma)>\delta_0>0$.
\end{itemize}
Then
$$\inj(z)>C.$$
\end{lema}

The proof is analogous to the smooth case, it relies on three comparison results that we next recall. First result is that injectivity
radius controls volume on curvature bounded above.

\begin{prop}[Günther-Bishop inequality, cf. {\cite[Thm III.4.2]{Chavel06}} ] \label{prop_GBineq}
Let $\mathcal M$ be a smooth Riemannian manifold with curvature $<\Lambda $, and let $x_0\in \mathcal M$ and $i_0=\inj x_0$. Then
$$\vol B(x_0,i_0) > \vol_\Lambda(i_0) .$$
\end{prop}

The same result holds easily for cone surfaces.
\begin{prop} \label{prop_GBineq_bis}
Let $\mathcal M$ be a cone surface with curvature $<\Lambda $, cone angles $>\alpha_0>0$, and let $x_0\in \mathcal M$ and
$i_0=\inj x_0$. Then
$$\vol B(x_0,i_0) > C(\alpha_0,i_0,\Lambda)$$
\end{prop}

Second result is that the volume of a ball controls the volume of a larger ball, on curvature bounded below.
\begin{prop}[Bishop-Gromov inequality {\cite[Thm 10.6.6]{BurBurIva}} ] \label{prop_BGineq}
Let $X$ be an $n$-dimensional Alexandrov space of curvature $\geq \Lambda$, then the map
$$r\mapsto \frac{\vol B(z,r)}{\vol_\Lambda(r)}$$
is decreasing, i.e., if $r_1\leq r_2$
$$\frac{\vol B(x,r_1)}{\vol_\Lambda(r_1)}\geq \frac{\vol B(z,r_2)}{\vol_\Lambda(r_2)}$$ 
\end{prop}

Third result is that volume controls injectivity radius on bounded curvature.

\begin{prop}[Cheeger-Gromov-Taylor inequality {\cite[Thm 4.3]{CheGroTay}} ] \label{prop_CGTineq}
Let $\mathcal M$ be a smooth $n$-manifold with bounded curvature $\Lambda_1<\sec<\Lambda_2$ and let
$r_0<\frac{1}{4}\frac{\pi}{\sqrt{\Lambda_2}}$ (no bound if
$\Lambda_2<0$). Let $z\in\mathcal M$. Then
$$\inj(z)>\frac{r_0}{2}\frac{1}{1+\frac{\vol_{\Lambda_1}(2r_0)}{\vol B(z,r_0)}}$$
\end{prop}

This proposition works with radii of balls under the conjugacy radius, $\frac{\pi}{\sqrt{\Lambda_2}}$ (cf. Klingenberg's theorem
\cite[Cor 5.7]{CheegerEbin}). That is, in $B(z,r_0)$ there are no conjugate points, and thus this result tells when the injectivity radius
can be realized as a geodesic loop based on $z$, given the topology of the ball.

We now prove our conic version of the injectivity radius decay.

\begin{proof}[Proof (of Lemma \ref{L:cptss:inj_smooth})]
We start with a bound on the injectivity radius of the basepoint, $i_0=\inj(x_0)$. By Proposition \ref{prop_GBineq} we can estimate
$\vol(B(x_0,i_0)$,

$$\vol B(x_0,i_0) \geq \vol_\Lambda (i_0) =C_1(\Lambda ,i_0)>0$$
and this bounds below $\vol B(z,2R)$ since $R>d(x_0,z)$ and then $B(x_0,i_0)\subseteq B(z,2R)$,
$$\vol B(z,2R)\geq \vol B(x_0,i_0) > C_1(\Lambda ,i_0) .$$
Now let $\delta = \frac{1}{2}\min\{\frac{1}{4}\frac{\pi}{\sqrt{\Lambda }},\delta_0\}$. We can use Proposition \ref{prop_BGineq} to bound
$\vol B(z,\delta)$,
$$\vol B(z,\delta) \geq C_2(\Lambda,R,\delta) \vol B(z,2R) >  C_3(R,\delta,\Lambda ,i_0) .$$
Finally we apply Proposition \ref{prop_CGTineq} to bound $\inj(z)$. In order to use this proposition on a cone surface, we need to
assume that a ball with radius under the conjugacy radius does not contain cone points (a pencil of geodesic rays through a cone point
converges instantaneously). We apply Proposition \ref{prop_CGTineq} with $r_0=\delta$, which justifies the choose of this $\delta$.
$$\inj(z)>\frac{\delta}{2}\ \frac{1}{1+\frac{\vol_\Lambda (2\delta)}{\vol B(z,\delta)}} \geq \frac{\delta}{2}\ \frac{1}{1+\frac{\vol_\Lambda
(2\delta)}
{C_3(R,\delta,\Lambda ,i_0)}} =: C(i_0,R,\Lambda ,\delta_0) .$$

\begin{figure}[ht]
\begin{center}
\def\svgwidth{0.6\textwidth}
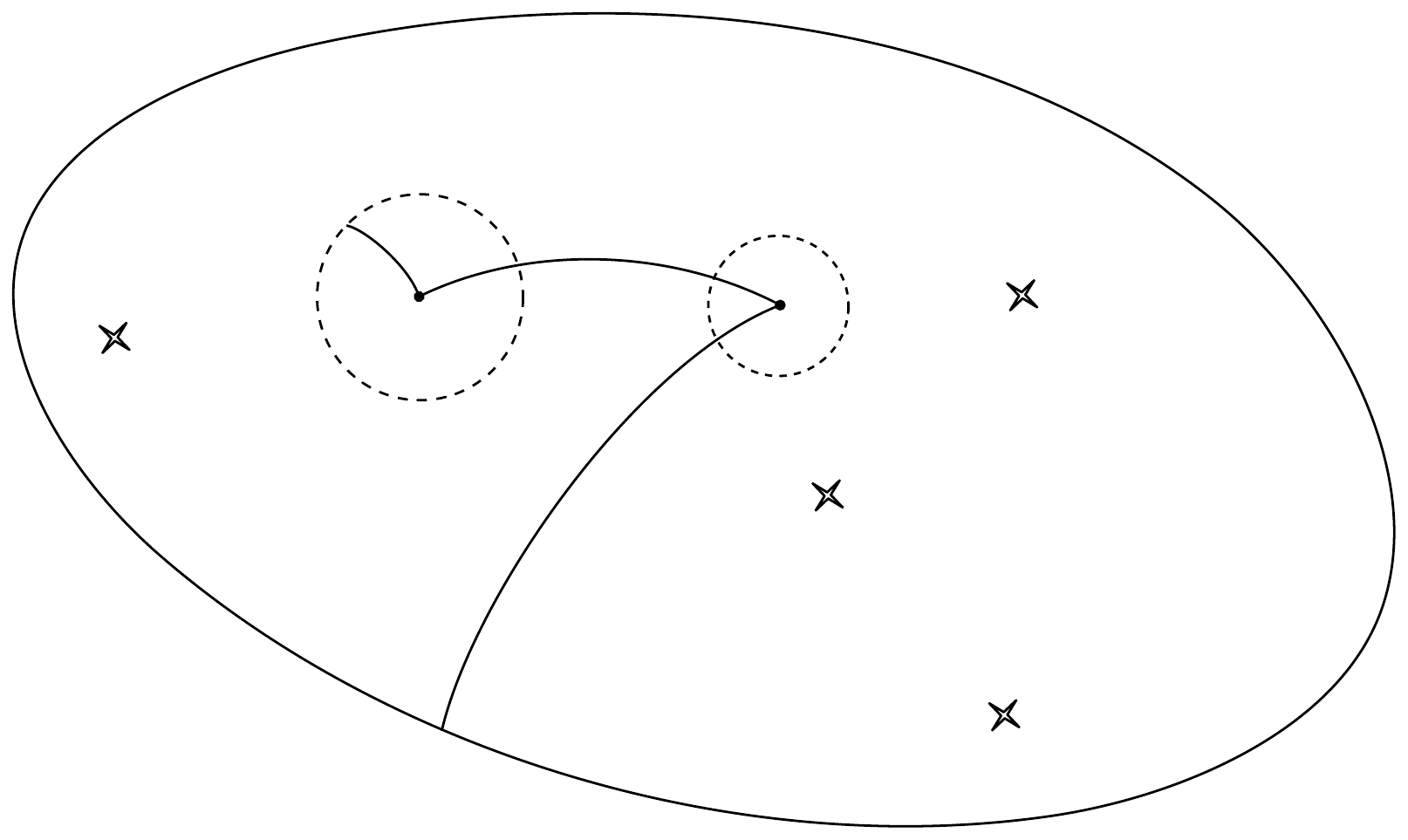
\caption{Controlled injectivity radius decay.}
\end{center}
\end{figure}
\end{proof}

\section{Proof of the Compactness Theorem for cone surfaces} \label{S:cptss:cptss_cone_proof}
\begin{proof}[Proof (of Theorem \ref{T:cptss:conesurfs})]
Let $\{(\mathcal M_k,x_k)\}_{k=1}^\infty$ a sequence of pointed cone surfaces inside the class $\mathfrak M$. By Gromov's compactness
theorem, there is a convergent subsequence to a pointed metric space $(\mathcal X_\infty, x_\infty)$ in the pointed Gromov-Hausdorff
topology. Since we are not assuming
compact
surfaces, the pointed convergence is relevant. Gromov's theorem states that for all $R>0$ there is a convergent subsequence
$$\mathcal M_k \cap B(x_k,R) \overset{GH}{\longrightarrow} \mathcal X_\infty \cap B(x_\infty,R)$$
in the Gromov-Hausdorff topology.

To simplify the notation, we will use $M_k=\mathcal M_k \cap B(x_k,R)$ and $X_k=\mathcal X_k \cap B(x_k,R)$. If $\mathcal M_k$ are compact,
we can take $\mathcal M_k =M_k$;
otherwise, $M_k$ are open noncomplete cone surfaces, and $M_k \rightarrow X_\infty$ in the Gromov-Hausdorff sense.

By Lemma \ref{L:cptss:conestruct}, there is a subsequence such that $M_k$ all have the same number of cone points, and their cone angles
form convergent sequences. By Lemma \ref{L:cptss:inj_cone}, the cone points are separated by a uniform distance. Hence we can pick
$\epsilon=\epsilon(R)>0$ not depending on $k$, such that
$$\epsilon < \frac{1}{10} \inf \{d(p,q) \ : \  p,q\in\Sigma_k, \ k=1\ldots\infty \} .$$

We remove on each surface a neighbourhood of radius $\epsilon$ of each cone point, and we form the sequence 
$$M_k^\epsilon := M_k\setminus \bigcup_{x\in\Sigma}B(x,\epsilon)$$
of noncomplete smooth Riemannian surfaces. By Lemma \ref{L:cptss:inj_smooth} and Theorem \ref{thm_compactness_manifolds}, 
$$M_k^\epsilon \underset{k\rightarrow \infty}{\longrightarrow} M_\infty^\epsilon$$
in the $\mathcal C^\infty$ sense, and $M_\infty^\epsilon \subset X_\infty$. If we reduce $\epsilon$ to $\epsilon/2$, we obtain another
$M_\infty^{\frac{\epsilon}{2}}$. Since $M_k^\epsilon \subset M_k^{\frac{\epsilon}{2}}$, this passes to the limit as $M_\infty^\epsilon
\subset M_\infty^{\frac{\epsilon}{2}}$, and both differential structures are compatible since an open covering of $M_\infty^\epsilon$ is
extended to an open covering of $M_\infty^{\frac{\epsilon}{2}}$, where the radius of the covering balls only depends of the distance to the
basepoint and the singular set. Therefore,
$$M_k^{\frac{\epsilon}{2^l}}\longrightarrow M_k$$
in the $\mathcal C^\infty$ sense as $l\rightarrow \infty$, because is a nested sequence, and hence
$$M_k^{\frac{\epsilon}{2^l}}\longrightarrow M_\infty^{\frac{\epsilon}{2^l}}$$
in the $\mathcal C^\infty$ sense as $k\rightarrow \infty$, for all $l$. Since $\{M_\infty^{\frac{\epsilon}{2^l}}\}_l$ is also a nested
sequence, we can take a subsequence of the diagonal sequence
$$M_k^{\frac{\epsilon}{2^k}}\longrightarrow M_\infty^{(0)}$$
in the $\mathcal C^\infty$ sense, where $M_\infty^{(0)} = \bigcup_{l>0} M_\infty^{\frac{\epsilon}{2^l}}$. This limit $M_\infty^{(0)}$ is a
smooth noncomplete surface inside the metric space $X_\infty$. We define $\Sigma_\infty = X_\infty \setminus M_\infty^{(0)}$.

\medskip

The only remaining issue is to check that $\Sigma_\infty$ consist of cone points with a model cone metric. Let $p_k\in M_k$ be a sequence
of cone points with a limit cone angle. As discussed in Chapter \ref{Ch:cone_rf}, in a neighbourhood of $p_k$ we can write the metric of
$M_k$ as
$$g_k=dr^2 + h_k^2(r,\theta)\ d\theta^2$$
with 
\begin{itemize}
 \item $\displaystyle h(0,\theta)= 0$,
 \item $\displaystyle \frac{\partial h_k}{\partial r} (0,\theta)= \frac{\alpha_k}{2\pi}$ where $\alpha_k$ is the cone angle,
 \item $\displaystyle \frac{\partial^2 h_k}{\partial r^2} (0,\theta)=0$
\end{itemize}
independently if $\theta$. Last condition follows from bounded curvature, since
$$\frac{\partial^2 h_k}{\partial r^2} = -K_k(r,\theta) \ h_k .$$
In other words, the function $h_k$ can be written as
$$h_k=\frac{\alpha_k}{2\pi} r + O(r^3)$$
and the function in $O(r^3)$ may depend on $\theta$. By the convergence on the smooth part seen above, for any $r>0$ we have
$$h_k \rightarrow h_\infty$$
as $k\rightarrow \infty$, and the convergence is uniform on compact sets not containing $r=0$. On the other hand, from Lemma
\ref{L:cptss:conestruct} we have
$$\alpha_k \rightarrow \alpha$$
as $k\rightarrow \infty$. We must check
\begin{itemize}
 \item $\displaystyle \lim_{r\rightarrow 0} h_\infty(r,\theta)=0$,
 \item $\displaystyle \lim_{r\rightarrow 0} \frac{\partial h_\infty}{\partial r} (r,\theta)=\frac{\alpha}{2\pi}$,
 \item $\displaystyle \lim_{r\rightarrow 0} \frac{\partial^2 h_\infty}{\partial r^2} (r,\theta)=0$.
\end{itemize}

For the third point, since the Riemannian curvature on the smooth part of the surfaces is uniformly bounded by hypothesis, also is the
Riemannian curvature in the limit surface, and hence the second derivative on $r$ is zero if and only if $h_\infty$ is zero. Hence it
follows from first point.

For the first point, the area element of the metrics $g_k$ is 
$$h_k(r,\theta)\ dr\wedge d\theta .$$
Since the curvature is bounded below by $-\Lambda$, by
comparison the volume element is less than the volume element of the hyperbolic space of curvature $-\Lambda$, namely
$$\frac{\sinh(\sqrt{\Lambda}r)}{\sqrt{\Lambda}} \ dr\wedge d\theta$$
thus, $$h_k(r,\theta)\leq \frac{\sinh(\sqrt{\Lambda}r)}{\sqrt{\Lambda}}.$$
This bound is uniform in $k$, and hence also applies in the limit $h_\infty$. Therefore
$$\lim_{r\rightarrow 0} |h_\infty(r,\theta)| \leq \lim_{r\rightarrow 0} \frac{\sinh(\sqrt{\Lambda}r)}{\sqrt{\Lambda}} =0 .$$

Finally we check the second point. From the fundamental theorem of calculus,
$$  \frac{\partial h_k}{\partial r}(r,\theta) = \int_0^r \frac{\partial^2 h_k}{\partial t^2} (t,\theta) \ dt + \frac{\partial h_k}{\partial
r} (0,\theta) $$
Then,
\begin{align*}
\left| \frac{\partial h_k}{\partial r}(r,\theta) - \frac{\alpha_k}{2\pi}\right|
&\leq \int_0^r \left| \frac{\partial^2 h_k}{\partial t^2} (t,\theta) \right| \ dt \\
&\leq \int_0^r |K (t,\theta) h_k(t,\theta)| \ dt \\
&\leq \Lambda \int_0^r |h_k(t,\theta)| \ dt \\
&\leq \Lambda \int_0^r \frac{\sinh(\sqrt{\Lambda}r)}{\sqrt{\Lambda}} \ dt \\
&= \Lambda \frac{-1+\cosh(\sqrt{\Lambda}r)}{\Lambda} = \Lambda \frac{1}{2}r^2 + O(r^4) \leq \Lambda r^2
\end{align*}
for $r$ small. This bound is uniform in $k$. Thus,
\begin{align*}
 \left| \frac{\partial h_\infty}{\partial r}(r,\theta) - \frac{\alpha}{2\pi}\right| 
&\leq \left| \frac{\partial h_\infty}{\partial r}(r,\theta) - \frac{\partial h_k}{\partial r}(r,\theta) \right| + 
 \left| \frac{\partial h_k}{\partial r}(r,\theta) - \frac{\alpha_k}{2\pi} \right| +
 \left| \frac{\alpha_k}{2\pi} - \frac{\alpha}{2\pi} \right| \\
&\leq \epsilon(r,k) + \Lambda r^2 + \epsilon'(k)
\end{align*}
Now, as $k\rightarrow\infty$, we have $\epsilon(r,k),\epsilon'(k) \rightarrow 0$, and hence
$$\left| \frac{\partial h_\infty}{\partial r}(r,\theta) - \frac{\alpha}{2\pi}\right| \leq  \Lambda r^2 $$
so
$$\lim_{k\rightarrow\infty}\left| \frac{\partial h_\infty}{\partial r}(r,\theta) - \frac{\alpha}{2\pi}\right| =0 .$$
This proves the convergence of the cone structure. The sequence $h_k$ converges in the $\mathcal C^1$ sense by Arzelà-Ascoli theorem, so the
metric tensors $g_k$ converge in $\mathcal C^1$ on that coordinate chart, and this implies Lipschitz convergence of the sequence of (metric)
surfaces. This finishes the proof of the theorem.
\end{proof}

\section{The compactness theorem for flows on cone surfaces} \label{S:cptss:cptss_flows}
We now move forwards to Ricci flows. Hamilton's compactness theorem for solutions of Ricci flow \cite{Hamilton_compactness} is the actual
result that is needed
for the analysis of singularities of the flow. Given a sequence of Riemannian manifolds, it induces a sequence of Ricci flows starting at
these
initial manifolds, and one expects it to converge to a limit Ricci flow starting at the limit of the original manifolds.
Let us recall the theorem in the smooth case.

\begin{tma}[Hamilton's compactness of flows]
Let $(M_k, O_k, g_k(t))$ be a sequence of complete pointed solutions to the Ricci flow, with $t\in(a,b)$, such that:
\begin{enumerate}
\item $|\Rm|\leq C_0 $ on $M_k \times (a,b) \ \forall k $,
\item $\inj_{g_k(0)}(O_k)\geq i_0 \ \forall k$.
\end{enumerate}
Then there exists a converging subsequence to a pointed flow $(M_\infty, O_\infty, g_\infty(t))$ with $t\in(a,b)$ in the $\mathcal C^\infty$
sense.
\end{tma}

The remarkable point is that the hypothesis can be relaxed from asking bounds on all the derivatives of the curvature to just
asking a bound on the curvature. This is due to the Bernstein-Bando-Shi estimates, \cite{Shi}, \cite[Ch 7]{ChowKnopf}, i.e. the Ricci flow
equation relates space and
time derivatives, and derivating the equation the regularity propagates to higher derivatives bounds. These bounds, however, depend on time
and get worse as $t\rightarrow 0$.

An outline of the proof of the compactness theorem for flows is the following. First one gets uniform bounds on the space and time
derivatives of the $g_k(t)$ metrics with respect to a background metric. This is contained in the following proposition, that includes
inside the Bernstein-Bando-Shi estimates.

\begin{prop}[{\cite[Lem 3.11]{TRFTA1}}]
Let $(\mathcal M,g)$ be a Riemannian $n$-manifold, and let $K\subset \mathcal M$ be a compact subset. Let $g_k(t)$ be a sequence of Ricci
flows defined on a neighbourhood of $K\times [T_1,T_2]$, where $t_0\in[T_1,T_2]$. Suppose
\begin{itemize}
 \item $C^{-1} g \leq g_k(t_0) \leq C g$ as quadratic forms, for some $C$ independent of $k$.
 \item $|\nabla^p g_k (t_0)| \leq C_p$, for all $k$ and $p$, for some $C_p$ independent of $k$, where $\nabla$ and the norm are respect to
$g$.
 \item $|\nabla_k^p \Rm_k|_k \leq C_p'$ on $K\times[T_1,T_2]$, for all $k$ and $p$, for some $C_p'$ independent of $k$, where $\nabla_k$ and
the norm are respect to $g_k$.
\end{itemize}
Then,
\begin{itemize}
 \item $B^{-1} g \leq g_k(t) \leq B g$, and
 \item $|\frac{\partial^q}{\partial t^q} \nabla^p g_k(t)| \leq \tilde C_{p,q}$,
\end{itemize}
on $K\times [T_1,T_2]$, where $B=B(t,t_0)=Ce^{2\sqrt{n-1}C_0'|t-t_0|}$ and $\tilde C_{p,q}$ independent of $k$.
\end{prop}

To finish the proof of Hamilton's compactness theorem for flows, one applies the $\mathcal C^\infty$ version of Arzelà-Ascoli theorem
(Theorem \ref{prop_ArzAsc_metrictensor}) to the
$(n+1)$-manifolds, using the obtained bounds on the metric and its derivatives in space and time.

So far it is the smooth case as done by Hamilton in \cite{Hamilton_compactness}, cf. \cite{TRFTA1}. We now turn into the case of flows on
cone surfaces.
The previous argument applies to neighbourhoods of smooth points (all bounds, such as Shi's estimates can be done locally). However, it
remains to
check the convergence for the functions defining the flow on a tubular neighbourhood of a cone point. Recall from Chapter \ref{Ch:cone_rf}
that the
existence of the angle preserving flow is given by the work in \cite{MazRubSes}. This flow preserves the cone angles and has, at least for
short time, bounded curvature and derivative of the curvature.

%

The second main theorem of this appendix is the following.

\begin{tma} \label{T:cptss:flows}
Let $(\mathcal M_k, g_k(t),O_k)$, with $t\in(a,b)$, be a sequence of pointed cone surfaces evolving according to the angle-preserving Ricci
flow. Assume that
\begin{enumerate}
 \item all cone angles are less than or equal to $\pi$,
 \item $|Rm_x| \leq C_0 $ for all $ x \in \mathcal M_k \times (a,b)$ , with $x$ not a cone point,
 \item $\inj_{g_k(0)} O_k \geq i_0$, if $O_k\notin \Sigma_k$.,
 \item Alternatively, $\inj_{g_k(0)} O_k \geq i_1$ and $\alpha>\alpha_0>0$, if $O_k$ is a cone point of angle $\alpha$.
\end{enumerate}
Then there exists a convergent subsequence to a pointed limit flow $(\mathcal M_\infty,g_\infty(t),O_\infty)$ in the $\mathcal C^\infty$
sense on the smooth part, and Lipschitz on the singular points.
\end{tma}

\begin{proof}
As noted before, the only issue we need to check is that, in a neighbourhood of a cone point of fixed angle $\alpha$, a sequence of
angle-preserving flows subconverges to an angle-preserving flow.

From the short-time existence theorem, for every initial metric $g_k(0)$ there exists an angle-preserving flow $g_k(t)$ for $t\in(a,b)$ a
uniform time interval since the curvature is bounded by hypothesis. By Proposition \ref{P:RF_in_gpolarcoords}, these metrics can be
written in geodesic polar coordinates around a cone point as
$$g_k(t)=dr^2 + h^k(r,\theta,t)^2 \ d\theta^2$$
and these functions satisfy the system
\begin{equation} \label{rf_polar_2}
 \left\{ \begin{array}{ll}
          h^k_{t}=h^k_{rr}-h^k_r \int{ \frac{h^k_{rr}}{h^k}} & \\
	  h^k = 0 &\mbox{ on } r=0 .\\
	  \end{array}
\right.
\end{equation}

We only need to apply the Arzelà-Ascoli theorem to the functions $h^k$ to get uniform convergence over compact sets. Note that these
compact sets may contain the cone point itself, since the functions are defined for $r\in[0,r_0)$. This is a very convenient property that
we would not have in conformal coordinates.

From the proof of Theorem \ref{T:cptss:conesurfs}, $h^k$ is uniformly bounded,
$$h^k \leq \frac{\sinh(\sqrt{C_0}r)}{\sqrt{C_0}} .$$
Also, the space derivatives of $h^k$ are also uniformly bounded since
$$\left| h^k_r - \frac{\alpha}{2\pi}\right| \leq C_0 r^2 .$$
Further, the second derivative $h^k_{rr}$ is also uniformly bounded since the curvature $K= - \frac{h^k_{rr}}{h^k}$ is uniformly bounded by
hypothesis.

Finally, the Ricci flow equation, written in the form of \eqref{rf_polar_2}, ensures that the time derivative of $h^k$ is uniformly
bounded. Hence, we have uniform bounds on space and time derivatives and therefore, by Arzelà-Ascoli theorem
there is a convergent subsequence of $h^k$ to a limit $h^\infty$. Up to now, we have that the limit function is $\mathcal C^1$ in space and
$\mathcal C^0$ in time. The regularity can be improved by checking second derivatives. From Proposition \ref{prop_fixedangle}, we know that
$h	^k_{rt}$ is
uniformly bounded and vanishes for $r=0$, hence the limit flow is also an angle-preserving flow and the cone angle is the same by the
Gromov-Hausdorff (and stronger) convergence of the time-slice surfaces. Similarly, higher derivatives (in space and time) of $h^k$ can be
derived from expressions on the space derivatives of $K=-h_{rr}/h$ (each time derivative is translated into two space derivatives by the
flow equation). But bounds on the derivatives of the curvature can be obtained from the bounds of the curvature itself by using
Bernstein-Bando-Shi estimates, cf. \cite[Thm 7.1]{ChowKnopf}. The proof of these estimates relies only on local computations in coordinates
and in the application of a maximum principle for functions. The maximum principle holds on cone surfaces by Theorem \ref{T:maxpple_cone},
and hence the
Bernstein-Bando-Shi estimates apply. This implies uniform bounds on the derivatives of $h^k$ and therefore the convergence of the metric
tensor on these coordinates is $\mathcal C^\infty$ both in space and time. 

From an intrinsic point of view, we can only claim Lipschitz
convergence at the singular points, since no tangent vectors exist at these points in the smooth Riemannian sense. The convergence is
$\mathcal C^\infty$ in the smooth part of the flow.
\end{proof}

%% file: comparison_triang.pdf_tex
\begingroup%
  \makeatletter%
  \providecommand\color[2][]{%
    \errmessage{(Inkscape) Color is used for the text in Inkscape, but the package 'color.sty' is not loaded}%
    \renewcommand\color[2][]{}%
  }%
  \providecommand\transparent[1]{%
    \errmessage{(Inkscape) Transparency is used (non-zero) for the text in Inkscape, but the package 'transparent.sty' is not loaded}%
    \renewcommand\transparent[1]{}%
  }%
  \providecommand\rotatebox[2]{#2}%
  \ifx\svgwidth\undefined%
    \setlength{\unitlength}{566.94956055bp}%
    \ifx\svgscale\undefined%
      \relax%
    \else%
      \setlength{\unitlength}{\unitlength * \real{\svgscale}}%
    \fi%
  \else%
    \setlength{\unitlength}{\svgwidth}%
  \fi%
  \global\let\svgwidth\undefined%
  \global\let\svgscale\undefined%
  \makeatother%
  \begin{picture}(1,0.77409669)%
    \put(0,0){\includegraphics[width=\unitlength]{comparison_triang.pdf}}%
    \put(0.71531457,0.72631819){\color[rgb]{0,0,0}\makebox(0,0)[lb]{\smash{$\mathbb H$}}}%
    \put(0.71600907,0.25916531){\color[rgb]{0,0,0}\makebox(0,0)[lb]{\smash{$\mathbb S$}}}%
    \put(0.01463716,0.72631819){\color[rgb]{0,0,0}\makebox(0,0)[lb]{\smash{$\mathcal M$}}}%
    \put(0.01606253,0.24633684){\color[rgb]{0,0,0}\makebox(0,0)[lb]{\smash{$\bar{\mathcal M}$}}}%
    \put(0.47416308,0.57815145){\color[rgb]{0,0,0}\makebox(0,0)[lb]{\smash{$l<\tilde l$}}}%
    \put(0.46133466,0.10492327){\color[rgb]{0,0,0}\makebox(0,0)[lb]{\smash{$\alpha' < \tilde{\alpha'}$}}}%
    \put(0.72173574,0.61895403){\color[rgb]{0,0,0}\makebox(0,0)[lb]{\smash{$d$}}}%
    \put(0.8556807,0.61895403){\color[rgb]{0,0,0}\makebox(0,0)[lb]{\smash{$d$}}}%
    \put(0.78776542,0.52101989){\color[rgb]{0,0,0}\makebox(0,0)[lb]{\smash{$\tilde l$}}}%
    \put(0.78332545,0.01731093){\color[rgb]{0,0,0}\makebox(0,0)[lb]{\smash{$\bar l$}}}%
    \put(0.1447538,0.03687854){\color[rgb]{0,0,0}\makebox(0,0)[lb]{\smash{$\bar l$}}}%
    \put(0.79397704,0.65939847){\color[rgb]{0,0,0}\makebox(0,0)[lb]{\smash{$\alpha$}}}%
    \put(0.14181867,0.67519448){\color[rgb]{0,0,0}\makebox(0,0)[lb]{\smash{$\alpha$}}}%
    \put(0.14435326,0.19921683){\color[rgb]{0,0,0}\makebox(0,0)[lb]{\smash{$\alpha'$}}}%
    \put(0.78590797,0.1840835){\color[rgb]{0,0,0}\makebox(0,0)[lb]{\smash{$\tilde{\alpha'}$}}}%
    \put(0.07723086,0.62440643){\color[rgb]{0,0,0}\makebox(0,0)[lb]{\smash{$d$}}}%
    \put(0.20627895,0.62440643){\color[rgb]{0,0,0}\makebox(0,0)[lb]{\smash{$d$}}}%
    \put(0.13820567,0.5047893){\color[rgb]{0,0,0}\makebox(0,0)[lb]{\smash{$l$}}}%
    \put(0.04052477,0.16486327){\color[rgb]{0,0,0}\makebox(0,0)[lb]{\smash{$\bar d_1$}}}%
    \put(0.23838321,0.17433703){\color[rgb]{0,0,0}\makebox(0,0)[lb]{\smash{$\bar d_2$}}}%
    \put(0.68516616,0.17184786){\color[rgb]{0,0,0}\makebox(0,0)[lb]{\smash{$\bar d_1$}}}%
    \put(0.87969086,0.16832021){\color[rgb]{0,0,0}\makebox(0,0)[lb]{\smash{$\bar d_2$}}}%
    \put(0.11490063,0.38415322){\color[rgb]{0,0,0}\makebox(0,0)[lb]{\smash{$f$}}}%
    \put(0.18961366,0.38213744){\color[rgb]{0,0,0}\makebox(0,0)[lb]{\smash{$2\delta$-iso}}}%
  \end{picture}%
\endgroup%

%% file: cone_noncollapse.pdf_tex
\begingroup%
  \makeatletter%
  \providecommand\color[2][]{%
    \errmessage{(Inkscape) Color is used for the text in Inkscape, but the package 'color.sty' is not loaded}%
    \renewcommand\color[2][]{}%
  }%
  \providecommand\transparent[1]{%
    \errmessage{(Inkscape) Transparency is used (non-zero) for the text in Inkscape, but the package 'transparent.sty' is not loaded}%
    \renewcommand\transparent[1]{}%
  }%
  \providecommand\rotatebox[2]{#2}%
  \ifx\svgwidth\undefined%
    \setlength{\unitlength}{740.52836914bp}%
    \ifx\svgscale\undefined%
      \relax%
    \else%
      \setlength{\unitlength}{\unitlength * \real{\svgscale}}%
    \fi%
  \else%
    \setlength{\unitlength}{\svgwidth}%
  \fi%
  \global\let\svgwidth\undefined%
  \global\let\svgscale\undefined%
  \makeatother%
  \begin{picture}(1,0.31096469)%
    \put(0,0){\includegraphics[width=\unitlength]{cone_noncollapse.pdf}}%
    \put(0.16718627,0.2851885){\color[rgb]{0,0,0}\makebox(0,0)[lb]{\smash{$p$}}}%
    \put(0.15845606,0.08439342){\color[rgb]{0,0,0}\makebox(0,0)[lb]{\smash{$q$}}}%
    \put(0.06788002,0.05001817){\color[rgb]{0,0,0}\makebox(0,0)[lb]{\smash{$B_R(\sigma)=N_R(\sigma)$}}}%
    \put(0.01016139,0.26336294){\color[rgb]{0,0,0}\makebox(0,0)[lb]{\smash{$\mathcal M$}}}%
    \put(0.39616357,0.25665924){\color[rgb]{0,0,0}\makebox(0,0)[lb]{\smash{$\mathbb H$}}}%
    \put(0.15736478,0.18151712){\color[rgb]{0,0,0}\makebox(0,0)[lb]{\smash{$\sigma$}}}%
    \put(0.66917401,0.14113986){\color[rgb]{0,0,0}\makebox(0,0)[lb]{\smash{$\tilde\sigma$}}}%
    \put(0.79903608,0.21971184){\color[rgb]{0,0,0}\makebox(0,0)[lb]{\smash{$V$}}}%
    \put(0.89070331,0.24153739){\color[rgb]{0,0,0}\makebox(0,0)[lb]{\smash{$\tilde N_R(\tilde\sigma)$}}}%
  \end{picture}%
\endgroup%

%% file: inj_decay.pdf_tex
\begingroup%
  \makeatletter%
  \providecommand\color[2][]{%
    \errmessage{(Inkscape) Color is used for the text in Inkscape, but the package 'color.sty' is not loaded}%
    \renewcommand\color[2][]{}%
  }%
  \providecommand\transparent[1]{%
    \errmessage{(Inkscape) Transparency is used (non-zero) for the text in Inkscape, but the package 'transparent.sty' is not loaded}%
    \renewcommand\transparent[1]{}%
  }%
  \providecommand\rotatebox[2]{#2}%
  \ifx\svgwidth\undefined%
    \setlength{\unitlength}{464.68525391bp}%
    \ifx\svgscale\undefined%
      \relax%
    \else%
      \setlength{\unitlength}{\unitlength * \real{\svgscale}}%
    \fi%
  \else%
    \setlength{\unitlength}{\svgwidth}%
  \fi%
  \global\let\svgwidth\undefined%
  \global\let\svgscale\undefined%
  \makeatother%
  \begin{picture}(1,0.5965637)%
    \put(0,0){\includegraphics[width=\unitlength]{inj_decay.pdf}}%
    \put(0.75219084,0.36276475){\color[rgb]{0,0,0}\makebox(0,0)[lb]{\smash{$\Sigma$}}}%
    \put(0.60559916,0.22263283){\color[rgb]{0,0,0}\makebox(0,0)[lb]{\smash{$\Sigma$}}}%
    \put(0.73102969,0.06522979){\color[rgb]{0,0,0}\makebox(0,0)[lb]{\smash{$\Sigma$}}}%
    \put(0.56401182,0.36412504){\color[rgb]{0,0,0}\makebox(0,0)[lb]{\smash{$z$}}}%
    \put(0.50421417,0.44817174){\color[rgb]{0,0,0}\makebox(0,0)[lb]{\smash{$B(z,\delta)$}}}%
    \put(0.28982208,0.35254824){\color[rgb]{0,0,0}\makebox(0,0)[lb]{\smash{$x_0$}}}%
    \put(0.27741197,0.42007376){\color[rgb]{0,0,0}\makebox(0,0)[lb]{\smash{$i_0$}}}%
    \put(0.41034036,0.4170922){\color[rgb]{0,0,0}\makebox(0,0)[lb]{\smash{$R$}}}%
    \put(0.40619089,0.20721203){\color[rgb]{0,0,0}\makebox(0,0)[lb]{\smash{$2R$}}}%
    \put(0.04895584,0.28535065){\color[rgb]{0,0,0}\makebox(0,0)[lb]{\smash{$\Sigma$}}}%
    \put(0.20424984,0.48351209){\color[rgb]{0,0,0}\makebox(0,0)[lb]{\smash{$B(x_0,i_0)$}}}%
    \put(0.79853595,0.51220682){\color[rgb]{0,0,0}\makebox(0,0)[lb]{\smash{$B(z,2R)$}}}%
  \end{picture}%
\endgroup%

%% file: tesis.bbl
\newcommand{\etalchar}[1]{$^{#1}$}
\begin{thebibliography}{CCG{\etalchar{+}}08}

\bibitem[Aub98]{Aubin}
Thierry Aubin.
\newblock {\em Some nonlinear problems in {R}iemannian geometry}.
\newblock Springer Monographs in Mathematics. Springer-Verlag, Berlin, 1998.

\bibitem[Bai09]{Baird}
Paul Baird.
\newblock A class of three-dimensional {R}icci solitons.
\newblock {\em Geom. Topol.}, 13(2):979--1015, 2009.

\bibitem[BBB{\etalchar{+}}10]{BBBMP}
Laurent Bessi{\`e}res, G{\'e}rard Besson, Michel Boileau, Sylvain Maillot, and
  Joan Porti.
\newblock {\em Geometrisation of 3-manifolds}, volume~13 of {\em EMS Tracts in
  Mathematics}.
\newblock European Mathematical Society (EMS), Z\"urich, 2010.

\bibitem[BBI01]{BurBurIva}
Dmitri Burago, Yuri Burago, and Sergei Ivanov.
\newblock {\em A course in metric geometry}, volume~33 of {\em Graduate Studies
  in Mathematics}.
\newblock American Mathematical Society, Providence, RI, 2001.

\bibitem[Bry]{Bryant}
Robert~L. Bryant.
\newblock {R}icci flow solitons in dimension three with {SO(3)}-symmetries.
\newblock {\em Unpublished}.
\newblock \url{http://www.math.duke.edu/~bryant/3DRotSymRicciSolitons.pdf}.

\bibitem[Cao96]{Cao_KRsolitons1}
Huai-Dong Cao.
\newblock Existence of gradient {K}\"ahler-{R}icci solitons.
\newblock In {\em Elliptic and parabolic methods in geometry ({M}inneapolis,
  {MN}, 1994)}, pages 1--16. A K Peters, Wellesley, MA, 1996.

\bibitem[Cao97]{Cao_KRsolitons2}
Huai-Dong Cao.
\newblock Limits of solutions to the {K}\"ahler-{R}icci flow.
\newblock {\em J. Differential Geom.}, 45(2):257--272, 1997.

\bibitem[CCC{\etalchar{+}}]{Caoetal}
Huai-Dong Cao, Giovanni Catino, Quiang Chen, Carlo Mantegazza, and Lorenzo
  Mazzieri.
\newblock Bach-flat gradient steady ricci solitons.
\newblock {\em Preprint (2011), arXiv:1107.4591 [math.DG]}.

\bibitem[CCCY03]{3CY}
H.~D. Cao, B.~Chow, S.~C. Chu, and S.~T. Yau, editors.
\newblock {\em Collected papers on {R}icci flow}, volume~37 of {\em Series in
  Geometry and Topology}.
\newblock International Press, Somerville, MA, 2003.

\bibitem[CCG{\etalchar{+}}07]{TRFTA1}
Bennett Chow, Sun-Chin Chu, David Glickenstein, Christine Guenther, James
  Isenberg, Tom Ivey, Dan Knopf, Peng Lu, Feng Luo, and Lei Ni.
\newblock {\em The {R}icci flow: techniques and applications. {P}art {I}},
  volume 135 of {\em Mathematical Surveys and Monographs}.
\newblock American Mathematical Society, Providence, RI, 2007.
\newblock Geometric aspects.

\bibitem[CCG{\etalchar{+}}08]{TRFTA2}
Bennett Chow, Sun-Chin Chu, David Glickenstein, Christine Guenther, James
  Isenberg, Tom Ivey, Dan Knopf, Peng Lu, Feng Luo, and Lei Ni.
\newblock {\em The {R}icci flow: techniques and applications. {P}art {II}},
  volume 144 of {\em Mathematical Surveys and Monographs}.
\newblock American Mathematical Society, Providence, RI, 2008.
\newblock Analytic aspects.

\bibitem[CE08]{CheegerEbin}
Jeff Cheeger and David~G. Ebin.
\newblock {\em Comparison theorems in {R}iemannian geometry}.
\newblock AMS Chelsea Publishing, Providence, RI, 2008.
\newblock Revised reprint of the 1975 original.

\bibitem[CGT82]{CheGroTay}
Jeff Cheeger, Mikhail Gromov, and Michael Taylor.
\newblock Finite propagation speed, kernel estimates for functions of the
  {L}aplace operator, and the geometry of complete {R}iemannian manifolds.
\newblock {\em J. Differential Geom.}, 17(1):15--53, 1982.

\bibitem[Cha06]{Chavel06}
Isaac Chavel.
\newblock {\em Riemannian geometry}, volume~98 of {\em Cambridge Studies in
  Advanced Mathematics}.
\newblock Cambridge University Press, Cambridge, second edition, 2006.
\newblock A modern introduction.

\bibitem[Che55]{Chern}
Shiing-shen Chern.
\newblock An elementary proof of the existence of isothermal parameters on a
  surface.
\newblock {\em Proc. Amer. Math. Soc.}, 6:771--782, 1955.

\bibitem[Che70]{Cheeger_finiteness}
Jeff Cheeger.
\newblock Finiteness theorems for {R}iemannian manifolds.
\newblock {\em Amer. J. Math.}, 92:61--74, 1970.

\bibitem[Cho91a]{Chow_orbifolds}
Bennett Chow.
\newblock On the entropy estimate for the {R}icci flow on compact
  {$2$}-orbifolds.
\newblock {\em J. Differential Geom.}, 33(2):597--600, 1991.

\bibitem[Cho91b]{Chow_sphere}
Bennett Chow.
\newblock The {R}icci flow on the {$2$}-sphere.
\newblock {\em J. Differential Geom.}, 33(2):325--334, 1991.

\bibitem[CK04]{ChowKnopf}
Bennett Chow and Dan Knopf.
\newblock {\em The {R}icci flow: an introduction}, volume 110 of {\em
  Mathematical Surveys and Monographs}.
\newblock American Mathematical Society, Providence, RI, 2004.

\bibitem[CLT06]{ChenLuTian}
Xiuxiong Chen, Peng Lu, and Gang Tian.
\newblock A note on uniformization of {R}iemann surfaces by {R}icci flow.
\newblock {\em Proc. Amer. Math. Soc.}, 134(11):3391--3393 (electronic), 2006.

\bibitem[CW91]{ChowWu}
Bennett Chow and Lang-Fang Wu.
\newblock The {R}icci flow on compact {$2$}-orbifolds with curvature negative
  somewhere.
\newblock {\em Comm. Pure Appl. Math.}, 44(3):275--286, 1991.

\bibitem[CZ06]{ChenZhu}
Bing-Long Chen and Xi-Ping Zhu.
\newblock Uniqueness of the {R}icci flow on complete noncompact manifolds.
\newblock {\em J. Differential Geom.}, 74(1):119--154, 2006.

\bibitem[DeT83]{DeTur}
Dennis~M. DeTurck.
\newblock Deforming metrics in the direction of their {R}icci tensors.
\newblock {\em J. Differential Geom.}, 18(1):157--162, 1983.

\bibitem[DLA06]{DumLliArt}
Freddy Dumortier, Jaume Llibre, and Joan~C. Art{\'e}s.
\newblock {\em Qualitative theory of planar differential systems}.
\newblock Universitext. Springer-Verlag, Berlin, 2006.

\bibitem[Gro07]{Gromov}
Misha Gromov.
\newblock {\em Metric structures for {R}iemannian and non-{R}iemannian spaces}.
\newblock Modern Birkh\"auser Classics. Birkh\"auser Boston Inc., Boston, MA,
  english edition, 2007.
\newblock Based on the 1981 French original, With appendices by M. Katz, P.
  Pansu and S. Semmes, Translated from the French by Sean Michael Bates.

\bibitem[GT]{Topgie}
Gregor Giesen and Peter Topping.
\newblock Existence of ricci flows of incomplete surfaces.
\newblock {\em Preprint (2010), to appear in Comm. P.D.E., arXiv:1007.3146v2
  [math.AP]}.

\bibitem[GW88]{GreeneWu}
R.~E. Greene and H.~Wu.
\newblock Lipschitz convergence of {R}iemannian manifolds.
\newblock {\em Pacific J. Math.}, 131(1):119--141, 1988.

\bibitem[Ham82]{Hamilton_3mfds}
Richard~S. Hamilton.
\newblock Three-manifolds with positive {R}icci curvature.
\newblock {\em J. Differential Geom.}, 17(2):255--306, 1982.

\bibitem[Ham88]{Hamilton_surfaces}
Richard~S. Hamilton.
\newblock The {R}icci flow on surfaces.
\newblock In {\em Mathematics and general relativity ({S}anta {C}ruz, {CA},
  1986)}, volume~71 of {\em Contemp. Math.}, pages 237--262. Amer. Math. Soc.,
  Providence, RI, 1988.

\bibitem[Ham93]{Hamilton_Harnack}
Richard~S. Hamilton.
\newblock The {H}arnack estimate for the {R}icci flow.
\newblock {\em J. Differential Geom.}, 37(1):225--243, 1993.

\bibitem[Ham95a]{Hamilton_compactness}
Richard~S. Hamilton.
\newblock A compactness property for solutions of the {R}icci flow.
\newblock {\em Amer. J. Math.}, 117(3):545--572, 1995.

\bibitem[Ham95b]{Hamilton_formsing}
Richard~S. Hamilton.
\newblock The formation of singularities in the {R}icci flow.
\newblock In {\em Surveys in differential geometry, {V}ol.\ {II} ({C}ambridge,
  {MA}, 1993)}, pages 7--136. Int. Press, Cambridge, MA, 1995.

\bibitem[Ive93]{Ivey}
Thomas Ivey.
\newblock Ricci solitons on compact three-manifolds.
\newblock {\em Differential Geom. Appl.}, 3(4):301--307, 1993.

\bibitem[Jef05]{Jeffres}
Thalia~D. Jeffres.
\newblock A maximum principle for parabolic equations on manifolds with cone
  singularities.
\newblock {\em Adv. Geom.}, 5(2):319--323, 2005.

\bibitem[KL08]{KleinerLott}
Bruce Kleiner and John Lott.
\newblock Notes on {P}erelman's papers.
\newblock {\em Geom. Topol.}, 12(5):2587--2855, 2008.

\bibitem[LY86]{LiYau}
Peter Li and Shing-Tung Yau.
\newblock On the parabolic kernel of the {S}chr\"odinger operator.
\newblock {\em Acta Math.}, 156(3-4):153--201, 1986.

\bibitem[Mel93]{Melrose}
Richard~B. Melrose.
\newblock {\em The {A}tiyah-{P}atodi-{S}inger index theorem}, volume~4 of {\em
  Research Notes in Mathematics}.
\newblock A K Peters Ltd., Wellesley, MA, 1993.

\bibitem[Mil63]{Mil}
J.~Milnor.
\newblock {\em Morse theory}.
\newblock Based on lecture notes by M. Spivak and R. Wells. Annals of
  Mathematics Studies, No. 51. Princeton University Press, Princeton, N.J.,
  1963.

\bibitem[MRS13]{MazRubSes}
Rafe Mazzeo, Yanir Rubinstein, and Natasa Sesum.
\newblock Ricci flow on surfaces with conic singularities.
\newblock {\em Preprint, arXiv:1306.6688 [math.DG]}, 2013.

\bibitem[Per02]{Perelman1}
Grisha Perelman.
\newblock The entropy formula for the {R}icci flow and its geometric
  applications.
\newblock {\em arXiv:math/0211159 [math.DG]}, 2002.

\bibitem[Per03a]{Perelman3}
Grisha Perelman.
\newblock Finite extinction time for the solutions to the {R}icci flow on
  certain three-manifolds.
\newblock {\em arXiv:math/0307245 [math.DG]}, 2003.

\bibitem[Per03b]{Perelman2}
Grisha Perelman.
\newblock {R}icci flow with surgery on three-manifolds.
\newblock {\em arXiv:math/0303109 [math.DG]}, 2003.

\bibitem[Pet87]{Peters}
Stefan Peters.
\newblock Convergence of {R}iemannian manifolds.
\newblock {\em Compositio Math.}, 62(1):3--16, 1987.

\bibitem[pro]{progP4}
P4: {P}olynomial {P}lanar {P}hase {P}ortraits.
\newblock Program available at \url{http://mat.uab.cat/~artes/p4/p4.htm}.

\bibitem[Ram12]{Ramos2}
Daniel Ramos.
\newblock An asymptotically cusped three dimensional expanding gradient {R}icci
  soliton.
\newblock {\em Preprint, arXiv:1211.4513 [math.DG]}, 2012.

\bibitem[Ram13]{Ramos3}
Daniel Ramos.
\newblock Gradient {R}icci solitons on surfaces.
\newblock {\em Preprint, arXiv:1304.6391 [math.DG]}, 2013.

\bibitem[Ram15]{Ramos1}
Daniel Ramos.
\newblock Smoothening cone points with {R}icci flow.
\newblock {\em Bull. Soc. Math. France}, 143(4):619--633, 2015.
\newblock Preprint arXiv:1107.4591 [math.DG] (2011).

\bibitem[Rat06]{Ratcliffe}
John~G. Ratcliffe.
\newblock {\em Foundations of hyperbolic manifolds}, volume 149 of {\em
  Graduate Texts in Mathematics}.
\newblock Springer, New York, second edition, 2006.

\bibitem[Ric12]{Richard}
Thomas Richard.
\newblock Canonical smoothing of compact {A}lexandrov surfaces via {R}icci
  flow.
\newblock {\em Preprint, arXiv:1204.5461 [math.DG]}, 2012.

\bibitem[Shi89]{Shi}
Wan-Xiong Shi.
\newblock Deforming the metric on complete {R}iemannian manifolds.
\newblock {\em J. Differential Geom.}, 30(1):223--301, 1989.

\bibitem[Sim12]{Simon}
Miles Simon.
\newblock Ricci flow of non-collapsed three manifolds whose {R}icci curvature
  is bounded from below.
\newblock {\em J. Reine Angew. Math.}, 662:59--94, 2012.

\bibitem[Top06]{Topping_lectures}
Peter Topping.
\newblock {\em Lectures on the {R}icci flow}, volume 325 of {\em London
  Mathematical Society Lecture Note Series}.
\newblock Cambridge University Press, Cambridge, 2006.

\bibitem[Top12]{Topping_revcusp}
Peter~M. Topping.
\newblock Uniqueness and nonuniqueness for {R}icci flow on surfaces: reverse
  cusp singularities.
\newblock {\em Int. Math. Res. Not. IMRN}, (10):2356--2376, 2012.

\bibitem[Tro91]{Troyanov}
Marc Troyanov.
\newblock Prescribing curvature on compact surfaces with conical singularities.
\newblock {\em Trans. Amer. Math. Soc.}, 324(2):793--821, 1991.

\bibitem[Wu91]{Wu}
Lang-Fang Wu.
\newblock The {R}icci flow on {$2$}-orbifolds with positive curvature.
\newblock {\em J. Differential Geom.}, 33(2):575--596, 1991.

\bibitem[Ye04]{Ye}
Rugang Ye.
\newblock On uniqueness of {$2$}-dimensional {$\kappa$}-solutions.
\newblock {\em Preprint}, 2004.
\newblock \url{http://www.math.ucsb.edu/~yer/2dkappa.pdf}.

\bibitem[Yin10]{Yin1}
Hao Yin.
\newblock {R}icci flow on surfaces with conical singularities.
\newblock {\em J. Geom. Anal.}, 20(4):970--995, 2010.

\bibitem[Yin13]{Yin2}
Hao Yin.
\newblock {R}icci flow on surfaces with conical singularities, {II}.
\newblock {\em Preprint, arXiv:1305.4355 [math.DG]}, 2013.

\end{thebibliography}
